\documentclass[final,leqno,onefignum,onetabnum]{siamltex1213}
\usepackage{amsfonts}
\usepackage{amsmath,booktabs,ctable,threeparttable}
\usepackage{amssymb,amsfonts,boxedminipage}
\newtheorem{remark}{Remark}[section]
\newtheorem{prop}{Proposition}[section]

\title{The Regularization Theory of the Krylov Iterative Solvers LSQR, CGLS,
LSMR and CGME For Linear Discrete Ill-Posed
Problems\thanks{This work was supported in part by
the National Science Foundation of China (No. 11371219)}}

\author{Zhongxiao Jia\thanks{Department of Mathematical Sciences, Tsinghua
University, 100084 Beijing, China. (\email{jiazx@tsinghua.edu.cn})}}

\begin{document}
\maketitle
\slugger{sirev}{xxxx}{xx}{x}{x--x}

\begin{abstract}
For the large-scale linear discrete ill-posed problem $\min\|Ax-b\|$ or $Ax=b$
with $b$ contaminated by a white noise, the Lanczos bidiagonalization based Krylov
solver LSQR and its mathematically equivalent CGLS are
most commonly used. They have intrinsic regularizing effects,
where the number $k$ of iterations plays the role
of regularization parameter. However, there has been no
answer to the long-standing fundamental concern: {\em for which kinds of
problems LSQR and CGLS can find best possible regularized solutions}?
The concern was actually expressed foresightedly by Bj\"{o}rck and
Eld\'{e}n in 1979.
Here a best possible regularized solution means that it is at least
as accurate as the best regularized
solution obtained by the truncated singular value decomposition (TSVD) method,
which and the best possible solution by standard-form Tikhonov regularization
are both of the same order of the worst-case error and cannot be improved
under the assumption that the solution to an underlying linear compact operator
equation is continuous or its derivative squares integrable.
In this paper we make a detailed analysis on the
regularization of LSQR for severely, moderately and mildly ill-posed problems.
We first consider the case that the singular values of $A$ are simple. We establish
accurate $\sin\Theta$ theorems for the 2-norm distance between the
underlying $k$-dimensional
Krylov subspace and the $k$-dimensional dominant right singular subspace
of $A$. Based on them and some follow-up results, for the first two kinds of
problems, we prove that LSQR finds a best possible regularized solution at
semi-convergence occurring at iteration $k_0$ and the following results hold for
$k=1,2,\ldots,k_0$: (i) the $k$-step Lanczos
bidiagonalization always generates a near best rank $k$ approximation to $A$;
(ii) the $k$ Ritz values always approximate the first $k$ large singular values
in natural order; (iii) the $k$-step LSQR
always captures the $k$ dominant SVD components of $A$. However,
for the third kind of problem, we prove that LSQR cannot find a best possible
regularized solution generally. We derive accurate estimates for
the diagonals and subdiagonals of the bidiagonal matrices generated by
Lanczos bidiagonalization, which can be used to decide if LSQR finds a best possible
regularized solution at semi-convergence. We also analyze the regularization
of the other two Krylov solvers LSMR and CGME that are
MINRES and the CG method applied to $A^TAx=A^Tb$ and
$\min\|AA^Ty-b\|$ with $x=A^Ty$, respectively, proving that the regularizing effects
of LSMR are similar to LSQR for each kind of problem and both are
superior to CGME.
We extend all the results to the case that $A$ has multiple singular values.
Numerical experiments confirm our theory on LSQR.
\end{abstract}



\pagestyle{myheadings}
\thispagestyle{plain}
\markboth{ZHONGXIAO JIA}{REGULARIZATION OF LSQR, CGLS, LSMR AND CGME}

\section{Introduction and Preliminaries}\label{intro}

Consider the linear discrete ill-posed problem
\begin{equation}
  \min\limits_{x\in \mathbb{R}^{n}}\|Ax-b\| \mbox{\,\ or \ $Ax=b$,}
  \ \ \ A\in \mathbb{R}^{m\times n}, \label{eq1}
  \ b\in \mathbb{R}^{m},
\end{equation}
where the norm $\|\cdot\|$ is the 2-norm of a vector or matrix, and
$A$ is extremely ill conditioned with its singular values decaying
to zero without a noticeable gap. \eqref{eq1} mainly arises
from the discretization of the first kind Fredholm integral equation
\begin{equation}\label{eq2}
Kx=(Kx)(t)=\int_{\Omega} k(s,t)x(t)dt=g(s)=g,\ s\in \Omega
\subset\mathbb{R}^q,
\end{equation}
where the kernel $k(s,t)\in L^2({\Omega\times\Omega})$ and
$g(s)$ are known functions, while $x(t)$ is the
unknown function to be sought. If $k(s,t)$ is non-degenerate
and $g(s)$ satisfies the Picard condition, there exists the unique squares
integrable solution
$x(t)$; see \cite{engl00,hansen98,hansen10,kirsch,mueller}. Here for brevity
we assume that $s$ and $t$ belong to the same set $\Omega\subset
\mathbb{R}^q$ with $q\geq 1$.
Applications include image deblurring, signal processing, geophysics,
computerized tomography, heat propagation, biomedical and optical imaging,
groundwater modeling, and many others; see, e.g.,
\cite{aster,engl93,engl00,hansen10,ito15,kaipio,kern,kirsch,mueller,
natterer,vogel02}.
The theory and numerical treatments of integral
equations can be found in \cite{kirsch,kythe}.
The right-hand side $b=\hat{b}+e$ is noisy and assumed to be
contaminated by a white noise $e$, caused by measurement, modeling
or discretization errors, where $\hat{b}$
is noise-free and $\|e\|<\|\hat{b}\|$.
Because of the presence of noise $e$ and the extreme
ill-conditioning of $A$, the naive
solution $x_{naive}=A^{\dagger}b$ of \eqref{eq1} bears no relation to
the true solution $x_{true}=A^{\dagger}\hat{b}$, where
$\dagger$ denotes the Moore-Penrose inverse of a matrix.
Therefore, one has to use regularization to extract a
best possible approximation to $x_{true}$.

The most common regularization, in its simplest form, is the
direct standard-form Tikhonov regularization
\begin{equation}\label{tikhonov}
  \min\limits_{x\in \mathbb{R}^{n}}{\|Ax-b\|^2+\lambda^2\|x\|^2}
\end{equation}
with $\lambda>0$ the regularization parameter
\cite{phillips,tikhonov63,tikhonov77}.
The solutions to \eqref{eq1} and \eqref{tikhonov} can be fully analyzed by
the singular value decomposition (SVD) of $A$. Let
\begin{equation}\label{eqsvd}
  A=U\left(\begin{array}{c} \Sigma \\ \mathbf{0} \end{array}\right) V^{T}
\end{equation}
be the SVD of $A$,
where $U = (u_1,u_2,\ldots,u_m)\in\mathbb{R}^{m\times m}$ and
$V = (v_1,v_2,\ldots,v_n)\in\mathbb{R}^{n\times n}$ are orthogonal,
$\Sigma = \diag (\sigma_1,\sigma_2,
\ldots,\sigma_n)\in\mathbb{R}^{n\times n}$ with the singular values
$\sigma_1>\sigma_2 >\cdots >\sigma_n>0$ assumed to be simple
throughout the paper except Section \ref{multiple}, and the superscript $T$ denotes
the transpose of a matrix or vector. Then
\begin{equation}\label{eq4}
  x_{naive}=\sum\limits_{i=1}^{n}\frac{u_i^{T}b}{\sigma_i}v_i =
  \sum\limits_{i=1}^{n}\frac{u_i^{T}\hat{b}}{\sigma_i}v_i +
  \sum\limits_{i=1}^{n}\frac{u_i^{T}e}{\sigma_i}v_i
  =x_{true}+\sum\limits_{i=1}^{n}\frac{u_i^{T}e}{\sigma_i}v_i
\end{equation}
with $\|x_{true}\|=\|A^{\dagger}\hat{b}\|=
\left(\sum_{k=1}^n\frac{|u_k^T\hat{b}|^2}{\sigma_k^2}\right)^{1/2}$.

Throughout the paper, we always assume that $\hat{b}$ satisfies the discrete
Picard condition $\|A^{\dagger}\hat{b}\|\leq C$ with some constant $C$ for $n$
arbitrarily large \cite{aster,gazzola15,hansen90,hansen90b,hansen98,hansen10,kern}.
It is an analog of the Picard condition in the  finite dimensional case;
see, e.g., \cite{hansen90}, \cite[p.9]{hansen98},
\cite[p.12]{hansen10} and \cite[p.63]{kern}. This
condition means that, on average, the Fourier coefficients
$|u_i^{T}\hat{b}|$ decay faster than $\sigma_i$ and enables
regularization to compute useful approximations to $x_{true}$, which
results in the following popular model that is used throughout Hansen's books
\cite{hansen98,hansen10} and the current paper:
\begin{equation}\label{picard}
  | u_i^T \hat b|=\sigma_i^{1+\beta},\ \ \beta>0,\ i=1,2,\ldots,n,
\end{equation}
where $\beta$ is a model parameter that controls the decay rates of
$| u_i^T \hat b|$.
Hansen \cite[p.68]{hansen10} points out, ``while this is a crude model,
it reflects the overall behavior often found in real problems."
One precise definition of the discrete Picard condition is
$| u_i^T \hat b|=\tau_i\sigma_i^{1+\zeta_i}$ with certain
constants $\tau_i\geq 0,\ \zeta_i>0,\ i=1,2,\ldots,n$. We remark that once
the $\tau_i>0$ and $\zeta_i$ do not differ greatly, such discrete Picard
condition does not affect our claims, rather it
complicates derivations and forms of the results.

The white noise $e$ has a number of attractive properties which
play a critical role in the regularization analysis: Its covariance matrix
is $\eta^2 I$, the expected values ${\cal E}(\|e\|^2)=m \eta^2$ and
${\cal E}(|u_i^Te|)=\eta,\,i=1,2,\ldots,n$, and
$\|e\|\approx \sqrt{m}\eta$ and $|u_i^Te|\approx \eta,\
i=1,2,\ldots,n$; see, e.g., \cite[p.70-1]{hansen98} and \cite[p.41-2]{hansen10}.
The noise $e$ thus affects $u_i^Tb,\ i=1,2,\ldots,n,$ {\em more or less equally}.
With \eqref{picard}, relation \eqref{eq4} shows that for large singular values
$|{u_i^{T}\hat{b}}|/{\sigma_i}$ is dominant relative to
$|u_i^{T}e|/{\sigma_i}$. Once
$| u_i^T \hat b| \leq | u_i^T e|$ from some $i$ onwards, the small singular
values magnify $|u_i^{T}e|/{\sigma_i}$, and the noise
$e$ dominates $| u_i^T b|/\sigma_i$ and must be suppressed. The
transition point $k_0$ is such that
\begin{equation}\label{picard1}
| u_{k_0}^T b|\approx | u_{k_0}^T \hat{b}|> | u_{k_0}^T e|\approx
\eta, \ | u_{k_0+1}^T b|
\approx | u_{k_0+1}^Te|
\approx \eta;
\end{equation}
see \cite[p.42, 98]{hansen10} and a similar description \cite[p.70-1]{hansen98}.
The $\sigma_k$ are then divided into the $k_0$ large ones and the $n-k_0$ small
ones. The truncated SVD (TSVD) method \cite{hansen98,hansen10} computes the TSVD
regularized solutions
\begin{equation}\label{solution}
  x^{tsvd}_k=\left\{\begin{array}{ll} \sum\limits_{i=1}^{k}\frac{u_i^{T}b}
  {\sigma_i}{v_i}\thickapprox
  \sum\limits_{i=1}^{k}\frac{u_i^{T}\hat{b}}
{\sigma_i}{v_i},\ \ \ &k\leq k_0;\\ \sum\limits_{i=1}^{k}\frac{u_i^{T}b}
{\sigma_i}{v_i}\thickapprox
\sum\limits_{i=1}^{k_0}\frac{u_i^{T}\hat{b}}{\sigma_i}{v_i}+
\sum\limits_{i=k_0+1}^{k}\frac{u_i^{T}e}{\sigma_i}{v_i},\ \ \ &k>k_0.
\end{array}\right.
\end{equation}
It is known from \cite[p.70-1]{hansen98} and
\cite[p.86-8,96]{hansen10} that $x_{k_0}^{tsvd}$ is
the best TSVD regularized solution to \eqref{eq1} and balances the
regularization and perturbation errors optimally. The parameter $k$ is a
regularization parameter that determines how many large SVD components of
$A$ are used to compute a regularized solution $x_k^{tsvd}$ to
\eqref{eq1}.

Let $U_k=(u_1,\ldots,u_k)$, $V_k=(v_1,\ldots,v_k)$ and $\Sigma_k=
{\rm diag}(\sigma_1,\ldots,\sigma_k)$, and define $A_k=U_k\Sigma_k V_k^T$.
Then $A_k$ is the best rank $k$ approximation to $A$ with
$\|A-A_k\|=\sigma_{k+1}$ (cf. \cite[p.12]{bjorck96}), and
$
x_{k}^{tsvd}=A_k^{\dagger}b
$
is the minimum-norm least squares solution to
$$
  \min\limits_{x\in \mathbb{R}^{n}}\|A_kx-b\|
$$
that perturbs $A$ to $A_k$ in \eqref{eq1}. This interpretation will be
often exploited later.

The solution $x_{\lambda}$ of the Tikhonov regularization
has a filtered SVD expansion
\begin{equation}\label{eqfilter}
  x_{\lambda} = \sum\limits_{i=1}^{n}f_i\frac{u_i^{T}b}{\sigma_i}v_i,
\end{equation}
where the $f_i=\frac{\sigma_i^2}{\sigma_i^2+\lambda^2}$ are called filters.
The TSVD method is a special parameter filtered method, where, in $x_k^{tsvd}$,
we take $f_i=1,\ i=1,2,\ldots,k$ and $f_i=0,\ i=k+1,\ldots,n$. The error
$x_{\lambda}-x_{true}$ can be written as the sum of the regularization and
perturbation errors, and an optimal $\lambda_{opt}$ aims to balance
these two errors and make the sum of their norms minimized
\cite{hansen98,hansen10,kirsch,vogel02}. The best possible regularized
solution $x_{\lambda_{opt}}$ retains the $k_0$ dominant SVD components
and dampens the other $n-k_0$ small SVD components as much as
possible \cite{hansen98,hansen10}. Apparently, the ability to acquire only
the largest SVD components of $A$ is fundamental in solving \eqref{eq1}.

A number of parameter-choice methods have been developed for finding
$\lambda_{opt}$ or $k_0$, such as the discrepancy principle \cite{morozov},
the L-curve criterion, whose use goes back to Miller \cite{miller} and
Lawson and Hanson \cite{lawson} and is termed much later and studied in detail
in \cite{hansen92,hansen93}, and the generalized cross validation
(GCV) \cite{golub79,wahba}; see, e.g.,
\cite{bauer11,hansen98,hansen10,kern,kilmer03,kindermann,neumaier98,reichel13,vogel02}
for numerous comparisons. All parameter-choice methods aim to
make $f_i/\sigma_i$ not small for $i=1,2,\ldots,k_0$
and $f_i/\sigma_i\approx 0$ for $i=k_0+1,\ldots,n$. Each of these methods
has its own merits and disadvantages, and
no one is absolutely reliable for all ill-posed problems.
For example, some of the mentioned parameter-choice methods
may fail to find accurate approximations to $\lambda_{opt}$;
see \cite{hanke96a,vogel96} for an analysis
on the L-curve method and \cite{hansen98} for some other parameter-choice
methods. A further investigation on paramater-choice methods is not our concern
in this paper.

The TSVD method is important in its own right.
It and the standard-form Tikhonov regularization
produce very similar solutions with essentially the minimum
2-norm error, i.e., the worst-case error \cite[p.13]{kirsch};
see \cite{varah79}, \cite{hansen90b}, \cite[p.109-11]{hansen98} and
\cite[Sections 4.2 and 4.4]{hansen10}. Indeed, for a linear compact
equation $Kx=g$ including \eqref{eq2} with the noisy $g$ and true solution
$x_{true}(t)$, under the source condition that its solution $x_{true}(t)
\in {\cal R}(K^*)$ or $x_{true}(t)\in {\cal R}(K^*K)$, the range of
the adjoint $K^*$ of $K$ or
that of $K^*K$, which amounts to assuming that $x_{true}(t)$ or its
derivative is squares integrable, the errors of
the best regularized solutions by the TSVD method and the Tikhonov
regularization are {\em order optimal, i.e., the same order
as the worst-case error} \cite[p.13,18,20,32-40]{kirsch},
\cite[p.90]{natterer} and \cite[p.7-12]{vogel02}. These conclusions
carries over to \eqref{eq1} \cite[p.8]{vogel02}.
Therefore, either of $x_{\lambda_{opt}}$ and $x_{k_0}^{tsvd}$
is a best possible solution to \eqref{eq1} under the above
assumptions and can be taken as standard reference
when assessing the regularizing effects of an iterative
solver. For the sake of clarity, we will take $x_{k_0}^{tsvd}$.

For \eqref{eq1} large, the TSVD method and the Tikhonov regularization
method are generally too demanding, and only iterative regularization
methods are computationally viable. A major class of methods has been
Krylov iterative solvers that project \eqref{eq1} onto a sequence of
low dimensional Krylov subspaces
and computes iterates to approximate $x_{true}$; see, e.g.,
\cite{aster,engl00,gilyazov,hanke95,hansen98,hansen10,kirsch}.
Of Krylov iterative solvers, the CGLS (or CGNR) method,
which implicitly applies the Conjugate Gradient (CG)
method \cite{golub89,hestenes} to the normal equations $A^TAx=A^Tb$
of \eqref{eq1}, and its mathematically equivalent LSQR algorithm \cite{paige82}
have been most commonly used. The Krylov solvers CGME
(or CGNE) \cite{bjorck96,bjorck15,craig,hanke95,hanke01,hps09} and
LSMR \cite{bjorck15,fong} are also choices, which amount to the
CG method applied to $\min\|AA^Ty-b\|$ with $x=A^Ty$ and MINRES \cite{paige75}
applied to $A^TAx=A^Tb$, respectively. These Krylov solvers have been
intensively studied and known to have regularizing
effects \cite{aster,eicke,gilyazov,hanke95,hanke01,hansen98,hansen10,hps09}
and exhibit semi-convergence \cite[p.89]{natterer};
see also \cite[p.314]{bjorck96}, \cite[p.733]{bjorck15},
\cite[p.135]{hansen98} and \cite[p.110]{hansen10}: The iterates
converge to $x_{true}$ and their norms increase steadily,
and the residual norms decrease in an initial stage; then afterwards the
noise $e$ starts to deteriorate the iterates so that they start to diverge
from $x_{true}$ and instead converge to $x_{naive}$,
while their norms increase considerably and the residual norms stabilize.
If we stop at the right time, then, in principle,
we have a regularization method, where the iteration number plays the
role of parameter regularization.
Semi-convergence is due to the fact that the projected problem starts to
inherit the ill-conditioning of \eqref{eq1} from some iteration
onwards, and the appearance of a small singular
value of the projected problem amplifies the noise considerably.

The regularizing effects of CG type methods were noticed by
Lanczos \cite{lanczos} and were rediscovered in \cite{johnsson,squire,tal}.
Based on these works and motivated by a heuristic explanation on good
numerical results with very few iterations using CGLS
in \cite{johnsson}, and realizing that such an excellent performance
can only be expected if convergence to the regular part
of the solution, i.e., $x_{k_0}^{tsvd}$, takes place before the effects of
ill-posedness show up, on page 13 of \cite{bjorck79},
Bj\"{o}rck and Eld\'{e}n in 1979 foresightedly expressed a fundamental concern
on CGLS (and LSQR): {\em More
research is needed to tell for which problems this approach will work, and
what stopping criterion to choose.} See also \cite[p.145]{hansen98}.
As remarked by Hanke and Hansen \cite{hanke93}, the paper \cite{bjorck79}
was the only extensive survey on algorithmic details until that time,
and a strict proof of the regularizing properties of conjugate gradients is
extremely difficult.
An enormous effort has long been made to the study of
regularizing effects of LSQR and CGLS (cf.
\cite{eicke,engl00,firro97,gilyza86,hanke95,hanke01,hansen98,hansen10,
hps09,kirsch,nemi,nolet,paige06,scales,vorst90}) in the Hilbert
or finite dimensional space setting, but a rigorous regularization theory
of LSQR and CGLS for \eqref{eq1} is still lacking,
and there has been no definitive answer to the above long-standing
fundamental question, and the same is for LSMR and CGME.

For $A$ symmetric, MINRES and MR-II applied to $Ax=b$
directly are alternatives and have been
shown to have regularizing effects \cite{calvetti01,hanke95,hanke96,
hansen10,jensen07,kilmer99}, but MR-II seems preferable since
the noisy $b$ is excluded in the underlying subspace \cite{huang15,jensen07}.
For $A$ nonsymmetric or multiplication with $A^{T}$ difficult to compute,
GMRES and RRGMRES are candidate methods
\cite{baglama07,calvetti02,calvetti02c,neuman12}, and the latter
may be better \cite{jensen07}. The hybrid
approaches based on the Arnoldi process have been first
proposed in \cite{calvetti00b} and studied in
\cite{calvetti01,calvetti03,lewis09,novati13}.
Gazzola and her coauthors \cite{gazzola14}--\cite{gazzola-online}
have described a general framework of the hybrid methods
and presented various Krylov-Tikhonov methods with different parameter-choice
strategies. Unfortunately, unlike LSQR and CGLS, these methods are
highly problem dependent and may not have regularizing effects for general
nonsymmetric ill-posed problems; see, e.g., \cite{jensen07} and \cite[p.126]{hansen10}.
The fundamental cause is that the underlying Krylov subspaces may not favor the
the dominant left and singular subspaces of $A$,
which are desired in solving \eqref{eq1}.

The behavior of ill-posed problems critically depends on the decay rate of
$\sigma_j$. The following characterization of the degree of ill-posedness
of \eqref{eq1} was introduced in \cite{hofmann86}
and has been widely used \cite{aster,engl00,hansen98,hansen10,mueller}:
If $\sigma_j=\mathcal{O}(j^{-\alpha})$, then \eqref{eq1}
is mildly or moderately ill-posed for $\frac{1}{2}<\alpha\le1$ or $\alpha>1$.
If $\sigma_j=\mathcal{O}(\rho^{-j})$ with $\rho>1$,
$j=1,2,\ldots,n$, then \eqref{eq1} is severely ill-posed.
Here for mildly ill-posed problems we add the requirement
$\alpha>\frac{1}{2}$, which does not appear in \cite{hofmann86}
but must be met for $k(s,t)\in L^2({\Omega\times\Omega})$
in \eqref{eq1} \cite{hanke93,hansen98}.
In the one-dimensional case, i.e., $q=1$, \eqref{eq1}
is severely ill-posed with $k(s,t)$ sufficiently smooth, and
it is moderately ill-posed with $\sigma_j=\mathcal{O}(j^{-p-1/2})$,
where $p$ is the highest order of continuous derivatives of
$k(s,t)$; see, e.g., \cite[p.8]{hansen98} and \cite[p.10-11]{hansen10}.
Clearly, the singular values $\sigma_j$ for a
severely ill-posed problem decay at the same rate $\rho^{-1}$,
while those of a moderately or mildly ill-posed problem decay
at the decreasing rate $\left(\frac{j}{j+1}\right)^{\alpha}$
that approaches one more quickly with $j$ for the mildly ill-posed problem than
for the moderately ill-posed problem.

If a regularized solution to \eqref{eq1} is at least as accurate as
$x_{k_0}^{tsvd}$, then it is called a best possible regularized solution.
Given \eqref{eq1}, if the regularized solution of an iterative regularization
solver at semi-convergence is a best possible one,
then, by the words of Bj\"{o}rck and Eld\'{e}n,
the solver {\em works} for the problem and is said to have
the {\em full} regularization. Otherwise, the solver is said to have
the {\em partial} regularization.

Because it has been unknown whether or not LSQR, CGLS, LSMR and CGME
have the full regularization for a given \eqref{eq1},
one commonly combines them with some explicit
regularization, so that the resulting hybrid variants
(hopefully) find best possible regularized solutions
\cite{aster,hansen98,hansen10}. A hybrid CGLS is to run CGLS
for several trial regularization parameters $\lambda$ and
picks up the best one among the candidates \cite{aster}. Its
disadvantages are
that regularized solutions cannot be updated with different $\lambda$
and there is no guarantee that the selected regularized solution
is a best possible one.
The hybrid LSQR variants have been advocated by Bj\"{o}rck and Eld\'{e}n
\cite{bjorck79} and O'Leary and Simmons \cite{oleary81}, and improved and
developed by Bj\"orck \cite{bjorck88} and Bj\"{o}rck, Grimme and
Van Dooren \cite{bjorck94}.
A hybrid LSQR first projects \eqref{eq1} onto Krylov
subspaces and then regularizes the projected problems explicitly.
It aims to remove the effects
of small Ritz values and expands a Krylov subspace until it
captures the $k_0$ dominant SVD components of $A$
\cite{bjorck88,bjorck94,hanke93,oleary81}.
The hybrid LSQR and CGME have been intensively studied in, e.g.,
\cite{bazan10,bazan14,berisha,chung08,hanke01,hanke93,lewis09,
neuman12,novati13,renaut} and \cite{aster,hansen10,hansen13}.
Within the framework of such hybrid solvers, it is hard to find a near-optimal
regularization parameter \cite{bjorck94,renaut}. More seriously,
as we will elaborate mathematically and numerically in the concluding
section of this paper, it may make no sense to speak of
the regularization of the projected problems and their optimal
regularization parameters since they may actually fail to satisfy
the discrete Picard conditions. In contrast,
if an iterative solver has the full regularization, we stop it after
semi-convergence. Obviously, we cannot emphasize too much
the importance of completely understanding the
regularization of LSQR, CGLS, LSMR and CGME. By the definition of
the full or partial regularization,
we now modify the concern of Bj\"{o}rck and Eld\'{e}n as:
{\em  Do LSQR, CGLS, LSMR and CGME have the full or partial regularization for
severely, moderately and mildly ill-posed problems? How to identify
their full or partial regularization in practice?}

In this paper, assuming exact arithmetic, we first focus on LSQR and make a
rigorous analysis on its regularization for severely,
moderately and mildly ill-posed problems. Due to the mathematical equivalence
of CGLS and LSQR, the assertions on the full or partial regularization of
LSQR apply to CGLS as well. We then analyze the regularizing
effects of LSMR and CGME and draw definitive conclusions.
We prove that LSQR has the full regularization for severely and
moderately ill-posed problems once $\rho>1$ and $\alpha>1$ suitably,
and it generally has only the partial regularization for mildly ill-posed
problems. In Section \ref{lsqr}, we describe the
Lanczos bidiagonalization process and LSQR, and make an
introductory analysis. In Section \ref{sine},
we establish accurate $\sin\Theta$ theorems for the 2-norm
distance between the underlying $k$-dimensional Krylov subspace and the
$k$-dimensional dominant right singular subspace of $A$. We then derive
some follow-up results that play a central role in analyzing
the regularization of LSQR. In Section \ref{rankapp},
we prove that a $k$-step Lanczos bidiagonalization always
generates a near best rank $k$ approximation to $A$, and the $k$ Ritz values
always approximate the first $k$ large singular values in natural order,
and no small Ritz value appears for $k=1,2,\ldots,k_0$.
This will show that LSQR has the full regularization.
For mildly ill-posed problems,
we prove that, for some $k\leq k_0$, the $k$ Ritz values generally
do not approximate the first $k$ large singular values in natural order
and LSQR generally has only the partial regularization.
In Section \ref{alphabeta}, we derive bounds for the entries of bidiagonal
matrices generated by Lanczos bidiagonalization, proving how fast they
decay and showing how to use them to reliably identify if
LSQR has the full regularization when
the degree of ill-posedness of \eqref{eq1} is unknown in advance.
Exploiting some of the results on LSQR, we analyze the regularization of
LSMR and CGME and prove that LSMR has similar regularizing effects to LSQR
for each kind of problem and both of them are superior to CGME.
In Section \ref{compare}, we present some perturbation results and
prove that LSQR resembles the TSVD method for severely and moderately
ill-posed problems. In Section \ref{multiple}, with a number of
nontrivial changes and reformulations, we extend all the results to
the case that $A$ has multiple singular values. In Section \ref{numer},
we report numerical experiments to confirm our theory on LSQR. Finally,
we summarize the paper with further remarks in Section \ref{concl}.

Throughout the paper, denote by $\mathcal{K}_{k}(C, w)=
span\{w,Cw,\ldots,C^{k-1}w\}$  the $k$-dimensional Krylov subspace generated
by the matrix $\mathit{C}$ and the vector $\mathit{w}$, and by $I$ and the
bold letter $\mathbf{0}$ the identity matrix
and the zero matrix with orders clear from the context, respectively.
For $B=(b_{ij})$, we define the nonnegative matrix $|B|=(|b_{ij}|)$,
and for $|C|=(|c_{ij}|)$, $|B|\leq |C|$ means
$|b_{ij}|\leq |c_{ij}|$ componentwise.

\section{The LSQR algorithm}\label{lsqr}

LSQR is based on the Lanczos bidiagonalization process, which
computes two orthonormal bases $\{q_1,q_2,\dots,q_k\}$ and
$\{p_1,p_2,\dots,p_{k+1}\}$  of $\mathcal{K}_{k}(A^{T}A,A^{T}b)$ and
$\mathcal{K}_{k+1}(A A^{T},b)$  for $k=1,2,\ldots,n$,
respectively. We describe the process as Algorithm 1.

{\bf Algorithm 1: \  $k$-step Lanczos bidiagonalization process}
\begin{remunerate}
\item Take $ p_1=b/\|b\| \in \mathbb{R}^{m}$, and define $\beta_1{q_0}=0$.

\item For $j=1,2,\ldots,k$
 \begin{romannum}
  \item
  $r = A^{T}p_j - \beta_j{q_{j-1}}$
  \item $\alpha_j = \|r\|;q_j = r/\alpha_j$
  \item
   $   z = Aq_j - \alpha_j{p_{j}}$
  \item
  $\beta_{j+1} = \|z\|;p_{j+1} = z/\beta_{j+1}.$
   \end{romannum}
\end{remunerate}

Algorithm 1 can be written in the matrix form
\begin{align}
  AQ_k&=P_{k+1}B_k,\label{eqmform1}\\
  A^{T}P_{k+1}&=Q_{k}B_k^T+\alpha_{k+1}q_{k+1}e_{k+1}^{T}.\label{eqmform2}
\end{align}
where $e_{k+1}$ denotes the $(k+1)$-th canonical basis vector of
$\mathbb{R}^{k+1}$, $P_{k+1}=(p_1,p_2,\ldots,p_{k+1})$,
$Q_k=(q_1,q_2,\ldots,q_k)$ and
\begin{equation}\label{bk}
  B_k = \left(\begin{array}{cccc} \alpha_1 & & &\\ \beta_2 & \alpha_2 & &\\ &
  \beta_3 &\ddots & \\& & \ddots & \alpha_{k} \\ & & & \beta_{k+1}
  \end{array}\right)\in \mathbb{R}^{(k+1)\times k}.
\end{equation}
It is known from \eqref{eqmform1} that
\begin{equation}\label{Bk}
B_k=P_{k+1}^TAQ_k.
\end{equation}
We remind that the singular values of $B_k$, called the Ritz values of $A$ with
respect to the left and right subspaces $span\{P_{k+1}\}$ and $span\{Q_k\}$,
are all simple. This basic fact will often be used later.

At iteration $k$, LSQR solves the problem
$\|Ax^{(k)}-b\|=\min_{x\in \mathcal{K}_k(A^TA,A^Tb)}
\|Ax-b\|$ and computes the iterates $x^{(k)}=Q_ky^{(k)}$ with
\begin{equation}\label{yk}
  y^{(k)}=\arg\min\limits_{y\in \mathbb{R}^{k}}\|B_ky-\|b\|e_1^{(k+1)}\|
  =\|b\| B_k^{\dagger} e_1^{(k+1)},
\end{equation}
where $e_1^{(k+1)}$ is the first canonical basis vector of $\mathbb{R}^{k+1}$,
and the residual norm $\|Ax^{(k)}-b\|$ decreases monotonically with respect to
$k$. We have $\|Ax^{(k)}-b\|=
\|B_k y^{(k)}-\|b\|e_1^{(k+1)}\|$ and $\|x^{(k)}\|=\|y^{(k)}\|$,
both of which can be cheaply computed.

Note that $\|b\|e_1^{(k+1)}=P_{k+1}^T b$. We have
\begin{equation}\label{xk}
x^{(k)}=Q_k B_k^{\dagger} P_{k+1}^Tb,
\end{equation}
that is, the iterate $x^{(k)}$ by LSQR is the minimum-norm least
squares solution to the perturbed
problem that replaces $A$ in \eqref{eq1} by its rank $k$ approximation
$P_{k+1}B_k Q_k^T$. Recall that the best rank $k$ approximation
$A_k$ to $A$ satisfies $\|A-A_k\|=\sigma_{k+1}$. We can relate
LSQR and the TSVD method from two perspectives. One of them is
to interpret LSQR as solving a nearby problem that perturbs $A_k$ to
$P_{k+1}B_k Q_k^T$, provided that
$P_{k+1}B_k Q_k^T$ is a near best rank $k$ approximation
to $A$ with an approximate accuracy $\sigma_{k+1}$. The other is
to interpret $x_k^{tsvd}$ and $x^{(k)}$ as the solutions to
the two perturbed problems of \eqref{eq1} that replace $A$ by
the rank $k$ approximations with the same quality
to $A$, respectively. Both perspectives lead to the consequence:
the LSQR iterate $x^{(k_0)}$ is as accurate as $x_{k_0}^{tsvd}$
and is thus a best possible regularized solution to \eqref{eq1}, provided that
$P_{k+1}B_k Q_k^T$ is a near best rank $k$
approximation to $A$ with the approximate accuracy $\sigma_{k+1}$
and the $k$ singular values of $B_k$ approximate the first $k$
large ones of $A$ in natural order for $k=1,2,\ldots,k_0$. Otherwise,
as will be clear later, $x^{(k_0)}$ cannot be as accurate as $x_{k_0}^{tsvd}$
if either $P_{k_0+1}B_{k_0} Q_{k_0}^T$ is not a near best rank $k_0$
approximation to $A$ or $B_{k_0}$ has at least one
singular value smaller than $\sigma_{k_0+1}$. We will give a precise
definition of a near best rank $k$ approximation later.

As stated in the introduction, the semi-convergence of LSQR must occur at some
iteration $k$. Under the discrete Picard condition
\eqref{picard}, if semi-convergence
occurs at iteration $k_0$, we are sure that LSQR has the full regularization
because $x^{(k_0)}$ has captured the $k_0$ dominant SVD components of $A$ and
effectively suppressed the other $n-k_0$ SVD components;
if semi-convergence occurs at some iteration
$k<k_0$, then LSQR has only the partial
regularization since it has not yet captured the needed $k_0$ dominant SVD
components of $A$.

\section{$\sin\Theta$ theorems for the distances between $\mathcal{K}_k(A^TA,A^Tb)$
and $span\{V_k\}$ as well as the others related} \label{sine}

Van der Sluis and Van der Vorst \cite{vorst86} prove the following result,
which has been used in Hansen \cite{hansen98} and the references therein
to illustrate and analyze the regularizing effects of LSQR and CGLS. We
will also investigate it further in our paper.

\begin{prop}\label{help}
LSQR with the starting vector $p_1=b/\|b\|$ and CGLS
applied to $A^TAx=A^Tb$ with the starting vector
$x^{(0)}=0$ generate the same iterates
\begin{equation}\label{eqfilter2}
  x^{(k)}=\sum\limits_{i=1}^nf_i^{(k)}\frac{u_i^{T}b}{\sigma_i}v_i,\
  k=1,2,\ldots,n,
\end{equation}
where
\begin{equation}\label{filter}
f_i^{(k)}=1-\prod\limits_{j=1}^k\frac{(\theta_j^{(k)})^2-\sigma_i^2}
{(\theta_j^{(k)})^2},\ i=1,2,\ldots,n,
\end{equation}
and the $\theta_j^{(k)}$ are the singular values of $B_k$
labeled as $\theta_1^{(k)}>\theta_2^{(k)}>\cdots>\theta_k^{(k)}$.
\end{prop}

\eqref{eqfilter2} shows that $x^{(k)}$ has a filtered SVD expansion of
form \eqref{eqfilter}. If all the Ritz values $\theta_j^{(k)}$
approximate the first $k$ singular values
$\sigma_j$ of $A$ in natural order, the filters $f_i^{(k)}\approx 1,\,
i=1,2,\ldots,k$ and the other $f_i^{(k)}$ monotonically decay to zero for
$i=k+1,k+2,\ldots,n$. If this is the case until $k=k_0$, the $k_0$-step LSQR
has the full regularization and computes a best possible regularized solution
$x^{(k_0)}$. However, if a small Ritz value appears before some
$k\leq k_0$, i.e., $\theta_{k-1}^{(k)}>\sigma_{k_0+1}$ and
$\sigma_{j^*}<\theta_k^{(k)}\leq\sigma_{k_0+1}$ with the smallest integer
$j^*>k_0+1$, then $f_i^{(k)}\in (0,1)$ tends to zero
monotonically for $i=j^*,j^*+1,\ldots,n$; on the other hand,
we have
$$
\prod\limits_{j=1}^k\frac{(\theta_j^{(k)})^2-\sigma_i^2}
{(\theta_j^{(k)})^2}=\frac{(\theta_k^{(k)})^2-\sigma_i^2}
{(\theta_k^{(k)})^2}\prod\limits_{j=1}^{k-1}
\frac{(\theta_j^{(k)})^2-\sigma_i^2}{(\theta_j^{(k)})^2}\leq 0,
\ i=k_0+1,\ldots,j^*-1
$$
since the first factor is non-positive and the second factor is positive.
Then we get $f_i^{(k)}\geq 1,\ i=k_0+1,\ldots,j^*-1$, so that $x^{(k)}$
is deteriorated and LSQR has only the partial regularization.
Hansen \cite[p.146-157]{hansen98} summarizes the known results on $f_i^{(k)}$,
where a bound for $|f_i^{(k)}-1|,\ i=1,2,\ldots,k$
is given in \cite[p.155]{hansen98} but there is no accurate estimate for
the bound. As we will see in Section \ref{compare},
the results to be established in this paper can be used for this purpose,
and, more importantly, we will show that the bound in \cite[p.155]{hansen98}
can be sharpened substantially.

The standard $k$-step Lanczos bidiagonalization method computes the
$k$ Ritz values $\theta_j^{(k)}$, which are used to approximate some of
the singular values of $A$, and is mathematically equivalent to the
symmetric Lanczos method for
the eigenvalue problem of $A^TA$ starting with $q_1=A^Tb/\|A^Tb\|$;
see \cite{bjorck96,bjorck15} or \cite{reichel05,jia03,jia10} for
several variations that are based on standard, harmonic,
refined projection \cite{bai,stewart01,vorst02}
or a combination of them. A general convergence
theory of harmonic and refined harmonic projection methods was lacking
in the books \cite{bai,stewart01,vorst02} and has later been established
in \cite{jia05}. As is known from \cite{bjorck96,meurant,parlett}, for a general
singular value distribution and a general vector $b$, some of the $k$ Ritz
values become good approximations to the
largest and smallest singular values of $A$ as $k$ increases.
If large singular values are well separated but
small singular values are clustered, large Ritz values
converge fast but small Ritz values converge very slowly.

For \eqref{eq1}, we see from \eqref{eqsvd} and
\eqref{picard} that $A^Tb$ contains more information on dominant right
singular vectors than on the ones corresponding to small singular values.
Therefore, $\mathcal{K}_k(A^TA,A^Tb)$ hopefully contains richer
information on the first $k$ right singular vectors $v_i$ than on the
other $n-k$ ones, at least for $k$ small.
Furthermore, note that $A$ has many small singular values clustered at zero.
Due to these two basic facts, all the Ritz values are expected to
approximate the large singular values of $A$ in natural order until some
iteration $k$, at which a small Ritz value shows up
and the regularized solutions then start to be contaminated by the
noise $e$ dramatically after that iteration. These qualitative
arguments are frequently used to analyze
and elaborate the regularizing effects of LSQR and CGLS; see,
e.g., \cite{aster,hansen98,hansen08,hansen10,hansen13} and the references
therein. Clearly, these arguments are not precise and cannot help us
draw any definitive conclusion on the full or partial
regularization of LSQR. For a severely ill-posed example from
seismic tomography, it is reported in \cite{vorst90} that the desired
convergence of the Ritz values actually holds as long as the discrete
Picard condition is satisfied and there is a good separation
among the large singular values of $A$. Unfortunately,
there has been no mathematical justification on these observations.

A complete understanding of the regularization
of LSQR includes accurate solutions of the following basic problems:
How well or accurately does $\mathcal{K}_{k}(A^{T}A, A^{T}b)$
approximate or capture the $k$-dimensional dominant right singular subspace of
$A$? How accurate is the rank $k$ approximation $P_{k+1}B_kQ_k^T$ to $A$?
Can it be a near best rank $k$ approximation to $A$?
How does the noise level $\|e\|$ affects the approximation accuracy of
$\mathcal{K}_{k}(A^{T}A, A^{T}b)$ and $P_{k+1}B_kQ_k^T$ for $k\leq k_0$ and $k>k_0$,
respectively? What sufficient conditions on $\rho$ and $\alpha$ are needed
to guarantee that $P_{k+1}B_kQ_k^T$ is a near best rank $k$ approximation to $A$?
When do the Ritz values $\theta_i^{(k)}$
approximate $\sigma_i,\ i=1,2,\ldots,k$ in natural order? When does at least
a small Ritz value appear, i.e., $\theta_k^{(k)}<\sigma_{k_0+1}$ before some
$k\leq k_0$? We will make a rigorous and detailed analysis on these problems
and some others related closely, present our results, and draw definitive
assertions on the regularization of LSQR for three kinds of ill-posed
problems.

In terms of the canonical angles $\Theta(\mathcal{X},\mathcal{Y})$ between
two subspaces $\mathcal{X}$ and $\mathcal{Y}$ of the same
dimension \cite[p.43]{stewartsun}, we first present the following $\sin\Theta$
theorem, showing how the $k$-dimensional Krylov subspace
$\mathcal{K}_{k}(A^{T}A, A^{T}b)$ captures or approximates the $k$-dimensional
dominant right singular subspace of $A$ for severely ill-posed problems.

\begin{theorem}\label{thm2}
Let the SVD of $A$ be as \eqref{eqsvd}. Assume that \eqref{eq1} is severely
ill-posed with $\sigma_j=\mathcal{O}(\rho^{-j})$ and $\rho>1$, $j=1,2,\ldots,n$,
and the discrete Picard condition \eqref{picard} is satisfied.
Let $\mathcal{V}_k=span\{V_k\}$ be the $k$-dimensional dominant right singular
subspace of $A$ spanned by the columns of $V_k=(v_1,v_2,\ldots,v_k)$
and $\mathcal{V}_k^R=\mathcal{K}_{k}(A^{T}A, A^{T}b)$. Then for
$k=1,2,\ldots,n-1$ we have
\begin{equation}\label{deltabound}
\|\sin\Theta(\mathcal{V}_k,\mathcal{V}_k^R)\|=
\frac{\|\Delta_k\|}{\sqrt{1+\|\Delta_k\|^2}}
\end{equation}
with $\Delta_k \in \mathbb{R}^{(n-k)\times k}$ to be defined by \eqref{defdelta}
and
\begin{equation}\label{k1}
\|\Delta_1\|\leq \frac{\sigma_{2}}{\sigma_1}\frac{|u_2^Tb|}{|u_1^Tb|}
\left(1+\mathcal{O}(\rho^{-2})\right),
\end{equation}
\begin{equation}\label{eqres1}
  \|\Delta_k\|\leq
  \frac{\sigma_{k+1}}{\sigma_k}\frac{|u_{k+1}^Tb|}{|u_{k}^Tb|}
  \left(1+\mathcal{O}(\rho^{-2})\right)
  |L_{k_1}^{(k)}(0)|,\ k=2,3,\ldots,n-1,
\end{equation}
where
\begin{equation}\label{lk}
|L_{k_1}^{(k)}(0)|=\max_{j=1,2,\ldots,k}|L_j^{(k)}(0)|,
\ |L_j^{(k)}(0)|=\prod\limits_{i=1,i\ne j}^k\frac{\sigma_i^2}{|\sigma_j^2-
\sigma_i^2|},\,j=1,2,\ldots,k.
\end{equation}
In particular, we have
\begin{align}
\|\Delta_1\|&\leq\frac{\sigma_2^{2+\beta}}{\sigma_1^{2+\beta}}
\left(1+\mathcal{O}(\rho^{-2})\right), \label{case5}\\
\|\Delta_k\|&\leq\frac{\sigma_{k+1}^{2+\beta}}{\sigma_k^{2+\beta}}
\left(1+\mathcal{O}(\rho^{-2})\right)|L_{k_1}^{(k)}(0)|,\
k=2,3,\ldots,k_0, \label{case1}\\
\|\Delta_k\|&\leq\frac{\sigma_{k+1}}{\sigma_k}\left(1+\mathcal{O}(\rho^{-2})\right)
|L_{k_1}^{(k)}(0)|, \ k=k_0+1,\ldots, n-1.\label{case2}
\end{align}
\end{theorem}

{\em Proof}.
Let $U_n=(u_1,u_2,\ldots,u_n)$ whose columns are the
first $n$ left singular vectors of $A$ defined by \eqref{eqsvd}.
Then the Krylov subspace $\mathcal{K}_{k}(\Sigma^2,
\Sigma U_n^Tb)=span\{DT_k\}$ with
\begin{equation*}\label{}
  D=\diag(\sigma_i u_i^Tb)\in\mathbb{R}^{n\times n},\ \
  T_k=\left(\begin{array}{cccc} 1 &
  \sigma_1^2&\ldots & \sigma_1^{2k-2}\\
1 &\sigma_2^2 &\ldots &\sigma_2^{2k-2} \\
\vdots & \vdots&&\vdots\\
1 &\sigma_n^2 &\ldots &\sigma_n^{2k-2}
\end{array}\right).
\end{equation*}
Partition the diagonal matrix $D$ and the matrix $T_k$ as follows:
\begin{equation*}\label{}
  D=\left(\begin{array}{cc} D_1 & 0 \\ 0 & D_2 \end{array}\right),\ \ \
  T_k=\left(\begin{array}{c} T_{k1} \\ T_{k2} \end{array}\right),
\end{equation*}
where $D_1, T_{k1}\in\mathbb{R}^{k\times k}$. Since $T_{k1}$ is
a Vandermonde matrix with $\sigma_j$ supposed to be distinct for $j=1,2,\ldots,k$,
it is nonsingular. Therefore, from $\mathcal{K}_{k}(A^{T}A, A^{T}b)=span\{VDT_k\}$
we have
\begin{equation}\label{kry}
\mathcal{V}_k^R=\mathcal{K}_{k}(A^{T}A, A^{T}b)=span
  \left\{V\left(\begin{array}{c} D_1T_{k1} \\ D_2T_{k2} \end{array}\right)\right\}
  =span\left\{V\left(\begin{array}{c} I \\ \Delta_k \end{array}\right)\right\},
\end{equation}
where
\begin{equation}\label{defdelta}
\Delta_k=D_2T_{k2}T_{k1}^{-1}D_1^{-1}\in \mathbb{R}^{(n-k)\times k}.
\end{equation}
Write $V=(V_k, V_k^{\perp})$, and define
\begin{equation}\label{zk}
Z_k=V\left(\begin{array}{c} I \\ \Delta_k \end{array}\right)
=V_k+V_k^{\perp}\Delta_k.
\end{equation}
Then $Z_k^TZ_k=I+\Delta_k^T\Delta_k$, and the columns of
$\hat{Z}_k=Z_k(Z_k^TZ_k)^{-\frac{1}{2}}$
form an orthonormal basis of $\mathcal{V}_k^R$. So we get an
orthogonal direct sum decomposition of $\hat{Z}_k$:
\begin{equation}\label{decomp}
\hat{Z}_k=(V_k+V_k^{\perp}\Delta_k)(I+\Delta_k^T\Delta_k)^{-\frac{1}{2}}.
\end{equation}
By definition and \eqref{decomp}, for the matrix 2-norm  we obtain
\begin{align}\label{sindef}
   \|\sin\Theta(\mathcal{V}_k,\mathcal{V}_k^R)\|
   &=\|(V_k^{\perp})^T\hat{Z}_k\|
   =\|\Delta_k(I+\Delta_k^T\Delta_k)^{-\frac{1}{2}}\|
   =\frac{\|\Delta_k\|}{\sqrt{1+\|\Delta_k\|^2}},
\end{align}
which is \eqref{deltabound}.

Next we estimate $\|\Delta_k\|$. For $k=2,3,\ldots,n-1$,
it is easily justified that the $j$-th column of $T_{k1}^{-1}$ consists of
the coefficients of the $j$-th Lagrange polynomial
\begin{equation*}\label{}
  L_j^{(k)}(\lambda)=\prod\limits_{i=1,i\neq j}^k
  \frac{\lambda-\sigma_i^2}{\sigma_j^2-\sigma_i^2}
\end{equation*}
that interpolates the elements of the $j$-th canonical basis vector
$e_j^{(k)}\in \mathbb{R}^{k}$ at the abscissas $\sigma_1^2,\sigma_2^2
\ldots, \sigma_k^2$. Consequently, the $j$-th column of $T_{k2}T_{k1}^{-1}$ is
\begin{equation}\label{tk12i}
  T_{k2}T_{k1}^{-1}e_j^{(k)}=(L_j^{(k)}(\sigma_{k+1}^2),\ldots,L_j^{(k)}
  (\sigma_{n}^2))^T, \ j=1,2,\ldots,k,
\end{equation}
from which we obtain
\begin{equation}\label{tk12}
  T_{k2}T_{k1}^{-1}=\left(\begin{array}{cccc} L_1^{(k)}(\sigma_{k+1}^2)&
  L_2^{(k)}(\sigma_{k+1}^2)&\ldots & L_k^{(k)}(\sigma_{k+1}^2)\\
L_1^{(k)}(\sigma_{k+2}^2)&L_2^{(k)}(\sigma_{k+2}^2) &\ldots &
L_k^{(k)}(\sigma_{k+2}^2) \\
\vdots & \vdots&&\vdots\\
L_1^{(k)}(\sigma_{n}^2)&L_2^{(k)}(\sigma_{n}^2) &\ldots &L_k^{(k)}(\sigma_{n}^2)
\end{array}\right)\in \mathbb{R}^{(n-k)\times k}.
\end{equation}
Since $|L_j^{(k)}(\lambda)|$ is monotonically
decreasing for $0\leq \lambda<\sigma_k^2$, it is bounded by $|L_j^{(k)}(0)|$.
With this property and the definition of $L_{k_1}^{(k)}(0)$
by \eqref{lk}, we get
\begin{align}
|\Delta_k|&=|D_2T_{k2}T_{k1}^{-1}D_1^{-1}| \notag \\
&\leq
\left(\begin{array}{cccc}
\frac{\sigma_{k+1}}{\sigma_1}\left|\frac{u_{k+1}^Tb}
{u_1^Tb}\right||L_{k_1}^{(k)}(0)|, &\frac{\sigma_{k+1}}{\sigma_2}
\left|\frac{u_{k+1}^Tb}
{u_2^Tb}\right||L_{k_1}^{(k)}(0)|, &\ldots&\frac{\sigma_{k+1}}{\sigma_k}
\left|\frac{u_{k+1}^Tb}{u_k^Tb}\right||L_{k_1}^{(k)}(0)| \\
\frac{\sigma_{k+2}}{\sigma_1}\left|\frac{ u_{k+2}^Tb}
{u_1^Tb}\right| |L_{k_1}^{(k)}(0)|, &\frac{\sigma_{k+2}}{\sigma_2}
\left|\frac{u_{k+2}^Tb}
{u_2^Tb}\right| |L_{k_1}^{(k)}(0)|,&
\ldots &\frac{\sigma_{k+2}}{\sigma_k}\left|\frac{u_{k+2}^Tb}
{u_k^Tb}\right| |L_{k_1}^{(k)}(0)| \\
\vdots &\vdots & &\vdots\\
\frac{\sigma_n}{\sigma_1}\left|\frac{u_n^Tb}
{u_1^Tb}\right| |L_{k_1}^{(k)}(0)|, &\frac{\sigma_n}{\sigma_2}\left|\frac{u_n^Tb}
{u_2^Tb}\right| |L_{k_1}^{(k)}(0)|,& \ldots &
\frac{\sigma_n}{\sigma_k}\left|\frac{u_n^Tb}{u_k^Tb}\right| |L_{k_1}^{(k)}(0)|
\end{array}
\right) \notag\\
&= |L_{k_1}^{(k)}(0)||\tilde\Delta_k|, \label{amplify}
\end{align}
where
\begin{equation}
|\tilde\Delta_k|=\left|(\sigma_{k+1} u_{k+1}^T b,\sigma_{k+2}u_{k+2}^Tb,
\ldots,\sigma_n u_n^T b)^T
\left(\frac{1}{\sigma_1 u_1^Tb},\frac{1}{\sigma_2 u_2^Tb},\ldots,
\frac{1}{\sigma_k u_k^Tb}\right)\right|  \label{delta1}
\end{equation}
is a rank one matrix. Therefore, by $\|C\|\leq \||C|\|$
(cf. \cite[p.53]{stewart98}), we get
\begin{align}
\|\Delta_k\| &\leq \||\Delta_k|\|\leq |L_{k_1}^{(k)}(0)|
\left\||\tilde\Delta_k|\right\| \notag\\
&=|L_{k_1}^{(k)}(0)|\left(\sum_{j=k+1}^n\sigma_j^2| u_j^Tb|^2\right)^{1/2}
\left(\sum_{j=1}^k \frac{1}{\sigma_j^2| u_j^Tb|^2}\right)^{1/2}.
\label{delta2}
\end{align}

By the discrete Picard condition \eqref{picard}, \eqref{picard1} and the
description between them,
for the white noise $e$, it is known from \cite[p.70-1]{hansen98} and
\cite[p.41-2]{hansen10} that
$
| u_j^T b|\approx | u_j^T \hat{b}|=\sigma_j^{1+\beta}
$
decrease as $j$ increases up to $k_0$ and then become stabilized as
$
| u_j^T b|\approx | u_j^T e |\approx \eta \approx \frac{\|e\|}{\sqrt{m}},
$
a small constant for $j>k_0$. In order to simplify the
derivation and present our results compactly, in terms of these
assumptions and properties, in later proofs we will use the following precise
equalities and inequalities:
\begin{align}
|u_j^T b|&= |u_j^T\hat{b}|=\sigma_j^{1+\beta},\ j=1,2,\ldots,k_0,\label{ideal3}\\
|u_j^T b|&=|u_j^T e|=\eta, \ j=k_0+1,\ldots,n,\label{ideal2}\\
|u_{j+1}^Tb|&\leq |u_j^Tb |,\ j=1,2,\ldots,n-1.\label{ideal}
\end{align}
From \eqref{ideal} and $\sigma_j=\mathcal{O}(\rho^{-j}),\ j=1,2,\ldots,n$,
for $k=1,2,\ldots,n-1$ we obtain
\begin{align}
\left(\sum_{j=k+1}^n\sigma_j^2| u_j^Tb|^2\right)^{1/2}
&= \sigma_{k+1}| u_{k+1}^Tb| \left(\sum_{j=k+1}^n
\frac{\sigma_j^2| u_j^Tb|^2}{\sigma_{k+1}^2| u_{k+1}^Tb|^2}\right)^{1/2}
\notag\\
&\leq \sigma_{k+1}| u_{k+1}^Tb| \left(\sum_{j=k+1}^n
\frac{\sigma_j^2}{\sigma_{k+1}^2}\right)^{1/2} \notag\\
&=\sigma_{k+1}| u_{k+1}^Tb|\left(1+\sum_{j=k+2}^n\mathcal{O}
(\rho^{2(k-j)+2})\right)^{1/2}
\notag \\
&=\sigma_{k+1}| u_{k+1}^Tb|\left(1+\mathcal{O}\left(\sum_{j=k+2}^n
\rho^{2(k-j)+2}\right)\right)^{1/2}
\notag \\
&=\sigma_{k+1}| u_{k+1}^Tb|\left(1+ \mathcal{O}\left(\frac{\rho^{-2}}
    {1-\rho^{-2}}\left(1-\rho^{-2(n-k-1)}\right)\right)\right)^{1/2}
\notag \\
&=\sigma_{k+1}| u_{k+1}^Tb| \left(1+\mathcal{O}(\rho^{-2})\right)^{1/2}\notag\\
&=\sigma_{k+1}| u_{k+1}^Tb| \left(1+\mathcal{O}(\rho^{-2})\right)
\label{severe1}
\end{align}
with $1+\mathcal{O}(\rho^{-2})$ replaced by one for $k=n-1$. In a similar manner,
for $k=2,3,\ldots,n-1$, from \eqref{ideal} we get
\begin{align*}
\left(\sum_{j=1}^k \frac{1}{\sigma_j^2| u_j^Tb|^2}\right)^{1/2}
&=\frac{1}{\sigma_k | u_k^T b|}\left(\sum_{j=1}^k\frac{\sigma_k^2| u_k^Tb|^2}
{\sigma_j^2| u_j^Tb|^2}\right)^{1/2}
\leq \frac{1}{\sigma_k | u_k^T b|}\left(\sum_{j=1}^k\frac{\sigma_k^2}
{\sigma_j^2}\right)^{1/2} \\
&=\frac{1}{\sigma_k | u_k^T b|}\left(1+\mathcal{O}\left(\sum_{j=1}^{k-1}
\rho^{2(j-k)}\right)\right)^{1/2} \\
&=\frac{1}{\sigma_k | u_k^T b|}\left(1+\mathcal{O}(\rho^{-2})\right).
\end{align*}
From the above and \eqref{delta2}, we finally obtain
$$
\|\Delta_k\|\leq \frac{\sigma_{k+1}}{\sigma_k}\frac{| u_{k+1}^T b|}
{| u_k^T b|}\left(1+\mathcal{O}(\rho^{-2})\right)|L_{k_1}^{(k)}(0)|,
\ k=2,3,\ldots,n-1,
$$
which proves \eqref{eqres1}.

Note that the Langrange polynomials $L_j^{(k)}(\lambda)$ require $k\geq 2$. So,
we need to treat the case $k=1$ independently:
from \eqref{defdelta} and \eqref{ideal}, observe that
$$
T_{k2}=(1,1,\ldots,1)^T,\ D_2T_{k2}=(\sigma_2u_2^Tb,\sigma_3 u_3^Tb,
\ldots,\sigma_n u_n^Tb)^T,\
T_{k1}^{-1}=1,\ D_1^{-1}=\frac{1}{\sigma_1 u_1^Tb}.
$$
Therefore, we have
\begin{equation}\label{deltaexp}
\Delta_1=(\sigma_2u_2^Tb,\sigma_3 u_3^Tb,
\ldots,\sigma_n u_n^Tb)^T\frac{1}{\sigma_1 u_1^Tb},
\end{equation}
from which and \eqref{severe1} for $k=1$ it is direct to get \eqref{k1}.

In terms of the discrete Picard condition \eqref{picard},
\eqref{picard1}, \eqref{ideal3} and \eqref{ideal2}, we have
\begin{equation}\label{ratio1}
\frac{|u_{k+1}^Tb|}{| u_k^T b|}= \frac{|u_{k+1}^T\hat{b}|}{| u_k^T\hat{b}|}
=\frac{\sigma_{k+1}^{1+\beta}}{\sigma_k^{1+\beta}}, \ k\leq k_0
\end{equation}
and
\begin{equation}\label{ratio2}
\frac{|u_{k+1}^Tb|}{| u_k^T b|}= \frac{|u_{k+1}^T e|}{| u_k^T e|}
=1,\  k>k_0.
\end{equation}
Applying them to \eqref{k1} and \eqref{eqres1} establishes \eqref{case5},
\eqref{case1} and \eqref{case2}, respectively.
\qquad\endproof

We next estimate the factor $|L_{k_1}^{(k)}(0)|$ accurately.

\begin{theorem}\label{estlk}
For the severely ill-posed problem and $k=2,3,\ldots,n-1$, we have
\begin{align}
|L_k^{(k)}(0)|&=1+\mathcal{O}(\rho^{-2}), \label{lkkest}\\
|L_j^{(k)}(0)|&=\frac{1+\mathcal{O}(\rho^{-2})}
{\prod\limits_{i=j+1}^k\left(\frac{\sigma_{j}}{\sigma_i}\right)^2}
=\frac{1+\mathcal{O}(\rho^{-2})}{\mathcal{O}(\rho^{(k-j)(k-j+1)})},
\ j=1,2,\ldots,k-1, \label{lj0}\\
|L_{k_1}^{(k)}(0)|&=\max_{j=1,2,\ldots,k}|L_j^{(k)}(0)|
=1+\mathcal{O}(\rho^{-2}). \label{lkk}
\end{align}
\end{theorem}

{\em Proof}.
Exploiting the Taylor series expansion and
$\sigma_i=\mathcal{O}(\rho^{-i})$ for $i=1,2,\ldots,n$,
by definition, for $j=1,2,\ldots,k-1$ we have
\begin{align}
|L_j^{(k)}(0)|&=\prod\limits_{i=1,i\neq j}^k
\left|\frac{\sigma_i^2}{\sigma_i^2-\sigma_j^2}\right|
  =\prod\limits_{i=1}^{j-1}\frac{\sigma_i^2}{\sigma_i^2-\sigma_j^2}
   \cdot\prod\limits_{i=j+1}^{k}\frac{\sigma_i^2}{\sigma_j^2-\sigma_{i}^2}
   \notag\\
& =\prod\limits_{i=1}^{j-1}\frac{1}
{1-\mathcal{O}(\rho^{-2(j-i)})}
\prod\limits_{i=j+1}^{k}\frac{1}
{1-\mathcal{O}(\rho^{-2(i-j)})}\frac{1}
{\prod\limits_{i=j+1}^{k}\mathcal{O}(\rho^{2(i-j)})} \notag\\
&=\frac{\left(1+\sum\limits_{i=1}^j \mathcal{O}(\rho^{-2i})\right)
\left(1+\sum\limits_{i=1}^{k-j+1} \mathcal{O}(\rho^{-2i})\right)}
{\prod\limits_{i=j+1}^{k}\mathcal{O}(\rho^{2(i-j)})} \label{lik}
\end{align}
by absorbing those higher order terms into the two $\mathcal{O}(\cdot)$
in the numerator. For $j=k$, we get
\begin{align*}
|L_k^{(k)}(0)|&=\prod\limits_{i=1}^{k-1}
\left|\frac{\sigma_i^2}{\sigma_i^2-\sigma_{k}^2}\right|
=\prod\limits_{i=1}^{k-1}\frac{1}
{1-\mathcal{O}(\rho^{-2(k-i)})}=
\prod\limits_{i=1}^{k-1}\frac{1}
{1-\mathcal{O}(\rho^{-2i})}\\
&=1+\sum\limits_{i=1}^k \mathcal{O}(\rho^{-2i})
=1+\mathcal{O}\left(\sum\limits_{i=1}^k\rho^{-2i}\right)\\
&=1+ \mathcal{O}\left(\frac{\rho^{-2}}
    {1-\rho^{-2}}(1-\rho^{-2k})\right)
=1+\mathcal{O}(\rho^{-2}),
\end{align*}
which is \eqref{lkkest}.

Note that for the numerator of \eqref{lik} we have
  $$
  1+\sum\limits_{i=1}^j \mathcal{O}(\rho^{-2i})
    =1+ \mathcal{O}\left(\sum\limits_{i=1}^j\rho^{-2i}\right)
    =1+ \mathcal{O}\left(\frac{\rho^{-2}}
    {1-\rho^{-2}}(1-\rho^{-2j})\right),
  $$
  and
  $$
    1+\sum\limits_{i=1}^{k-j+1} \mathcal{O}(\rho^{-2i})
    =1+ \mathcal{O}\left(\sum\limits_{i=1}^{k-j+1}\rho^{-2i}\right)
    =1+ \mathcal{O}\left(\frac{\rho^{-2}}{1-\rho^{-2}}
    (1-\rho^{-2(k-j+1)})\right),
  $$
whose product for any $k$ is
  $$
  1+ \mathcal{O}\left(\frac{2\rho^{-2}}{1-\rho^{-2}}\right)
  +\mathcal{O}\left(\left(\frac{\rho^{-2}}{1-\rho^{-2}}\right)^2\right)=
  1+ \mathcal{O}\left(\frac{2\rho^{-2}}{1-\rho^{-2}}\right)
  =  1+\mathcal{O}(\rho^{-2}).
  $$
On the other hand, note that the denominator of \eqref{lik} is defined by
$$
\prod\limits_{i=j+1}^k\left(\frac{\sigma_{j}}{\sigma_i}\right)^2
=\prod\limits_{i=j+1}^{k}\mathcal{O}(\rho^{2(i-j)})
=\mathcal{O}((\rho\cdot\rho^2\cdots\rho^{k-j})^2)
=\mathcal{O}(\rho^{(k-j)(k-j+1)}),
$$
which, together with the above estimate
for the numerator of \eqref{lik}, proves \eqref{lj0}.
Notice that
$
\prod\limits_{i=j+1}^k\left(\frac{\sigma_j}{\sigma_i}\right)^2
$
is always {\em bigger than one} for $j=1,2,\ldots,k-1$.
Therefore, for any $k$, combining \eqref{lkkest} and \eqref{lj0}
gives \eqref{lkk}.
\qquad\endproof

\begin{remark}\label{severerem}
\eqref{lkkest} and \eqref{lkk} have essentially been shown in {\rm \cite{huangjia}}.
Here we have given a general and complete proof. From \eqref{lkk}, we get
\begin{equation}\label{rhoL}
\left(1+\mathcal{O}(\rho^{-2})\right)|L_{k_1}^{(k)}(0)|= 1+
\mathcal{O}(\rho^{-2}), \ k=2,3,\ldots,n-1,
\end{equation}
so the results in Theorem~\ref{thm2} are simplified as
\begin{align}
\|\Delta_k\|&\leq\frac{\sigma_{k+1}^{2+\beta}}{\sigma_k^{2+\beta}}
\left(1+\mathcal{O}(\rho^{-2})\right),\ k=1,2,\ldots, k_0,
\label{case3}\\
\|\Delta_k\|&\leq\frac{\sigma_{k+1}}{\sigma_k}\left(1+\mathcal{O}(\rho^{-2})\right),
\ k=k_0+1,\ldots,n-1. \label{case4}
\end{align}
\end{remark}

\begin{remark}
\eqref{lj0} illustrates that $|L_j^{(k)}(0)|$ increases fast
with $j$ increasing and the smaller $j$, the smaller $|L_j^{(k)}(0)|$.
\eqref{case3} and \eqref{case4} indicate that $\mathcal{V}_k^R$ captures
$\mathcal{V}_k$ better for $k\leq k_0$ than for $k>k_0$. That is,
after the transition point $k_0$, the noise $e$ starts to
deteriorate $\mathcal{V}_k^R$ and impairs its ability to capture
$\mathcal{V}_k$.
\end{remark}

In what follows we establish accurate estimates for
$\|\sin\Theta(\mathcal{V}_k,\mathcal{V}_k^R)\|$ for
moderately and mildly ill-posed problems.

\begin{theorem}\label{moderate}
Assume that \eqref{eq1} is moderately ill-posed with $\sigma_j=
\zeta j^{-\alpha},\ j=1,2,\ldots,n$,
where $\alpha>\frac{1}{2}$ and $\zeta>0$ is some constant,
and the other assumptions and notation are the same as in Theorem~\ref{thm2}.
Then \eqref{deltabound} holds with
\begin{align}
\|\Delta_1\|&\leq\frac{|u_2^Tb|}{|u_1^Tb|}\sqrt{\frac{1}{2\alpha-1}},\label{k2}\\
  \|\Delta_k\|&\leq\frac{|u_{k+1}^Tb|}{|u_{k}^Tb|}
\sqrt{\frac{k^2}{4\alpha^2-1}+\frac{k}{2\alpha-1}}|L_{k_1}^{(k)}(0)|,\
k=2,3,\ldots,n-1. \label{modera1}
\end{align}
Particularly, we have
\begin{align}
\|\Delta_1\|&\leq \frac{\sigma_2^{1+\beta}}{\sigma_1^{1+\beta}}
\sqrt{\frac{1}{2\alpha-1}},\label{mod1}\\
\|\Delta_k\|&\leq \frac{\sigma_{k+1}^{1+\beta}}{\sigma_k^{1+\beta}}
\sqrt{\frac{k^2}{4\alpha^2-1}+\frac{k}{2\alpha-1}}|L_{k_1}^{(k)}(0)|, \
k=2,3,\ldots, k_0, \label{modera2} \\
\|\Delta_k\|&\leq
\sqrt{\frac{k^2}{4\alpha^2-1}+\frac{k}{2\alpha-1}}|L_{k_1}^{(k)}(0)|,\
k=k_0+1,\ldots, n-1.
\label{modera3}
\end{align}
\end{theorem}

{\em Proof}.
Following the proof of Theorem~\ref{thm2}, we know that $|\Delta_k|\leq
|L_{k_1}^{(k)}(0)||\tilde\Delta_k|$ still holds with
$\tilde{\Delta}_k$ defined by \eqref{delta1}. So we only need to
bound the right-hand side of \eqref{delta2}. For $k=1,2,\ldots, n-1$, from
\eqref{ideal} we get
\begin{align}
\left(\sum_{j=k+1}^n\sigma_j^2| u_j^Tb|^2\right)^{1/2}
&= \sigma_{k+1}| u_{k+1}^Tb| \left(\sum_{j=k+1}^n
\frac{\sigma_j^2| u_j^Tb|^2}{\sigma_{k+1}^2| u_{k+1}^Tb|^2}\right)^{1/2}
\notag\\
&\leq \sigma_{k+1}| u_{k+1}^Tb| \left(\sum_{j=k+1}^n
\frac{\sigma_j^2}{\sigma_{k+1}^2}\right)^{1/2} \notag\\
&= \sigma_{k+1}| u_{k+1}^Tb| \left(\sum_{j=k+1}^n \left(\frac{j}{k+1}
\right)^{-2\alpha}\right)^{1/2} \notag \\
&=\sigma_{k+1}| u_{k+1}^Tb|
\left((k+1)^{2\alpha}\sum_{j=k+1}^n \frac{1}{j^{2\alpha}}\right)^{1/2}
\notag\\
&< \sigma_{k+1}| u_{k+1}^Tb| (k+1)^{\alpha}\left(\int_k^{\infty}
\frac{1}{x^{2\alpha}} dx\right)^{1/2}
\notag \\
&= \sigma_{k+1}| u_{k+1}^Tb|\left(\frac{k+1}{k}\right)^{\alpha}
\sqrt{\frac{k}{2\alpha-1}} \notag\\
&=\sigma_{k+1}| u_{k+1}^Tb|\frac{\sigma_k}{\sigma_{k+1}} \sqrt{\frac{k}
{2\alpha-1}} \notag\\
&=\sigma_k | u_{k+1}^Tb|\sqrt{\frac{k}
{2\alpha-1}}.
\label{modeest}
\end{align}
Since the function $x^{2\alpha}$ with any $\alpha> \frac{1}{2}$
is convex over the interval $[0,1]$, for $k=2,3,\ldots, n-1$, from \eqref{ideal}
we obtain
\begin{align}
\left(\sum_{j=1}^k \frac{1}{\sigma_j^2| u_j^Tb|^2}\right)^{1/2}
&=\frac{1}{\sigma_k | u_k^T b|}\left(\sum_{j=1}^k
\frac{\sigma_k^2| u_k^Tb|^2}
{\sigma_j^2| u_j^Tb|^2}\right)^{1/2}
\leq\frac{1}{\sigma_k | u_k^T b|}\left(\sum_{j=1}^k\frac{\sigma_k^2}
{\sigma_j^2}\right)^2 \notag \\
&=\frac{1}{\sigma_k | u_k^T b|}
\left(\sum_{j=1}^k \left(\frac{j}{k}
\right)^{2\alpha}\right)^{1/2} \notag \\
&=\frac{1}{\sigma_k | u_k^T b|}
\left(k\sum_{j=1}^{k} \frac{1}{k}\left(\frac{j-1}{k}
\right)^{2\alpha}+1\right)^{1/2} \label{sum1} \\
&< \frac{1}{\sigma_k | u_k^T b|} \left(k\int_0^1
x^{2\alpha}dx+1\right)^{1/2}\notag \\
&=\frac{1}{\sigma_k | u_k^T b|} \sqrt{\frac{k}{2\alpha+1}+1}. \label{estimate2}
\end{align}
Substituting the above and \eqref{modeest} into \eqref{delta2} establishes
\eqref{modera1}, from which and \eqref{ratio1}, \eqref{ratio2}
it follows that \eqref{modera2} and \eqref{modera3} hold. For $k=1$,
we still have \eqref{deltaexp}, from which and \eqref{modeest}
we obtain \eqref{k2}. From \eqref{ratio1} and \eqref{k2}
we get \eqref{mod1}.
\qquad\endproof

\begin{remark}
For a purely technical reason and for the sake of precise presentation,
we have used the simplified singular value
model $\sigma_j=\zeta j^{-\alpha}$ to replace the general form
$\sigma_j=\mathcal{O}(j^{-\alpha})$, where the constant in each
$\mathcal{O}(\cdot)$ is implicit. This model,
though simple, reflects the essence of moderately and mildly ill-posed
problems and avoids some troublesome derivations and
non-transparent formulations.
\end{remark}

\begin{remark}
In the spirit of the proof of Theorem~\ref{estlk}, exploiting
the first order Taylor expansion, we have an
estimate
\begin{align}
  |L_{k_1}^{(k)}(0)|& \approx |L_k^{(k)}(0)|=
\prod\limits_{i=1}^{k-1}\frac{\sigma_i^2}{\sigma_i^2-\sigma_k^2}
=\prod\limits_{i=1}^{k-1}\frac{1}
{1-(\frac{i}{k})^{2\alpha}} \notag\\
&\approx 1+\sum\limits_{i=1}^{k-1}\left(\frac{i}{k}\right)^{2\alpha}
=1+k\sum\limits_{i=1}^k\frac{1}{k}\left(\frac{i-1}{k}\right)^{2\alpha}
\notag\\
&< 1+k\int_0^1 x^{2\alpha}dx=1+\frac{k}{2\alpha+1},
\label{lkkmoderate}
\end{align}
where the right-hand side of \eqref{lkkmoderate}
increases linearly with respect to $k$.
\end{remark}

\begin{remark}\label{rem3.5}
\eqref{k2} and \eqref{mod1} indicate that $\|\Delta_1\|<1$ is
guaranteed for moderately ill-posed problems with $\alpha>1$.
One might worry that the upper bounds \eqref{k2} and
\eqref{modera1} overestimate $\|\Delta_k\|$ and thus
$\|\sin\Theta(\mathcal{V}_k,\mathcal{V}_k^R)\|$ considerably
because, in the proof, we have bounded the opaque $\sum$ in \eqref{sum1}
from above by the compact integral \eqref{estimate2} nearest to it, which
can be overestimates for $k=1,2,\ldots,k_0$. It is not the case provided that
$k_0$ is not very small. In fact, since $\alpha>\frac{1}{2}$,
we can bound \eqref{sum1} from below by the integral nearest to it:
\begin{align*}
\frac{1}{\sigma_k | u_k^T b|}\left(k\sum_{j=1}^{k} \frac{1}{k}\left(\frac{j-1}{k}
\right)^{2\alpha}+1\right)^{1/2}
&> \frac{1}{\sigma_k | u_k^T b|} \left(k\int_0^{\frac{k-1}{k}}
x^{2\alpha}dx+1\right)^{1/2}\notag \\
&=\frac{1}{\sigma_k | u_k^T b|} \sqrt{\frac{k}{2\alpha+1}
\left(\frac{k-1}{k}\right)^{2\alpha+1}+1},
\end{align*}
which is near to \eqref{estimate2} once  $k\leq k_0$ is not very small.
The smaller $\alpha$, the smaller the difference between
the upper and lower bounds, i.e., the sharper \eqref{estimate2}.
\end{remark}

\begin{remark}
It is easily seen from \eqref{deltabound} that
$\|\sin\Theta(\mathcal{V}_k,\mathcal{V}_k^R)\|$ increases monotonically
with respect to $\|\Delta_k\|$. For $\|\Delta_k\|$ reasonably small and
$\|\Delta_k\|$ large we have
$$
\|\sin\Theta(\mathcal{V}_k,\mathcal{V}_k^R)\|\approx \|\Delta_k\|
\ \mbox{ and } \
\|\sin\Theta(\mathcal{V}_k,\mathcal{V}_k^R)\|\approx 1,
$$
respectively. From \eqref{picard} and \eqref{picard1}, we obtain
$k_0=\lfloor\eta^{-\frac{1}{\alpha(1+\beta)}}\rfloor-1$,
where $\lfloor\cdot\rfloor$ is the Gaussian function. For the white noise
$e$, we have $\eta\approx \frac{\|e\|}{\sqrt{m}}$.
As a result, for moderately ill-posed problems with $\alpha>1$,
$k_0$ is typically small and at most modest for a practical noise $e$, whose
relative size $\frac{\|e\|}{\|\hat{b}\|}$ typically ranges from
$10^{-4}$ to $10^{-2}$. This means that for a moderately ill-posed problem
$\|\Delta_k\|$ is at most modest and cannot be large, so that
$\|\sin\Theta(\mathcal{V}_k,\mathcal{V}_k^R)\|<1$ fairly.
\end{remark}

\begin{remark}
For severely ill-posed problems, since all the $\frac{\sigma_{k+1}}{\sigma_k}
\sim \rho^{-1}$, a constant, \eqref{case3} and \eqref{case4}
indicate that $\|\sin\Theta(\mathcal{V}_k,\mathcal{V}_k^R)\|$
is essentially unchanged for $k=1,2,\ldots,k_0$ and $k=k_0+1,\ldots,n-1$,
respectively, that is, $\mathcal{V}_k^R$ captures $\mathcal{V}_k$ with almost
the same accuracy for $k\leq k_0$ and $k>k_0$, respectively. However,
the situation is different for moderately ill-posed
problems. For them,  $\frac{\sigma_{k+1}}{\sigma_k}=
\left(\frac{k}{k+1}\right)^{\alpha}$ increases slowly as $k$ increases, and
the factor $\sqrt{\frac{k^2}{4\alpha^2-1}+\frac{k}{2\alpha-1}}
|L_{k_1}^{(k)}(0)|$ increases as $k$ grows. Therefore, \eqref{modera2} and
\eqref{modera3} illustrate that $\|\sin\Theta(\mathcal{V}_k,\mathcal{V}_k^R)\|$
increases slowly with $k\leq k_0$ and $k> k_0$, respectively.
This means that $\mathcal{V}_k^R$ may not capture
$\mathcal{V}_k$ so well as it does for severely ill-posed problems as $k$
increases. In particular, starting with some $k>k_0$,
$\|\sin\Theta(\mathcal{V}_k,\mathcal{V}_k^R)\|$ starts to approach one,
which indicates that, for $k$ big, $\mathcal{V}_k^R$ will contain substantial
information on the right singular vectors corresponding to the $n-k$ small
singular values of $A$.
\end{remark}

\begin{remark}\label{mildrem}
For mildly ill-posed problems with $\frac{1}{2}<
\alpha\leq 1$, there are some distinctive features. Note
from \eqref{picard} and \eqref{picard1} that
$k_0$ is now considerably bigger than
that for a severely or moderately ill-posed problem
with the same noise level $\|e\|$ and $\beta$. As a result,
firstly, for $\alpha\leq 1$ and the same $k$, the factor $\frac{\sigma_{k+1}}
{\sigma_{k}}=\left(\frac{k}{k+1}\right)^{\alpha}$ is
bigger than that for the moderately ill-posed problem;
secondly, $\sqrt{\frac{k^2}{4\alpha^2-1}+\frac{k}{2\alpha-1}} \sim k$
if $\alpha\approx 1$ and is much bigger than $k$ and can be
arbitrarily large if $\alpha\approx \frac{1}{2}$; thirdly,
since $\frac{1}{2}<\alpha\leq 1$, for $k\geq 3$ that ensures
$\frac{2\alpha+1}{k}\leq 1$, we have
\begin{align}
 |L_{k_1}^{(k)}(0)|&\geq |L_{k}^{(k)}(0)|=
\prod\limits_{i=1}^{k-1}\frac{\sigma_i^2}{\sigma_i^2-\sigma_{k}^2}
=\prod\limits_{i=1}^{k-1}\frac{1}
{1-(\frac{i}{k})^{2\alpha}} \notag\\
&> 1+\sum\limits_{i=1}^{k-1}\left(\frac{i}{k}\right)^{2\alpha}
>1+k\int_0^{\frac{k-1}{k}} x^{2\alpha}dx \notag\\
&=1+\frac{k\left(\frac{k-1}{k}\right)^{2\alpha+1}}{2\alpha+1}
\approx 1+\frac{k}{2\alpha+1}\left(1-\frac{2\alpha+1}{k}\right)=
\frac{k}{2\alpha+1},\label{lowerbound}
\end{align}
which also holds for moderately ill-posed problems and
is bigger than one considerably for $\frac{1}{2}<\alpha\leq 1$ as
$k$ increases up to $k_0$.
Our accurate bound \eqref{modera2} thus becomes increasingly large as $k$
increases up to $k_0$ for mildly ill-posed
problems, causing that $\|\Delta_k\|$ is large and
$\|\sin\Theta(\mathcal{V}_k,\mathcal{V}_k^R)\|\approx 1$ starting with
some $k\leq k_0$. Consequently, $\mathcal{V}_{k_0}^R$ cannot
effectively capture the $k_0$ dominant right singular vectors and
contains substantial information on the right singular vectors
corresponding to the $n-k_0$ small singular values.
\end{remark}

\begin{remark}
In \cite[Thm 2.1]{huangjia}, the authors derived some
bounds for $\|\Delta_k\|$ and
$\|\sin\Theta(\mathcal{V}_k,\mathcal{V}_k^R)\|$. There, without
realizing the crucial fact that $|\Delta_k|$ can be effectively
bounded by a rank one matrix and the key point that
$D_2T_{k2}T_{k1}^{-1}D_1^{-1}$ must be treated as a whole other
than separately, by \eqref{defdelta} the authors made use of
$$
\|\Delta_k\|\leq \|\Delta\|_F=\left\|D_2T_{k2}T_{k1}^{-1}D_1^{-1}\right\|_F
   \leq \|D_2\|\left\|T_{k2}T_{k1}^{-1}\right\|_F\left\|D_1^{-1}\right\|
$$
and
$\left\|T_{k2}T_{k1}^{-1}\right\|_F\leq |L_{k_1}^{(k)}(0)|\sqrt{k(n-k)}$
(cf. \eqref{tk12}) to obtain bounds for $\|\Delta_k\|$ and
$\|\sin\Theta(\mathcal{V}_k,\mathcal{V}_k^R)\|$. These bounds are too
pessimistic because of the appearance of the fatal factor $\sqrt{k(n-k)}$,
which ranges from
$\sqrt{2(n-2)}$ to $\frac{n}{2}$ for $k=2,3,\ldots,n-1$, too large
amplification for $n$ large. In contrast, our new estimates, which hold for
both $\|\Delta_k\|$ and $\|\Delta_k\|_F$, are much more accurate and
$\sqrt{k(n-k)}$ has been removed.
\end{remark}

Before proceeding, we tentatively investigate how
$\|\sin\Theta(\mathcal{V}_k,\mathcal{V}_k^R)\|$ affects the smallest Ritz
value $\theta_k^{(k)}$. This problem is of central importance for
understanding the regularizing effects of LSQR. We aim to lead the reader
to a first manifestation that (i) we may have $\theta_k^{(k)}>\sigma_{k+1}$
when $\|\sin\Theta(\mathcal{V}_k,\mathcal{V}_k^R)\|<1$ fairly, that is,
no small Ritz value may appear provided that
$\mathcal{V}_k^R$ captures $\mathcal{V}_k$ with only {\em some}
other than high accuracy, and (ii) we must have
$\theta_k^{(k)}\leq\sigma_{k+1}$, that is, $\theta_k^{(k)}$ cannot approximate
$\sigma_k$ in natural order£¬ meaning that $\theta_k^{(k)}\leq\sigma_{k_0+1}$
no later than iteration $k_0$, once
$\|\sin\Theta(\mathcal{V}_k,\mathcal{V}_k^R)\|$ is sufficiently close to one.

\begin{theorem}\label{initial}
Let $\|\sin\Theta(\mathcal{V}_k,\mathcal{V}_k^R)\|^2=1-\varepsilon_k^2$ with
$0< \varepsilon_k< 1$, $k=1,2,\ldots,n-1$, and let the unit-length
$\tilde{q}_k\in\mathcal{V}_k^R$
be a vector that has the smallest acute angle with $span\{V_k^{\perp}\}$, i.e.,
the closest to $span\{V_k^{\perp}\}$, where $V_k^{\perp}$ is the matrix consisting
of the last $n-k$ columns of $V$ defined by \eqref{eqsvd}. Then it holds that
\begin{equation}\label{rqi}
\varepsilon_k^2\sigma_k^2+
(1-\varepsilon_k^2)\sigma_n^2< \tilde{q}_k^TA^TA\tilde{q}_k<
\varepsilon_k^2\sigma_{k+1}^2+
(1-\varepsilon_k^2)\sigma_1^2.
\end{equation}
If $\varepsilon_k\geq \frac{\sigma_{k+1}}{\sigma_k}$,
then
\begin{equation}
\sqrt{\tilde{q}_k^TA^TA\tilde{q}_k}>\sigma_{k+1};
\label{est1}
\end{equation}
if $\varepsilon_k^2\leq\frac{\delta}
{(\frac{\sigma_1}{\sigma_{k+1}})^2-1}$ for a given arbitrarily small
$\delta>0$, then
\begin{equation}\label{thetasigma}
\theta_k^{(k)}<(1+\delta)^{1/2}\sigma_{k+1}.
\end{equation}
\end{theorem}

{\em Proof}.
Since the columns of $Q_k$ generated by Lanczos bidiagonalization form an
orthonormal basis of $\mathcal{V}_k^R$, by definition and the assumption on
$\tilde{q}_k$ we have
\begin{align}
\|\sin\Theta(\mathcal{V}_k,\mathcal{V}_k^R)\|&=\|(V_k^{\perp})^TQ_k\|
=\|V_k^{\perp}(V_k^{\perp})^TQ_k\| \notag\\
&=\max_{\|c\|=1}\|V_k^{\perp}(V_k^{\perp})^TQ_kc\|
=\|V_k^{\perp}(V_k^{\perp})^T Q_kc_k\| \notag\\
&=\|V_k^{\perp}(V_k^{\perp})^T\tilde{q}_k\|=\|(V_k^{\perp})^T\tilde{q}_k\|
=\sqrt{1-\varepsilon_k^2}
\label{qktilde}
\end{align}
with $\tilde{q}_k=Q_kc_k\in\mathcal{V}_k^R$ and $\|c_k\|=1$.
Since $\mathcal{V}_k$ is the orthogonal complement of $span\{V_k^{\perp}\}$,
by definition we know that $\tilde{q}_k\in \mathcal{V}_k^R$ has the largest acute
angle with $\mathcal{V}_k$, that is, it is the vector in $\mathcal{V}_k^R$
that contains the least information on $\mathcal{V}_k$.

Expand $\tilde{q}_k$ as the following orthogonal direct sum decomposition:
\begin{equation}\label{decompqk}
\tilde{q}_k=V_k^{\perp}(V_k^{\perp})^T\tilde{q}_k+V_kV_k^T\tilde{q}_k.
\end{equation}
Then from $\|\tilde{q}_k\|=1$ and \eqref{qktilde} we obtain
\begin{align}\label{angle2}
\|V_k^T\tilde{q}_k\|&=\|V_kV_k^T\tilde{q}_k\|=
\sqrt{1-\|V_k^{\perp}(V_k^{\perp})^T\tilde{q}_k\|^2}=\sqrt{1-(1-\varepsilon_k^2)}
=\varepsilon_k.
\end{align}
From \eqref{decompqk}, we next bound the Rayleigh quotient of $\tilde{q}_k$
with respect to $A^TA$ from below. By the SVD \eqref{eqsvd} of $A$ and
$V=(V_k,V_k^{\perp})$, we partition
$$
\Sigma=\left(\begin{array}{cc}
\Sigma_k &\\
&\Sigma_k^{\perp}
\end{array}
\right),
$$
where $\Sigma_k=\diag(\sigma_1,\sigma_2,\ldots,\sigma_k)$ and
$\Sigma_k^{\perp}=\diag(\sigma_{k+1},\sigma_{k+2},\ldots,\sigma_n)$.
Making use of $A^TAV_k=V_k\Sigma_k^2$ and $A^TAV_k^{\perp}=
V_k^{\perp}(\Sigma_k^{\perp})^2$ as well as $V_k^TV_k^{\perp}=\mathbf{0}$,
we obtain
\begin{align}
\tilde{q}_k^TA^TA\tilde{q}_k&=\left(V_k^{\perp}(V_k^{\perp})^T\tilde{q}_k+V_kV_k^T
\tilde{q}_k\right)^TA^TA \left(V_k^{\perp}(V_k^{\perp})^T\tilde{q}_k+
V_kV_k^T\tilde{q}_k\right) \notag\\
&=\left(\tilde{q}_k^TV_k^{\perp}(V_k^{\perp})^T+\tilde{q}_k^TV_kV_k^T\right)
\left(V_k^{\perp}(\Sigma_k^{\perp})^2(V_k^{\perp})^T\tilde{q}_k+V_k\Sigma_k^2V_k^T
\tilde{q}_k\right) \notag\\
&=\tilde{q}_k^TV_k^{\perp}(\Sigma_k^{\perp})^2(V_k^{\perp})^T\tilde{q}_k
+\tilde{q}_k^TV_k\Sigma_k^2V_k^T\tilde{q}_k. \label{expansion}
\end{align}
Observe that it is impossible for $(V_k^{\perp})^T\tilde{q}_k$ and
$V_k^T\tilde{q}_k$ to be the eigenvectors of $(\Sigma_k^{\perp})^2$
and $\Sigma_k^2$ associated with their respective smallest eigenvalues
$\sigma_n^2$ and $\sigma_k^2$ simultaneously, which are
the $(n-k)$-th canonical vector $e_{n-k}$ of $\mathbb{R}^{n-k}$ and
the $k$-th canonical vector $e_k$ of $\mathbb{R}^{k}$, respectively;
otherwise, we have $\tilde{q}_k=v_n$ and
$\tilde{q}_k=v_k$ simultaneously, which are impossible as $k<n$. Therefore,
from \eqref{expansion} and \eqref{qktilde}, \eqref{angle2}
we obtain the strict inequality
\begin{align*}
\tilde{q}_k^TA^TA\tilde{q}_k&> \|(V_k^{\perp})^T\tilde{q}_k\|^2
\sigma_n^2+\|V_k^T\tilde{q}_k\|^2\sigma_k^2
=(1-\varepsilon_k^2)\sigma_n^2+\varepsilon_k^2 \sigma_k^2,
\end{align*}
from which it follows that the lower bound of \eqref{rqi} holds. Similarly,
from \eqref{expansion}  and \eqref{qktilde}, \eqref{angle2}
we obtain the upper bound of \eqref{rqi}:
$$
\tilde{q}_k^TA^TA\tilde{q}_k <\|(V_k^{\perp})^T\tilde{q}_k\|^2
\|(\Sigma_k^{\perp})^2\|+\|V_k^T\tilde{q}_k\|^2\|\Sigma_k^2\|
=(1-\varepsilon_k^2)\sigma_{k+1}^2+\varepsilon_k^2 \sigma_1^2.
$$

From the lower bound of \eqref{rqi}, we see that if
$\varepsilon_k$ satisfies $\varepsilon_k^2 \sigma_k^2\geq \sigma_{k+1}^2$,
i.e., $\varepsilon_k\geq \frac{\sigma_{k+1}}{\sigma_k}$,
then $\sqrt{\tilde{q}_k^TA^TA\tilde{q}_k}>\sigma_{k+1}$, i.e.,
\eqref{est1} holds.

From \eqref{Bk}, we obtain $B_k^TB_k=Q_k^TA^TAQ_k$.
Note that $(\theta_k^{(k)})^2$ is the smallest eigenvalue
of the symmetric positive definite matrix $B_k^TB_k$.
Therefore, we have
\begin{equation}\label{rqi2}
(\theta_k^{(k)})^2=\min_{\|c\|=1} c^TQ_k^TA^TAQ_kc=
\min_{q\in \mathcal{V}_k^R,\ \|q\|=1} q^TA^TAq
=\hat{q}_k^TA^TA\hat{q}_k,
\end{equation}
where $\hat{q}_k$ is, in fact, the Ritz vector of $A^TA$ from
$\mathcal{V}_k^R$ corresponding to the smallest Ritz value $(\theta_k^{(k)})^2$.
Therefore, for $\tilde{q}_k$ defined in Theorem~\ref{initial} we have
$$
\theta_k^{(k)}\leq \sqrt{\tilde{q}_k^TA^TA\tilde{q}_k},
$$
from which it follows from \eqref{rqi} that
$(\theta_k^{(k)})^2<(1-\varepsilon_k^2)\sigma_{k+1}^2+\varepsilon_k^2 \sigma_1^2$.
As a result, for any $\delta>0$, we can choose $\varepsilon_k\geq 0$ such that
$$
(\theta_k^{(k)})^2<(1-\varepsilon_k^2)\sigma_{k+1}^2+\varepsilon_k^2 \sigma_1^2
\leq (1+\delta)\sigma_{k+1}^2,
$$
i.e., \eqref{thetasigma} holds,
solving which for $\varepsilon_k^2$ gives $\varepsilon_k^2\leq\frac{\delta}
{(\frac{\sigma_1}{\sigma_{k+1}})^2-1}$.
\qquad\endproof

\begin{remark}
We analyze $\theta_k^{(k)}$ for $\varepsilon_k\geq \frac{\sigma_{k+1}}{\sigma_k}$.
A key observation and interpretation is that, in the sense of $\min$
in \eqref{rqi2}, $\hat{q}_k\in \mathcal{V}_k^R$ is the optimal vector that
extracts the least information from $\mathcal{V}_k$ and the richest information
from $span\{V_k^{\perp}\}$.
From Theorem~\ref{initial}, since $\mathcal{V}_k$ is the orthogonal complement
of $span\{V_k^{\perp}\}$, we know that $\tilde{q}_k\in \mathcal{V}_k^R$
has the largest acute angle with $\mathcal{V}_k$, that is,
it contains the least information from $\mathcal{V}_k$ and the richest information
from $span\{V_k^{\perp}\}$. Therefore, $\hat{q}_k$ and $\tilde{q}_k$ have a similar
optimality, so that we have
\begin{equation}\label{approxeq}
\theta_k^{(k)}\approx \sqrt{\tilde{q}_k^TA^TA\tilde{q}_k}.
\end{equation}
Combining this estimate with \eqref{est1}, we may have
$\theta_k^{(k)}>\sigma_{k+1}$ when
$\varepsilon_k\geq \frac{\sigma_{k+1}}{\sigma_k}$.
\end{remark}

\begin{remark}
We inspect the condition $\varepsilon_k\geq\frac{\sigma_{k+1}}
{\sigma_k}$ for \eqref{est1} and get insight into whether or not the true
$\varepsilon_k$ resulting from three kinds of ill-posed problems satisfies
it. For severely ill-posed problems, the lower bound $\frac{\sigma_{k+1}}
{\sigma_k}$ is basically the constant $\rho^{-1}$; for moderately ill-posed
problems with $\alpha>1$, the bound increases with increasing
$k\leq k_0$, and it cannot be close to one provided that $\alpha>1$ suitably
or $k_0$ not big; for mildly ill-posed problems with $\alpha<1$, the bound
increases faster than it does for moderately ill-posed problems, and
it may well approach one for $k\leq k_0$. Therefore,
the condition for \eqref{est1} requires that
$\|\sin\Theta(\mathcal{V}_k,\mathcal{V}_k^R)\|$
be not close to one for severely and moderately ill-posed problems,
but $\|\sin\Theta(\mathcal{V}_k,\mathcal{V}_k^R)\|$
must be close to zero for mildly ill-posed problems.
In view of \eqref{deltabound} and
$\|\sin\Theta(\mathcal{V}_k,\mathcal{V}_k^R)\|^2=1-\varepsilon_k^2$,
we have $\|\Delta_k\|^2=\frac{1-\varepsilon_k^2}{\varepsilon_k^2}$.
Thus, the condition $\varepsilon_k\geq\frac{\sigma_{k+1}}
{\sigma_k}$ for \eqref{est1} amounts to requiring
that $\|\Delta_k\|$ be at most modest and cannot
be large for severely and moderately ill-posed problems but
it must be fairly small for mildly ill-posed problems. Unfortunately,
Theorems~\ref{thm2}--\ref{moderate} and the remarks followed
indicate that $\|\Delta_k\|$ increases with $k$ increasing and is generally
large for a mildly ill-posed problem, while it
increases slowly with $k\leq k_0$ for a moderately
ill-posed problem with $\alpha>1$ suitably, and  by \eqref{case3}
it is approximately a constant $\rho^{-(2+\beta)}$, which is smaller
than one considerably for a
severely ill-posed problem with $\rho>1$ not close to one.
Consequently, for mildly ill-posed problems,
because the actual $\|\Delta_k\|$ can hardly be small and is generally
large, the true $\varepsilon_k$ is small and may well be
close to one, so that the condition $\varepsilon_k\geq\frac{\sigma_{k+1}}{\sigma_k}$
generally fails to meet as $k$ increases, while it is satisfied for severely
or moderately ill-posed problems with $\rho>1$ or $\alpha>1$ suitably.
\end{remark}

\begin{remark}\label{appear}
\eqref{thetasigma} shows that there is
at least one Ritz value $\theta_k^{(k)}\leq\sigma_{k+1}$ when
$\|\sin\Theta(\mathcal{V}_k,\mathcal{V}_k^R)\|$ is sufficiently close
to one since we can choose $\delta$ small enough such that
$(1+\delta)^{1/2}\sigma_{k+1}$ is close to $\sigma_{k+1}$ arbitrarily.
As we have shown, $\|\sin\Theta(\mathcal{V}_k,\mathcal{V}_k^R)\|$ cannot be
close to one for severely or moderately ill-posed problems with $\rho>1$ or
$\alpha>1$ suitably, but it is generally so for mildly ill-posed problems.
This means that for some $k\leq k_0$ it is very possible to have
$\theta_k^{(k)}\leq\sigma_{k+1}$ for mildly ill-posed problems.
\end{remark}

We must be aware that our above analysis on
$\theta_k^{(k)}>\sigma_{k+1}$ is not rigorous because
we cannot quantify {\em how small}
$\sqrt{\tilde{q}_k^TA^TA\tilde{q}_k}-\theta_k^{(k)}$ is.
From $\theta_k^{(k)}\leq\sqrt{\tilde{q}_k^TA^TA\tilde{q}_k}$, it is apparent
that the condition $\varepsilon_k\geq\frac{\sigma_{k+1}}{\sigma_k}$
may not be sufficient for $\theta_k^{(k)}>\sigma_{k+1}$ though it is so for
$\sqrt{\tilde{q}_k^TA^TA\tilde{q}_k}>\sigma_{k+1}$.
We delay our detailed and rigorous analysis
to Section \ref{rankapp}, where we present a number of deep-going and accurate
results on the key problems stated in the last second paragraph before
Theorem~\ref{thm2}, including the precise behavior of $\theta_k^{(k)}$.
One of the results will be on the sufficient conditions
for $\theta_k^{(k)}>\sigma_{k+1}$, which are satisfied when
certain deterministic and mild restrictions on $\rho$ or $\alpha$
are imposed for severely or moderately ill-posed problems.
However, we will see that $\alpha<1$ for mildly ill-posed problems
never meets the sufficient conditions to be presented there.

Theorems~\ref{thm2}--\ref{moderate} establish necessary
background for answering the fundamental
concern by Bj\"{o}rck and Eld\'{e}n, and their proof approaches
also provide key ingredients for some of the later results.
We next present the following results,
which will play a central role in our later analysis.

\begin{theorem}\label{thm3}
Assume that the dicrete Picard condition \eqref{picard} is satisfied,
let $\Delta_k\in \mathbb{R}^{(n-k)\times k}$ be defined as
\eqref{defdelta} and $L_j^{(k)}(0)$ and $L_{k_1}^{(k)}(0)$
be defined as \eqref{lk}, and write
$\Delta_k=(\delta_1,\delta_2,\ldots,\delta_k)$.  Then for severely
ill-posed problems and $k=1,2,\ldots,n-1$ we have
\begin{align}
\|\delta_j\|&\leq \frac{\sigma_{k+1}}{\sigma_j}\frac{| u_{k+1}^Tb|}
{| u_j^T b|}\left(1+\mathcal{O}(\rho^{-2})\right)| L_j^{(k)}(0)|,
\ k>1,\ j=1,2,\ldots,k, \label{columndelta} \\
\|\delta_1\|& \leq \frac{\sigma_{2}}{\sigma_1}\frac{| u_2^Tb|}
{| u_1^T b|}\left(1+\mathcal{O}(\rho^{-2})\right)|,\ k=1
\label{columndelta1}
\end{align}
and
\begin{equation} \label{prodnorm}
\|\Sigma_k\Delta_k^T\|\leq \left\{\begin{array}{ll}
\sigma_{k+1}\frac{| u_{k+1}^Tb|}{| u_k^T b|}
\left(1+\mathcal{O}(\rho^{-2})\right)
& \mbox{ for } 1\leq k\leq k_0,\\
\sigma_{k+1}\sqrt{k-k_0+1}\left(1+\mathcal{O}(\rho^{-2})\right)
& \mbox{ for } k_0<k\leq n-1;
\end{array}
\right.
\end{equation}
for moderately or mild ill-posed problems with the singular values
$\sigma_j=\zeta j^{-\alpha}$ and $\zeta$ a positive constant we have
\begin{align}
\|\delta_j\|&\leq \frac{\sigma_k}{\sigma_j}\frac{| u_{k+1}^Tb|}
{| u_j^T b|} \sqrt{\frac{k}{2\alpha-1}}|L_j^{(k)}(0)|,
\ k>1,\ j=1,2,\ldots,k, \label{columndelta2} \\
\|\delta_1\|&\leq \frac{| u_2^Tb|}
{| u_1^T b|} \sqrt{\frac{1}{2\alpha-1}},\ k=1 \label{columnnorm}
\end{align}
and
\begin{equation}\label{prodnorm2}
\|\Sigma_k\Delta_k^T\|\leq \left\{\begin{array}{ll}
\sigma_1\frac{| u_2^Tb|}{| u_1^T b|}\sqrt{\frac{1}{2\alpha-1}} &
 \mbox{ for } k=1,\\
\sigma_k\frac{| u_{k+1}^Tb|}{| u_k^T b|}\sqrt{\frac{k^2}{4\alpha^2-1}+
\frac{k}{2\alpha-1}}
|L_{k_1}^{(k)}(0)|& \mbox{ for } 1<k\leq k_0, \\
\sigma_k \sqrt{\frac{k k_0}{4\alpha^2-1}+
\frac{k(k-k_0+1)}{2\alpha-1}}|L_{k_1}^{(k)}(0)|
& \mbox{ for } k_0<k\leq n-1.
\end{array}
\right.
\end{equation}
\end{theorem}

{\em Proof}.
From \eqref{defdelta} and \eqref{amplify}, for $j=1,2,\ldots,k$ and $k>1$
we have
\begin{equation}\label{deltaj}
\|\delta_j\|^2\leq |L_j^{(k)}(0)|^2 \sum_{i=k+1}^n\frac{\sigma_{i}^2}
{\sigma_j^2}\frac{| u_i^T b|^2}{| u_j^T b|^2}
\end{equation}
and from \eqref{deltaexp}, for $k=1$ we have
\begin{equation}\label{deltas}
\|\delta_1\|^2=\sum_{i=2}^n\frac{\sigma_{i}^2}
{\sigma_1^2}\frac{| u_i^T b|^2}{| u_1^T b|^2}.
\end{equation}
For severely ill-posed problems, $k=1,2,\ldots,n-1$ and
$j=1,2,\ldots,k$, from \eqref{severe1} we obtain
\begin{align*}
\sum_{i=k+1}^n\frac{\sigma_{i}^2}
{\sigma_j^2}\frac{| u_i^T b|^2}{| u_j^T b|^2}&=
\frac{1}{\sigma_j^2| u_j^Tb|^2}
\sum_{i=k+1}^n\sigma_{i}^2 |u_i^T b|^2 \\
&\leq \frac{\sigma_{k+1}^2}
{\sigma_j^2}\frac{| u_{k+1}^T b|^2}{| u_j^T b|^2}
\left(1+\mathcal{O}(\rho^{-2})\right).
\end{align*}
For moderately or mildly ill-posed problems, $k=1,2,\ldots,n-1$ and
$j=1,2,\ldots,k$, from \eqref{modeest} we obtain
\begin{align*}
\sum_{i=k+1}^n\frac{\sigma_{i}^2}
{\sigma_j^2}\frac{| u_i^T b|^2}{| u_j^T b|^2}
&=\frac{1}{\sigma_j^2| u_j^Tb|^2}
\sum_{i=k+1}^n\sigma_{i}^2| u_i^T b|^2\\
&\leq \frac{\sigma_k^2}
{\sigma_j^2}\frac{| u_{k+1}^T b|^2}{| u_j^T b|^2}
\frac{k}{2\alpha-1}.
\end{align*}
Combining the above with \eqref{deltaj}, \eqref{lkk} and \eqref{rhoL}, we
obtain \eqref{columndelta}, while \eqref{columndelta2} follows
from the above and \eqref{deltaj} directly. For $k=1$, from \eqref{deltas}
and the above we get \eqref{columndelta1} and \eqref{columnnorm}, respectively.

By \eqref{delta1}, for $k>1$ we have
$$
|\Delta_k\Sigma_k|\leq |L_{k_1}^{(k)}(0)|\left|(\sigma_{k+1} u_{k+1}^T b,
\sigma_{k+2}u_{k+2}^Tb,\ldots,\sigma_n u_n^T b)^T
\left(\frac{1}{u_1^Tb},\frac{1}{u_2^Tb},\ldots,
\frac{1}{u_k^Tb}\right)\right|.
$$
Therefore, we get
\begin{align}
\|\Sigma_k\Delta_k^T\|&=\|\Delta_k\Sigma_k\|\leq \left\||\Delta_k\Sigma_k|\right\|
\notag\\
&\leq |L_{k_1}^{(k)}(0)|\left(\sum_{j=k+1}^n\sigma_j^2| u_j^Tb|^2\right)^{1/2}
\left(\sum_{j=1}^k \frac{1}{| u_j^Tb|^2}\right)^{1/2}. \label{sigdel}
\end{align}
By \eqref{deltaexp}, for $k=1$ we have
$$
\|\Delta_1\Sigma_1\|=\left(\sum_{j=2}^n\sigma_j^2| u_j^Tb|^2\right)^{1/2}
\frac{1}{| u_1^Tb|}.
$$
We have derived the bounds \eqref{severe1} and
\eqref{modeest} for $\left(\sum_{j=k+1}^n\sigma_j^2| u_j^Tb|^2\right)^{1/2}$
for severely and moderately or mildly ill-posed problems, respectively, from
which we obtain \eqref{prodnorm} and \eqref{prodnorm2} for $k=1$. In order
to bound $\|\Sigma_k\Delta_k^T\|$ for $k>1$, we need to estimate
$\left(\sum_{j=1}^k\frac{1}{| u_j^Tb|^2}\right)^{1/2}$.
We next carry out this task for severely and moderately or mildly ill-posed
problems, respectively, for each kind of which we consider the cases of
$k\leq k_0$ and $k>k_0$ separately.

Case of $k\leq k_0$ for severely ill-posed problems: From the discrete
Picard condition \eqref{picard} and \eqref{ideal3}, we obtain
\begin{align*}
\sum_{j=1}^k
\frac{1}{| u_j^Tb|^2}&= \frac{1}{| u_k^Tb|^2} \sum_{j=1}^k
\frac{| u_k^Tb|^2}{| u_j^Tb|^2}
=\frac{1}{| u_k^Tb|^2} \left(1+\mathcal{O}\left(\sum_{j=1}^{k-1}\rho^{2(j-k)
(1+\beta)}\right)\right)\\
&=\frac{1}{| u_k^Tb|^2} \left(1+\mathcal{O}(\rho^{-2(1+\beta)}) \right).
\end{align*}

Case of $k> k_0$ for severely ill-posed problems: From \eqref{ideal3} and
\eqref{ideal2}, we obtain
\begin{align*}
\sum_{j=1}^k
\frac{1}{| u_j^Tb|^2}&= \frac{1}{| u_k^Tb|^2}\left( \sum_{j=1}^{k_0}
\frac{| u_k^Tb|^2}{| u_j^Tb|^2}+\sum_{j=k_0+1}^{k}
\frac{| u_k^Tb|^2}{| u_j^Tb|^2}\right)\\
&=\frac{1}{| u_k^Tb|^2}\left(1+\mathcal{O}\left(\sum_{j=1}^{k_0-1}
\rho^{2(j-k_0)(1+\beta)}\right)+k-k_0\right)\\
&=\frac{1}{| u_k^Tb|^2}\left(1+\mathcal{O}(\rho^{-2(1+\beta)})+k-k_0\right).
\end{align*}
Substituting the above two relations for the two cases into \eqref{sigdel} and
combining them with \eqref{severe1} and \eqref{lkk}, we
get \eqref{prodnorm}.

Case of $k\leq k_0$ for moderately or mildly ill-posed problems: From \eqref{ideal3}
we have
\begin{align*}
\sum_{j=1}^k
\frac{1}{| u_j^Tb|^2}&= \frac{1}{| u_k^Tb|^2} \sum_{j=1}^k
\frac{| u_k^Tb|^2}{| u_j^Tb|^2}= \frac{1}{| u_k^Tb|^2}\sum_{j=1}^k
\left(\frac{j}{k}\right)^{2\alpha (1+\beta)}\\
&<\frac{1}{| u_k^Tb|^2} \sum_{j=1}^k
\left(\frac{j}{k}\right)^{2\alpha}
=\frac{1}{| u_k^T b|^2} k\sum_{j=1}^k \frac{1}{k}\left(\frac{j}{k}
\right)^{2\alpha} \notag \\
&< \frac{1}{| u_k^T b|^2} \left(k \int_0^1
x^{2\alpha}dx+1 \right)=\frac{1}{| u_k^Tb|^2}\left(\frac{k}{2\alpha+1}+1\right).
\end{align*}

Case of $k> k_0$ for moderately or mildly ill-posed problems: From \eqref{ideal3} and
\eqref{ideal2} we have
\begin{align*}
\sum_{j=1}^k\frac{1}{| u_j^Tb|^2}
&= \frac{1}{| u_k^Tb|^2} \left(\sum_{j=1}^{k_0}
\frac{| u_k^Tb|^2}{| u_j^Tb|^2}
+\sum_{j=k_0+1}^{k}
\frac{| u_k^Tb|^2}{| u_j^Tb|^2}\right)\\
&= \frac{1}{| u_k^Tb|^2} \left(\sum_{j=1}^{k_0}
\left(\frac{j}{k_0}\right)^{2\alpha (1+\beta)}+k-k_0\right)\\
&<\frac{1}{| u_k^Tb|^2} \left(\sum_{j=1}^{k_0}
\left(\frac{j}{k_0}\right)^{2\alpha}+k-k_0\right)\\
&\leq \frac{1}{| u_k^Tb|^2} \left(\frac{k_0}{2\alpha+1}+1+k-k_0\right).
\end{align*}
Substituting the above two bounds for the two cases into \eqref{sigdel} and
combining them with \eqref{modeest}, we get \eqref{prodnorm2}.
\qquad\endproof

\eqref{prodnorm} and \eqref{prodnorm2} indicate that $\|\Sigma_k\Delta_k^T\|$
decays swiftly as $k$ increases. As has been seen, we must take some cares to
accurately bound $\|\Sigma_k\Delta_k^T\|$. Indeed,
for $1<k\leq k_0$, if we had simply bounded it by
\begin{equation}\label{rough}
\|\Sigma_k\Delta_k^T\|\leq \|\Sigma_k\|\|\Delta_k^T\|=\sigma_1\|\Delta_k\|,
\end{equation}
the factors $\sigma_{k+1}$
in \eqref{prodnorm} and $\sigma_k$ in \eqref{prodnorm2} would have been replaced by
$\frac{\sigma_1\sigma_{k+1}}{\sigma_k}\approx \sigma_1\rho^{-1}$ and
$\sigma_1$, respectively, by substituting the estimates \eqref{eqres1}
and \eqref{modera1} for $\|\Delta_k\|$ into the above. Such bounds
overestimate $\|\Sigma_k\Delta_k^T\|$ too much as $k$ increases, and
are useless to precisely analyze the regularization of
LSQR, CGME and LSMR for ill-posed problems since they make us
impossible to get those predictively accurate results to be presented
in Sections \ref{rankapp}--\ref{compare}.

As a byproduct, we consider an interesting problem that has its own right,
though its solution will not be used in this paper:
How close to the Krylov subspace $\mathcal{V}_k^R$ is the individual
right singular vector $v_j$ for $j\leq k$ and $k=1,2,\ldots,n-1$? Denote
by $\sin\angle(v_j,\mathcal{V}_k^R)$ the distance between
$v_j$ and $\mathcal{V}_k^R$, which is defined as
$$
\sin\angle(v_j,\mathcal{V}_k^R)=\|(I-\Pi_k)v_j\|=\min_{w\in \mathcal{V}_k^R}
\|v_j-w\|
$$
with $\Pi_k$ the orthogonal projector onto $\mathcal{V}_k^R$. Then
we present the following result.

\begin{theorem}\label{indiv}
Let $\Delta_k=(\delta_1,\delta_2,\ldots,\delta_k)$ be defined by
\eqref{defdelta}. Then for $k=1,2,\ldots,n-1$ and $j=1,2,\ldots,k$ we have
\begin{align}
\frac{\sigma_{\min}(\Delta_k)}{\sqrt{1+\sigma_{\min}^2(\Delta_k)}}
\leq
\sin\angle(v_j,\mathcal{V}_k^R)&\leq
\min\{\|\sin\Theta(\mathcal{V}_k,\mathcal{V}_k^R)\|,\|\delta_j\|\},
\label{errorvj}
\end{align}
where $\sigma_{\min}(\cdot)$ denotes the smallest singular value of
a matrix.
\end{theorem}

{\em Proof}.
We first prove the upper bound of \eqref{errorvj}.
Since the columns of $Z_k$ defined by \eqref{zk} form a basis of
$\mathcal{V}_k^R$, its $j$-th column
$Z_k e_j\in \mathcal{V}_k^R$. As a result, we get
\begin{align*}
\sin\angle(v_j,\mathcal{V}_k^R)&=\min_{w\in \mathcal{V}_k^R}
\|v_j-w\|\leq \|v_j-Z_ke_j\|\\
&=\|v_j-(V_k+V_k^{\perp}\Delta_k)e_j\|=\|v_j-v_j-V_k^{\perp}\delta_j\|\\
&=\|V_k^{\perp}\delta_j\|=\|\delta_j\|.
\end{align*}
Recall from \eqref{decomp} that the columns of $\hat{Z}_k$
form an orthonormal basis of $\mathcal{V}_k^R$, and suppose that
$(\hat{Z}_k,\hat{Z}_k^{\perp})$ is orthogonal. Then the columns of
$\hat{Z}_k^{\perp}$ are an orthonormal basis of the orthogonal complement
of $\mathcal{V}_k^R$ with respect to $\mathbb{R}^n$. Particularly,
$$
\hat{Z}_k^{\perp}=(V_k^{\perp}-V_k\Delta_k^T)(I+\Delta_k\Delta_k^T)^{-\frac{1}{2}}
$$
meets the requirement. By definition, we obtain
\begin{align*}
\|\sin\Theta(\mathcal{V}_k,\mathcal{V}_k^R)\|&=\|(\hat{Z}_k^{\perp})^TV_k\|
=\|\hat{Z}_k^{\perp}(\hat{Z}_k^{\perp})^TV_k\|=\max_{\|c\|=1}
\|\hat{Z}_k^{\perp}(\hat{Z}_k^{\perp})^TV_kc\|\\
&=\max_{\|c\|=1}\|(I-\hat{Z}_k\hat{Z}_k^T)V_kc\|=\max_{\|c\|=1}\|(I-\Pi_k)V_kc\|,
\end{align*}
from which and $v_j=V_ke_j$ it follows that
$$
\sin\angle(v_j,\mathcal{V}_k^R)
\leq \|\sin\Theta(\mathcal{V}_k,\mathcal{V}_k^R)\|
$$
by taking $c=e_j,\ j=1,2,\ldots,k$. So the upper bound of \eqref{errorvj} holds.

We next derive the lower bound of \eqref{errorvj}. We obtain from above
that
\begin{align*}
\sin\angle(v_j,\mathcal{V}_k^R)&=\|(I-\Pi_k)v_j\|=\|(\hat{Z}_k^{\perp})^Tv_j\|\\
&=\|(I+\Delta_k\Delta_k^T)^{-\frac{1}{2}}\left((V_k^{\perp})^T-
\Delta_kV_k^T\right)v_j\|\\
&=\|(I+\Delta_k\Delta_k^T)^{-\frac{1}{2}}\Delta_k e_j\|\\
&\geq \sigma_{\min}\left((I+\Delta_k\Delta_k^T)^{-\frac{1}{2}}\Delta_k\right)
=\frac{\sigma_{\min}(\Delta_k)}{\sqrt{1+\sigma_{\min}^2(\Delta_k)}}. \qquad\endproof
\end{align*}

We remark that the lower bound in \eqref{errorvj} is just the sine of the
smallest canonical angle of $\mathcal{V}_k$ and $\mathcal{V}_k^R$. Since
$v_j\in \mathcal{V}_k$, it is natural that $\angle(v_j,\mathcal{V}_k^R)$ lies
between the smallest and largest angles of $\mathcal{V}_k$ and $\mathcal{V}_k^R$,
as \eqref{errorvj} indicates. The nontrivial point of the upper bound in
\eqref{errorvj} is that $\sin\angle(v_j,\mathcal{V}_k^R)$ can be
much smaller than $\|\sin\Theta(\mathcal{V}_k,\mathcal{V}_k^R)\|$,
as indicated by the bounds \eqref{columndelta} and \eqref{columndelta2},
especially for $j$ not close to $k$.
Combining \eqref{errorvj} with \eqref{columndelta} and \eqref{columndelta2},
we see that the smaller $j$, the closer $v_j$ is to $\mathcal{V}_k^R$.

\section{The rank $k$ approximation $P_{k+1}B_kQ_k^T$ to $A$, the Ritz values
$\theta_i^{(k)}$ and the regularization of LSQR}\label{rankapp}

Making use of Theorems~\ref{thm2}--\ref{thm3}, we are able to solve those
key problems stated before Theorem~\ref{thm2} and give definitive
answers to the fundamental concern by Bj\"{o}rck and Eld\'{e}n, proving
that LSQR has the full regularization
for severely or moderately ill-posed problems with $\rho>1$ or
$\alpha>1$ suitably and it, in general, has only the partial regularization
for mildly ill-posed problems.

Define
\begin{equation}\label{gammak}
\gamma_k = \|A-P_{k+1}B_kQ_k^T\|,
\end{equation}
which measures the accuracy of the rank $k$ approximation $P_{k+1}B_kQ_k^T$ to $A$
generated by Lanczos bidiagonalization. Recall \eqref{xk} and
the comments followed. It is known that the full
or partial regularization of LSQR uniquely depends on whether or not
$\gamma_k\approx \sigma_{k+1}$ holds, where we will make the precise meaning `$\approx$'
clear by introducing the definition of near best
rank $k$ approximation to $A$, and on whether or not the $k$ singular values
$\theta_i^{(k)}$, i.e., Ritz values, of $B_k$, approximate
the $k$ large singular values $\sigma_i$ of $A$ in natural order
for $k=1,2,\ldots, k_0$. If both of them hold, LSQR has the full regularization;
if either of them is not satisfied, LSQR has only the partial
regularization.

This section consists of three subsections. In Section \ref{rankaccur}, we
present accurate estimates for $\gamma_k$ for the three kinds of ill-posed
problems under consideration. We prove that, under some reasonable conditions
on $\rho$ or $\alpha$, the matrix $P_{k+1}B_kQ_k^T$
is a near best rank $k$ approximation to $A$. In Section \ref{ritzapprox},
we deepen the results in Section \ref{rankaccur} and show
how the $k$ Ritz values $\theta_i^{(k)}$ behave. We derive the sufficient
conditions on $\rho$ and $\alpha$ for which
they approximate the first $k$ large singular values $\sigma_i$ of $A$
in natural order. In Section \ref{morerank}, we consider general best and near
best rank approximations to $A$ with respect to the 2-norm. For
$A$ with $\sigma_i=\zeta i^{-\alpha},\ i=1,2,\ldots,n$,
we analyze the nonzero singular values of such a rank
$k$ approximation, and prove that they approximate the first $k$
large singular values of $A$ for $\alpha>1$ suitably but can
fail to do so for $\frac{1}{2}<\alpha \leq 1$. These results will help
understand the regularizing effects of LSQR.

\subsection{The accuracy of rank $k$ approximation $P_{k+1}B_kQ_k^T$ to $A$
and more related}\label{rankaccur}

We first present one of the main results in this paper.

\begin{theorem}\label{main1}
Assume that the discrete Picard condition \eqref{picard} is
satisfied. Then for $k=1,2,\ldots,n-1$ we have
\begin{equation}\label{final}
  \sigma_{k+1}\leq \gamma_k\leq \sqrt{1+\eta_k^2}\sigma_{k+1}
\end{equation}
with
\begin{equation} \label{const1}
\eta_k\leq \left\{\begin{array}{ll}
\xi_k\frac{| u_{k+1}^Tb|}{| u_k^T b|}
\left(1+\mathcal{O}(\rho^{-2})\right)
& \mbox{ for } 1\leq k\leq k_0,\\
\xi_k\sqrt{k-k_0+1}\left(1+\mathcal{O}(\rho^{-2})\right)
& \mbox{ for } k_0<k \leq n-1
\end{array}
\right.
\end{equation}
for severely ill-posed problems and
\begin{equation}\label{const2}
\eta_k\leq \left\{\begin{array}{ll}
\xi_1\frac{\sigma_1}{\sigma_2}\frac{| u_2^Tb|}{| u_1^Tb|}
\sqrt{\frac{1}{2\alpha-1}} & \mbox{ for } k=1, \\
\xi_k\frac{\sigma_k}{\sigma_{k+1}}\frac{|u_{k+1}^T b|}{|u_k^T b|}
\sqrt{\frac{k^2}{4\alpha^2-1}+\frac{k}{2\alpha-1}}
|L_{k_1}^{(k)}(0)|& \mbox{ for } 1< k\leq k_0, \\
\xi_k\frac{\sigma_k}{\sigma_{k+1}}\sqrt{\frac{k k_0}{4\alpha^2-1}+
\frac{k(k-k_0+1)}{2\alpha-1}}|L_{k_1}^{(k)}(0)|
& \mbox{ for } k_0<k\leq n-1
\end{array}
\right.
\end{equation}
for moderately or mildly ill-posed problems with $\sigma_j=\zeta j^{-\alpha},\
j=1,2,\ldots,n$, where
$\xi_k=\sqrt{\left(\frac{\|\Delta_k\|}{1+\|\Delta_k\|^2}\right)^2+1}$ for
$\|\Delta_k\|<1$ and $\xi_k\leq\frac{\sqrt{5}}{2}$ for $\|\Delta_k\|\geq 1$
with $\Delta_k$ defined
by \eqref{defdelta}.
\end{theorem}

{\em Proof}.
Since $A_k$ is the best rank $k$ approximation to
$A$ with respect to the 2-norm and $\|A-A_k\|=\sigma_{k+1}$,
the lower bound in \eqref{final} holds. Next we prove the upper bound.

From \eqref{eqmform1}, we obtain
\begin{align}
\gamma_k
&= \|A-P_{k+1}B_kQ_k^T\|= \|A-AQ_kQ_k^T\|= \|A(I-Q_kQ_k^T)\|. \label{gamma2}
\end{align}
From Algorithm 1, \eqref{kry}, \eqref{zk} and \eqref{decomp}, we obtain
$$
\mathcal{V}_k^R
=\mathcal{K}_{k}(A^{T}A,A^{T}b)=span\{Q_k\}=span\{\hat{Z}_k\}
$$
with $Q_k$ and $\hat{Z}_k$ being orthonormal, and the
orthogonal projector onto $\mathcal{V}_k^R$ is thus
\begin{equation}\label{twobasis}
Q_kQ_k^T=\hat{Z}_k\hat{Z}_k^T.
\end{equation}
Keep in mind that $A_k=U_k\Sigma_k V_k^T$. It is direct to justify that
$(U_k\Sigma_k V_k^T)^T(A-U_k\Sigma_k V_k^T)=\mathbf{0}$ for $k=1,2,\ldots,n-1$.
Therefore, exploiting this and noting that $\|I-\hat{Z}_k\hat{Z}_k^T\|=1$ and
$V_k^TV_k^{\perp}=\mathbf{0}$  for $k=1,2,\ldots,n-1$,
we get from \eqref{gamma2}, \eqref{twobasis} and \eqref{decomp} that
\begin{align}
\gamma_k^2 &= \|(A-U_k\Sigma_kV_k^T+U_k\Sigma_kV_k^T)(I-\hat{Z}_k\hat{Z}_k^T)\|^2
  \notag\\
  &=\max_{\|y\|=1}\|(A-U_k\Sigma_kV_k^T+U_k\Sigma_kV_k^T)
  (I-\hat{Z}_k\hat{Z}_k^T)y\|^2 \notag\\
  &=\max_{\|y\|=1}\|(A-U_k\Sigma_kV_k^T)(I-\hat{Z}_k\hat{Z}_k^T)y+
 U_k\Sigma_kV_k^T(I-\hat{Z}_k\hat{Z}_k^T)y\|^2\notag\\
 &=\max_{\|y\|=1}\left(\|(A-U_k\Sigma_kV_k^T)(I-\hat{Z}_k\hat{Z}_k^T)y\|^2+
 \| U_k\Sigma_kV_k^T(I-\hat{Z}_k\hat{Z}_k^T)y\|^2\right)\notag\\
 &\leq \|(A-U_k\Sigma_kV_k^T)(I-\hat{Z}_k\hat{Z}_k^T)\|^2+
 \| U_k\Sigma_kV_k^T(I-\hat{Z}_k\hat{Z}_k^T)\|^2 \notag \\
 &\leq \sigma_{k+1}^2+\| \Sigma_kV_k^T(I-\hat{Z}_k\hat{Z}_k^T)\|^2 \notag\\
 &\leq \sigma_{k+1}^2+\|\Sigma_kV_k^T\left(I-(V_k+V_k^{\perp}\Delta_k)(I+
 \Delta_k^T\Delta_k)^{-1}(V_k+V_k^{\perp}\Delta_k)^T\right)\|^2\notag\\
 &= \sigma_{k+1}^2 + \left\|\Sigma_k\left(V_k^T-(I+
 \Delta_k^T\Delta_k)^{-1}(V_k+V_k^{\perp}\Delta_k)^T\right)\right\|^2 \notag\\
  &= \sigma_{k+1}^2 + \left\|\Sigma_k(I+
 \Delta_k^T\Delta_k)^{-1}\left((I+\Delta_k^T\Delta_k)V_k^T-
 \left(V_k+V_k^{\perp}\Delta_k\right)^T\right)\right\|^2 \notag\\
  &= \sigma_{k+1}^2+ \|\Sigma_k(I+
 \Delta_k^T\Delta_k)^{-1}\left(\Delta_k^T\Delta_kV_k^T-\Delta_k^T
 (V_k^{\perp})^T\right)\|^2\notag\\
  &= \sigma_{k+1}^2 + \|\Sigma_k(I+
 \Delta_k^T\Delta_k)^{-1}\Delta_k^T\Delta_kV_k^T-\Sigma_k(I+
 \Delta_k^T\Delta_k)^{-1}\Delta_k^T(V_k^{\perp})^T\|^2 \label{twomatrix}\\
    &\leq \sigma_{k+1}^2 +\|\Sigma_k(I+
 \Delta_k^T\Delta_k)^{-1}\Delta_k^T\Delta_k\|^2+\|\Sigma_k(I+
 \Delta_k^T\Delta_k)^{-1}\Delta_k^T\|^2\notag\\
 &=\sigma_{k+1}^2+\epsilon_k^2,
 \label{estimate1}
\end{align}
where the last inequality follows by using $V_k^T V_k^{\perp}=\mathbf{0}$ and
the definition of the induced matrix 2-norm to amplify the second term
in \eqref{twomatrix}.

We estimate $\epsilon_k$ accurately below. To this end, we need to use two
key identities and some results related. By the SVD of $\Delta_k$, it is direct
to justify that
\begin{equation}\label{inden1}
(I+
 \Delta_k^T\Delta_k)^{-1}\Delta_k^T\Delta_k=\Delta_k^T\Delta_k(I+
 \Delta_k^T\Delta_k)^{-1}
\end{equation}
and
\begin{equation}\label{inden2}
(I+
 \Delta_k^T\Delta_k)^{-1}\Delta_k^T=\Delta_k^T(I+
 \Delta_k\Delta_k^T)^{-1}.
\end{equation}
Define the function $f(\lambda)=\frac{\lambda}{1+\lambda^2}$ with
$\lambda\in [0,\infty)$. Since the derivative
$f^{\prime}(\lambda)=\frac{1-\lambda^2}{(1+\lambda^2)^2}$,
$f(\lambda)$ is monotonically increasing for $\lambda\in [0,1]$
and decreasing for $\lambda\in [1,\infty)$, and the maximum of
$f(\lambda)$ over $\lambda\in [0,\infty)$ is $\frac{1}{2}$, which attains at
$\lambda=1$. Based on
these properties and exploiting the SVD of $\Delta_k$, for the matrix 2-norm
we get
\begin{equation}\label{compact}
\|\Delta_k(I+\Delta_k^T\Delta_k)^{-1}\|=\frac{\|\Delta_k\|}{1+\|\Delta_k\|^2}
\end{equation}
for $\|\Delta_k\|<1$ and
\begin{equation}\label{noncomp}
\|\Delta_k(I+\Delta_k^T\Delta_k)^{-1}\|\leq\frac{1}{2}
\end{equation}
for $\|\Delta_k\|\geq 1$ (Note: in this case, since $\Delta_k$ may have at least
one singular value smaller than one, we do not
have an expression like \eqref{compact}). It then follows
from \eqref{estimate1}, \eqref{compact}, \eqref{noncomp}
and $\|(1+\Delta_k\Delta_k^T)^{-1}\|\leq 1$ that
\begin{align}
\epsilon_k^2&=\|\Sigma_k \Delta_k^T\Delta_k(I+
 \Delta_k^T\Delta_k)^{-1}\|^2+\|\Sigma_k \Delta_k^T (I+
 \Delta_k\Delta_k^T)^{-1}\|^2 \label{separa}\\
 &\leq \|\Sigma_k\Delta_k^T\|^2\|\Delta_k(I+\Delta_k^T\Delta_k)^{-1}\|^2+
 \|\Sigma_k\Delta_k^T\|^2 \|(1+\Delta_k\Delta_k^T)^{-1}\|^2 \notag\\
 &\leq \|\Sigma_k\Delta_k^T\|^2\left(\|\Delta_k
 (I+\Delta_k^T\Delta_k)^{-1}\|^2+1\right)\notag\\
  &=\|\Sigma_k\Delta_k^T\|^2\left(\left(\frac{\|\Delta_k\|}
  {1+\|\Delta_k\|^2}\right)^2+1\right)
  =\xi_k^2\|\Sigma_k\Delta_k^T\|^2 \notag
\end{align}
for $\|\Delta_k\|<1$ and
$$
\epsilon_k\leq \|\Sigma_k\Delta_k^T\|\sqrt{\|\Delta_k
(I+\Delta_k^T\Delta_k)^{-1}\|^2+1}=\xi_k\|\Sigma_k\Delta_k^T\|
\leq \frac{\sqrt{5}}{2}\|\Sigma_k\Delta_k^T\|
$$
for $\|\Delta_k\|\geq 1$. Replace $\|\Sigma_k\Delta_k^T\|$ by
its bounds \eqref{prodnorm} and \eqref{prodnorm2} in the
above, insert the resulting bounds for $\epsilon_k$
into \eqref{estimate1}, and let $\epsilon_k=\eta_k\sigma_{k+1}$.
Then we obtain the upper bound in \eqref{final} with $\eta_k$
satisfying \eqref{const1} and \eqref{const2} for severely and moderately
or mildly ill-posed problems, respectively.
\qquad\endproof

Note from \eqref{ideal3} that
$$
\frac{|u_{k+1}^T b|}{| u_k^T b|}=
\frac{\sigma_{k+1}^{1+\beta}}{\sigma_k^{1+\beta}},\ k\leq k_0.
$$
Therefore, for the right-hand side of \eqref{const2} and $k\leq k_0$
we have
$$
\frac{\sigma_k}{\sigma_{k+1}}\frac{| u_{k+1}^T b|}{| u_k^T b|}=
\left(\frac{\sigma_{k+1}}{\sigma_k}\right)^{\beta}<1.
$$

\begin{remark}\label{decayrate}
For severely ill-posed problems, from \eqref{case3}, \eqref{case4} and
the definition of $\xi_k$ we know that
$$
\xi_k(1+\mathcal{O}(\rho^{-2}))
=1+\mathcal{O}(\rho^{-2})
$$
for both $k\leq k_0$ and $k>k_0$. Therefore, from
\eqref{const1} and \eqref{ideal3}, for $k\leq k_0$ we have
\begin{equation}\label{etak0}
\eta_k\leq \xi_k
\frac{|u_{k+1}^Tb |}{|u_k^T b |}\left(1+\mathcal{O}(\rho^{-2})\right)
=\frac{|u_{k+1}^Tb |}{|u_k^T b |}=
\frac{\sigma_{k+1}^{1+\beta}}{\sigma_k^{1+\beta}}=\mathcal{O}(\rho^{-1-\beta})<1
\end{equation}
by ignoring the smaller term $\mathcal{O}(\rho^{-1-\beta})
\mathcal{O}(\rho^{-2})
=\mathcal{O}(\rho^{-3-\beta})$, and for $k>k_0$ we have
\begin{equation}\label{incres}
\eta_k\leq \xi_k\sqrt{k-k_0+1}\left(1+\mathcal{O}(\rho^{-2})\right)
=\sqrt{k-k_0+1}
\end{equation}
by ignoring the smaller term $\sqrt{k-k_0+1}\mathcal{O}(\rho^{-2})$,
which increases slowly with $k$.
\end{remark}

\begin{remark}
For the moderately or mildly ill-posed problems with
$\sigma_j=\zeta j^{-\alpha}$,
from the derivation on $\eta_k$ and its estimate \eqref{const2},
by comparing \eqref{k2} and \eqref{modera1} with \eqref{const2}, for $k\leq k_0$
we approximately have
\begin{equation}\label{etadelta}
\frac{\sigma_k}{\sigma_{k+1}}\|\Delta_k\|\leq
\eta_k\leq \frac{\sqrt{5}}{2}\frac{\sigma_k}{\sigma_{k+1}}\|\Delta_k\|,
\end{equation}
and for $k>k_0$, from \eqref{lkkmoderate} and \eqref{lowerbound}
we approximately have
\begin{align}
\eta_k&< \frac{\sigma_k}{\sigma_{k+1}}\sqrt{\frac{k k_0}{4\alpha^2-1}+
\frac{k(k-k_0+1)}{2\alpha-1}}|L_{k_1}^{(k)}(0)| \notag\\
&\sim \frac{k^{3/2}\sqrt{k_0}}{(2\alpha+1)\sqrt{{4\alpha^2-1}}}+
\frac{k^{3/2}\sqrt{k-k_0+1}}{(2\alpha+1)\sqrt{{2\alpha-1}}}, \label{asym}
\end{align}
which increases faster than the right-hand side of \eqref{incres} with
respect to $k$.
\end{remark}

\begin{remark}\label{decayrate2}
From \eqref{final}, \eqref{const1} and \eqref{etak0},
for severely ill-posed problems we have
$$
1<\sqrt{1+\eta_k^2}<1+\frac{1}{2}{\eta_k^2}\leq
1+\frac{1}{2}\frac{\sigma_{k+1}^{2(1+\beta)}}{\sigma_k^{2(1+\beta)}}
\sim 1+\frac{1}{2}\rho^{-2(1+\beta)},
$$
and $\gamma_k$ is an accurate approximation to
$\sigma_{k+1}$ for $k\leq k_0$ and marginally less accurate for $k>k_0$.
Thus, the rank $k$ approximation $P_{k+1}B_kQ_k^T$ is as accurate as
the best rank $k$ approximation $A_k$ within the
factor $\sqrt{1+\eta_k^2}\approx 1$ for $k\leq k_0$ and $\rho>1$ suitably.
For moderately ill-posed problems, $\gamma_k$ is still an excellent
approximation to $\sigma_{k+1}$, and the rank $k$ approximation
$P_{k+1}B_kQ_k^T$ is almost as accurate as the best rank $k$ approximation
$A_k$ for $k\leq k_0$. Therefore, $P_{k+1}B_kQ_k^T$ plays the same role as
$A_k$ for these two kinds of ill-posed problems and $k\leq k_0$, it is known
from the clarification in Section \ref{lsqr} that LSQR may have the full
regularization. We will, afterwards, deepen this theorem and derive
more results, proving that LSQR must have the full regularization for
these two kinds of problems provided that $\rho>1$ and $\alpha>1$ suitably.

For both severely and moderately ill-posed problems, we note that the
situation is not so satisfying for increasing $k>k_0$. But at that time,
a possibly big $\eta_k$ does not do harm to our regularization purpose
since we will prove that, provided that $\rho>1$ and $\alpha>1$ suitably,
LSQR has the full regularization and has already found
a best possible regularized solution at semi-convergence occurring at
iteration $k_0$. If it is the case, we will simply stop performing it
after semi-convergence.
\end{remark}

\begin{remark}\label{mildre}
For mildly ill-posed problems, the situation is fundamentally different.
As clarified in Remark~\ref{mildrem}, we have
$\sqrt{\frac{k^2}{4\alpha^2-1}+\frac{k}{2\alpha-1}}>1$ and
$|L_k^{(k)}(0)|>1$ considerably as $k$ increases up to $k_0$
because of $\frac{1}{2}<\alpha\leq 1$,
leading to $\eta_k>1$ substantially. This means that
$\gamma_{k_0}$ is substantially bigger than $\sigma_{k_0+1}$ and can
well lie between $\sigma_{k_0}$ and $\sigma_1$, so that
the rank $k_0$ approximation $P_{k_0+1}B_{k_0}Q_{k_0}^T$ is much less accurate
than the best rank $k_0$ approximation $A_{k_0}$ and LSQR has only the partial
regularization.
\end{remark}

\begin{remark}
For a given ill-posed problem, the noise level $\|e\|$ only affects $k_0$ but
has no effect on the overall decay rate of $\gamma_k$.
\end{remark}

\begin{remark}
There are several subtle treatments in the proof of Theorem~\ref{main1}, each of
which turns out to be absolutely necessary. Ignoring or missing any one of them
would be fatal and make us fail to obtain accurate estimates for $\epsilon_k$
defined by \eqref{estimate1}: The first is the treatment of
$\|U_k\Sigma_kV_k^T(I-\hat{Z}_k\hat{Z}_k^T)\|$.
By the definition of $\|\sin\Theta(\mathcal{V}_k,\mathcal{V}_k^R)\|$, if we had
amplified it by
$$
\|U_k\Sigma_kV_k^T(I-\hat{Z}_k\hat{Z}_k^T)\|
\leq \|\Sigma_k\|\|V_k^T(I-\hat{Z}_k\hat{Z}_k^T)\|=
\sigma_1\|\sin\Theta(\mathcal{V}_k,\mathcal{V}_k^R)\|,
$$
we would have obtained a too large overestimate, which is almost a fixed
constant for severely
ill-posed problems and $k=1,2,\ldots,k_0$ and increases with
$k=1,2,\ldots,k_0$ for moderately and mildly ill-posed problems. Such rough
estimates are useless to get a meaningful bound for $\gamma_k$. The key is to treat
$U_k\Sigma_kV_k^T(I-\hat{Z}_k\hat{Z}_k^T)$ as a
whole other rather separate it in the above way, so that we can
bound its norm accurately. The second is the use of
\eqref{inden1} and \eqref{inden2}. The third is the extraction of
$\|\Sigma_k\Delta_k^T\|$ from \eqref{separa} as a whole other than
amplify it to $\|\Sigma_k\|\|\Delta_k\|=\sigma_1\|\Delta_k\|$,
i.e., the fatal overestimate \eqref{rough}.
The fourth is accurate estimates for it;
see \eqref{prodnorm} and \eqref{prodnorm2} in Theorem~\ref{thm3}.
For example, without using \eqref{inden1} and \eqref{inden2}, we would have
no way but to obtain
\begin{align*}
\epsilon_k^2 &\leq \|\Sigma_k\|^2\|(I+
 \Delta_k^T\Delta_k)^{-1}\Delta_k^T\Delta_k\|^2+\|\Sigma_k\|^2\|(I+
 \Delta_k^T\Delta_k)^{-1}\Delta_k^T\|^2\\
 &=\sigma_1^2\left(\frac{\|\Delta_k\|^2}{1+\|\Delta_k\|^2}\right)^2+\sigma_1^2
 \|(I+\Delta_k^T\Delta_k)^{-1}\Delta_k^T\|^2\\
 &=\sigma_1^2\left(\frac{\|\Delta_k\|^2}{1+\|\Delta_k\|^2}\right)^2+\sigma_1^2
 \|\Delta_k(I+\Delta_k^T\Delta_k)^{-1}\|^2.
\end{align*}
From \eqref{compact}, \eqref{noncomp} and the previous
estimates for $\|\Delta_k\|$, such bound is too pessimistic
and completely useless in our context, and
it even does not decrease and could not be small as $k$ increases, while
our estimates for $\epsilon_k=\eta_k\sigma_{k+1}$
in Theorem~\ref{main1} are much more accurate and
decay swiftly as $k$ increases, as indicated by \eqref{const1}
and \eqref{const2}.
\end{remark}

In order to prove the full or partial regularization of LSQR for
\eqref{eq1} completely and rigorously, besides Theorem~\ref{main1},
it appears that we need to introduce a precise definition of the near best
rank $k$ approximation
$P_{k+1}B_kQ_k^T$ to $A$, i.e., the precise meaning of
$\gamma_k\approx \sigma_{k+1}$. By definition \eqref{gammak}, the rank
$k$ matrix $P_{k+1}B_kQ_k^T$ is called a near best rank $k$ approximation
to $A$ if it satisfies
\begin{equation}\label{near}
\sigma_{k+1}\leq \gamma_k<\sigma_k \mbox{ and } \gamma_k-\sigma_{k+1}
<\sigma_k-\gamma_k,\mbox{ i.e., } \gamma_k<\frac{\sigma_k+\sigma_{k+1}}{2},
\end{equation}
that is, $\gamma_k$ lies between $\sigma_k$ and $\sigma_{k+1}$ and is closer to
$\sigma_{k+1}$. This definition is natural.
For an ill-posed problem \eqref{eq1}, since there is no
considerable gap of $\sigma_k$ and $\sigma_{k+1}$, the definition means that
$\gamma_k$ must approximate $\sigma_{k+1}$ more accurately as $k$ increases.
We mention in passing that a near best rank $k$ approximation to $A$ from
an ill-posed problem is much more stringent than it is for a matrix from
a (numerically ) rank-deficient problem where the large singular values are
very well separated from the small ones and there is a substantial gap
between two groups of singular values. In addition, we point out that it may
be much harder to computationally obtain a near best $k$ rank approximation to
the large $A$ from the ill-posed problem than for a
numerically rank deficient matrix of the same order.

Based on Theorem~\ref{main1}, for the severely and moderately or mildly ill-posed
problems with the singular value models $\sigma_k=\zeta\rho^{-k}$ and
$\sigma_k=\zeta k^{-\alpha}$, we next derive the sufficient conditions on
$\rho$ and $\alpha$ that guarantee that $P_{k+1}B_kQ_k^T$ is a near best rank $k$
approximation to $A$ for $k=1,2,\ldots,k_0$. We analyze if and how the
sufficient conditions are satisfied for three kinds of ill-posed problems.

\begin{theorem}\label{nearapprox}
For a given \eqref{eq1}, assume that the discrete Picard condition
\eqref{picard} is satisfied. Then, in the sense of \eqref{near},
$P_{k+1}B_kQ_k^T$ is a near best rank $k$ approximation to $A$
for $k=1,2,\ldots,k_0$ if
\begin{equation}\label{condition}
\sqrt{1+\eta_k^2}<\frac{1}{2}\frac{\sigma_k}{\sigma_{k+1}}+\frac{1}{2}.
\end{equation}
For the severely ill-posed problems with $\sigma_k=\zeta\rho^{-k}$ and
the moderately or mildly ill-posed problems with $\sigma_k=\zeta k^{-\alpha}$,
$P_{k+1}B_kQ_k^T$ is a near best rank $k$ approximation to $A$
for $k=1,2,\ldots,k_0$ if $\rho>2$ and $\alpha$ satisfies
\begin{equation}\label{condition1}
2\sqrt{1+\eta_k^2}-1<\left(\frac{k_0+1}{k_0}\right)^{\alpha},
\end{equation}
respectively.
\end{theorem}

{\em Proof}.
By \eqref{final}, we see that $\gamma_k\leq \sqrt{1+\eta_k^2}\sigma_{k+1}$.
Therefore, $P_{k+1}B_kQ_k^T$ is a near best rank $k$ approximation to $A$ in
the sense of \eqref{near} provided that
$$
\sqrt{1+\eta_k^2}\sigma_{k+1}<\sigma_k
$$
and
$$
\sqrt{1+\eta_k^2}\sigma_{k+1}<\frac{\sigma_k+\sigma_{k+1}}{2},
$$
from which \eqref{condition} follows.

From \eqref{etak0}, for the severely ill-posed problems with
$\sigma_k=\zeta\rho^{-k}$ and $\rho>1$ we have
\begin{equation}\label{simp}
\sqrt{1+\eta_k^2}<1+\frac{1}{2}\eta_k^2\leq 1+\frac{1}{2}\rho^{-2(1+\beta)}
<1+\rho^{-1}, \ k=1,2,\ldots,k_0,
\end{equation}
from which it follows that
\begin{align}\label{ampli}
\sqrt{1+\eta_k^2}\sigma_{k+1}
&<(1+\rho^{-1})\sigma_{k+1}.
\end{align}
Since $\sigma_k/\sigma_{k+1}=\rho$, \eqref{condition} holds provided that
$$
1+\rho^{-1}<\frac{1}{2}\rho+\frac{1}{2},
$$
i.e., $\rho^2-\rho-2>0$, solving which for $\rho$ we get $\rho>2$. For the
moderately or mildly ill-posed problems with $\sigma_k=\zeta k^{-\alpha}$,
it is direct from \eqref{condition} to
get
$$
2\sqrt{1+\eta_k^2}-1<\left(\frac{k+1}{k}\right)^{\alpha}.
$$
Since $\left(\frac{k+1}{k}\right)^{\alpha}$ decreases
monotonically as $k$ increases, its minimum over $k=1,2,\ldots,k_0$
is $\left(\frac{k_0+1}{k_0}\right)^{\alpha}$.
Therefore, we obtain \eqref{condition1}.
\qquad\endproof

\begin{remark}
Given the noise level $\|e\|$, the discrete Picard condition \eqref{picard}
and \eqref{picard1}, from the bound \eqref{const2} for
$\eta_k,\,k=1,2,\ldots,k_0$, we see that the bigger $\alpha>1$ is, the smaller
$k_0$ and $\eta_k$ are. Therefore,
there must be $\alpha>1$ such that \eqref{condition1} holds.
Here we should remind that it is more suitable to
regard the conditions on $\rho$ and $\alpha$ as an indication that
$\rho$ and $\alpha$ must not be close to one other than precise requirements
since we have used the bigger \eqref{simp} and simplified
models $\sigma_k=\zeta \rho^{-k}$ and $\sigma_k=\zeta k^{-\alpha}$.
\end{remark}

\begin{remark}
For the mildly ill-posed problems with $\sigma_k=\zeta k^{-\alpha}$,
Theorem~\ref{moderate} has shown that
$\|\Delta_k\|$ is generally not small and can be arbitrarily large
for $k=1,2,\ldots,k_0$. From \eqref{etadelta}, we see
that $\eta_k$ has comparable size to $\|\Delta_k\|$. Note that the right-hand
side $\left(\frac{k_0+1}{k_0}\right)^{\alpha}\leq 2$ for
$\frac{1}{2}<\alpha\leq 1$ and any $k_0\geq 1$. Consequently, \eqref{condition1}
cannot be met generally for mildly ill-posed problems. The rare possible
exceptions are that $k_0$ is only very few and $\alpha$ is close to one since,
in such case, $\eta_k$ is not large for $k=1,2,\ldots,k_0$.
So, $P_{k+1}B_kQ_k^T$ is generally not a near best rank $k$ approximation
to $A$ for $k=1,2,\ldots, k_0$ for this kind of problem.
\end{remark}

\subsection{The approximation behavior of the Ritz values $\theta_i^{(k)}$}
\label{ritzapprox}

In this subsection, starting with Theorem~\ref{main1}, we prove that, under
certain sufficient conditions on $\rho$ and $\alpha$ for the severely
and moderately ill-posed problems with the models $\sigma_i=\zeta\rho^{-i}$ and
$\sigma_i=\zeta i^{-\alpha}$, respectively,
the $k$ Ritz values $\theta_i^{(k)}$ approximate the first
$k$ large singular values $\sigma_i$ in natural order  for $k=1,2,\ldots,k_0$,
which means that no Ritz value smaller than $\sigma_{k_0+1}$ appears.
Combining this result with Theorem~\ref{nearapprox},
we can draw the definite conclusion that LSQR must have the full
regularization for these two kinds of problems provided that $\rho>1$ and
$\alpha>1$ suitably.

\begin{theorem}\label{ritzvalue}
Assume that \eqref{eq1} is severely ill-posed with
$\sigma_i=\zeta\rho^{-i}$ and $\rho>1$ or moderately ill-posed with
$\sigma_i=\zeta i^{-\alpha}$ and $\alpha>1$,
and the discrete Picard condition \eqref{picard} is
satisfied. Let the Ritz values $\theta_i^{(k)}$ be labeled
as $\theta_1^{(k)}>\theta_2^{(k)}>\cdots>\theta_{k}^{(k)}$.
Then
\begin{align}
0<\sigma_i-\theta_i^{(k)} &\leq \sqrt{1+\eta_k^2}\sigma_{k+1},\
i=1,2,\ldots,k.\label{error}
\end{align}
If $\rho\geq 1+\sqrt{2}$ or $\alpha>1$ satisfies
\begin{equation}\label{condm}
1+\sqrt{1+\eta_{k}^2}<\left(\frac{k_0+1}{k_0}\right)^{\alpha},\
k=1,2,\ldots,k_0,
\end{equation}
then the $k$ Ritz values $\theta_i^{(k)}$ strictly interlace
the first large $k+1$ singular values of $A$ and approximate
the first $k$ large ones in natural order for $k=1,2,\ldots,k_0$:
\begin{align}
\sigma_{i+1}&<\theta_i^{(k)}<\sigma_i,\,i=1,2,\ldots,k,
\label{error2}
\end{align}
meaning that there is no Ritz value $\theta_i^{(k)}$ smaller than $\sigma_{k_0+1}$
for $k=1,2,\ldots, k_0$.
\end{theorem}

{\em Proof}.
Note that for $k=1,2,\ldots,k_0$ the $\theta_i^{(k)},\ i=1,2,\ldots,k$ are
just the nonzero singular values of $P_{k+1}B_kQ_k^T$, whose other $n-k$
singular values are zeros. We write
$$
A=P_{k+1}B_k Q_k^T+(A-P_{k+1}B_k Q_k^T)
$$
with $\|A-P_{k+1}B_k Q_k^T\|=\gamma_k$
by definition \eqref{gammak}. Then by the Mirsky's theorem of
singular values \cite[p.204, Thm 4.11]{stewartsun}, we have
\begin{equation}\label{errbound}
| \sigma_i-\theta_i^{(k)}|\leq \gamma_k\leq
\sqrt{1+\eta_k^2}\sigma_{k+1},\ i=1,2,\ldots,k.
\end{equation}
Since the singular values of $A$ are simple and $b$ has components in all the
left singular vectors $u_1,u_2,\ldots, u_n$ of $A$, Lanczos bidiagonalization,
i.e., Algorithm 1, can be run to completion, producing $P_{n+1},\ Q_n$ and
the lower bidiagonal $B_n\in \mathbb{R}^{(n+1)\times n}$ such that
\begin{equation}\label{fulllb}
P^TAQ_n=\left(\begin{array}{c}
B_n\\
\mathbf{0}
\end{array}
\right)
\end{equation}
with the $m\times m$ matrix $P=(P_{n+1},\hat{P})$ and $n\times n$ matrix $Q_n$
orthogonal and all
the $\alpha_i$ and $\beta_{i+1}$, $i=1,2,\ldots,n$, of $B_n$ being positive.
Note that the singular values of $B_k,\ k=1,2,\ldots,n,$
are all simple and that $B_k$ consists of the first $k$ columns of $B_n$
with the last $n-k$ {\em zero} rows deleted. Applying the Cauchy's {\em strict}
interlacing theorem \cite[p.198, Corollary 4.4]{stewartsun} to the singular
values of $B_k$ and $B_n$, we have
\begin{align}
\sigma_{n-k+i}< \theta_i^{(k)}&< \sigma_i,\ i=1,2,\ldots,k.
\label{interlace}
\end{align}
Therefore, \eqref{errbound} becomes
\begin{equation}\label{ritzapp}
0< \sigma_i-\theta_i^{(k)}\leq\gamma_{k}\leq
\sqrt{1+\eta_k^2}\sigma_{k+1},\ i=1,2,\ldots,k,
\end{equation}
which proves \eqref{error}.
That is, the $\theta_i^{(k)}$ approximate $\sigma_i$ from below
for $i=1,2,\ldots,k$ with the errors no more than
$\gamma_k\leq \sqrt{1+\eta_k^2}\sigma_{k+1}$.
For $i=1,2,\ldots,k$, notice that $\rho^{-k+i}\leq 1$. Then
from \eqref{ritzapp}, \eqref{simp} and $\sigma_i=\zeta\rho^{-i}$ we obtain
\begin{align*}
\theta_i^{(k)}&\geq \sigma_i-\gamma_k>\sigma_i-
(1+\rho^{-1})\sigma_{k+1}\\
&=\zeta\rho^{-i}-\zeta (1+\rho^{-1})\rho^{-(k+1)}\\
&=\zeta\rho^{-(i+1)}(\rho-(1+\rho^{-1})\rho^{-k+i})\\
&\geq \zeta\rho^{-(i+1)}(\rho-\rho^{-1}-1)\\
&\geq\zeta\rho^{-(i+1)}=\sigma_{i+1},
\end{align*}
provided that $\rho-\rho^{-1}\geq 2$, solving which we get $\rho\geq 1+\sqrt{2}$.
Together with the upper bound of \eqref{interlace}, we have proved \eqref{error2}.

For the moderately ill-posed problems with $\sigma_i=\zeta i^{-\alpha},\
i=1,2,\ldots,k$ and $k=1,2,\ldots,k_0$,
we get
\begin{align*}
\theta_i^{(k)}&\geq \sigma_i-\gamma_k\geq\sigma_i-\sqrt{1+\eta_k^2}
\sigma_{k+1}\\
&=\zeta i^{-\alpha}-\zeta \sqrt{1+\eta_k^2}(k+1)^{-\alpha}\\
&=\zeta (i+1)^{-\alpha}\left(\left(\frac{i+1}{i}\right)^{\alpha}
-\sqrt{1+\eta_k^2}\left(\frac{i+1}{k+1}\right)^{\alpha}\right)\\
&>\zeta (i+1)^{-\alpha}=\sigma_{i+1},
\end{align*}
i.e., \eqref{error2} holds, provided that $\eta_k>0$ and $\alpha>1$ are such that
$$
\left(\frac{i+1}{i}\right)^{\alpha}
-\sqrt{1+\eta_k^2}\left(\frac{i+1}{k+1}\right)^{\alpha}>1,
$$
which means that
$$
\sqrt{1+\eta_k^2}<\left(\left(\frac{i+1}{i}\right)^{\alpha}-1\right)
\left(\frac{k+1}{i+1}\right)^{\alpha}=
\left(\frac{k+1}{i}\right)^{\alpha}-\left(\frac{k+1}{i+1}\right)^{\alpha},\
i=1,2,\ldots,k.
$$
It is easily justified that the above right-hand side monotonically
decreases with respect to $i=1,2,\ldots,k$, whose minimum
attains at $i=k$ and equals $\left(\frac{k+1}{k}\right)^{\alpha}-1$.
Furthermore, since $\left(\frac{k+1}{k}\right)^{\alpha}-1$ decreases
monotonically as $k$ increases, its minimum over $k=1,2,\ldots,k_0$
is $\left(\frac{k_0+1}{k_0}\right)^{\alpha}-1$,
which is just the condition \eqref{condm}.
\qquad\endproof

\begin{remark}
Similar to \eqref{condition1},
there must be $\alpha>1$ such that \eqref{condm} holds. Again, we stress
that the conditions on $\rho$ and $\alpha$ should be regarded as an indicator that
$\rho$ and $\alpha$ must not be close to one other than precise requirements
since we have used the amplified \eqref{simp} and the simplified
models $\sigma_i=\zeta \rho^{-i}$ and $\sigma_i=\zeta i^{-\alpha}$.
Comparing Theorem~\ref{nearapprox} with Theorem~\ref{ritzvalue}, we
find out that, as far as the severely or moderately ill-posed problems are
concerned, for $k=1,2,\ldots,k_0$ the near best rank approximation
$P_{k+1}B_kQ_k^T$ essentially means that the singular values
$\theta_i^{(k)}$ of $B_k$ approximate the first $k$
large singular values $\sigma_i$ of $A$ in natural order, provided that
$\rho>1$ or $\alpha>1$ suitably.
\end{remark}

\begin{remark}
Under the conditions of Theorems~\ref{nearapprox}--\ref{ritzvalue},
let us explore how the results in them
depend on $\|\sin\Theta(\mathcal{V}_k,\mathcal{V}_k^R)\|$.
\eqref{etak0} and \eqref{case1} indicate that,
for the severely ill-posed problems with $\sigma_k=\zeta\rho^{-k}$,
ignoring higher order small terms,
we have $\eta_k\leq\rho^{-1-\beta}$ and $\|\Delta_k\|\leq \rho^{-2-\beta}<1$
for $k\leq k_0$; for the moderately ill-posed problems with
$\sigma_k=\zeta k^{-\alpha}$, \eqref{etadelta} indicates that $\eta_k$
and $\|\Delta_k\|$ are comparable in size for $k\leq k_0$,
while \eqref{modera2} shows that $\|\Delta_k\|$ is
at most of modest size for $k\leq k_0$. As a result, Theorem~\ref{thm2}
and Theorem~\ref{moderate} demonstrate that
$\|\sin\Theta(\mathcal{V}_k,\mathcal{V}_k^R)\|<\frac{1}{\sqrt{2}}$ and
$\|\sin\Theta(\mathcal{V}_k,\mathcal{V}_k^R)\|<1$ fairly for severely and
moderately ill-posed problems, respectively. In other words,
the largest canonical angle between
$\mathcal{V}_k^R$ and $\mathcal{V}_k$ does not exceed $\frac{\pi}{4}$
and is considerably smaller than $\frac{\pi}{2}$
for these two kinds of problems and $k\leq k_0$, respectively.
\end{remark}

\begin{remark}\label{extract}
Theorems~\ref{main1}--\ref{ritzvalue} show that, for $k=1,2,\ldots,k_0$,
the $k$-step Lanczos bidiagonalization is guaranteed to extract or acquire
the first $k$ dominant SVD components for the severely
or moderately ill-posed problems with $\rho>1$ or $\alpha>1$ suitably,
so that LSQR has the full regularization for these two kinds of ill-posed
problems and can obtain best possible regularized solutions $x^{(k_0)}$ at
semi-convergence.
\end{remark}

Let us have a closer look at the regularization of LSQR for mildly ill-posed
problems. We observe that the sufficient condition \eqref{condm} for \eqref{error2}
is never met for this kind of problem because
$
\left(\frac{k_0+1}{k_0}\right)^{\alpha}\leq 2
$
for any $k_0$ and $\frac{1}{2}< \alpha\leq 1$. This indicates that,
for $k=1,2,\ldots,k_0$, the $k$ Ritz values $\theta_i^{(k)}$ may not approximate
the the first $k$ large singular values $\sigma_i$ in natural order
and particularly there is at least one Ritz value
$\theta_{k_0}^{(k_0)}<\sigma_{k_0+1}$, causing that $x^{(k_0)}$ is already
deteriorated and cannot be as accurate as the best TSVD solution $x_{k_0}^{tsvd}$,
so that LSQR has only the partial regularization.
We can also make use of Theorem~\ref{initial} to explain the partial
regularization of LSQR: Theorem~\ref{moderate} has shown that
$\|\Delta_k\|$ is generally not small and
may become arbitrarily large as $k$ increases up to $k_0$ for mildly ill-posed
problems, meaning that $\|\sin\Theta(\mathcal{V}_k,\mathcal{V}_k^R)\|\approx 1$,
as the sharp bound \eqref{modera2} indicates, from which it follows that
a small Ritz value $\theta_{k_0}^{(k_0)}<\sigma_{k_0+1}$ generally appears.

\subsection{General best or best rank $k$ approximations to $A$ and their
implications on LSQR}
\label{morerank}

We investigate the general best or near best rank $k$ approximations to $A$ with
$\sigma_k=\zeta k^{-\alpha}$ and $\alpha>\frac{1}{2}$.
We aim to show that, for each of such rank
$k$ approximations, its smallest
nonzero singular value may be smaller than $\sigma_{k+1}$ for
$\frac{1}{2}<\alpha\leq 1$, that is, its nonzero singular values
may not approximate the $k$ large singular values of $A$ in natural
order, while the smallest nonzero singular value
of such a rank $k$ approximation is guaranteed to be bigger
than $\sigma_{k+1}$ if
only $\alpha>1$ suitably. As it will turn out, this can help
us further understand the regularization of LSQR for mildly and moderately
ill-posed problems. Finally, we investigate the behavior of
the Ritz values $\theta_i^{(k)},\ i=1,2,\ldots,k$ when $P_{k+1}B_kQ_k^T$
is not a near best rank $k$ approximation to $A$ for mildly ill-posed problems.

First of all, we point out an intrinsic fact that both the best and near best
rank $k$ approximations to $A$ with respect to the 2-norm are {\em not} unique.
This fact is important for further understanding Theorem~\ref{ritzvalue}.

Let $C_k$ be a best or near best rank
$k$ approximation to $A$ with $\|A-C_k\|=(1+\epsilon)\sigma_{k+1}$ with any
$\epsilon\geq 0$ satisfying $(1+\epsilon)\sigma_{k+1}<\frac{\sigma_k+
\sigma_{k+1}}{2}$ (Note: $\epsilon=0$ corresponds to a best rank $k$ approximation),
i.e., $(1+\epsilon)\sigma_{k+1}$ is between $\sigma_{k+1}$
and $\sigma_k$ and closer to $\sigma_{k+1}$, by which we get
$$
1+2\epsilon<\frac{\sigma_k}{\sigma_{k+1}}.
$$
It is remarkable to note that $C_k$ is not unique. For example, among others,
all the
$$
C_k=A_k(\theta,j)=A_k-\sigma_{k+1}U_k
{\rm diag}(\theta(1+\epsilon),\ldots,\theta(1+\epsilon),
\underbrace{(1+\epsilon)}_j,\theta(1+\epsilon),\ldots,\theta(1+\epsilon))
V_k^T
$$
with any $0\leq\theta\leq 1$ and $1\leq j\leq k-1$ is a family of
best or near best rank $k$ approximations to $A$.
The smallest nonzero singular value of $A_k(\theta,j)$ is
$\sigma_k-\theta(1+\epsilon)\sigma_{k+1}$. Since $\sigma_k=\zeta k^{-\alpha}$ and
$
\left(\frac{k+1}{k}\right)^{\alpha}<2
$
for any $k>1$ and $\frac{1}{2}<\alpha\leq 1$,
we obtain
\begin{equation}\label{mildmoderate2}
\sigma_k-\theta(1+\epsilon)\sigma_{k+1}=\sigma_{k+1}
\left(\left(\frac{k+1}{k}\right)^{\alpha}-\theta(1+\epsilon)\right)<\sigma_{k+1}
\end{equation}
for $\theta$ sufficiently close to one.
This shows that $\sigma_k-\theta(1+\epsilon)\sigma_{k+1}$ does not lie between
$\sigma_{k+1}$ and $\sigma_k$ and interlace them for $k>1$. In this case, for a
given $\alpha\in (\frac{1}{2},1]$, the bigger $k$ is, the
smaller $\left(\frac{k+1}{k}\right)^{\alpha}-\theta(1+\epsilon)$ is, and the
further is $\sigma_k-\theta(1+\epsilon)\sigma_{k+1}$ away from $\sigma_{k+1}$.
On the other hand, for $\theta$ sufficiently small we
have
\begin{equation}\label{mildmoderate}
\left(\frac{k+1}{k}\right)^{\alpha}-\theta(1+\epsilon)>1,
\end{equation}
that is, $\sigma_k-\theta(1+\epsilon)\sigma_{k+1}$ interlaces
$\sigma_{k+1}$ and $\sigma_k$ for $\theta$
sufficiently small.

For $A$ with $\sigma_k=\zeta k^{-\alpha}$ and $\alpha>1$, the situation is
much better since, for any $k$, the requirement \eqref{mildmoderate} is met
for {\em any} $0\leq\theta\leq 1$ provided that $\alpha>1$ suitably,
leading to $\sigma_k-\theta(1+\epsilon)\sigma_{k+1}>\sigma_{k+1}$, meaning
that the smallest singular value
$\sigma_k-\theta(1+\epsilon)\sigma_{k+1}$ of a near best rank approximation
$A_k(\theta,j)$ interlaces $\sigma_{k+1}$ and $\sigma_k$.

However, we should be aware that the above analysis is made
for the {\em worst-case}: For any best or a
near best rank $k$ approximation $C_k$ to $A$,
the minimum of the smallest nonzero singular values of all the $C_k$ is exactly
$\sigma_k-(1+\epsilon)\sigma_{k+1}$. We now prove this. Suppose that
$\sigma_k(C_k)$ is the smallest nonzero singular value of a given such $C_k$.
Then from $\|A-C_k\|=(1+\epsilon)\sigma_{k+1}$,
by the standard perturbation theory we have
$$
|\sigma_k-\sigma_k(C_k)|\leq (1+\epsilon)\sigma_{k+1}.
$$
Clearly, the minimum of all the $\sigma_k(C_k)$ is attained if and
only if the above equality holds, which is exactly
$\sigma_k-(1+\epsilon)\sigma_{k+1}$. On the other side, by
construction, we also see that the smallest singular value
$\sigma_k-\theta(1+\epsilon)\sigma_{k+1}$
of $C_k$ is arbitrarily close to or equal to $\sigma_k$ by
taking $\theta$ arbitrarily small or zero, which means that
\eqref{mildmoderate} holds. In this case, we observe from the equality
in \eqref{mildmoderate2} that $\sigma_k-\theta(1+\epsilon)\sigma_{k+1}
>\sigma_{k+1}$ and interlaces $\sigma_{k+1}$ and $\sigma_k$.

As far as LSQR is concerned, notice that the condition \eqref{condm} for the
interlacing property \eqref{error2} is derived by assuming the worst case that
$\sigma_k-\theta_k^{(k)}=\gamma_k\leq\sqrt{1+\eta_k^2}\sigma_{k+1}$, i.e.,
$\theta_k^{(k)}$ is supposed to be the smallest possible nonzero one among
all the $\sigma_k(C_k)$, where $C_k$ belongs to the set of near best
$k$ approximations that satisfy $\|A-C_k\|=\gamma_k\leq\sqrt{1+\eta_k^2}
\sigma_{k+1}$. For mildly ill-posed problems, the above arguments indicate
that although in the worst case some of the $k$ Ritz values $\theta_i^{(k)}$
may not approximate the first $k$ large singular values $\sigma_i$
of $A$ in natural order, it is possible so in practice
in case $P_{k+1}B_kQ_k^T$ is {\em occasionally} a near best
rank $k$ approximation to $A$ for some small $k\leq k_0$.

Unfortunately, as we have shown previously, $P_{k+1}B_kQ_k^T$ is
rarely a near best rank $k$ approximation to $A$ for
mildly ill-posed problems, i.e., $\gamma_k>\sigma_k$ generally.
Recall the second part of Theorem~\ref{initial} and Remark~\ref{appear},
which have shown rigorously that there is at least one Ritz value
$\theta_k^{(k)}<\sigma_{k+1}$ if $\varepsilon_k$
is sufficiently small there, that is, $\eta_k$ or equivalently $\|\Delta_k\|$
is large. This is exactly the case that $P_{k+1}B_kQ_k^T$ is not a near best
rank $k$ approximation to $A$, causing that LSQR has only the partial
regularization.

We can make a further analysis on the behavior of
$\theta_i^{(k)},\ i=1,2,\ldots,k$ when \eqref{eq1} is mildly ill-posed.
Suppose that $\gamma_k\in [\sigma_{j+1},\sigma_j]$ for some $j\leq k$, which
means that $P_{k+1}B_kQ_k^T$ is definitely not a near best rank $k$
approximation to $A$ when $j<k$. Below we derive the smallest
upper bound for $\sigma_j-\gamma_k$ and obtain the biggest
lower bound for $\theta_j^{(k)}$. For $\sigma_j=\zeta j^{-\alpha}$
with $\frac{1}{2}<\alpha\leq 1$ we have
\begin{align*}
\sigma_j-\gamma_k&\leq \sigma_j-\sigma_{j+1}
=\sigma_{j+1}\left(\left(\frac{j+1}{j}\right)^{\alpha}-1\right)\\
&\leq \sigma_{j+1}\left(1+\frac{\alpha}{j}-1\right)\\
&=\frac{\alpha}{j}\sigma_{j+1}=\frac{\alpha}{\zeta j^{1-\alpha}}\zeta
j^{-\alpha}\sigma_{j+1}
=\frac{\alpha}{\zeta j^{1-\alpha}}\zeta^2 j^{-\alpha} (j+1)^{-\alpha}\\
&=
\frac{\alpha}{\zeta j^{1-\alpha}}\zeta \sigma_{j(j+1)}=
\frac{\alpha}{j^{1-\alpha}}\sigma_{j(j+1)},
\end{align*}
in which $\frac{\alpha}{j^{1-\alpha}}<1$
decreases with increasing $j$ for $\alpha<1$ and is one for $\alpha=1$.
Therefore, the smallest
upper bound for $\sigma_j-\gamma_k$ is no more than $\sigma_{j(j+1)}$,
which is smaller than $\sigma_{k+1}$ once $j(j+1)>k$. In view of the above
and \eqref{ritzapp}, for $\gamma_k\in [\sigma_{j+1},\sigma_j]$, since
$\theta_j^{(k)}\geq\sigma_j-\gamma_k$ and has the biggest lower bound
$\sigma_{j(j+1)}$, we may have $\theta_j^{(k)}<\sigma_{k+1}$ provided that
$j(j+1)>k$. Moreover, when $\theta_j^{(k)}<\sigma_{k+1}$,
by the labeling rule, there are $k-j+1$ Ritz values
$\theta_j^{(k)},\theta_{j+1}^{(k)},\ldots,\theta_k^{(k)}$ smaller
than $\sigma_{k+1}$. As a result, for $k=k_0$, there are $k_0-j+1$
Ritz values smaller than $\sigma_{k_0+1}$ that deteriorate the LSQR
iterate $x^{(k_0)}$, so that LSQR has only the partial regularization.

\section{Decay rates of $\alpha_k$ and $\beta_{k+1}$ and the regularization of
LSMR and CGME} \label{alphabeta}

In this section, we will present a number of results on the decay rates of
$\alpha_k,\ \beta_{k+1}$ and $\gamma_k$ and on certain other rank
$k$ approximations to $A$ and $A^TA$ constructed by Lanczos bidiagonalization.
The decay rates of $\alpha_k$ and $\beta_{k+1}$ are
particularly useful for practically detecting the degree of
ill-posedness of \eqref{eq1} and identifying the full or partial regularization
of LSQR and LSMR. The results on the new rank $k$
approximations critically determine the full or partial regularization of the
Krylov iterative regularization solvers LSMR \cite{fong} and CGME
\cite{craig,hanke95,hanke01,hps09}. In Section \ref{decay}, we
prove how $\alpha_k$ and $\beta_{k+1}$ decay by relating them to
$\gamma_k$ and the estimates established for it. Then we show how to
exploit the decay rate of $\alpha_k+\beta_{k+1}$ to identify
the degree of ill-posedness of \eqref{eq1} and the regularization of LSQR.
In Section \ref{lsmr}, we prove that the regularization of LSMR resembles LSQR for each
of the three kinds of ill-posed problems. In Section \ref{cgme}, we prove that
the regularizing effects of CGME have intrinsic indeterminacy and are inferior
to those of LSQR and LSMR. In Section \ref{rrqr},
we compare LSQR with some standard randomized algorithms \cite{halko11} and
strong rank-revealing QR, i.e., RRQR, factorizations \cite{gu96,hong92},
and show that the former solves
ill-posed problems more accurately than the latter two ones at no more cost.

\subsection{Decay rates of $\alpha_k$ and $\beta_{k+1}$ and their practical
use}\label{decay}

We consider how $\alpha_k$ and $\beta_{k+1}$ decay in certain pronounced
manners and show how to use them to identify the full or partial regularization
of LSQR in practice.

\begin{theorem}\label{main2}
With the notation defined previously, the following results hold:
\begin{eqnarray}
  \alpha_{k+1}&<&\gamma_k\leq \sqrt{1+\eta_k^2}\sigma_{k+1},
  \ k=1,2,\ldots,n-1,\label{alpha}\\
 \beta_{k+2}&<& \gamma_k\leq\sqrt{1+\eta_k^2}\sigma_{k+1}, \ k=1,2,\ldots,n-1,
 \label{beta}\\
\alpha_{k+1}\beta_{k+2}&\leq &
\frac{\gamma_k^2}{2}\leq
   \frac{(1+\eta_k^2)\sigma_{k+1}^2}{2}, \ k=1,2,\ldots,n-1,
  \label{prod2}\\
\gamma_{k+1}&<&\gamma_k,\  \ k=1,2,\ldots,n-2. \label{gammamono}
\end{eqnarray}
\end{theorem}

{\em Proof}.
From \eqref{fulllb}, since $P$ and $Q_n$ are orthogonal matrices,
we have
\begin{align}
\gamma_k &=\|A-P_{k+1}B_kQ_k^T\|=\|P^T(A-P_{k+1}B_kQ_k^T)Q_n\| \label{invar}\\
&=
\left\| \left(\begin{array}{c}
B_n \\
\mathbf{0}
\end{array}
\right)-(I,\mathbf{0} )^TB_k (I,\mathbf{0} )\right\|=\|G_k\| \label{gk}
\end{align}
with
\begin{align}\label{gk1}
G_k&=\left(\begin{array}{cccc}
\alpha_{k+1} & & & \\
\beta_{k+2}& \alpha_{k+2} & &\\
&
  \beta_{k+3} &\ddots & \\& & \ddots & \alpha_{n} \\
  & & & \beta_{n+1}
  \end{array}\right)\in \mathbb{R}^{(n-k+1)\times (n-k)}
\end{align}
resulting from deleting the $(k+1)\times k$ leading principal matrix of $B_n$
and the first $k$ zero rows and columns of the resulting matrix.
From the above, for $k=1,2,\ldots,n-1$ we have
\begin{align}
\alpha_{k+1}^2+\beta_{k+2}^2&=\|G_ke_1\|^2\leq \|G_k\|^2=\gamma_k^2,
\label{alphabetasum1}
\end{align}
which shows that $\alpha_{k+1}< \gamma_k$ and $\beta_{k+2}<\gamma_k$
since $\alpha_{k+1}>0$ and $\beta_{k+2}>0$.
So from \eqref{final}, we get \eqref{alpha} and \eqref{beta}. On the other
hand, noting that
\begin{align*}
2\alpha_{k+1}\beta_{k+2}&\leq \alpha_{k+1}^2+\beta_{k+2}^2\leq \gamma_k^2,
\end{align*}
we get \eqref{prod2}.

Note that $\alpha_k>0$ and $\beta_{k+1}>0,\ k=1,2,\ldots,n$.
By $\gamma_k=\|G_k\|$ and \eqref{gk1}, note that $\gamma_{k+1}=\|G_{k+1}\|$
equals the 2-norm of the submatrix deleting the first column of $G_k$.
Applying the Cauchy's strict interlacing theorem to the singular values
of this submatrix and $G_k$, we obtain \eqref{gammamono}.
\qquad\endproof

\begin{remark}
For severely and moderately ill-posed problems, based on
the results in the last section,
\eqref{alpha} and \eqref{beta} show that $\alpha_{k+1}$ and $\beta_{k+2}$
decay as fast as $\sigma_{k+1}$ for $k\leq k_0$ and their decays may become slow
for $k>k_0$. For mildly ill-posed problems,
since $\eta_k$ are generally bigger than one considerably for $k\leq k_0$,
$\alpha_{k+1}$ and $\beta_{k+2}$ cannot generally decay as fast as
$\sigma_{k+1}$, and their decays become slower for $k>k_0$.

Gazzola and his
coauthors \cite{gazzola14,gazzola-online} claim without rigorous proofs
that $\alpha_{k+1}\beta_{k+1}=\mathcal{O}(k\sigma_k^2)$ and
$\alpha_{k+1}\beta_{k+2}=\mathcal{O}(k\sigma_{k+1}^2)$  for severely ill-posed
problems with the constants in $\mathcal{O}(\cdot)$ unknown
(see Proposition 4 of \cite{gazzola-online}), but they do not show how
fast each of them decays; see Proposition 6 of \cite{gazzola-online}.
In contrast, our \eqref{alpha}, \eqref{beta} and \eqref{prod2} are
rigorous and quantitative for all three kinds of ill-posed problems.
In \cite[Corollary 3.1]{gazzola16}, the authors have derived
the product inequality
$$
\prod_{k=1}^l\alpha_{k+1}\beta_{k+1}\leq \prod_{k=1}^l\sigma_k^2,\
l=1,2,\ldots,n-1.
$$
Whether or not this inequality is sharp is unknown, as they point out.
By it, they empirically claim that $\alpha_{k+1}\beta_{k+1}$
may decay as fast as $\sigma_k^2$ when the inequality is sharp; conversely,
if it is not sharp, nothing can be said on how fast
$\alpha_{k+1}\beta_{k+1}$ decays.
\end{remark}

We now shed light on \eqref{alpha} and \eqref{beta}.
For a given \eqref{eq1}, its degree of ill-posedness
is either known or unknown. If it is unknown, \eqref{alpha} is of
practical importance and can be exploited to identify whether or not LSQR has
the full regularization without extra cost in an automatic and
reliable way, so is \eqref{beta}.
From the proofs of \eqref{alpha} and \eqref{beta}, we find that
$\alpha_{k+1}$ and $\beta_{k+2}$ are as small as $\gamma_k$. Since our theory
and analysis in Section \ref{rankapp} have proved that $\gamma_k$
decays as fast as $\sigma_{k+1}$ for severely or moderately ill-posed problems
with $\rho>1$ or $\alpha>1$ suitably and it decays more slowly than
$\sigma_{k+1}$ for mildly il-posed problems, the decay rate
of $\sigma_k$ can be judged by that of $\alpha_k$
or $\beta_{k+1}$ or better judged by that of $\alpha_k+\beta_{k+1}$
reliably, as shown below.

Given \eqref{eq1}, run LSQR until semi-convergence
occurs at iteration $k^*$. Check how $\alpha_k+\beta_{k+1}$ decays as $k$
increases during the process. If, on average, it decays in
an obviously exponential way, then \eqref{eq1} is a severely ill-posed problem.
In this case, LSQR has the full regularization, and semi-convergence means
that we have found a best possible regularized solution. If, on average,
$\alpha_k$ decays as fast as $k^{-\alpha}$ with $\alpha>1$ considerably, then
\eqref{eq1} is surely a moderately ill-posed problem, and LSQR also has found a
best possible regularized solution at semi-convergence. If, on average,
it decays at most as fast as or more slowly than $k^{-\alpha}$
with $\alpha$ no more than one, \eqref{eq1} is a mildly ill-posed problem.
Notice that the noise $e$ does not deteriorate regularized solutions until
semi-convergence. Therefore, if a hybrid LSQR is used, then it is more reasonable
and also cheaper to apply regularization to projected problems only
from iteration $k^*+1$ onwards
other than from the very start, i.e., the first iteration, as done in the hybrid
Lanczos bidiagonalization/Tikhonov regularization scheme \cite{berisha}, until
a best possible regularized solution
is found. For a hybrid LSMR, regularization is applied to the projected
problems generated in LSMR in the same way.

\subsection{The regularization of LSMR}\label{lsmr}

Based on the previous results, we can rigorously analyze the
regularizing effects of of LSMR \cite{fong,bjorck15} and
draw definitive conclusions on its regularization for three kinds of
ill-posed problems.

LSMR is mathematically equivalent to MINRES applied to
$A^TAx=A^Tb$, and its iterate $x_k^{lsmr}$ minimizes
$\|A^T(b-Ax)\|$ over $x\in \mathcal{V}_k^R$,
and the residual norm $\|A^T(b-Ax_k^{lsmr})\|$ decreases
monotonically with respect to $k$. In our notation, noting from Algorithm 1
that $Q_{k+1}^TA^TAQ_k=(B_k^TB_k,\alpha_{k+1}\beta_{k+1}e_k)^T$ with rank $k$,
it is known from Section 2.2 of \cite{fong} that
\begin{equation}\label{lsmrsolution}
x_k^{lsmr}=Q_ky_k^{lsmr}=Q_k(Q_{k+1}^TA^TAQ_k)^{\dagger}Q_{k+1}^TA^Tb,
\end{equation}
which can be efficiently computed and updated.
So LSMR amounts to solving the modified problem that perturbs the matrix $A^TA$
in $A^TAx=A^Tb$ to its rank $k$ approximation $Q_{k+1}Q_{k+1}^TA^TAQ_kQ_k^T$,
and the iterate $x_k^{lsmr}$ is the minimum-norm least squares solution to
the modified problem
\begin{equation}\label{lsmrrank}
\min\|Q_{k+1}Q_{k+1}^TA^TAQ_kQ_k^Tx-A^Tb\|.
\end{equation}
It is direct to verify that the TSVD solution $x_k^{tsvd}$ is exactly
the minimum-norm least squares solution to the modified
problem $\min\|A_k^TA_kx-A^Tb\|$ that replaces
$A^TA$ by its 2-norm best rank $k$ approximation $A_k^TA_k$ in $A^TAx=A^Tb$.
As a result, the regularization problem for LSMR now becomes that of accurately
estimating $\|A^TA-Q_{k+1}Q_{k+1}^TA^TAQ_kQ_k^T\|$,
investigates how close it is to $\sigma_{k+1}^2=\|A^TA-A_k^TA_k\|$ and analyzes
whether or not the singular values of $Q_{k+1}^TA^TAQ_k$ approximate
the $k$ large singular values $\sigma_i^2,\ i=1,2,\ldots,k$ of $A^TA$ in natural
order.

\begin{theorem}\label{aprod}
For LSMR and $k=1,2,\ldots,n-1$, we have
\begin{equation}\label{aproderror}
\gamma_k^2\leq \|A^TA-Q_{k+1}Q_{k+1}^TA^TAQ_kQ_k^T\|\leq
\sqrt{1+m_k(\gamma_{k-1}/\gamma_k)^2}\gamma_k^2
\end{equation}
with $0\leq m_k<1$.
\end{theorem}

{\em Proof}.
For the orthogonal matrix $Q_n$ generated by Algorithm 1, noticing that
$\alpha_{n+1}=0$, from \eqref{eqmform1} and \eqref{eqmform2} we obtain
$Q_n^TA^TAQ_n=B_n^TB_n$ and
\begin{align*}
\|A^TA-Q_{k+1}Q_{k+1}^TA^TAQ_kQ_k^T\|&=
\|Q_n^T(A^TA-Q_{k+1}Q_{k+1}^TA^TAQ_kQ_k^T)Q_n\|\\
&=\|B_n^TB_n-(I,\mathbf{0})^T(B_k^TB_k,\alpha_{k+1}
\beta_{k+1}e_k)^T(I,\mathbf{0})\|\\
&=\|F_k\|,
\end{align*}
where $F_k$ is the $(n-k+1)\times (n-k)$ matrix that is generated by deleting
the $(k+1)\times k$ leading principal matrix of the symmetric tridiagonal
matrix $B_n^TB_n$ and the first $k-1$ zero rows and $k$ zero columns
of the resulting matrix. Note that $B_n^TB_n$ has the diagonals
$\alpha_k^2+\beta_{k+1}^2$, $k=1,2,\ldots,n$ and the super- and sub-diagonals
$\alpha_k\beta_k,\ k=2,3,\ldots,n$. We have
\begin{align}\label{fk}
F_k&=\left(\begin{array}{ccccc}
\alpha_{k+1}\beta_{k+1} & & & &\\
\alpha_{k+1}^2+\beta_{k+2}^2 &\alpha_{k+2}\beta_{k+2} &&&\\
\alpha_{k+2}\beta_{k+2} &\alpha_{k+2}^2+\beta_{k+3}^2&\ddots & &\\
& \alpha_{k+3}\beta_{k+3} & \ddots& & \\
& & &\alpha_{n-1}\beta_{n-1}& \\
 & & \ddots& \alpha_{n-1}^2+\beta_n^2 &\alpha_n\beta_n\\
  & & &\alpha_n\beta_n&\alpha_n^2+\beta_{n+1}^2\\
\end{array}
\right)
\end{align}
According to \eqref{gk1}, it is direct to check that
the $(n-k)\times (n-k)$ symmetric tridiagonal matrix $G_k^TG_k$ is
a submatrix of $F_k$ that deletes its first row $\alpha_{k+1}\beta_{k+1}e_1^T$.
Therefore, we have $\|F_k\|\geq\|G_k^TG_k\|=\|G_k\|^2=\gamma_k^2$
with the last equality being from \eqref{invar} and \eqref{gk}, which proves
the lower bound in \eqref{aproderror}.

On the other hand, noting the strict inequalities
in \eqref{alpha} and \eqref{beta}, since
$$
F_k^TF_k=(G_k^TG_k)^2+\alpha_{k+1}^2\beta_{k+1}^2e_1e_1^T,
$$
from \cite[p.98]{wilkinson} we obtain
$$
\|F_k\|^2=\|G_k\|^4+m^{\prime}_k\alpha_{k+1}^2\beta_{k+1}^2\leq \gamma_k^4+
m_k\gamma_{k-1}^2\gamma_k^2
$$
with $0\leq m^{\prime}_k\leq 1$ and $0\leq m_k<m^{\prime}_k$ if
$m^{\prime}_k>0$, from which the upper bound of \eqref{aproderror} follows
directly.
\qquad\endproof

Recall that LSQR is mathematically equivalent to CGLS that implicitly applies
the CG method to $A^TAx=A^Tb$. By \eqref{xk}, \eqref{eqmform1},
\eqref{eqmform2} and \eqref{Bk}, noting that $P_{k+1}P_{k+1}^Tb=b$,
we obtain the LSQR iterates
\begin{align*}
x^{(k)}&=Q_kB_k^{\dagger}P_{k+1}^Tb=Q_k(B_k^TB_k)^{-1}B_k^TP_{k+1}^Tb\\
&=Q_k(Q_k^TA^TAQ_k)^{-1}Q_k^TA^TP_{k+1}P_{k+1}^Tb=
Q_k(Q_k^TA^TAQ_k)^{-1}Q_k^TA^Tb,
\end{align*}
which is the minimum-norm least squares solution to the modified problem
\begin{equation}\label{lsqrrank}
\min\|Q_kQ_k^TA^TAQ_kQ_k^T x-A^Tb\|
\end{equation}
that replaces $A^TA$ by its rank $k$ approximation $Q_kQ_k^TA^TAQ_kQ_k^T
=Q_kB_k^TB_kQ_k^T$ in $A^TAx=A^Tb$. As a result, in the sense of solving
$A^TAx=A^Tb$, for LSQR, the accuracy of such rank $k$ approximation
is $\|A^TA-Q_kQ_k^TA^TAQ_kQ_k^T\|$. We can establish the following result which
relates the corresponding approximation accuracy concerning LSMR to
that concerning LSQR.

\begin{theorem}\label{lsqrmr}
For the rank $k$ approximations to $A^TA$ defined in \eqref{lsmrrank}
and \eqref{lsqrrank} involved in LSMR and LSQR, we have
\begin{align}\label{lsqrmrest}
\|A^TA-Q_{k+1}Q_{k+1}^TA^TAQ_kQ_k^T\|&\leq
\|A^TA-Q_kQ_k^TA^TAQ_kQ_k^T\|.
\end{align}
\end{theorem}

{\proof} Similar to the proof of Theorem~\ref{aprod}, it is direct to verify
that
\begin{align*}
\|A^TA-Q_kQ_k^TA^TAQ_kQ_k^T\|&=
\|Q_n^T(A^TA-Q_kQ_k^TA^TAQ_kQ_k^T)Q_n\|\\
&=\|B_n^TB_n-(I,\mathbf{0})^TB_k^TB_k(I,\mathbf{0})\|\\
&=\|F_k^{\prime}\|,
\end{align*}
where $F_k^{\prime}$ is an $(n-k+1)\times (n-k+1)$ matrix whose
first column is $(0,\alpha_{k+1}\beta_{k+1},\mathbf{0})^T$ and last $n-k$
columns are just the matrix $F_k$ defined by \eqref{fk}. Therefore, we have
$\|F_k\|\leq \|F_k^{\prime}\|$, which is just \eqref{lsqrmr}.
\qquad\endproof

This theorem indicates that, as far as solving $A^TAx=A^Tb$
is concerned, the rank $k$ approximation in LSMR is at least as accurate as
that in LSQR. However, regarding LSQR applied to \eqref{eq1} directly,
Theorem~\ref{main1} is much more attractive since it not only deals with
the rank $k$ approximation to $A$ directly but also estimates the accuracy of
the rank $k$ approximation in terms
of $\sigma_{k+1}$ more compactly and informatively.

\begin{remark}
According to the results and analysis in Section \ref{rankapp},
we have $\gamma_{k-1}/\gamma_k\sim \rho$ for severely ill-posed problems,
and $\gamma_{k-1}/\gamma_k\sim (k/(k-1))^{\alpha}$ at most for
moderately and mildly ill-posed problems.
In comparison with Theorem~\ref{main1}, noting the form of
the lower and upper bounds of \eqref{aproderror}, we see that
$Q_{k+1}Q_{k+1}^TA^TAQ_kQ_k^T=Q_{k+1}(B_k^TB_k,\alpha_{k+1}\beta_{k+1}e_k)^TQ_k^T$
as a rank $k$ approximation to $A^TA$
is basically as accurate as $P_{k+1}B_kQ_k^T$ as a rank $k$ approximation
to $A$.
\end{remark}

\begin{remark}
From \cite[p.33]{stewartsun}, the singular values of
$(B_k^TB_k,\alpha_{k+1}\beta_{k+1}e_k)^T$ are correspondingly bigger
than those of $B_k^TB_k$, i.e., $(\theta_i^{(k)})^2$. Therefore,
the smallest singular value of $(B_k^TB_k,\alpha_{k+1}\beta_{k+1}e_k)^T$
is no less than $(\theta_k^{(k)})^2$.
As a result, $(B_k^TB_k,\alpha_{k+1}\beta_{k+1}e_k)^T$ has no
singular values smaller than $\sigma_{k_0+1}^2$ before $k\leq k_0$,
provided that $\theta_k^{(k)}>\sigma_{k_0+1}$ for
$k\leq k_0$. This means that the noise deteriorates the iterates
$x_k^{lsmr}$ no sooner than it does for the LSQR iterates $x^{(k)}$.
\end{remark}

\begin{remark}
A combination of Theorem~\ref{lsqrmr} and the above two remarks means that the
regularizing effects of LSMR are highly competitive with and not inferior to
those of LSQR for each kind of ill-posed problem under consideration. Consequently,
from the theory of LSQR in Section \ref{rankapp}, we conclude that LSMR has
the full regularization for severely or moderately ill-posed problems with
$\rho>1$ or $\alpha>1$ suitably. However, Theorem~\ref{aprod} indicates
that LSMR generally has only the partial regularization for mildly ill-posed
problems since $\gamma_{k_0}$ is generally bigger than $\sigma_{k_0+1}$
considerably; see Remark~\ref{mildre}.
\end{remark}

\begin{remark}
We can define a near best rank $k$ approximation to $A^TA$ similar to
\eqref{near}. Based on \eqref{aproderror}, if, in LSMR, we simply take
$\|A^TA-Q_{k+1}Q_{k+1}^TA^TAQ_kQ_k^T\|=\gamma_k^2$ for the ease of presentation,
we can establish an analog of Theorem~\ref{nearapprox} for LSMR. In the
meantime, completely parallel to the proof of Theorem~\ref{ritzvalue}, we can
also derive an analog of Theorem~\ref{ritzvalue} for LSMR, in which
the sufficient conditions on $\eta_k$ that ensure that the singular values
of $(B_k^TB_k,\alpha_{k+1}\beta_{k+1}e_k)^T$
approximate the first $k$ large singular values of $A^TA$ in natural order
are found to be
$$
2+\eta_k^2<\left(\frac{k_0+1}{k_0}\right)^{2\alpha},\ k=1,2,\ldots,k_0
$$
for $A$ with $\sigma_i=\zeta i^{-\alpha}$.
\end{remark}

\begin{remark}
Since LSMR and LSQR have similar regularizing effects for each kind of ill-posed
problem, we can judge the full or partial regularization
of LSMR by inspecting the decay rate of $\alpha_k+\beta_{k+1}$ with respect
to $k$, as has been done for LSQR,
\end{remark}

In Section \ref{rankapp} we have interpreted LSQR
as solving the modified problem that perturbs $A$ in \eqref{eq1}
to its rank $k$ approximation $P_{k+1}B_kQ_k^T$. The
regularization of LSQR then is up to the accuracy of such rank $k$
approximation to $A$ and how the $k$ large singular values of $A$ are
approximated by the nonzero singular values of $B_k$. We will treat CGME
in the same way later. It might be hopeful to treat LSMR in this preferable
and more direct way. From \eqref{lsmrsolution}, LSMR is also equivalent to
computing the minimum-norm least squares solution to the modified
problem
$$
\min\|\left(Q_k(Q_{k+1}^TA^TAQ_k)^{\dagger}
Q_{k+1}^TA^T\right)^{\dagger} x-b\|,
$$
which perturbs $A$ in \eqref{eq1} to its rank $k$ approximation
$\left(Q_k(Q_{k+1}^TA^TAQ_k)^{\dagger}Q_{k+1}^TA^T\right)^{\dagger}$.
However, an analysis of such formulation appears intractable because
there is no explicit way to
remove two generalized inverses $\dagger$ in such rank $k$
approximation, which makes it impossible to accurately estimate
$\|A-\left(Q_k(Q_{k+1}^TA^TAQ_k)^{\dagger}Q_{k+1}^TA^T\right)^{\dagger}\|$
in terms of $\sigma_{k+1}$.

\subsection{The other rank $k$ approximations to $A$ generated by Lanczos
bidiagonalization and the regularization of CGME}\label{cgme}

By \eqref{eqmform1} and \eqref{eqmform2}, we get
\begin{align}
P_{k+1}P_{k+1}^TA &= P_{k+1}(B_kQ_k^T+\alpha_{k+1}e_{k+1}q_{k+1}^T)\notag\\
&=P_{k+1}(B_k, \alpha_{k+1}e_{k+1})Q_{k+1}^T \notag\\
&=P_{k+1}\bar{B}_kQ_{k+1}^T, \label{leftbidiag}
\end{align}
where $Q_{k+1}=(Q_k,q_{k+1})$, and $\bar{B}_k=(B_k,\alpha_{k+1}e_{k+1})\in
\mathbb{R}^{(k+1)\times (k+1)}$ is lower bidiagonal with rank $k+1$.
Thus, it follows from \cite{hanke95,hanke01,hps09} that
CGME is the CG method applied to $\min\|AA^Ty-b\|$ and $x=A^Ty$, where
the $k$-th iterate $x_k^{cgme}$ minimizes the error $\|A^{\dagger}b-x\|$, i.e.,
$\|x_{naive}-x\|$, over $x\in \mathcal{V}_k^R$, and the error norm
$\|x_{naive}-x_k^{cgme}\|$ decreases monotonically with respect to $k$.
By Lanczos bidiagonalization, it is known from \cite{hanke95,hanke01,hps09} that
$x_k^{cgme}=Q_ky_k^{cgme}$ with $y_k^{cgme}=\|b\|\bar{B}_{k-1}^{-1}e_1^{(k)}$
and the residual norm $\|Ax_k^{cgme}-b\|=\beta_{k+1}|e_k^T y_k^{cgme}|$
with $e_k$ the $k$-th canonical vector of dimension $k$.
Noting that $\|b\|e_1^{(k)}=P_k^Tb$, we have
\begin{equation}\label{cgmesolution}
x_k^{cgme}=Q_k\bar{B}_{k-1}^{-1}P_k^Tb.
\end{equation}
Therefore, $x_k^{cgme}$ is the minimum-norm least squares solution to the
modified problem that replaces $A$ in \eqref{eq1} by its rank $k$
approximation $P_k\bar{B}_{k-1}Q_k^T=P_kP_k^TA$.

\begin{theorem}\label{cgmeappr}
For the rank $k+1$ approximation $P_{k+1}P_{k+1}^TA$ and the rank $k$ approximation
in CGME, we have
\begin{eqnarray}
  \|(I-P_{k+1}P_{k+1}^T)A\|
  &\leq& \gamma_k\leq\sqrt{1+\eta_k^2}\sigma_{k+1}, \label{leftsub}\\
\gamma_k<\|A-P_k\bar{B}_{k-1}Q_k^T\|&\leq& \gamma_{k-1}. \label{cgmelowup}
\end{eqnarray}
\end{theorem}

{\proof}
Since $P_{k+1}P_{k+1}^T(I-P_{k+1}P_{k+1}^T)=0$, we obtain
\begin{align}
  \gamma_k^2 &=\|A-P_{k+1}B_kQ_k^T\|^2 \notag\\
  &=\|P_{k+1}P_{k+1}^TA-P_{k+1}B_kQ_k^T+(I-P_{k+1}P_{k+1}^T)A\|^2\notag\\
  &=\max_{\|y\|=1}\|\left(
  (P_{k+1}P_{k+1}^TA-P_{k+1}B_kQ_k^T)+(I-P_{k+1}P_{k+1}^T)A\right)y\|^2\notag\\
   &=\max_{\|y\|=1} \|P_{k+1}P_{k+1}^T(P_{k+1}P_{k+1}^TA-
           P_{k+1}B_kQ_k^T)y+(I-P_{k+1}P_{k+1}^T)Ay\|^2 \notag\\
           &=\max_{\|y\|=1}\left(\|P_{k+1}P_{k+1}^T(P_{k+1}P_{k+1}^TA-
           P_{k+1}B_kQ_k^T)y\|^2+\|(I-P_{k+1}P_{k+1}^T)Ay\|^2\right)\notag\\
            &=\max_{\|y\|=1}\left(\|P_{k+1}(P_{k+1}^TA-B_kQ_k^T)y\|^2+
           \|(I-P_{k+1}P_{k+1}^T)Ay\|^2\right)\notag\\
           &=\max_{\|y\|=1}\left(\|(P_{k+1}^TA-B_kQ_k^T)y\|^2+
           \|(I-P_{k+1}P_{k+1}^T)Ay\|^2\right)\notag\\
            &\geq \max_{\|y\|=1} \|(I-P_{k+1}P_{k+1}^T)Ay\|^2 \notag\\
            &=\|(I-P_{k+1}P_{k+1}^T)A\|^2,\notag
\end{align}
which, together with \eqref{final}, establishes \eqref{leftsub}.

From \eqref{leftbidiag} and \eqref{leftsub} we obtain
\begin{equation}\label{leftapp}
\|(I-P_{k+1}P_{k+1}^T)A\|=\|A-P_{k+1}\bar{B}_kQ_{k+1}^T\|\leq \gamma_k.
\end{equation}

The upper bound of \eqref{cgmelowup} is direct \eqref{leftapp} and
\eqref{gammamono} by noting that
$$
\|A-P_k\bar{B}_{k-1}Q_k^T\|=\|(I-P_kP_k^T)A\|\leq\gamma_{k-1}.
$$
Along the proof path of Theorem~\ref{main2}, we obtain
$$
\|A-P_k\bar{B}_{k-1}Q_k^T\|=\|(\beta_{k+1}e_1,G_k)\|
$$
with $G_k$ defined by \eqref{gk1}. It is straightforward to justify
that the singular values of $G_k\in \mathbb{R}^{(n-k+1)\times (n-k)}$ strictly
interlace those of $(\beta_k e_1,G_k)\in \mathbb{R}^{(n-k+1)\times (n-k+1)}$
by noting that $(\beta_{k+1}e_1,G_k)^T (\beta_{k+1}e_1,G_k)$ is
an {\em unreduced} symmetric tridiagonal matrix, from which and
$\|G_k\|=\gamma_k$ (cf. \eqref{invar} and \eqref{gk}) the
lower bound of \eqref{cgmelowup} follows.
\qquad\endproof

By the definition \eqref{gammak} of $\gamma_k$, this theorem indicates that
$P_k\bar{B}_{k-1}Q_k^T$ is definitely a less accurate
rank $k$ approximation to $A$ than $P_{k+1}B_kQ_k^T$ in LSQR.
Moreover, a combination of it and Theorem~\ref{main1} indicates that
$P_k\bar{B}_{k-1}Q_k^T$ may never be a near best rank $k$ approximation
to $A$ even for severely and moderately ill-posed problems because,
unlike LSQR, there do not exist sufficient conditions on
$\rho>1$ and $\alpha>1$ to meet this requirement.
For mildly ill-posed problems, CGME generally has only the partial
regularization since $\gamma_k$ has been proved to be generally bigger than
$\sigma_{k+1}$ substantially and is rarely close to $\sigma_{k+1}$.

Next we consider the other issue that is as equally important as
the rank $k$ approximation in CGME: the behavior
of the singular values of $\bar{B}_{k-1}$, which are denoted by
$\bar{\theta}_i^{(k-1)},\ i=1,2,\ldots,k$
labeled in the decreasing order. Observe that
$\bar{B}_{k-1}$ consists of the first $k$ rows of $B_k$. Since $B_k B_k^T$ is
an $(k+1)\times (k+1)$ unreduced symmetric tridiagonal matrix, whose eigenvalues
are $(\theta_1^{(k)})^2,(\theta_2^{(k)})^2,\ldots,(\theta_k^{(k)})^2,0$,
and $\bar{B}_{k-1}\bar{B}_{k-1}^T$ is the $k\times k$ leading principal
submatrix of $B_k B_k^T$, whose eigenvalues are
$(\bar{\theta}_1^{(k-1)})^2,(\bar{\theta}_2^{(k-1)})^2, \ldots,
(\bar{\theta}_k^{(k-1)})^2$, by the strict interlacing property of eigenvalues,
we obtain
\begin{equation}\label{secondinter}
\theta_1^{(k)}>\bar{\theta}_1^{(k-1)}>\theta_2^{(k)}>
\bar{\theta}_2^{(k-1)}> \cdots >\theta_k^{(k)}> \bar{\theta}_k^{(k-1)}>0, \
k=1,2,\ldots,n.
\end{equation}
On the other hand, note that $\alpha_{n+1}=0$ and
$\bar{\theta}_i^{(n)}=\theta_i^{(n)}=\sigma_i,\ i=1,2,\ldots,n$, i.e.,
the singular values of $\bar{B}_n$ are $\sigma_1,\sigma_2,\ldots,\sigma_n$
and zero, which is denoted by the dummy $\sigma_{n+1}=0$.  Since
$\bar{\theta}_{n+1}^{(n)}=\sigma_{n+1}=0$ and the first $k$ rows of $\bar{B}_n$
are $(\bar{B}_{k-1},\mathbf{0})\in \mathbb{R}^{k\times n}$, whose singular
values are $\bar{\theta}_1^{(k-1)},\ldots,\bar{\theta}_k^{(k-1)}$, by
applying the strict interlacing property
of singular values to $(\bar{B}_{k-1},\mathbf{0})$ and $\bar{B}_n$, for
$k=1,2,\ldots,n-1$ we have
\begin{equation}\label{interbar}
\sigma_{n+1-k+i}<\bar{\theta}_i^{(k-1)}<\sigma_i,\ i=1,2,\ldots,k,
\end{equation}
from which it follows that
\begin{equation}\label{thetak}
0<\bar{\theta}_k^{(k-1)}<\sigma_k.
\end{equation}

\eqref{secondinter} and \eqref{thetak} indicate that,
unlike $\theta_k^{(k)}$ that lies between $\sigma_{k+1}$ and $\sigma_k$
and approximates $\sigma_k$ for severely or moderately ill-posed problems
with $\rho>1$ or $\alpha>1$ suitably (cf. \eqref{error2}), the
lower bound for $\bar{\theta}_k^{(k-1)}$ is simply zero,
and there does not exist a better one for it. This means that
$\bar{\theta}_k^{(k-1)}$ may be much smaller than $\sigma_{k_0+1}$ and actually
it can be arbitrarily small, independently of the degree $\rho$ or $\alpha$ of
ill-posedness. In other words, the size of $\rho$ or $\alpha$ does not have any
intrinsic effects on the lower bound of $\bar{\theta}_k^{(k-1)}$, and
one thus cannot control $\bar{\theta}_k^{(k-1)}$ from below by choosing $\rho$ or
$\alpha$. In the meantime, \eqref{secondinter} tells us that
$\bar{\theta}_k^{(k-1)}<\theta_k^{(k)}$. These facts, together with
Theorem~\ref{cgmeappr}, show that the regularization of
CGME is inferior to that of LSQR and LSMR for each kind of problem.
On the one hand, they mean that
CGME has the partial regularization for mildly ill-posed problems; on the
other hand, the regularizing effects of CGME have indeterminacy
for severely and moderately ill-posed problems, that is,
it may or may not have the full regularization for these
two kinds of problems. Clearly, CGME has the full regularization
only when $P_k\bar{B}_{k-1}Q_k^T$
is as accurate as the rank $k$ approximation $P_{k+1}B_kQ_k^T$ and
$\bar{\theta}_k^{(k-1)}\approx \theta_k^{(k)},\
k=1,2,\ldots,k_0$ for these two kinds of problems with $\rho>1$ and
$\alpha>1$ considerably, but unfortunately there is no
guarantee that these requirements are satisfied mathematically.

The above analysis indicates that CGME itself is not reliable and cannot be
trusted to compute best possible regularized solutions.
In principle, one can detect the full or partial regularization
of CGME as follows: One first exploits the decay rate of
$\alpha_k+\beta_{k+1}$ to identify the degree of ill-posedness of \eqref{eq1}.
If \eqref{eq1} is mildly ill-posed, CGME has only
the partial regularization. If \eqref{eq1} is recognized as severely or moderately
ill-posed, one then needs to do two things to identify the regularization
of CGME: check if $\|A-P_k\bar{B}_{k-1}Q_k\|
\approx \|A-P_{k+1}B_kQ_k\|$, and compute the singular values of both $B_k$ and
$\bar{B}_{k-1}$ and check if $\bar{\theta}_k^{(k-1)}\approx \theta_k^{(k)}$.
If both hold, CGME has the full regularization; if either of them does not hold,
it has only the partial regularization.

We can informally deduce more features on CGME. For the LSQR iterate $x^{(k)}$,
note that the optimality requirement of CGME means that
$\|x_{naive}-x_k^{cgme}\|\leq \|x_{naive}-x^{(k)}\|$. Since
$$
\|x_{naive}-x_k^{cgme}\|=\|x_{naive}-x_{true}+x_{true}-x_k^{cgme}\|
\leq \|x_{naive}-x_{true}\|+\|x_{true}-x_k^{cgme}\|
$$
and
$$
\|x_{naive}-x^{(k)}\|=\|x_{naive}-x_{true}+x_{true}-x^{(k)}\|
\leq \|x_{naive}-x_{true}\|+\|x_{true}-x^{(k)}\|
$$
with the first terms in the right-hand sides being the same constant,
not rigorously speaking, we should have
\begin{equation}\label{cgmelsqr}
\|x_{true}-x_k^{cgme}\|\leq \|x_{true}-x^{(k)}\|
\end{equation}
until the semi-convergence of CGME. Keep in mind that
the regularization of CGME is inferior to or are at most
as good as that of LSQR for each kind of ill-posed problem.
Both $\|x_{true}-x_k^{cgme}\|$ and
$\|x_{true}-x^{(k)}\|$ first decrease until their respective
semi-convergence and then become increasingly large as $k$ increases.
As a result, we deduce that (i) $x_k^{cgme}$ is at least as accurate as
$x^{(k)}$ until the semi-convergence of CGME and (ii)
CGME reaches semi-convergence no later than LSQR;
otherwise, \eqref{cgmelsqr} indicates
that the optimal regularized solution by CGME at semi-convergence
would be more accurate than that by LSQR at semi-convergence, which
contradicts the property that LSQR has better
regularization than CGME. The experiments in \cite{hanke01} justify this
assertion; see Figure 3.1 and Figure 5.2 there.

Next let us return to \eqref{leftapp} and show
how to extract a rank $k$ approximation to $A$ from the
rank $k+1$ approximation $P_{k+1}\bar{B}_kQ_{k+1}^T$ as best as possible.

\begin{theorem}\label{approx}
Let $\bar{C}_k$ be the best rank $k$ approximation to $\bar{B}_k$ with respect
to the 2-norm. Then
\begin{align}
\|A-P_{k+1}\bar{C}_kQ_{k+1}^T\|&\leq \sigma_{k+1}+\gamma_k,\label{lowrank}\\
\|A-P_{k+1}\bar{C}_kQ_{k+1}^T\|&\leq \bar{\theta}_{k+1}^{(k)}+\gamma_k,
\label{better}
\end{align}
where $\bar{\theta}_{k+1}^{(k)}$ is the smallest singular value of $\bar{B}_k$.
\end{theorem}

{\em Proof}. Write $A-P_{k+1}\bar{C}_kQ_{k+1}^T=A-P_{k+1}\bar{B}_kQ_{k+1}^T+
P_{k+1}(\bar{B}_k-\bar{C}_k)Q_{k+1}^T$. Then from \eqref{leftbidiag} we obtain
\begin{align}
\|A-P_{k+1}\bar{C}_kQ_{k+1}^T\| & \leq \|A-P_{k+1}\bar{B}_kQ_{k+1}^T\|+
\|P_{k+1}(\bar{B}_k-\bar{C}_k)Q_{k+1}^T\|  \label{barbk}\\
&= \|A-P_{k+1}\bar{B}_kQ_{k+1}^T\|+
\|P_{k+1}P_{k+1}^TA-P_{k+1}\bar{C}_kQ_{k+1}^T\|. \label{decom2}
\end{align}
By the assumption on $C_k$ and \eqref{leftbidiag},
$P_{k+1}\bar{C}_kQ_{k+1}^T$ is the best rank $k$ approximation to
$P_{k+1}\bar{B}_kQ_{k+1}^T=P_{k+1}P_{k+1}^TA$. Keep in mind that $A_k$ is
the best rank $k$ approximation to $A$. Since $P_{k+1}P_{k+1}^TA_k$ is a rank
$k$ approximation to $P_{k+1}P_{k+1}^TA$, we get
\begin{align*}
\|P_{k+1}P_{k+1}^TA-P_{k+1}\bar{C}_kQ_{k+1}^T\| &\leq \|P_{k+1}P_{k+1}^T(A-A_k)\|\\
&\leq \|A-A_k\|=\sigma_{k+1},
\end{align*}
from which, \eqref{leftapp} and \eqref{decom2} it follows
that \eqref{lowrank} holds.

Since $P_{k+1}$ and $Q_{k+1}$ are orthonormal, by the 2-norm invariance, we
obtain
$$
\|P_{k+1}(\bar{B}_k-\bar{C}_k)Q_{k+1}^T\|=\|\bar{B}_k-\bar{C}_k\|
=\bar{\theta}_{k+1}^{(k)},
$$
from which and \eqref{barbk} it follows that \eqref{better} holds.
\qquad\endproof

We point out that \eqref{lowrank} may be conservative since
we have amplified $\|P_{k+1}(\bar{B}_k-\bar{C}_k)Q_{k+1}^T\|$
twice and obtained its bound $\sigma_{k+1}$, which can be a considerable
overestimate. In comparison with \eqref{gammak} and \eqref{final},
the bound \eqref{lowrank} indicates that $P_{k+1}\bar{C}_kQ_{k+1}^T$
may not be as accurate as $P_{k+1}B_kQ_k^T$, but
\eqref{better} illustrates that $P_{k+1}\bar{C}_kQ_{k+1}^T$
can be as accurate as $P_{k+1}B_kQ_k^T$ because
$\bar{\theta}_{k+1}^{(k)}<\theta_{k+1}^{(k+1)}<\sigma_{k+1}$
from \eqref{secondinter} and \eqref{error2}. Moreover,
as we have explained, $\bar{\theta}_{k+1}^{(k)}$ can be arbitrarily small.
If so, $\bar{\theta}_{k+1}^{(k)}$ is negligible in \eqref{better} and
$P_{k+1}\bar{C}_kQ_{k+1}^T$ is at least as accurate as
$P_{k+1}B_kQ_k^T$.%

We now present a new but informal analysis to show why
$P_{k+1}\bar{C}_kQ_{k+1}^T$ may be at least as accurate as
$P_{k+1}B_kQ_k^T$ as a rank $k$ approximation
to $A$. Keep in mind that $\bar{\theta}_i^{(k)},\
i=1,2,\ldots,k+1$ be the singular values of $\bar{B}_k$.
Then the singular values of $\bar{C}_k$ are
$\bar{\theta}_i^{(k)},\ i=1,2,\ldots,k$. Since $\alpha_{k+1}>0$ for
$k\leq n-1$, applying the strict interlacing property of singular
values to $B_k$ and $\bar{B}_k$, we have
\begin{equation}\label{bkbarbk}
\bar{\theta}_1^{(k)}>\theta_1^{(k)}>\bar{\theta}_2^{(k)}>\cdots>
\bar{\theta}_k^{(k)}>\theta_k^{(k)}>\bar{\theta}_{k+1}^{(k)}>0, \
k=1,2,\ldots,n-1.
\end{equation}
The above relationships, together with \eqref{ritzapp}, prove that
\begin{equation}\label{ck}
\sigma_i-\bar{\theta}_i^{(k)}<\sigma_i-\theta_i^{(k)}\leq \gamma_k,\,
i=1,2,\ldots,k,
\end{equation}
that is, the $\bar{\theta}_i^{(k)}$ are more accurate than
$\theta_i^{(k)}$ as approximations to $\sigma_i,\ i=1,2,\ldots,k$.
By the standard perturbation theory, note from \eqref{leftapp} that
$$
\sigma_i-\bar{\theta}_i^{(k)}\leq \|A-P_{k+1}\bar{B}_kQ_{k+1}^T\|\leq \gamma_k,
\,i=1,2,\ldots,k+1,
$$
while the singular value differences between $A$ and $P_{k+1}\bar{C}_kQ_{k+1}^T$
are $\sigma_i-\bar{\theta}_i^{(k)},\ i=1,2,\ldots,k$ and
$\sigma_i,\ i=k+1,\ldots,n$, all of which, from \eqref{ck} and \eqref{final},
are no more than $\gamma_k$. Based on these rigorous facts and
the relationship between $C_k$ and $\bar{B}_k$,
it is possible that $\|A-P_{k+1}\bar{C}_kQ_{k+1}^T\|\leq \gamma_k$,
and if it is so, then by definition \eqref{gammak} $P_{k+1}\bar{C}_kQ_{k+1}^T$
is a more accurate rank $k$ approximation to $A$ than $P_{k+1}B_kQ_k^T$ is.

\subsection{A comparison with standard randomized algorithms and RRQR
factorizations}\label{rrqr}

We compare the rank approximations $P_{k+1}B_kQ_k^T$ and
$P_{k+1}\bar{C}_kQ_{k+1}^T$ by Lanczos bidiagonalization with those by some
standard randomized algorithms and RRQR factorizations, and demonstrate
that the former ones are much more accurate than the latter ones
for severely and moderately ill-posed problems.

Note \eqref{gamma2}. Compare \eqref{final}, \eqref{leftsub} and \eqref{lowrank}
or \eqref{better} with the corresponding results
(1.9), (5.6), (6.3) and Theorem 9.3 in \cite{halko11}
for standard randomized algorithms and those on the strong
RRQR factorization \cite{gu96}, where the
constants in front of $\sigma_{k+1}$ are like $\sqrt{kn}$
and $\sqrt{1+4k(n-k)}$, respectively, which are far bigger than one.  Within
the framework of the RRQR factorizations, it is known from \cite{hong92} that
the optimal factor of such kind is $\sqrt{k(n-k)+\min\{k,n-k\}}$
but to find corresponding permutations is
an NP-hard problem, whose cost increases exponentially
with $n$; see also \cite[p.298]{bjorck15}. Clearly, the strong
RRQR factorizations are near-optimal within the framework,
and they suit well for finding a high quality low rank $k$ approximation
to a matrix whose $k$ large singular values are much bigger than
the $n-k$ small ones.

Unfortunately, the standard randomized algorithms and RRQR factorization do
not very nicely fit into solving ill-posed problems: they have regularizing
effects but, in general, cannot find best possible regularized solutions.
We argue as follows: Since there are no considerable gaps of singular values,
the RRQR factorization techniques can hardly find a near best rank $k$
approximation to $A$ in the sense of \eqref{near}, which is vital to
solve \eqref{eq1} to find a best possible regularized solution.
In contrast, for a severely or moderately ill-posed problem with $\rho>1$
or $\alpha>1$ suitably, the rank $k$ approximations $P_{k+1}B_kQ_k^T$ are near
best ones for $k=1,2,\ldots, k_0$ and no singular value smaller than $\sigma_{k_0+1}$
appears. Besides, it is easy to check that the $k$-step
Lanczos bidiagonalization costs fewer flops than the standard randomized
algorithms do for a
sparse $A$, and it is more efficient than the strong RRQR factorization for
a dense $A$, which includes $\mathcal{O}(mnk)$ flops and the overhead cost of
searching permutations.

For further developments and recent advances on
randomized algorithms, we refer to Gu's work \cite{gu2015},
where he has considered randomized
algorithms within the subspace iteration framework proposed in
\cite{halko11}, presented a number of
new results and improved the error bounds for
the rank $k$ approximations that are iteratively extracted. Such
approaches may be promising to solve ill-posed problems.

\section{The filters $f_i^{(k)}$ and
a comparison of LSQR and the TSVD method}\label{compare}

Based on Proposition~\ref{help}, exploiting Theorem~\ref{main1},
Theorem~\ref{ritzvalue} and
Theorem~\ref{main2}, we present the following results, which,
from the viewpoint of Tikhonov regularization, explain why
LSQR has the full regularization for severely and moderately ill-posed problems
with $\alpha>1$ and $\alpha>1$ suitably and why it generally has the partial
regularization for mildly ill-posed problems.

\begin{theorem}\label{main3}
For the severely or moderately ill-posed problems with $\rho>1$ or $\alpha>1$,
under the assumptions of Theorem~\ref{ritzvalue}, let
$f_i^{(k)}$ be defined by \eqref{filter}. Then for $k=1,2,\ldots,k_0$ we have
\begin{align}
|f_i^{(k)}-1| &\approx \frac{2\sigma_{k+1}}
{\sigma_i}\left|\prod\limits_{j=1,j\not= i}^{k}
\left(1-\left(\frac{\sigma_i}{\sigma_j}\right)^2\right)\right|,
\ i=1,2,\ldots,k, \label{fi1}\\
f_i^{(k)}&\approx  \sigma_i^2\sum\limits_{j=1}^{k}\frac{1}{\sigma_j^2},
\ i=k+1,k+2,\ldots,n. \label{fi2}
\end{align}
\end{theorem}

{\em Proof}.
For $k=1,2,\ldots,k_0$,
it follows from \eqref{filter} that
\begin{equation*}
|f_i^{(k)}-1|=\left|\frac{(\theta_i^{(k)})^2-\sigma_i^2}
{(\theta_i^{(k)})^2}
\prod\limits_{j=1,j\not=i}^{k}\frac{(\theta_j^{(k)})^2-\sigma_i^2}
{(\theta_j^{(k)})^2}\right|, \ i=1,2,\ldots,k.
\end{equation*}
To simplify presentations and illuminate the essence,
for the severely and moderately ill-posed problems with $\rho>1$ and
$\alpha>1$ suitably, we simply replace $\sqrt{1+\eta_k^2}$ in \eqref{error}
by one.  On the other hand, we replace the denominator of
$\frac{(\theta_i^{(k)})^2-\sigma_i^2}{(\theta_i^{(k)})^2}$
by $\sigma_i^2$. Then by \eqref{error} we approximately have
$$
\sigma_i^2-(\theta_i^{(k)})^2=(\sigma_i-\theta_i^{(k)})
(\sigma_i+\theta_i^{(k)})\approx  2\sigma_{k+1}\sigma_i.
$$
For $j=1,2,\ldots,k$ but $i$, replace $\theta_j^{(k)}$
by $\sigma_j$ approximately. Then \eqref{fi1} follows.

By \eqref{error2}, since $\theta_{k}^{(k)}>
\sigma_i$ for $i=k+1,\ldots,n$,
the factors $\sigma_i/\theta_j^{(k)}<1$ and decay to
zero with increasing $i$ for each fixed $j\leq k$. Therefore,
for $i=k+1,\ldots,n$ we get
\begin{align*}
f_i^{(k)}&=1-\prod_{j=1}^{k} \left(1-\left(\frac{\sigma_i}{\theta_j^{(k)}}
\right)^2\right) \\
&=1-\left(1-\sum_{j=1}^{k} \left(\frac{\sigma_i}{\theta_j^{(k)}}\right)^2
\right)+\mathcal{O} \left(\frac{\sigma_i^4}{\left(\theta_{k-1}^{(k)}
\theta_{k}^{(k)}\right)^2}\right)\\
&=\sum_{j=1}^{k} \left(\frac{\sigma_i}{\theta_j^{(k)}}\right)^2
+\mathcal{O} \left(\frac{\sigma_i^4}{\left(\theta_{k-1}^{(k)}
\theta_{k}^{(k)}\right)^2}\right)
\end{align*}
Replace $\theta_j^{(k)}$ by its upper bound $\sigma_j,\ j=1,2,\ldots,k$ in
the above, and note that the second term is higher order small relative to
the first term. Then \eqref{fi2} follows.
\qquad\endproof

\begin{remark}
For $k=1,2,\ldots,k_0$,
since $\sigma_i,\, i=1,2,\ldots,k$ are dominant singular values, the factors
$$
\left|\prod\limits_{j=1,j\not= i}^{k}
\left(1-\left(\frac{\sigma_i}{\sigma_j}\right)^2\right)\right|,
\ i=1,2,\ldots,k,
$$
are modest. Consequently, \eqref{fi1} indicates that
the $f_i^{(k)}\approx 1$ with the errors
$\mathcal{O}(\sigma_{k+1}/\sigma_i)$ for
$i=1,2,\ldots,k$, while
\eqref{fi2} shows that the $f_i^{(k)}$ are at least as small as
$\sigma_i^2/\sigma_{k}^2$ for $i=k+1,\ldots,n$
and decrease with increasing $i$.
\end{remark}

\begin{remark}
For mildly ill-posed problems and $k=1,2,\ldots,k_0$,
as we have shown in Remark~\ref{appear}
and Section \ref{rankapp}, it is generally the case that
$\theta_{k}^{(k)}<\sigma_{k_0+1}$. Suppose $\theta_{k}^{(k)}>\sigma_{j^*}$
with the smallest integer $j^*>k_0+1$. Then we have shown in the paragraph
after Proposition~\ref{help} that $f_i^{(k)}\geq 1, \ i=k_0+1,\ldots, j^*-1$.
As a result, LSQR has only the partial regularization.
\end{remark}

Recall that $\Delta_k=(\delta_1,\delta_2,\ldots,\delta_k)$,
and define $U_k^{o}=(u_{k+1},\ldots,u_n)$.
In terms of \eqref{ideal3}--\eqref{ideal},
Hansen \cite[p.151,155, Theorems 6.4.1-2]{hansen98}
presents the following bounds
\begin{equation}\label{ferror1}
|f_i^{(k)}-1|\leq \frac{\sigma_{k+1}}{\sigma_i}
\frac{(\|(U_k^{o})^Tb\||+\sigma_{k+1}\|\Delta_k\|\|x^{(k)}\|}
{| u_i^Tb |}\|\delta_i\|,\ i=1,2,\ldots,k,
\end{equation}
\begin{equation}\label{deltainf}
\|\delta_i\|_{\infty}\leq \frac{\sigma_{k+1}}{\sigma_i}
\frac{| u_{k+1}^T b|}{| u_i^T b|} |L_i^{(k)}(0) |,\ i=1,2,\ldots,k,
\end{equation}
and
\begin{equation}\label{ferror2}
0 \leq f_i^{(k)}\leq \frac{\sigma_i^2}{\sigma_k^2} |L_k^{(k)}(0)|
\sum_{j=1}^k f_j^{(k)},\ \ i=k+1,\ldots,n,
\end{equation}
where $\|\cdot\|_{\infty}$ is the infinity norm of a vector.

We now address a few points on the bounds \eqref{ferror1} and \eqref{ferror2}.
First, there had no estimates for $\|\Delta_k\|$ and $|L_i^{(k)}(0)|,\
i=1,2,\ldots,k$; second, what we need is $\|\delta_k\|$ other than
$\|\delta_k\|_{\infty}$, and as is seen from its proof in
\cite[p.151]{hansen98}, it is relatively easy to obtain the accurate bound
\eqref{deltainf} for $\|\delta_i\|_{\infty}$, whereas it is
hard to derive an accurate one for $\|\delta_i\|$.
Because of lacking accurate estimates, it is unclear how small or large
the bound \eqref{ferror1} and \eqref{ferror2} are. Moreover, as it will
appear soon, the
factor $\sigma_{k+1}\|(U_k^{o})^T b\|$ in the numerator of \eqref{ferror1}
may be a too crude overestimate, such that the bound \eqref{ferror1} is
pessimistic and is useless to estimate $|f_i^{(k)}-1|,\
i=1,2,\ldots,k$ and $f_i^{(k)},\ i=k+1,\ldots,n$.
Let us have a closer look at these points.
Obviously, exploiting \eqref{deltainf}, we can only obtain the bounds
\begin{align*}
\|\delta_i\|_{\infty} &\leq \|\delta_i\|\leq \sqrt{n-k}\|\delta_i\|_{\infty},\\
\max_{i=1,2,\ldots,k}\|\delta_i\|&\leq \|\Delta_k\|\leq
\sqrt{k}\max_{i=1,2,\ldots,k}\|\delta_i\|,
\end{align*}
from which it follows that
$$
\max_{i=1,2,\ldots,k}\|\delta_i\|_{\infty}\leq \|\Delta_k\|\leq
\sqrt{k(n-k)}\max_{i=1,2,\ldots,k}\|\delta_i\|_{\infty}.
$$
As a result, the estimates for both $\|\delta_i\|$ and $\|\Delta_k\|$ are too
crude for $n$ large and $k$ small. Indeed, as we have seen previously,
their accurate estimates are much involved and complicated.
In Theorem~\ref{thm3}, we have derived accurate estimates for
$\|\delta_i\|,\ i=1,2,\ldots,k$; see \eqref{columndelta}, \eqref{columndelta1}
for severely ill-posed problems and \eqref{columndelta2}, \eqref{columnnorm}
for moderately and mildly ill-posed problems. Theorems~\ref{thm2}--\ref{moderate}
have given sharp estimates for $\|\Delta_k\|$ for three kinds of problems,
respectively.

The factor $\sigma_{k+1}\|(U_k^{o})^T b\|$ itself in the numerator
of \eqref{ferror1}, though simple and elegant in form, does not give clear and
quantitative information on its size. As a matter of fact, one must
analyze its size carefully for the two cases
$k\leq k_0$ and $k>k_0$, respectively,
for each kind of ill-posed problem; see the discrete Picard condition
\eqref{picard} and \eqref{picard1}. For each of these two cases, using our proof
approach used for Theorems~\ref{thm2}, \ref{moderate} and \ref{thm3},
we can obtain accurate estimates for $\|(U_k^{o})^T b\|=
\left(\sum_{j=k+1}^n |u_j^T b|^2\right)^{1/2}$ for three kinds of
ill-posed problems, respectively. However, the point is that the
factor $\sigma_{k+1}\|(U_k^{o})^T b\|$ results from a substantial
amplification in the derivation. It is seen from the last line
of \cite[p.155]{hansen98} that this factor results from simply bounding it by
$$
\|\Sigma_k^{\perp} (U_k^{(o)})^T b\|\leq \sigma_{k+1} \|(U_k^{(o)})^Tb\|,
$$
where $\Sigma_k^{\perp}={\rm diag} (\sigma_{k+1},\ldots,\sigma_n)$.
For our context, this amplification is fatal, and it is subtle
to obtain sharp bounds for $\|\Sigma_k^{\perp}(U_k^{(o)})^T b\|$.
We observe that $\|\Sigma_k^{\perp} (U_k^{(o)})^T b\|$ is nothing but
the first square root factor in \eqref{delta2}, for which we have
established the accurate estimates \eqref{severe1} and \eqref{modeest}
for severely, moderately and mildly ill-posed problems, respectively,
which hold for $k=1,2,\ldots,n-1$ and are {\em independent} of $n$.
It can be checked that these bounds for $\|\Sigma_k^{\perp} (U_k^{(o)})^T b\|$ are
substantially smaller than $\sigma_{k+1} \|(U_k^{(o)})^Tb\|$.

After the above substantial improvements on \eqref{ferror1} and
\eqref{ferror2}, we can exploit
the accurate bounds for $\|\Delta_k\|$ and $\|\delta_i\|$ in
Theorems~\ref{thm2}--\ref{moderate} and Theorem~\ref{thm3}, as
well as the remarks on them,
to accurately estimate the bounds \eqref{ferror1}
and \eqref{ferror2}. From them we can draw the full regularization
LSQR for severely and moderately ill-posed problems with
$\rho>1$ and $\alpha>1$ suitably and its partial regularization
for mildly ill-posed problems.

Making use of some standard perturbation results from Hansen~\cite{hansen10},
we can quantitatively relate LSQR to the TSVD method and analyze the
differences between their corresponding regularized solutions
and differences between $P_{k+1}B_kQ_k^T x^{(k)}$ and $A_k x_k^{tsvd}$ predicting
the right-hand side $b$ for $k=1,2,\ldots,k_0$.

\begin{theorem}\label{lsqrtsvd}
For the severely or moderately ill-posed problem~\eqref{eq1},
let $A_k$ be the rank $k$ best approximation to $A$, and assume that
$\|E_k\|=\|P_{k+1}B_kQ_k^T-A_k\|\leq \sigma_k-\sigma_{k+1}$. Then
for $k=1,2,\ldots,k_0$
we have
\begin{align}
\frac{\|x^{(k)}-x_k^{tsvd}\|}{\|x_k^{tsvd}\|}&\leq
\frac{\kappa(A_k)}{1-\epsilon_k}\left(\frac{\|E_k\|}{\|A_k\|}
+\frac{\epsilon_k}{1-\epsilon_k-\hat{\epsilon}_k}
\frac{\|A_kx_k^{tsvd}-b\|}{\|A_k x_k^{tsvd}\|}\right), \label{comperror}
\end{align}
\begin{align}
\frac{\|P_{k+1}B_kQ_k^T x^{(k)}-A_k x_k^{tsvd}\|}{\|b\|}&\leq
\frac{\epsilon_k}{1-\epsilon_k}, \label{compres}
\end{align}
where
$$
\kappa(A_k)=\frac{\sigma_1}{\sigma_k},\ \ \epsilon_k=\frac{\|E_k\|}{\sigma_k},\ \
\hat{\epsilon}_k=\frac{\sigma_{k+1}}{\sigma_k}.
$$
\end{theorem}

{\em Proof}.
For the problem $\min\|A_k x-b\|$ that replaces $A$ by $A_k$
in \eqref{eq1}, we regard the rank $k$ matrix $P_{k+1}B_kQ_k^T$
as a perturbed $A_k$ with the perturbation matrix
$E_k=P_{k+1}B_kQ_k^T-A_k$. Then by the standard perturbation
results on the TSVD solutions \cite[p.65-6]{hansen10}, we obtain
\eqref{comperror} and \eqref{compres} directly.
\qquad\endproof

\begin{remark}
Write $\|E_k\|=\|P_{k+1}B_kQ_k^T-A+A-A_k\|$.
Since the rank $k$ matrices $P_{k+1}B_kQ_k^T$ and
$A_k$ have the $k$ nonzero singular values $\theta_i^{(k)}$ and $\sigma_i$,
$i=1,2,\ldots,k$, respectively, from Mirsky's theorem
\cite[p.204, Theorem 4.11]{stewartsun} we get the bounds
\begin{align}
\max_{i=1,\ldots,k}|\sigma_i-\theta_i^{(k)}|&\leq \|E_k\|=\|A_k-P_{k+1}B_kQ_k^T\|,
\label{first}\\
\max\{\max_{i=1,\ldots,k}|\sigma_i-\theta_i^{(k)}|,
\sigma_{k+1}\}&\leq \|A-P_{k+1}B_kQ_k^T\|=\gamma_k, \label{second}
\end{align}
where the lower bound in \eqref{first} is no more than the one in \eqref{second}.
It is then expected that
$\|E_k\|\leq \gamma_k\approx \sigma_{k+1}$ for severely and
moderately ill-posed problems.
Therefore, we have $\epsilon_k\approx \hat{\epsilon}_k<1$, and
\eqref{compres} indicates that $\|P_{k+1}B_kQ_k^T x^{(k)}-A_k x_k^{tsvd}\|$,
is basically no more than $\epsilon_k$, $k=1,2,\ldots,k_0$.
\end{remark}

\begin{remark}
From \eqref{comperror}, since the possibly not small factor
$$
\frac{\|A_kx_k^{tsvd}-b\|}{\|A_k x_k^{tsvd}\|}
$$
enters the bound \eqref{comperror}, two regularized solutions $x^{(k)}$ and
$x_k^{tsvd}$ may differ considerably even though $P_{k+1}B_kQ_k^T x^{(k)}$
and $A_k x_k^{tsvd}$ predict the right-hand side $b$ with similar accuracy
for $k=1,2,\ldots,k_0$. This is the case for the inconsistent ill-posed
problem $\min\|Ax-b\|$ with $m>n$, where $\|A_kx_k^{tsvd}-b\|$ decreases
with respect to $k$ until
$$
\|A_{k_0}x_{k_0}^{tsvd}-b\|^2=\|Ax_{k_0}^{tsvd}-b\|^2\approx
\frac{n-k_0}{m}\|e\|^2+\|(I-U_nU_n^T)b\|^2,
$$
with $U_n$ the first $n$ columns of the $m\times m$ left singular vector
matrix $U$ and $\|(I-U_nU_n^T)b\|$
the incompatible part of $b$ lying outside of the range of $A$
(cf. \cite[p.71,88]{hansen10}). Here we remark that the term
$\|(I-U_nU_n^T)b\|$ appears in the relation (4.17) of \cite[p.71]{hansen10}
but is missing in the above right-hand side \cite[p.88]{hansen10}.
For the consistent $Ax=b$, since $\|(I-U_nU_n^T)b\|=0$,
the right-hand side of \eqref{comperror} is approximately 
$$
\frac{\sigma_{k_0+1}}{\sigma_{k_0}}\left(1+\frac{\sigma_{1}}{\sigma_{k_0}}
\sqrt{\frac{n-k_0}{m}}\frac{\|e\|}{\|b\|}\right).
$$
We see from the above and \eqref{compres} that two different
regularized solutions can be quite different even if their residual norms
are of similar very sizes, as addressed by Hansen
\cite[p.123-4, Theorem 5.7.1]{hansen98}. However, we point out that
the accuracy of different regularized solutions as approximations to
$x_{true}$ can be compared. If the norms of errors of them and $x_{true}$ have
very comparable sizes, they are equally accurate regularized solutions
to \eqref{eq1}.
\end{remark}

\begin{remark}
Note that $\|E_k\|=\|P_{k+1}B_kQ_k^T-A_k\|\leq \sigma_k-\sigma_{k+1}$
is assumed only for severely and moderately ill-posed problems. As the previous
analysis has indicated, we have $\|E_k\|\approx \gamma_k\approx\sigma_{k+1}$.
As a result, it is easily justified
that this assumption is valid for these two kinds of problems provided
that $\rho>1$ and $\alpha>1$ suitably. However, the assumption fails
to hold for the mildly ill-posed problems with
$\sigma_i=\zeta i^{-\alpha},\ i=1,2,\ldots,n$
and $\frac{1}{2}<\alpha\leq 1$ since, for $k>1$, we have
$$
\sigma_k-\sigma_{k+1}=\sigma_{k+1}
\left(\left(\frac{k+1}{k}\right)^{\alpha}-1\right)<\sigma_{k+1}
\leq \gamma_k\approx \|E_k\|.
$$
\end{remark}

\section{The extension to the case that $A$ has multiple singular values}
\label{multiple}

Previously, under the assumption that the singular values of $A$ are simple,
we have proved the results and made a detailed analysis on them. Recall
the basic fact that the singular values $\theta_i^{(k)},\ i=1,2,\ldots,k$
of $B_k$ are always simple mathematically, independent of whether
the singular values of $A$ are simple or multiple. In other words,
the Lanczos bidiagonalization process works as if the singular values
of $A$ are simple, and the Ritz values $\theta_i^{(k)},\ i=1,2,\ldots,k$,
are the approximations to some of the {\em distinct} singular values of $A$.
In this section, we will show that, by making a number of suitable and nontrivial
changes and reformulations, our previous results and analysis can
be extended to the case that $A$ has multiple singular values.

Assume that $A$ has $s$ distinct singular values
$\sigma_1>\sigma_2>\cdots>\sigma_s>0$ with $\sigma_i$ being $c_i$ multiple
and $s\leq n$. In order to treat this case, we need to make a number of
preliminary preparations and necessary modifications or reformulations.
Below let us show the detail.

First of all, we need to take $b$ into consideration and present a new form
SVD of $A$ by selecting a specific set of left and right singular vectors
corresponding to a multiple singular value $\sigma_i$ of $A$,
so that the discrete Picard condition \eqref{picard} holds for one
particularly chosen left singular vector associated with $\sigma_i$. Specifically,
for the $c_i$ multiple $\sigma_i$, the orthonormal basis of the
corresponding left singular subspace can be chosen so that
$b$ has a nonzero orthogonal projection on just
one unit length left singular vector $u_i$ in the
singular subspace and no components in the remaining $c_i-1$ ones. Precisely,
let the columns of $F_i$ form an orthonormal basis of the left singular
subspace associated with $\sigma_i$, each of which satisfies \eqref{picard}.
Then we take
\begin{equation}\label{newuk}
u_i=\frac{F_iF_i^Tb}{\|F_i^Tb\|},
\end{equation}
where $F_iF_i^T$ is the orthogonal projector onto the
left singular subspace with $\sigma_i$,
and define the corresponding unit length right singular vector by
$v_i = A^T u_i/\sigma_i$. We select the other $c_i-1$ orthonormal
left singular vectors which are orthogonal to $u_i$ and, together with $u_i$,
form the left singular subspace associated with $\sigma_i$, and define the
corresponding unit length right singular vectors in the same way as $v_i$,
which and $v_i$ form an orthonormal basis of the unique right singular
subspace with $\sigma_i$.
After such treatment, we get the desired SVD of $A$. We stress that
$u_i$ defined above is unique since the orthogonal projection of $b$
onto the left singular subspace with $\sigma_i$ is unique
and equal to $F_iF_i^Tb$ for a given orthonormal $F_i$.

Now we need to prove that $u_i$ satisfies the discrete Picard
condition \eqref{picard} essentially. To see this, by the Cauchy--Schwarz
inequality, \eqref{newuk} and the assumption that each column of $F_i$
satisfies the discrete Picard condition \eqref{picard}, we get
\begin{equation}\label{estpicard}
|u_i^T \hat{b}|=\frac{| b^TF_iF_i^T \hat{b}|}{\|F_i^Tb\|}\leq \|F_i^T\hat{b}\|
=\sqrt{c_i}\sigma_i^{1+\beta},\ i=1,2,\ldots,s.
\end{equation}
Therefore, the Fourier coefficients $|u_i^T \hat{b}|$, on average, decay faster
than the singular values $\sigma_i,\,i=1,2,\ldots,s$. This is exactly
what the discrete Picard condition means; see the description before
\eqref{picard}. Recall that \eqref{picard} is a simplified model of this
condition. Based on the estimate \eqref{estpicard},
we recover \eqref{picard} by simply resetting \eqref{estpicard} as
\begin{equation}\label{newpicard}
|u_i^T \hat{b}|=\sigma_i^{1+\beta},\ i=1,2,\ldots,s.
\end{equation}

With help of the SVD of $A$ described above, it is crucial to
observe that $x_k^{tsvd}$ in \eqref{solution} is now the sum consisting of
the first $k$ {\em distinct} dominant SVD components of $A$. Furthermore,
for \eqref{eq1} and such reformulation of \eqref{solution},
the matrix $A$ in them can be equivalently replaced by the new $m\times n$
matrix
\begin{equation}\label{aprime}
A^{\prime}=U\Sigma^{\prime} V^T,
\end{equation}
where $\Sigma^{\prime}={\rm diag}(\sigma_1,\sigma_2,
\ldots,\sigma_s,\mathbf{0})$, $U_s=(u_1,u_2,\ldots,u_s)$ and
$V_s=(v_1,v_2,\ldots,v_s)$ are the first $s$ columns of
$U$ and $V$, respectively, the last $n-s$ columns of $U$ are
the other left singular vectors of $A$ that are orthogonal to $b$
by the construction stated above, and the last $n-s$ columns of $V$ are
the other corresponding right singular vectors of $A$.
Obviously, for the new SVD of $A$ defined above, $A^{\prime}$ is of rank $s$
with the $s$ simple nonzero singular values $\sigma_1,\sigma_2,\ldots,\sigma_s$,
its left and right singular vector matrices $U$ and $V$ are the corresponding
ones of $A$ with proper column exchanges, respectively.
We have $x_{true}=A^{\dagger}\hat{b}=(A^{\prime})^{\dagger}\hat{b}$
and the TSVD regularized solutions $x_k^{tsvd}=(A_k^{\prime})^{\dagger}b$,
where $A_k^{\prime}$ is the best rank $k$ approximation to $A^{\prime}$
with respect to the 2-norm. In addition, we comment that from the discrete Picard
condition
$$
\|A^{\dagger}\hat{b}\|=\|(A^{\prime})^{\dagger}\hat{b}\|
=\left(\sum_{k=1}^s\frac{|u_k^T\hat{b}|^2}{\sigma_k^2}\right)^{1/2}\leq C,
$$
independently of $n$ and $s$,
we can obtain \eqref{newpicard} directly in the same way as done in
the introduction for \eqref{picard}.

Another fundamental change is that the $k$-dimensional dominant right
singular space of $A$ now becomes that of $A^{\prime}$, i.e.,
$\mathcal{V}_k=span\{V_k\}$ with $V_k=(v_1,v_2,\ldots,v_k)$ associated with
the first $k$ large singular values of $A^{\prime}$. It is the subspace of
concern in the case that $A$ has multiple singular values.
We will also denote $U_k=(u_1,u_2,\ldots,u_k)$ and $\mathcal{U}_k=span\{U_k\}$.
As for Krylov subspaces, by the SVD of $A$ and that of
$A^{\prime}$, expanding $b$ as $b=\sum_{j=1}^s\xi_j u_j+(I-U_sU_s^T)b$,
we easily justify
\begin{equation}\label{ata}
\mathcal{K}_k(A^TA, A^Tb)=\mathcal{K}_k((A^{\prime})^TA^{\prime},(A^{\prime})^Tb)
\end{equation}
and
\begin{equation}\label{aat}
\mathcal{K}_k(AA^T, b)=\mathcal{K}_k(A^{\prime}(A^{\prime})^T,b)
\end{equation}
by noting that
\begin{equation}\label{aprimea}
(A^TA)^i A^Tb=\left((A^{\prime})^TA^{\prime}\right)^i(A^{\prime})^Tb
=\sum_{j=1}^s \xi_j \sigma_j^{2i+1}v_j
\end{equation}
for any integer $i\geq 0$ and
\begin{equation}\label{aaprime}
(AA^T)^i b=\left(A^{\prime}(A^{\prime})^T\right)^ib
=\sum_{j=1}^s \xi_j \sigma_j^{2i}u_j
\end{equation}
for any integer $i\geq 1$. Thus, for the given $b$, Lanczos
bidiagonalization works on $A$ exactly as if it does on $A^{\prime}$.
That is, \eqref{eqmform1}--\eqref{Bk} generated by Algorithm 1 hold when
$A$ is replaced by $A^{\prime}$, and the $k$ Ritz values $\theta_i^{(k)}$
approximate $k$ nonzero singular values of $A^{\prime}$.
Moreover, \eqref{aprimea} and \eqref{aaprime} indicate
$$
\mathcal{K}_{s+1}((A^{\prime})^TA^{\prime},
(A^{\prime})^Tb)=\mathcal{K}_s((A^{\prime})^TA^{\prime},
(A^{\prime})^Tb),\
\mathcal{K}_{s+2}(A^{\prime}(A^{\prime})^T,b)=
\mathcal{K}_{s+1}(A^{\prime}(A^{\prime})^T,b).
$$
As a result, since $(A^{\prime})^Tb$ has nonzero components in all the
eigenvectors $v_1,v_2,\ldots,v_s$ of $(A^{\prime})^TA^{\prime}$ associated
with its nonzero distinct eigenvalues $\sigma_1^2,\sigma_2^2,\ldots,
\sigma_s^2$, Lanczos bidiagonalization cannot break down until step $s+1$,
and the singular values $\theta_i^{(s)}$ of $B_s$ are exactly the singular
values $\sigma_1,\sigma_2,\ldots,\sigma_s$ of $A^{\prime}$.
At step $s$, Lanczos bidiagonalization on $A$
generates the $(s+1)\times s$ lower bidiagonal matrix
\begin{align}
P_{s+1}^T AQ_s&=P_{s+1}^T A^{\prime} Q_s=B_s \label{lbaprime}
\end{align}
and
\begin{align}
\mathcal{V}_s&=span\{Q_s\},\ \mathcal{U}_s\subset
span\{P_{s+1}\}. \label{uv}
\end{align}

Having done the above, what we need is to estimate
how $\mathcal{K}_k(A^TA, A^Tb)=\mathcal{K}_k((A^{\prime})^TA^{\prime},
(A^{\prime})^Tb)$ approximates or captures the $k$-dimensional
dominant right subspaces $\mathcal{V}_k$, $k=1,2,\ldots,s-1$.
This is a crucial step and the starting point of all the later
analysis. In what follows let us show how to adapt the beginning part
of the proof of Theorem~\ref{thm2} to the case that $A$ has multiple singular
values.

Observe the Krylov subspace $\mathcal{K}_{k}((\Sigma^{\prime})^2,
\Sigma^{\prime} U^Tb)=span\{\hat{D}\hat{T}_k\}$ with
$$
  \hat{D}=\diag(\sigma_1 u_1^Tb,\ldots,\sigma_s u_s^Tb,
  \mathbf{0})=\left(\begin{array}{cc}
  D& \\
  & \mathbf{0}
  \end{array}\right)\in\mathbb{R}^{n\times n}
$$
and
$$
  \hat{T}_k=\left(\begin{array}{cccc} 1 &
  \sigma_1^2&\ldots & \sigma_1^{2k-2}\\
1 &\sigma_2^2 &\ldots &\sigma_2^{2k-2} \\
\vdots & \vdots&&\vdots\\
1 &\sigma_s^2 &\ldots &\sigma_s^{2k-2}\\
0 & 0 &\ldots& 0\\
\vdots & \vdots&&\vdots\\
0 & 0 &\ldots& 0
\end{array}\right)=\left(\begin{array}{c}
T_k\\
\mathbf{0}
\end{array}\right)\in \mathbb{R}^{n\times k}.
$$
Partition the diagonal matrix $D$ and the matrix $T_k$ into the forms
\begin{equation*}
  D=\left(\begin{array}{cc} D_1 & 0 \\ 0 & D_2 \end{array}\right)
  \in \mathbb{R}^{s\times s},\ \ \
  T_k=\left(\begin{array}{c} T_{k1} \\ T_{k2} \end{array}\right)
  \in \mathbb{R}^{s\times k},
\end{equation*}
where $D_1, T_{k1}\in\mathbb{R}^{k\times k}$ and
$D_2=\diag(u_{k+1}^Tb,\ldots,u_s^Tb)$. Since $T_{k1}$ is
a Vandermonde matrix with $\sigma_j$ distinct for $j=1,2,\ldots,k$, it is
nonsingular. Therefore, from
$$
\mathcal{K}_k((A^{\prime})^TA^{\prime},
(A^{\prime})^Tb)=span\{V\hat{D}\hat{T}_k\}
$$
and the structures of $\hat{D}$ and $\hat{T}_k$, we obtain
\begin{equation*}
\mathcal{K}_k((A^{\prime})^TA^{\prime},
(A^{\prime})^Tb)=span\{V_sD T_k\}=span \left\{V_s\left
(\begin{array}{c} D_1T_{k1} \\ D_2T_{k2} \end{array}\right)\right\}
=span\left\{V_s\left(\begin{array}{c} I \\ \Delta_k \end{array}\right)\right\},
\end{equation*}
with
\begin{equation*}
\Delta_k=D_2T_{k2}T_{k1}^{-1}D_1^{-1},
\end{equation*}
meaning that $\mathcal{K}_k((A^{\prime})^TA^{\prime},
(A^{\prime})^Tb)$ is orthogonal to the last $n-s$ columns of $V$.

Write
\begin{equation}\label{parti}
V_s=(V_k, V_k^{\perp}),\ \ V=(V_s,\hat{V}_s),
\end{equation}
and define
\begin{equation*}
Z_k=V_s\left(\begin{array}{c} I \\ \Delta_k \end{array}\right)
=V_k+V_k^{\perp}\Delta_k.
\end{equation*}
Then $Z_k^TZ_k=I+\Delta_k^T\Delta_k$, the
columns of $\hat{Z}_k=Z_k(Z_k^TZ_k)^{-\frac{1}{2}}$
form an orthonormal basis of $\mathcal{K}_k((A^{\prime})^TA^{\prime},
(A^{\prime})^Tb)$, and we get the orthogonal
direct sum decomposition
\begin{equation*}
\hat{Z}_k=(V_k+V_k^{\perp}\Delta_k)(I+\Delta_k^T\Delta_k)^{-\frac{1}{2}}.
\end{equation*}

Denote $\mathcal{V}_k^R=\mathcal{K}_k((A^{\prime})^TA^{\prime},
(A^{\prime})^Tb)$. For $\|\sin\Theta(\mathcal{V}_k, \mathcal{V}_k^R)\|$,
based on the above, we get \eqref{deltabound} by
replacing $V_k^{\perp}$ in \eqref{sindef} by
$(V_k^{\perp},\hat{V}_s)$ defined as \eqref{parti} and noting
that $\hat{Z}_k$ is orthogonal to $\hat{V}_s$. Then it is direct
to derive the same bounds for $\|\sin\Theta(\mathcal{V}_k, \mathcal{V}_k^R)\|$
as those established previously in completely the same way.

As for the extension of Theorem~\ref{initial}, by definition and \eqref{parti},
we need to replace $V_k^{\perp}$ in \eqref{qktilde}
by $(V_k^{\perp},\hat{V}_s)$ defined as \eqref{parti}.
The unit-length $\tilde{q}_k\in\mathcal{V}_k^R$ is now a vector that
has the smallest acute angle with $span\{(V_k^{\perp},\hat{V}_s)\}$,
and we modify \eqref{decompqk} as
$$
\tilde{q}_k=\hat{V}_s\hat{V}_s^T\tilde{q}_k
+V_k^{\perp}(V_k^{\perp})^T\tilde{q}_k+V_kV_k^T\tilde{q}_k.
$$
Recall that the columns of $\hat{V}_s$ are the right singular vectors
of $A^{\prime}$ corresponding to zero singular values. It disappears
when forming the Rayleigh quotient of $(A^{\prime})^T A^{\prime}$ with respect
to $\tilde{q}_k$. The proof of Theorem~\ref{initial} then
carries over to $A^{\prime}$, and the
results hold for the case that $A$ has multiple singular values.

Another fundamental change is that, when speaking of a rank $k$ approximation,
we now mean that for $A^{\prime}$. Note that the best rank $k$ approximation
$A^{\prime}_k$ to $A^{\prime}$ is $A^{\prime}_k=U_k\Sigma_k^{\prime}V_k^T$,
$k=1,2,\ldots,s$, where $U_k$ and $V_k$ are defined as before,
and $\Sigma_k^{\prime}={\rm diag}(\sigma_1,\sigma_2,\ldots,\sigma_k)$.
The $k$-step Lanczos bidiagonalzation process on $A$ now generates
rank $k$ approximations $P_{k+1}B_kQ_k^T$ in LSQR and $P_k\bar{B}_{k-1}Q_k$
in CGME to $A^{\prime}$ and the rank $k$ approximation
$Q_{k+1}Q_{k+1}^T(A^{\prime})^TA^{\prime}Q_kQ_k^T$ in LSMR
to $(A^{\prime})^TA^{\prime}$, where $P_{k+1}$ and $Q_k$ are
the first $k+1$ and $k$ columns of $P_{s+1}$ and $Q_s$ in \eqref{lbaprime}.
We then need to estimate the approximation accuracy of these rank
$k$ approximations and compare them with that of the best rank $k$
approximations $A_k^{\prime}$ and
$(A^{\prime}_k)^TA^{\prime}_k$, respectively. Meanwhile, for each of these
three rank $k$ approximation matrices, we need to
analyze how its $k$ nonzero singular values approximate $k$ singular
values of $A^{\prime}$ or $(A^{\prime})^T A^{\prime}$.

For the rank $k$ approximation $P_{k+1}B_kQ_k^T$ to $A^{\prime}$
in LSQR, similar to \eqref{gammak}, we define
$$
\gamma_k^{\prime}=\|A^{\prime}-P_{k+1}B_kQ_k^T\|,\ k=1,2,\ldots,s-1.
$$
Then, without any changes but the replacement of the index $n$ by $s$,
all the results in Section \ref{rankapp} and \eqref{leftsub}
in Theorem~\ref{main2} carry over to the multiple singular
value case.

The final important note is how to extend the results presented
Section \ref{decay}--\ref{cgme} to the multiple singular
value case. We have to derive the three key relations similar
to \eqref{invar}, \eqref{gk} and \eqref{gk1}, where the fact that Lanczos
bidiagonalization can be run to $n$ steps without breakdown is exploited.
In the case that $A$ has multiple singular values,
since Lanczos diagonalization on
$A$ must break down at step $s+1$, there are no $P_{n+1}$
and $Q_n$ as in \eqref{invar}. To this end, from \eqref{lbaprime} we augment
$P_{s+1}$ and $Q_s$ to the $m\times m$ and $n\times n$
orthogonal matrices $P=(P_{s+1},\hat{P})$ and $Q=(Q_s,\hat{Q})$,
respectively, from which and \eqref{lbaprime} we obtain
$$
P^TAQ=P^TA^{\prime}Q=\left(\begin{array}{cc}
B_s&\mathbf{0}\\
\mathbf{0}& \mathbf{0}
\end{array}
\right).
$$
Having this relation, like \eqref{invar} and \eqref{gk}, we get
$$
\gamma_k^{\prime}=\|A^{\prime}-P_{k+1}B_kQ_k^T\|
=\|P^T\left(A^{\prime}-P_{k+1}B_kQ_k^T\right)Q\|=\|G_k\|,
$$
where $G_k$ is the right bottom $(s-k+1)\times (s-k)$ matrix of
$B_s$, similar to $G_k$ in \eqref{gk1}. Then Theorem~\ref{main2}
extends naturally to the multiple singular value case without
any change but the replacement of the index $n$ by $s$, and all the
other results and analysis in Section \ref{decay}--\ref{cgme} carry over
to this case as well. The results in Section \ref{compare} hold without any change
whenever $A$, $A_k$ and the index $n$ are replaced by $A^{\prime}$,
$A^{\prime}_k$ and $s$, respectively.

In summary, based on the above reformulations, changes and preliminary work,
except Section \ref{rrqr}, we have extended all the results and analysis
in Sections \ref{lsqr}--\ref{compare}
to the case that $A$ has multiple singular values, just as we have done
for the simple singular value case. In the analysis, derivation and results,
the index $n$ is often replaced by $s$ whenever needed, and when this is
necessary is clear from the context related.

\section{Numerical experiments}\label{numer}

For a number of problems from Hansen's regularization toolbox \cite{hansen07},
Huang and Jia \cite{huangjia} have numerically justified the full
regularization of LSQR for severely and moderately ill-posed problems and its
partial regularization for mildly ill-posed problems, where each $A$ is
$1,024\times 1,024$.
In this section, we report numerical experiments to
confirm our theory and illustrate the full or partial regularization
of LSQR in much more detail. For the first two kinds of problems,
we demonstrate that $\gamma_k,\ \alpha_{k+1}$
and $\beta_{k+2}$ decay as fast as $\sigma_{k+1}$. We compare LSQR
and the hybrid LSQR with the TSVD method applied to projected
problems after semi-convergence. In the experiments, we use the L-curve
criterion, the function $\mathsf{lcorner}$ in \cite{hansen07}, to determine
an actually optimal regularization parameter.
For each of severely and moderately ill-posed problems, we show that the
regularized solution obtained by LSQR at semi-convergence is at least as
accurate as the best TSVD regularized solution, indicating
that LSQR has the full regularization. In the meantime, we
show that the regularized solution obtained by LSQR at semi-convergence
is considerably less accurate than that by the hybrid LSQR for mildly ill-posed
problems, demonstrating that LSQR has only the partial regularization.
As a byproduct, we compare LSQR with GMRES and RRGMRES and illustrate that the
latter ones have no regularizing effects for general nonsymmetric
ill-posed problems.

We choose several ill-posed problems from Hansen's regularization
toolbox \cite{hansen07}, which include
the severely ill-posed problems $\mathsf{shaw,\ wing,\ i\_laplace}$,
the moderately ill-posed problems $\mathsf{heat,\ phillips}$, and
the mildly ill-posed problem $\mathsf{deriv2}$.
All the codes are from \cite{hansen07}, and the problems
arise from discretizations of \eqref{eq2}.
We remind that, as far as solving \eqref{eq1} is concerned,
our primary goal consists in justifying the regularizing effects
of iterative solvers for \eqref{eq1}, which are {\em unaffected by the size}
of \eqref{eq1} and only depends on the degree of
ill-posedness, the noise level $\|e\|$ and the actual discrete Picard
condition, provided that the condition number of \eqref{eq1},
measured by the ratio between the largest and smallest singular values
of each $A$, is large enough.
Therefore, for this purpose, as extensively done in the
literature (see, e.g., \cite{hansen98,hansen10} and the references therein
as well as many other papers), it is enough to report the results on small
and/or medium sized discrete ill-posed problems since the condition numbers of
these $A$ are already huge or large, which, in finite precision arithmetic,
are roughly $10^{16}, 10^{8}$ and $10^{6}$ for severely, moderately and
mildly ill-posed problems with $n=256$, respectively.
Indeed, for $n$ large, say, 10,000 or more, we have observed that
LSQR and the hybrid LSQR have the same behavior as for small $n$,
e.g., $n=256$ used in this paper. Also,
an important reason is that such choice enables us to fully
justify the regularization effects of LSQR by comparing it with
the TSVD method, which suits only for small and/or medium sized problems
because of its computational complexity for $n$ large. For each example,
we generate a $256\times 256$ matrix $A$, the true
solution $x_{true}$ and noise-free right-hand side $\hat{b}$.
In order to simulate the noisy data, we generate white noise
vectors $e$ such that the relative noise levels
$\varepsilon=\frac{\|e\|}{\|\hat{b}\|}=10^{-2}, 10^{-3}, 10^{-4}$, respectively.
We mention that, to better illustrate the behavior of the hybrid LSQR,
we, in the concluding section, will report some important
observations on $\mathsf{phillips}$ and $\mathsf{deriv2}$ of $n=1,024$ and $10,240$,
whose condition numbers are as large as $1.7\times 10^{15}$ and $1.3\times 10^8$
for $n=10,240$, respectively. To simulate exact arithmetic,
LSQR uses full reorthogonalization in Lanczos bidiagonalization.
All the computations are carried out in Matlab 7.8 with the machine precision
$\epsilon_{\rm mach}= 2.22\times10^{-16}$ under the Miscrosoft
Windows 7 64-bit system.

\subsection{The accuracy of rank $k$ approximations}

{\bf Example 1}.
This problem $\mathsf{shaw}$ arises from one-dimensional image restoration and is
obtained by discretizing \eqref{eq2} with $[-\frac{\pi}{2}, \frac{\pi}{2}]$ as
the domains of $s$ and $t$, where
\begin{align*}
  k(s,t)&=(\cos(s)+\cos(t))^2\left(\frac{\sin(u)}{u}\right)^2,\
  u=\pi(\sin(s)+\sin(t)),\\
x(t)&=2\exp(-6(t-0.8)^2)+\exp(-2(t+0.5)^2).
\end{align*}

{\bf Example 2}. This problem $\mathsf{wing}$ has  a discontinuous solution
and is obtained by discretizing \eqref{eq2} with $[0, 1]$ as
the domains of $s$ and $t$, where
\begin{align*}
  k(s,t)&=t\exp(-st^2),\ \ \ g(s)=\frac{\exp(-\frac{1}{9}s)-\exp(-\frac{4}{9}s)}{2s},\\
  x(t)&=\left\{\begin{array}{ll} 1,\ \ \ &\frac{1}{3}<t<\frac{2}{3};\\ 0,\ \ \
  &elsewhere. \end{array}\right.
\end{align*}

The problems $\mathsf{shaw}$ and $\mathsf{wing}$ are severely ill-posed
with the singular values $\sigma_k=\mathcal{O}(e^{-4 k})$
for $\mathsf{shaw}$ and $\sigma_k=\mathcal{O}(e^{-9k})$ for $\mathsf{wing}$,
respectively.

\begin{figure}
\begin{minipage}{0.48\linewidth}
  \centerline{\includegraphics[width=6.0cm,height=5cm]{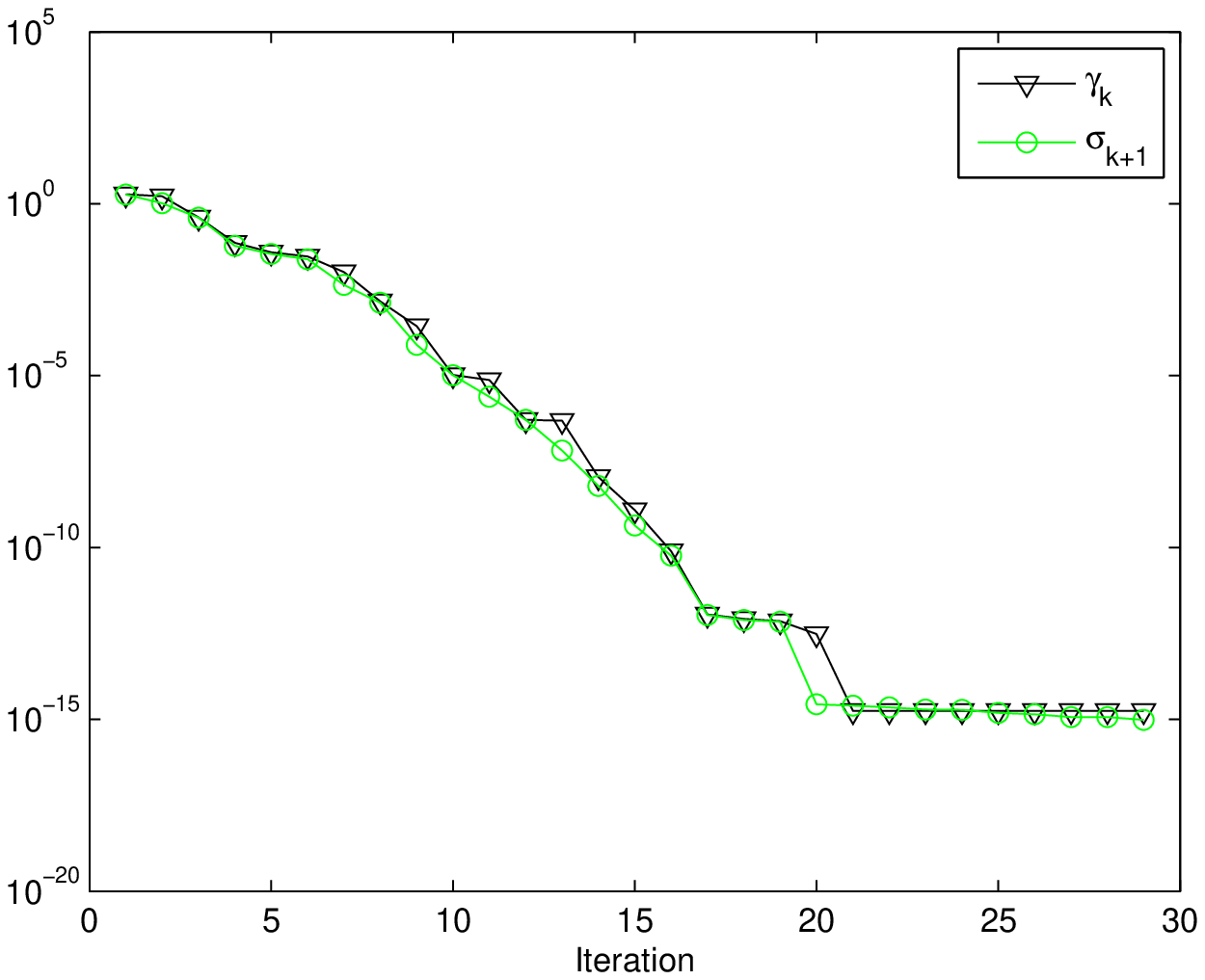}}
  \centerline{(a)}
\end{minipage}
\hfill
\begin{minipage}{0.48\linewidth}
  \centerline{\includegraphics[width=6.0cm,height=5cm]{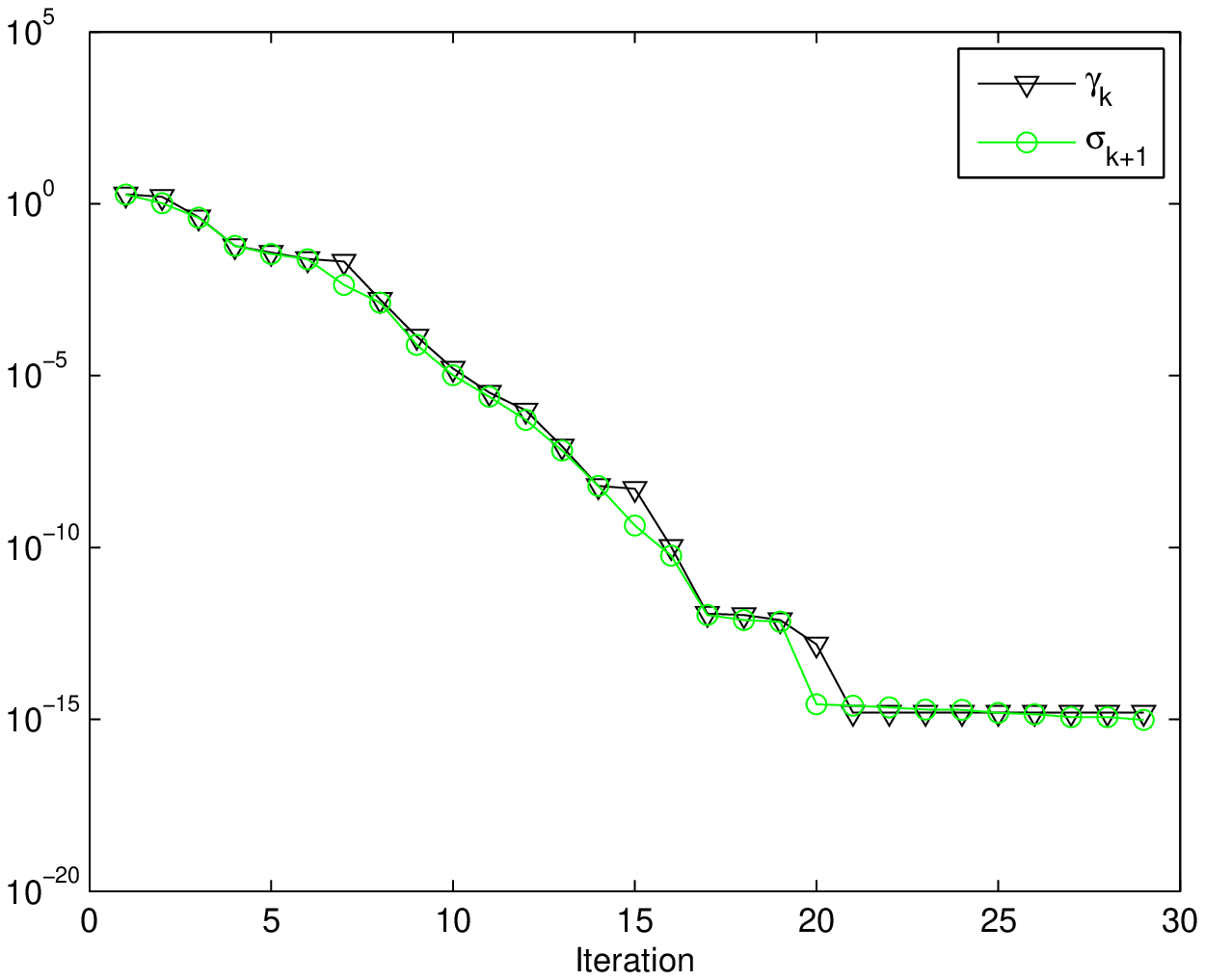}}
  \centerline{(b)}
\end{minipage}
\vfill
\begin{minipage}{0.48\linewidth}
  \centerline{\includegraphics[width=6.0cm,height=5cm]{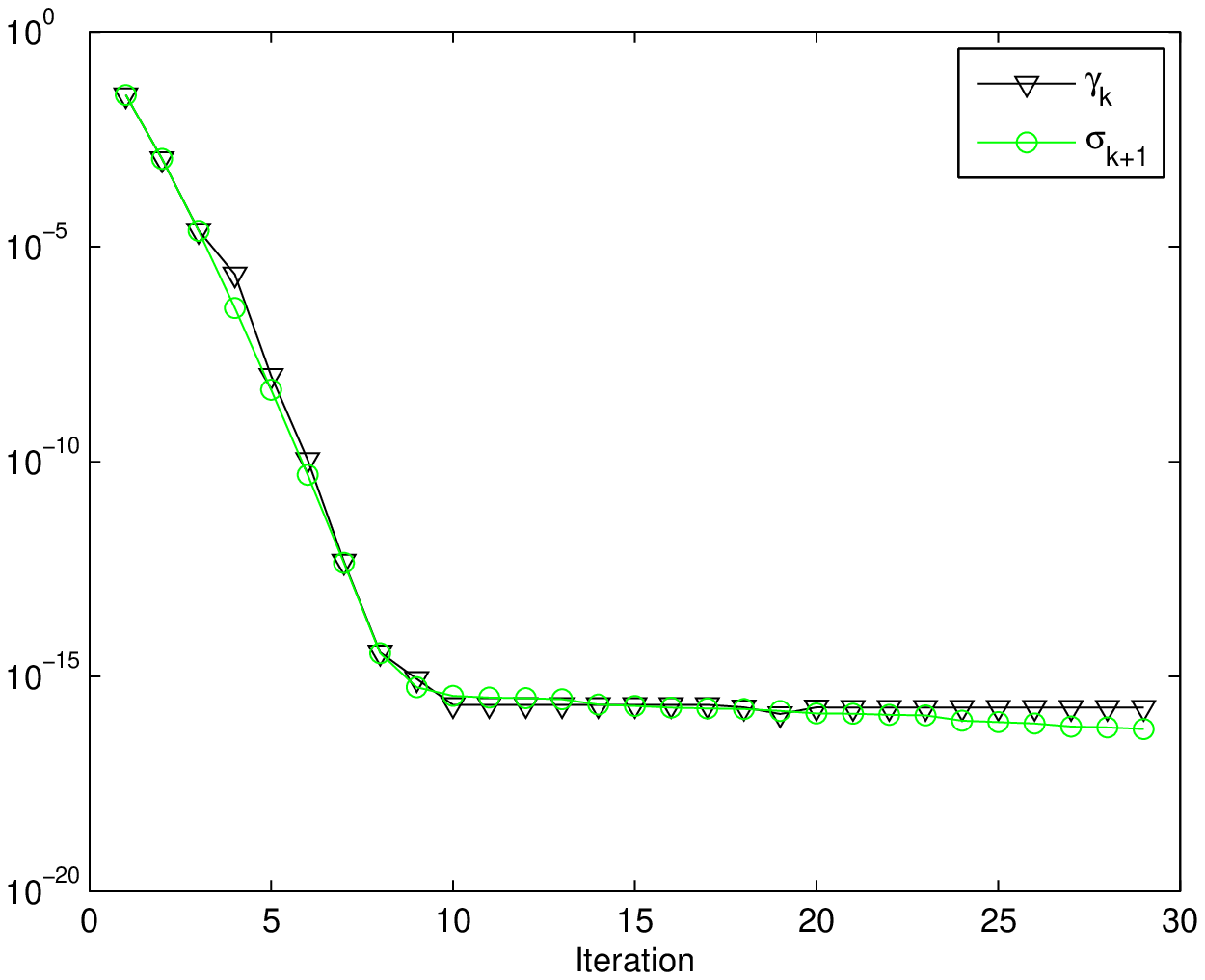}}
  \centerline{(c)}
\end{minipage}
\hfill
\begin{minipage}{0.48\linewidth}
  \centerline{\includegraphics[width=6.0cm,height=5cm]{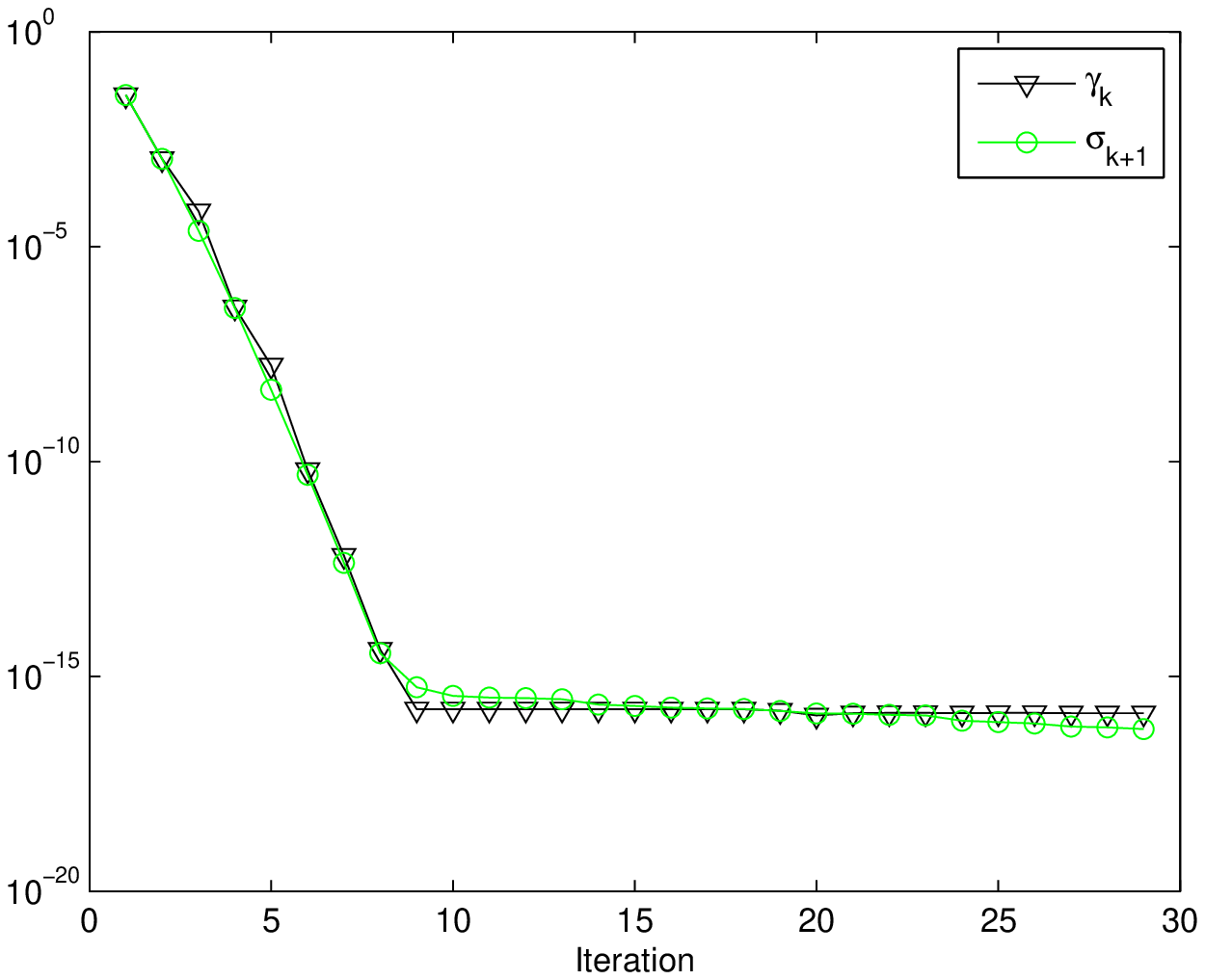}}
  \centerline{(d)}
\end{minipage}
\caption{(a)-(b): Decay curves of the sequences $\gamma_k$ and
$\sigma_{k+1}$ for $\mathsf{shaw}$ with $\varepsilon=10^{-2}$
(left) and $\varepsilon=10^{-3}$ (right); (c)-(d): Decay curves of the
sequences $\gamma_k$ and $\sigma_{k+1}$ for $\mathsf{wing}$ with
$\varepsilon=10^{-3}$ (left) and $\varepsilon=10^{-4}$ (right).} \label{fig1}
\end{figure}

In Figure~\ref{fig1}, we display the decay curves of the $\gamma_k$ for
$\mathsf{shaw}$ with $\varepsilon=10^{-2}, 10^{-3}$ and for $\mathsf{wing}$
with $\varepsilon=10^{-3}, 10^{-4}$, respectively.
We observe that the three curves with different $\varepsilon$ are almost
unchanged. This is in accordance with our Remark~\ref{decayrate}, where
it is stated that the decay rate of $\gamma_k$ is little affected
by noise levels for severely ill-posed problems, since
$\gamma_k$ primarily depends on the decay rate of $\sigma_{k+1}$
and different noise levels only affect the value of
$k_0$ other than the decay rate of $\gamma_k$. In addition,
we have observed that $\gamma_k$ and $\sigma_{k+1}$ decay
until they level off at $\epsilon_{\rm mach}$ due to round-off errors.
Most importantly, the results have clearly confirmed the theory
that $\gamma_k$ decreases as fast as $\sigma_{k+1}$, and we have
$\gamma_k\approx\sigma_{k+1}$, whose decay curves
are almost indistinguishable.

\begin{figure}
\begin{minipage}{0.48\linewidth}
  \centerline{\includegraphics[width=6.0cm,height=5cm]{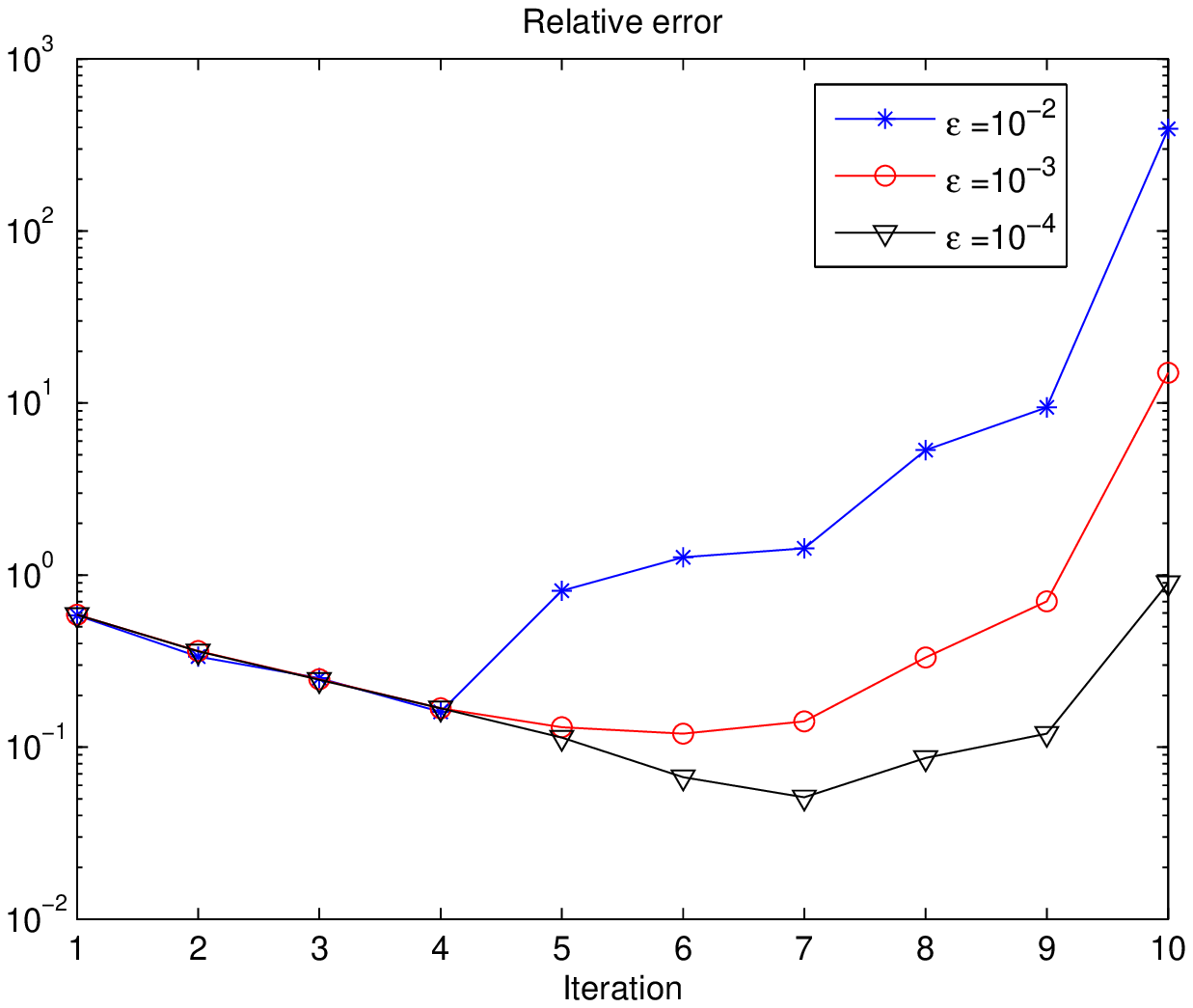}}
  \centerline{(a)}
\end{minipage}
\hfill
\begin{minipage}{0.48\linewidth}
  \centerline{\includegraphics[width=6.0cm,height=5cm]{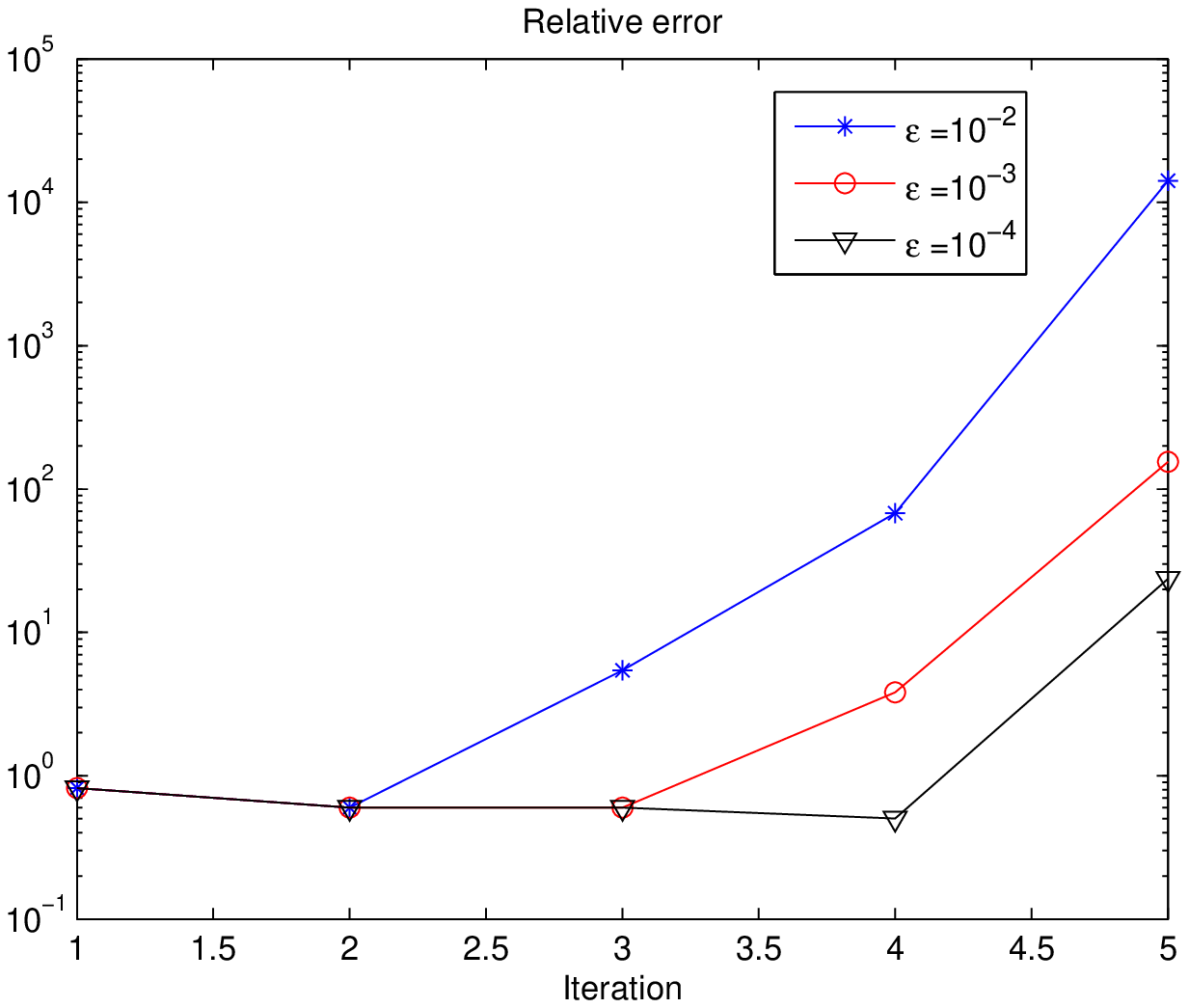}}
  \centerline{(b)}
\end{minipage}
\caption{ The relative errors $\|x^{(k)}-x_{true}\|/
\|x_{true}\|$ with $\varepsilon=10^{-2}, 10^{-3}, 10^{-4}$
for $\mathsf{shaw}$ (left) and $\mathsf{wing}$ (right).}
\label{fig2}
\end{figure}

In Figure~\ref{fig2}, we plot the relative errors $\|x^{(k)}
-x_{true}\|/\|x_{true}\|$ with different $\varepsilon$ for
these two problems. As we have seen, LSQR exhibits clear semi-convergence.
Moreover, for a smaller $\varepsilon$, we get a more accurate regularized
solution at cost of more iterations, as $k_0$ is bigger from \eqref{picard}
and \eqref{picard1}.

{\bf Example 3}.
This problem $\mathsf{heat}$ is moderately ill-posed, arises from the inverse
heat equation, and is obtained by discretizing \eqref{eq2}
with $[0, 1]$ as integration interval, where the kernel $k(s,t)=k(s-t)$ with
\begin{equation*}
  k(t)=\frac{t^{-3/2}}{2\sqrt{\pi}}\exp\left(-\frac{1}{4t}\right).
\end{equation*}

{\bf Example 4}.
This is the $\mathsf{phillips}$ famous problem, a moderately ill-posed one.
It is obtained by discretizing  \eqref{eq2}
with $[-6, 6]$ as the domains of $s$ and $t$, where
\begin{align*}
  k(s,t)&=\left\{\begin{array}{ll} 1+\cos\left(\frac{\pi(s-t)}{3}\right),
  \ \ \ &|s-t|<3,\\ 0,\ \ \ &|s-t|\geq 3, \end{array}\right.
\\
 g(s)&=(6-|s|)\left(1+\frac{1}{2}\cos\left(\frac{\pi s}{3}\right)\right)+
  \frac{9}{2\pi}\sin\left(\frac{\pi|s|}{3}\right),\\
 x(t)&=\left\{\begin{array}{ll} 1+\cos\left(\frac{\pi t}{3}\right),\ \ \ &|t|<3,\\
  0,\ \ \ &|t|\geq 3. \end{array}\right.
\end{align*}

\begin{figure}
\begin{minipage}{0.48\linewidth}
  \centerline{\includegraphics[width=6.0cm,height=5cm]{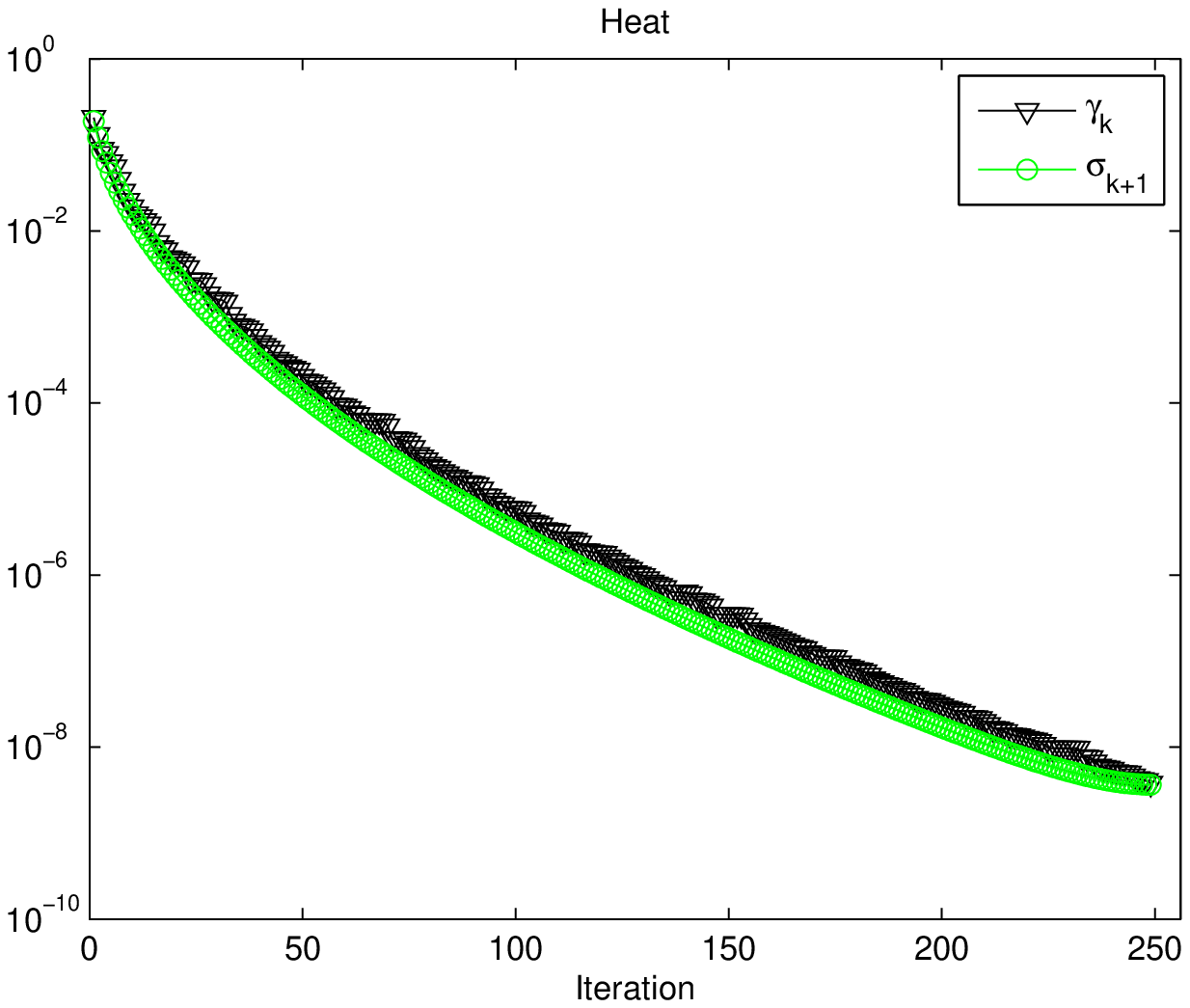}}
  \centerline{(a)}
\end{minipage}
\hfill
\begin{minipage}{0.48\linewidth}
  \centerline{\includegraphics[width=6.0cm,height=5cm]{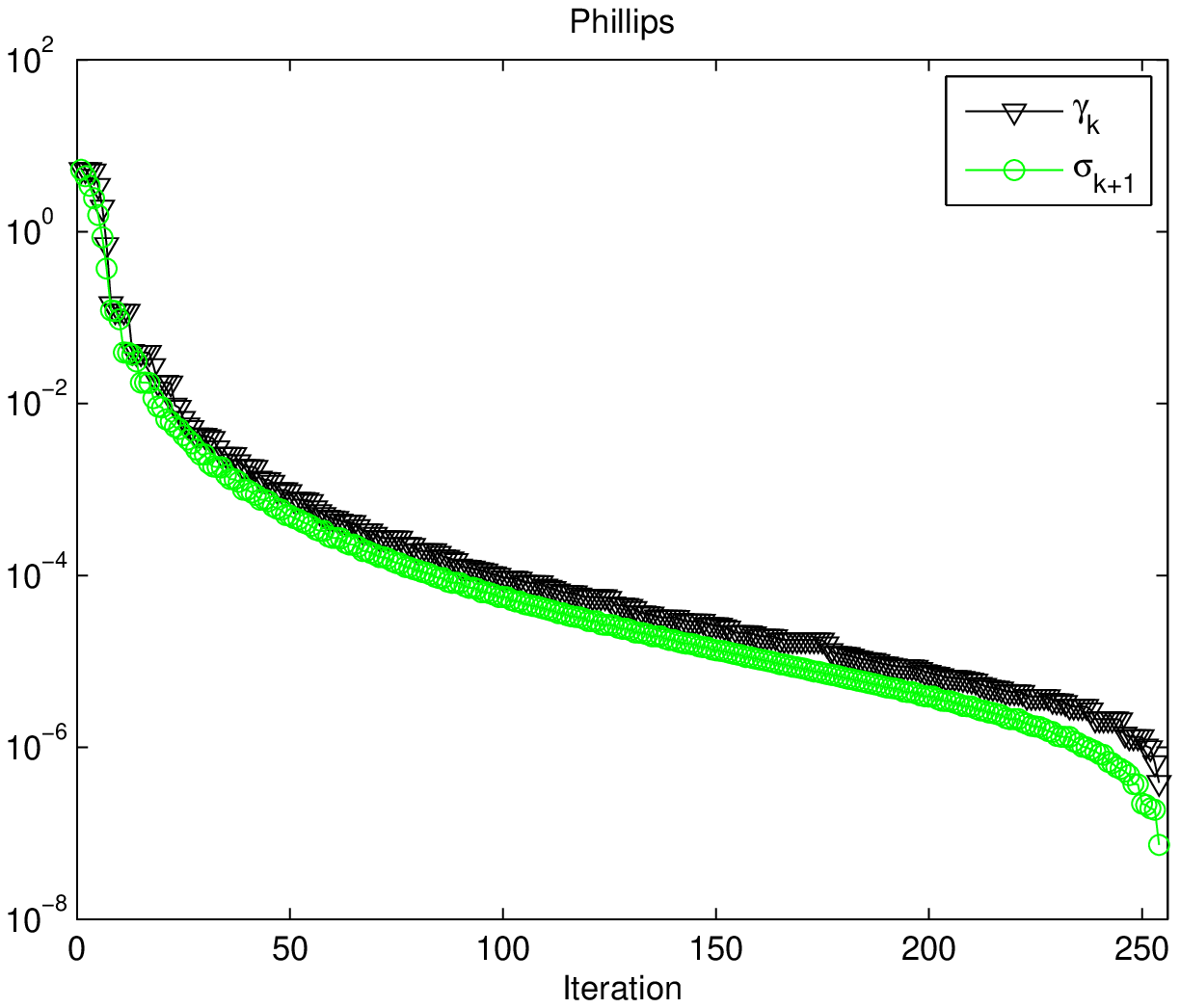}}
  \centerline{(b)}
\end{minipage}
\caption{(a): Decay curves of the sequences $\gamma_k$ and $\sigma_{k+1}$
for $\mathsf{heat}$ with (left) and (b): Decay curves of the
sequences $\gamma_k$ and $\sigma_{k+1}$ for
$\mathsf{phillips}$ with $\varepsilon=10^{-3}$ (right).}\label{fig3}
\end{figure}

From Figure~\ref{fig3}, we see that $\gamma_k$ decreases almost as fast as
$\sigma_{k+1}$ for the moderately ill-posed problems $\mathsf{heat}$ and
$\mathsf{phillips}$. However, slightly different
from severely ill-posed problems, $\gamma_k$, though
excellent approximations to $\sigma_{k+1}$, may not be so very accurate.
This is expected, as the constants $\eta_k$ in \eqref{const2} are generally
bigger than those in \eqref{const1} for severely ill-posed problems. Also,
different from Figure~\ref{fig1}, we observe from Figure~\ref{fig3} that
$\gamma_k$ deviates more from $\sigma_{k+1}$ with $k$ increasing, especially
for the problem $\mathsf{phillips}$. This confirms
Remarks~\ref{decayrate}--\ref{decayrate2} on moderately ill-posed problems.

\begin{figure}
\begin{minipage}{0.48\linewidth}
  \centerline{\includegraphics[width=6.0cm,height=5cm]{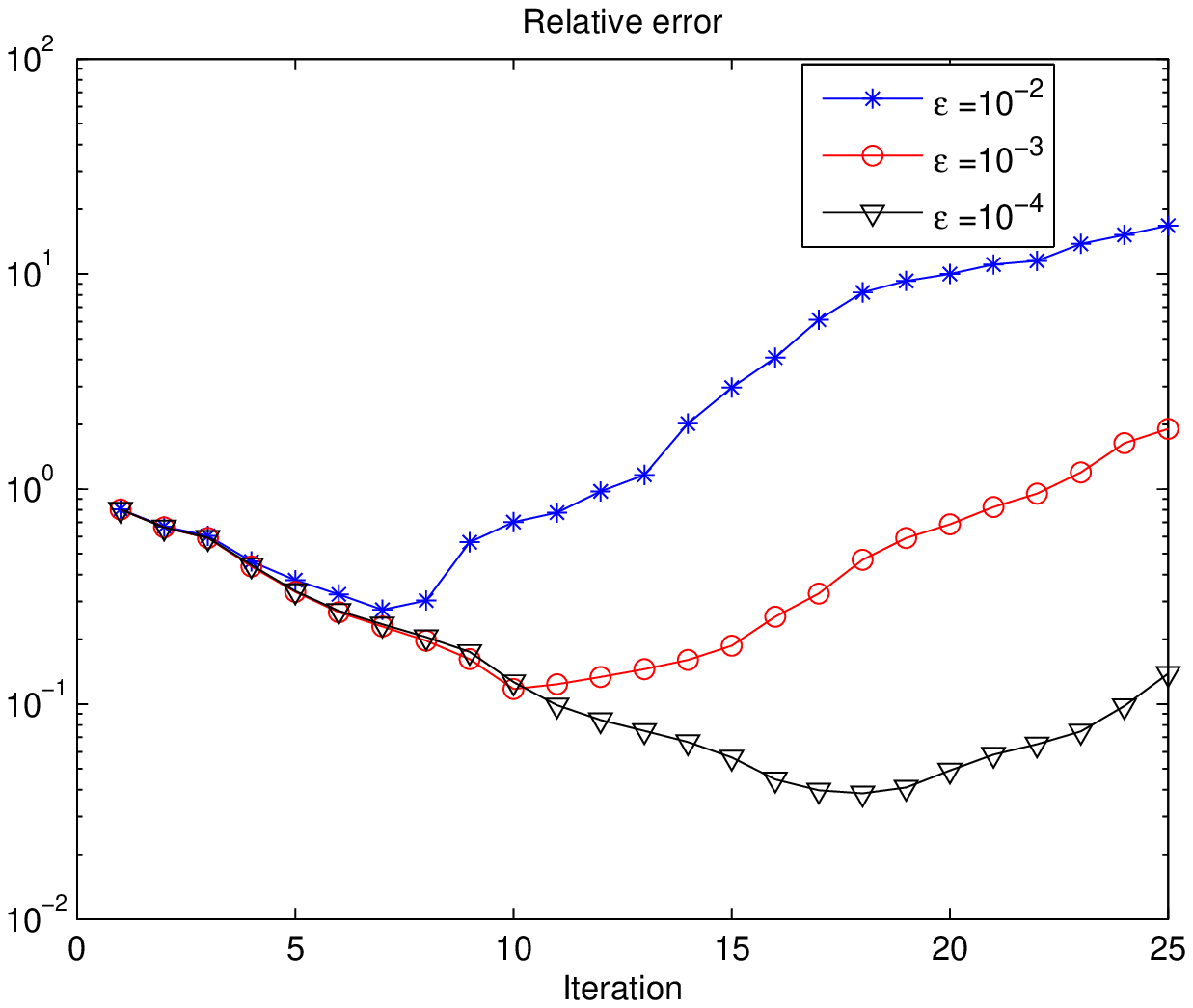}}
  \centerline{(a)}
\end{minipage}
\hfill
\begin{minipage}{0.48\linewidth}
  \centerline{\includegraphics[width=6.0cm,height=5cm]{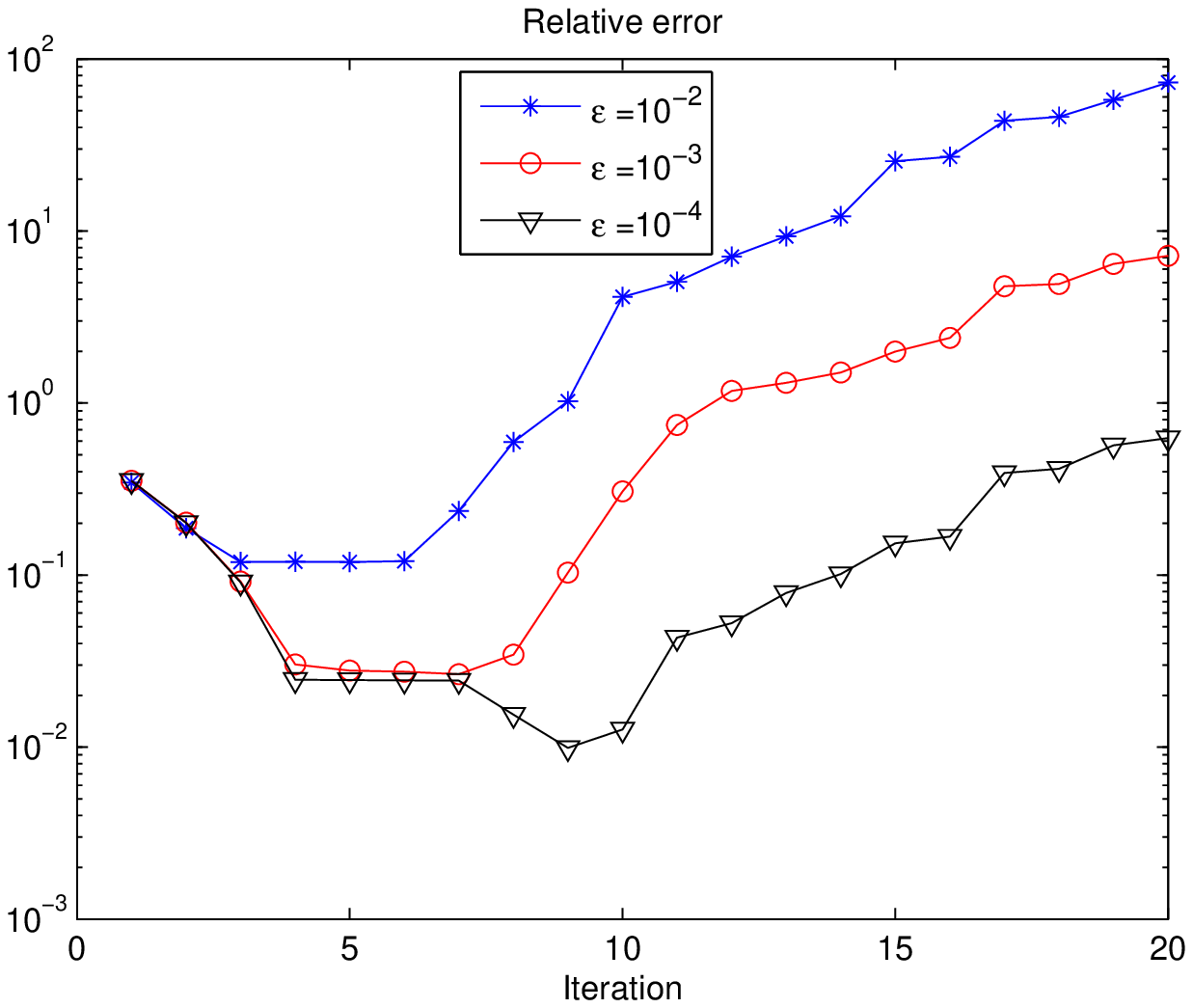}}
  \centerline{(b)}
\end{minipage}
\caption{ The relative errors $\|x^{(k)}-x_{true}\|/\|x_{true}\|$
with $\varepsilon=10^{-2}, 10^{-3}, 10^{-4}$ for
$\mathsf{heat}$ (left) and $\mathsf{phillips}$ (right).} \label{fig4}
\end{figure}

In Figure~\ref{fig4}, we depict the relative errors of $x^{(k)}$, and from them
we observe analogous phenomena to those for severely ill-posed problems.
The only distinction is that LSQR now needs more iterations, i.e.,
a bigger $k_0$ is needed for moderately ill-posed problems with the same
$\varepsilon$, as is seen from \eqref{picard} and \eqref{picard1}.

{\bf Example 5}.
The mildly ill-posed problem $\mathsf{deriv2}$ is obtained by
discretizing \eqref{eq2} with $[0, 1]$ as the domains of $s$ and $t$,
where the kernel $k(s,t)$ is the
Green's function for the second derivative:
\begin{equation*}\label{}
  k(s,t)=\left\{\begin{array}{ll} s(t-1),\ \ \ &s<t;\\ t(s-1),\ \ \
  &s\geq t, \end{array}\right.
\end{equation*}
and the solution $x(t)$ and
the right-hand side $g(s)$ are given by
\begin{equation*}\label{}
  x(t)=\left\{\begin{array}{ll} t,
  \ \ \ &t<\frac{1}{2};\\ 1-t,\ \ \ &t\geq\frac{1}{2}, \end{array}\right. \ \ \
  g(s)=\left\{\begin{array}{ll} (4s^3-3s)/24,\ \ \ &s<\frac{1}{2};
  \\ (-4s^3+12s^2-9s+1)/24,\ \ \ &s\geq\frac{1}{2}. \end{array}\right.
\end{equation*}

\begin{figure}
\begin{minipage}{0.48\linewidth}
  \centerline{\includegraphics[width=6.0cm,height=5cm]{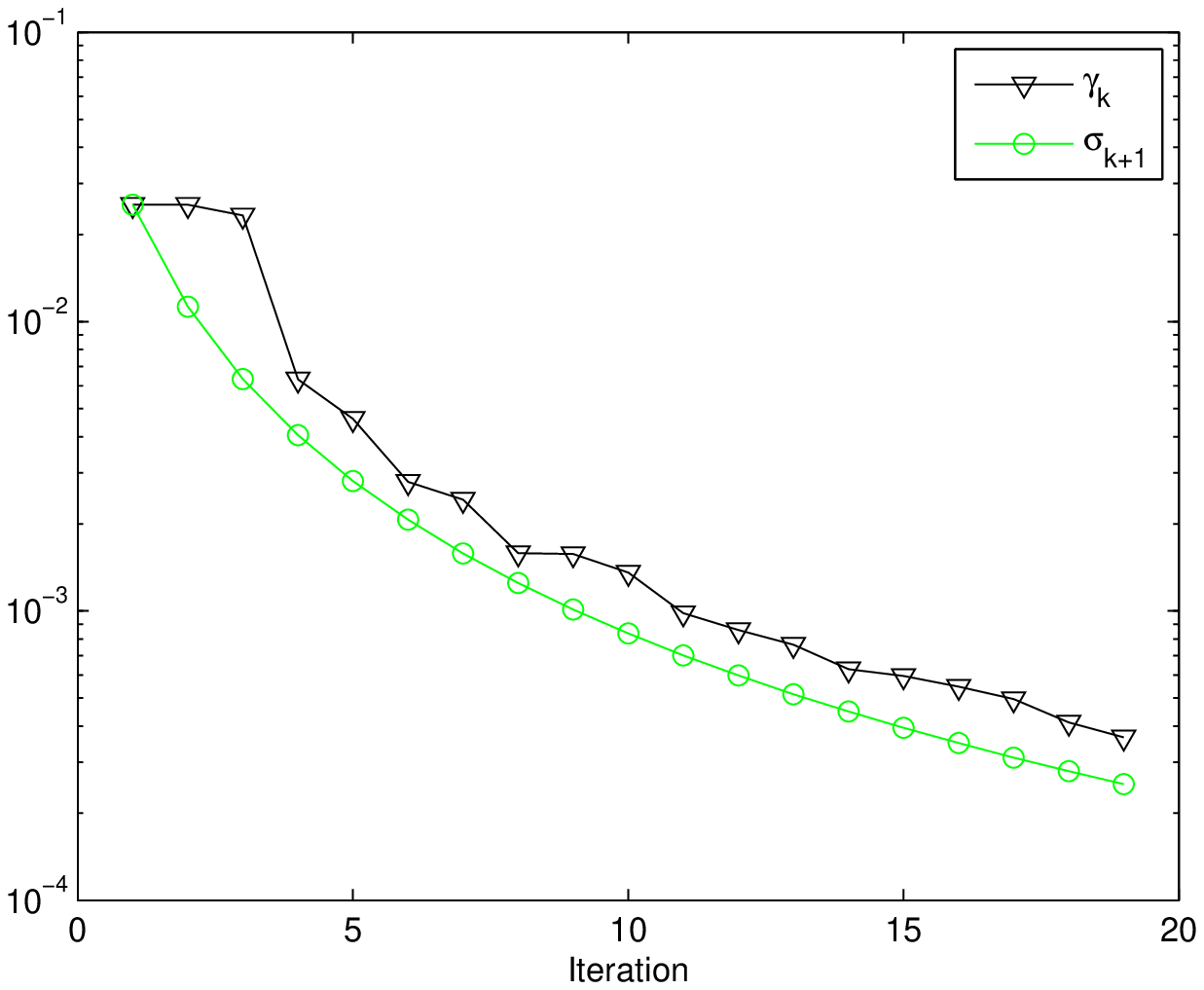}}
  \centerline{(a)}
\end{minipage}
\hfill
\begin{minipage}{0.48\linewidth}
  \centerline{\includegraphics[width=6.0cm,height=5cm]{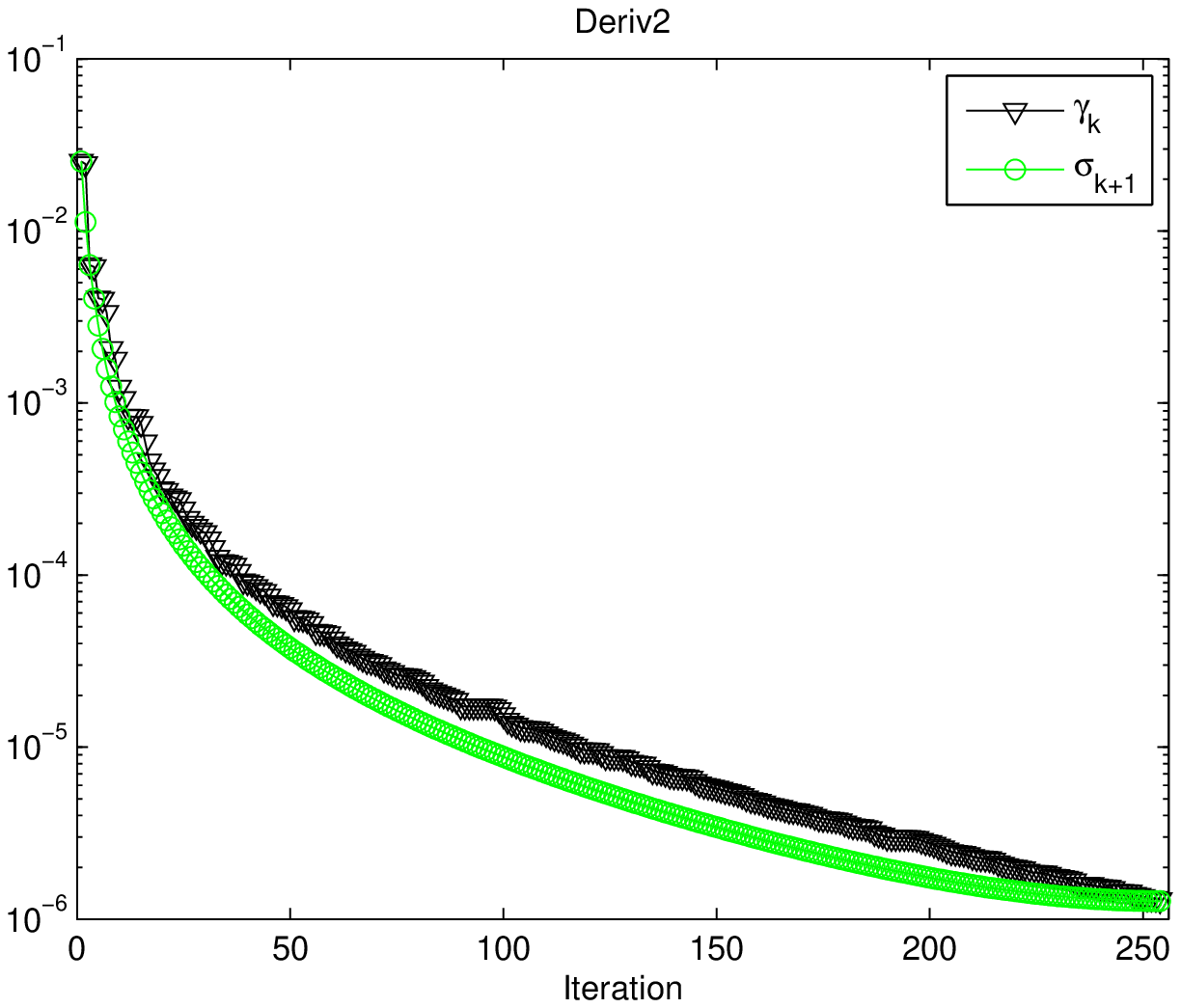}}
  \centerline{(b)}
\end{minipage}
\caption{(a)-(b): Decay curves of the partial and complete sequences
$\gamma_k$ and $\sigma_{k+1}$ for $\mathsf{deriv2}$ with
$\varepsilon=10^{-3}$}
\label{figmild}
\end{figure}

Figure~\ref{figmild} (a)-(b) display
the decay curves of the partial and complete sequences $\gamma_k$ and
$\sigma_{k+1}$, respectively. We see that, different from severely and
moderately ill-posed problems, $\gamma_k$ does not decay so fast as
$\sigma_{k+1}$ and deviates from $\sigma_{k+1}$ significantly.
Recall that Theorem~\ref{main1} holds for mildly
ill-posed problems, where $\eta_k$ defined by \eqref{const2} is considerably
bigger than one. These observations justify our theory and confirm
that the rank $k$ approximations to $A$ generated by Lanczos
bidiagonalization are not as accurate as those for severely and moderately
problems.

\subsection{A comparison of LSQR and the hybrid LSQR}

For the severely ill-posed $\mathsf{shaw}$, $\mathsf{wing}$
and the moderately ill-posed $\mathsf{heat}$, $\mathsf{phillips}$,
we compare the regularizing effects of LSQR and the hybrid LSQR with the TSVD
method applied to the projected problems after semi-convergence,
and demonstrate that they compute the same best possible regularized
solution for each problem and LSQR thus
has the full regularization. For the mildly ill-posed problem
$\mathsf{deriv2}$, we show that LSQR has only the partial
regularization and the hybrid LSQR can compute a best possible regularized
solution.

In the sequel, we report the results only for
the noise level $\varepsilon=10^{-3}$. Results for the other two
$\varepsilon$ are analogous and thus omitted unless stated otherwise.

\begin{figure}
\begin{minipage}{0.48\linewidth}
  \centerline{\includegraphics[width=6.0cm,height=5cm]{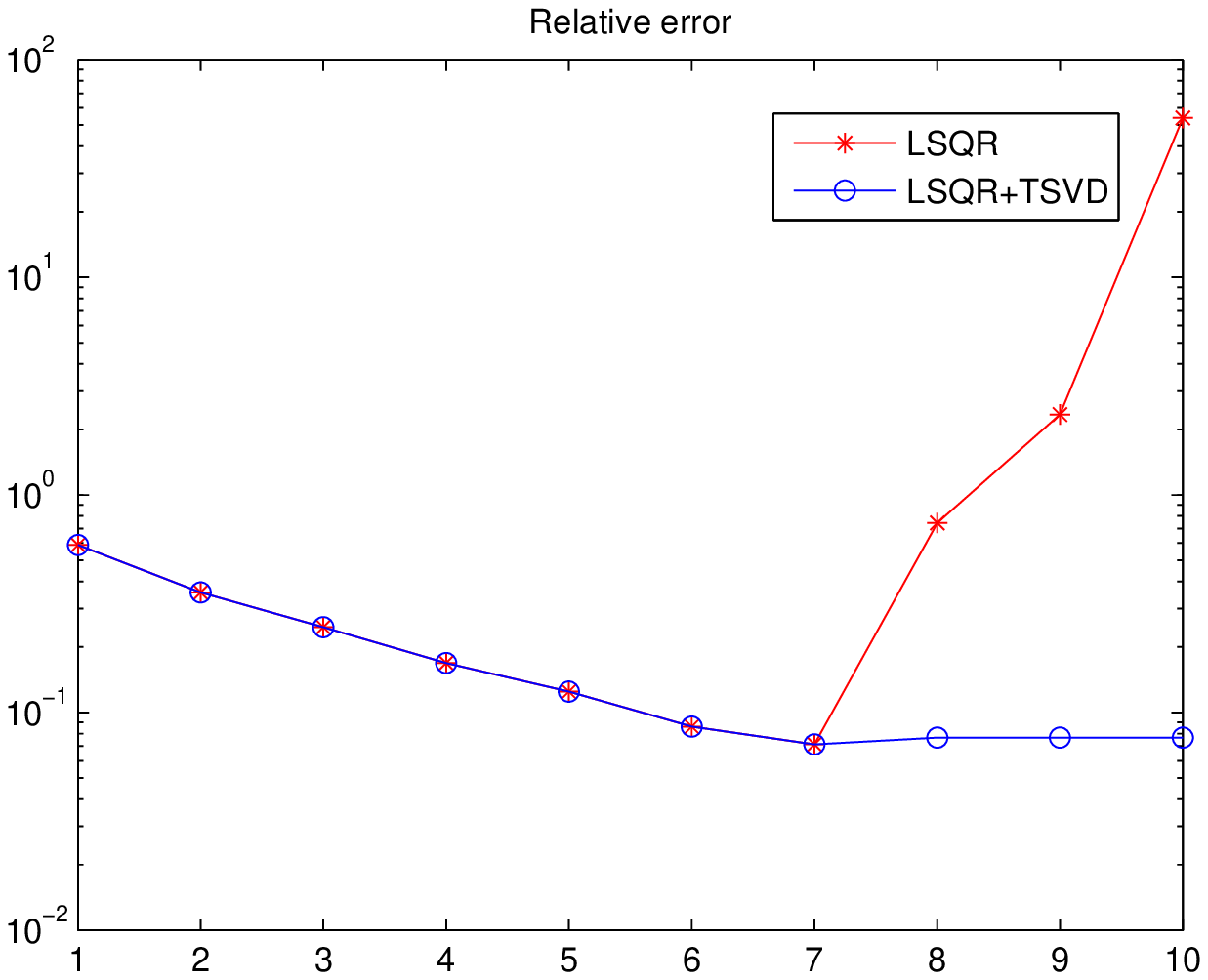}}
  \centerline{(a)}
\end{minipage}
\hfill
\begin{minipage}{0.48\linewidth}
  \centerline{\includegraphics[width=6.0cm,height=5cm]{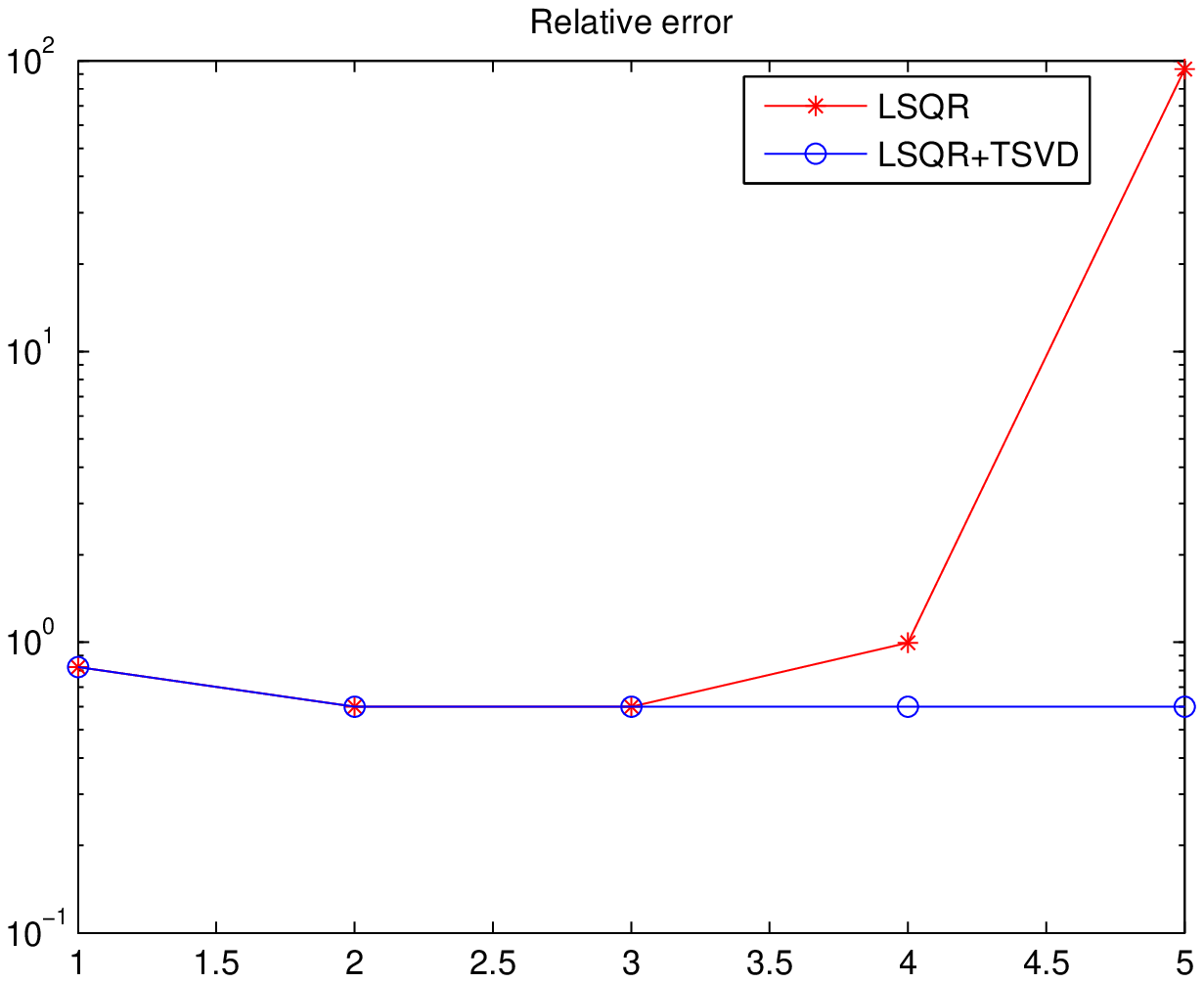}}
  \centerline{(b)}
\end{minipage}
\vfill
\begin{minipage}{0.48\linewidth}
  \centerline{\includegraphics[width=6.0cm,height=5cm]{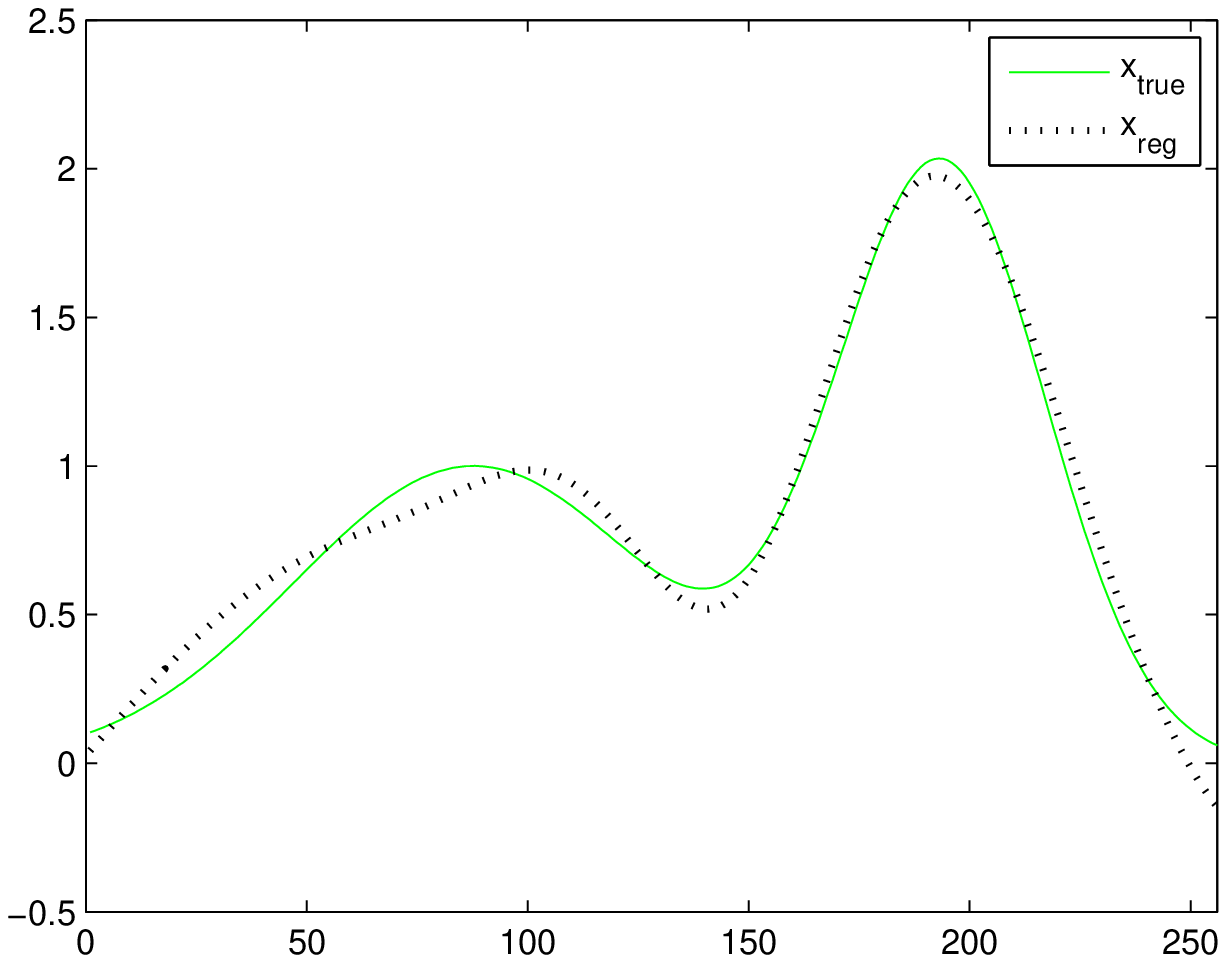}}
  \centerline{(c)}
\end{minipage}
\hfill
\begin{minipage}{0.48\linewidth}
  \centerline{\includegraphics[width=6.0cm,height=5cm]{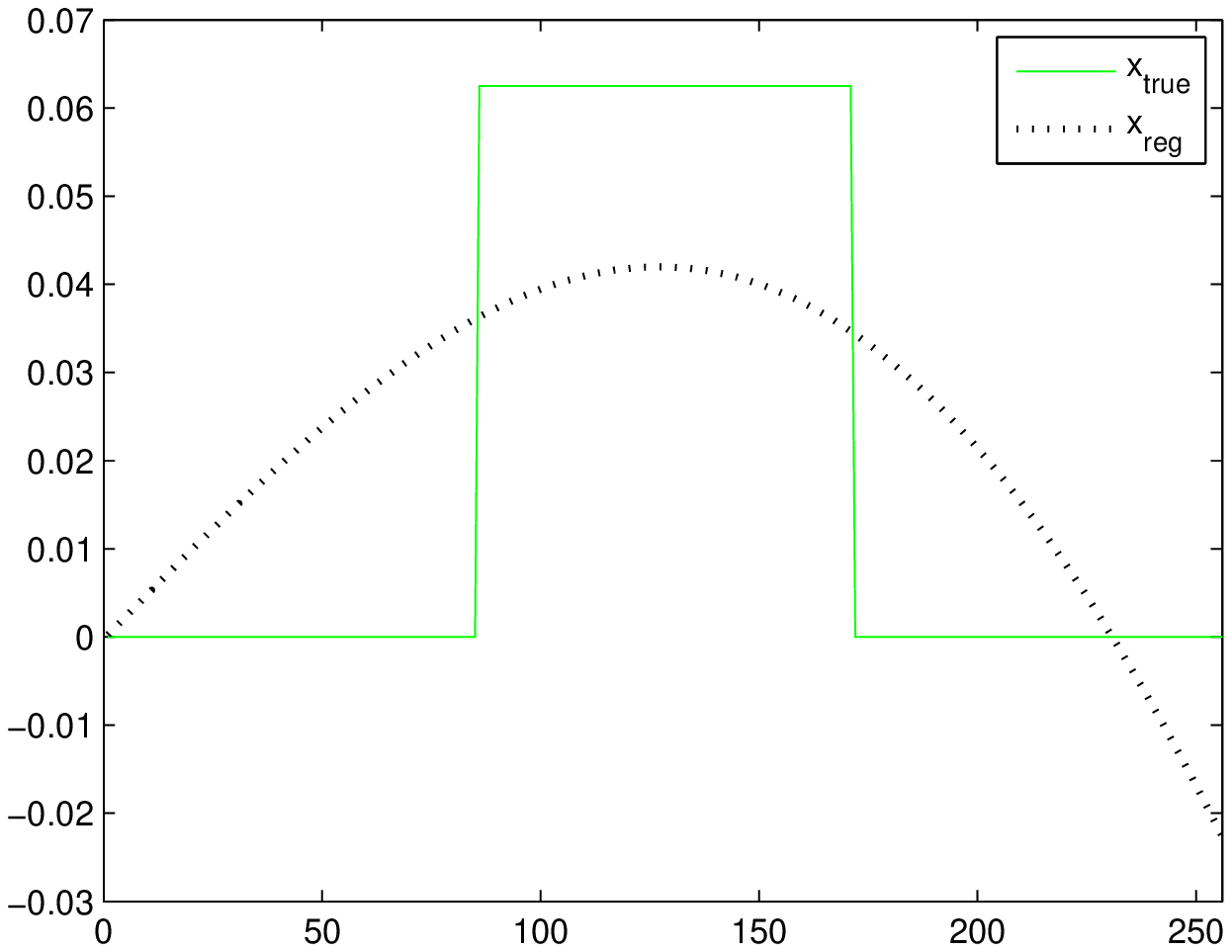}}
  \centerline{(d)}
\end{minipage}
\caption{(a)-(b): The relative errors $\|x^{(k)}-x_{true}\|/\|x_{true}\|$
by LSQR and the hybrid LSQR for $\varepsilon=10^{-3}$; (c)-(d):
The best regularized solutions $x_{reg}$ by LSQR
for $\mathsf{shaw}$ (left) and $\mathsf{wing}$ (right).}
\label{fig5a}
\end{figure}

\begin{figure}
\begin{minipage}{0.48\linewidth}
  \centerline{\includegraphics[width=6.0cm,height=5cm]{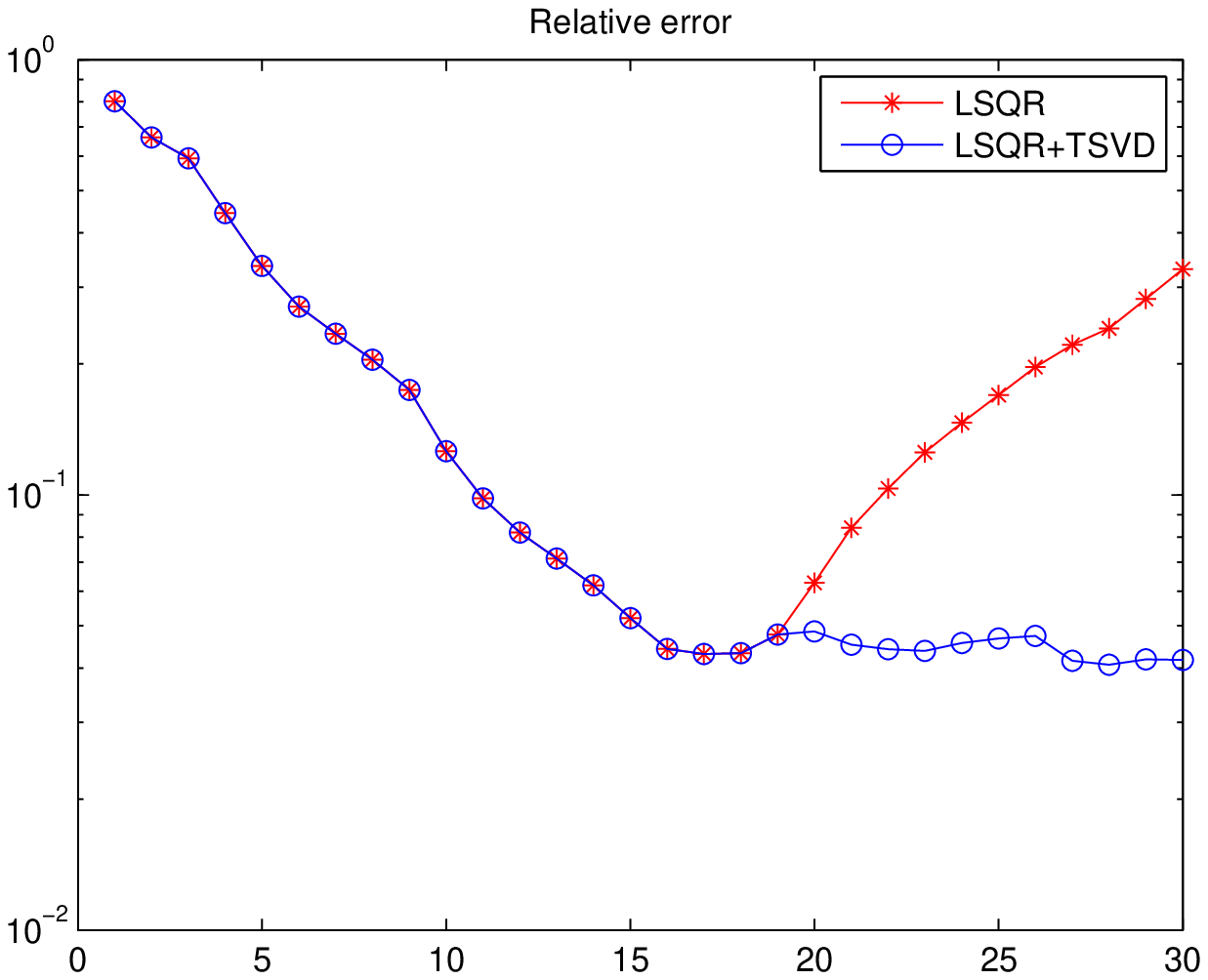}}
  \centerline{(a)}
\end{minipage}
\hfill
\begin{minipage}{0.48\linewidth}
  \centerline{\includegraphics[width=6.0cm,height=5cm]{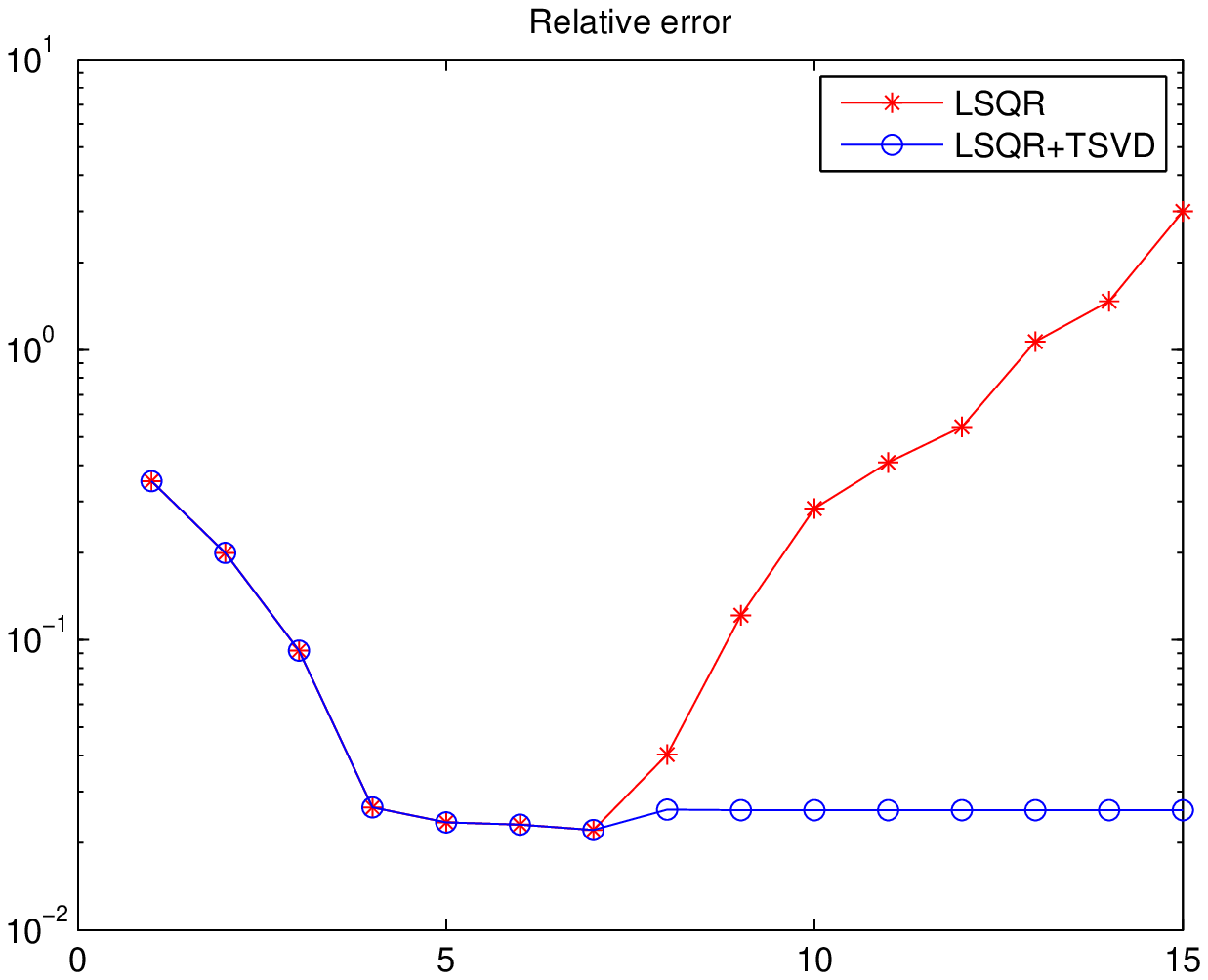}}
  \centerline{(b)}
\end{minipage}
\vfill
\begin{minipage}{0.48\linewidth}
  \centerline{\includegraphics[width=6.0cm,height=5cm]{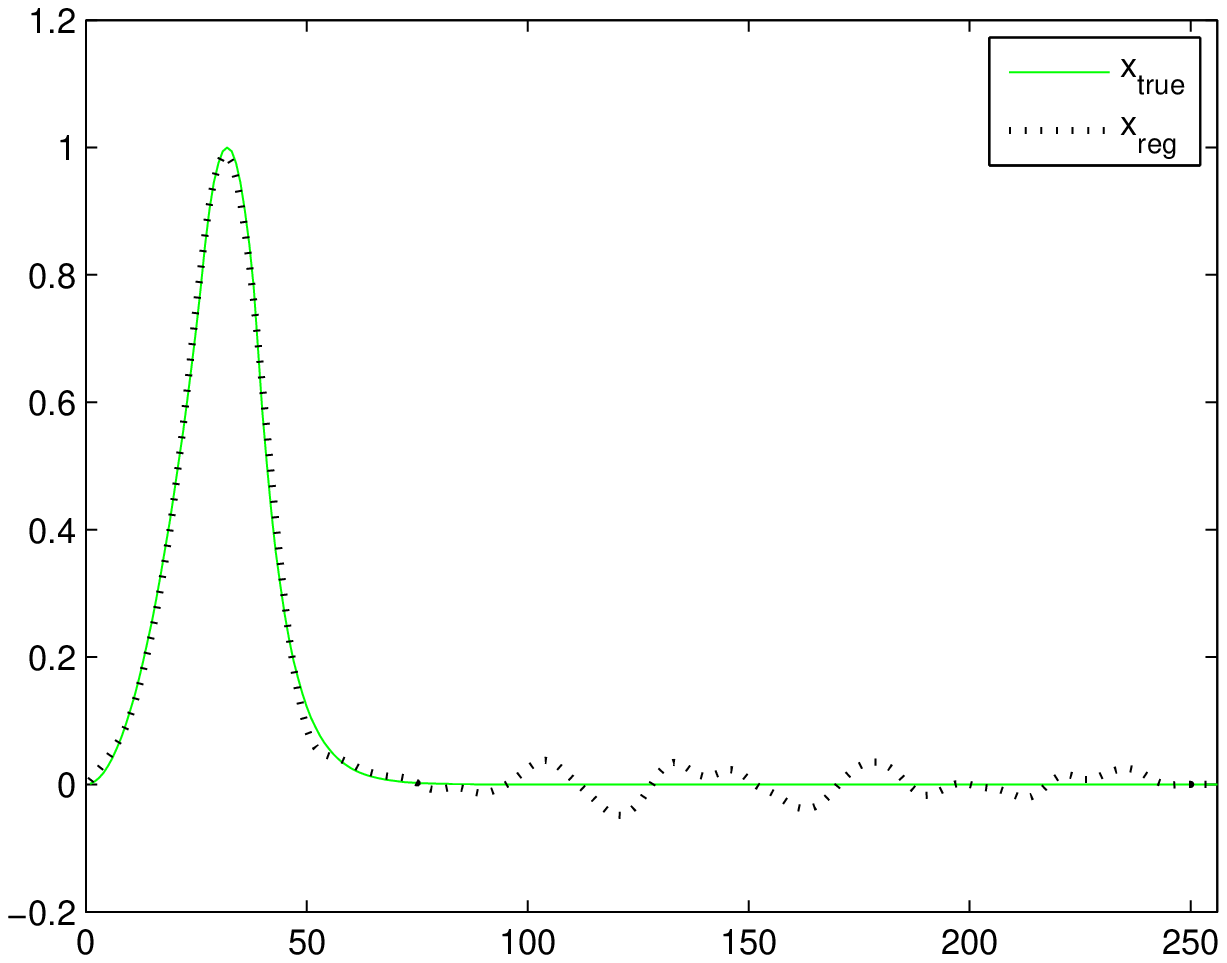}}
  \centerline{(c)}
\end{minipage}
\hfill
\begin{minipage}{0.48\linewidth}
  \centerline{\includegraphics[width=6.0cm,height=5cm]{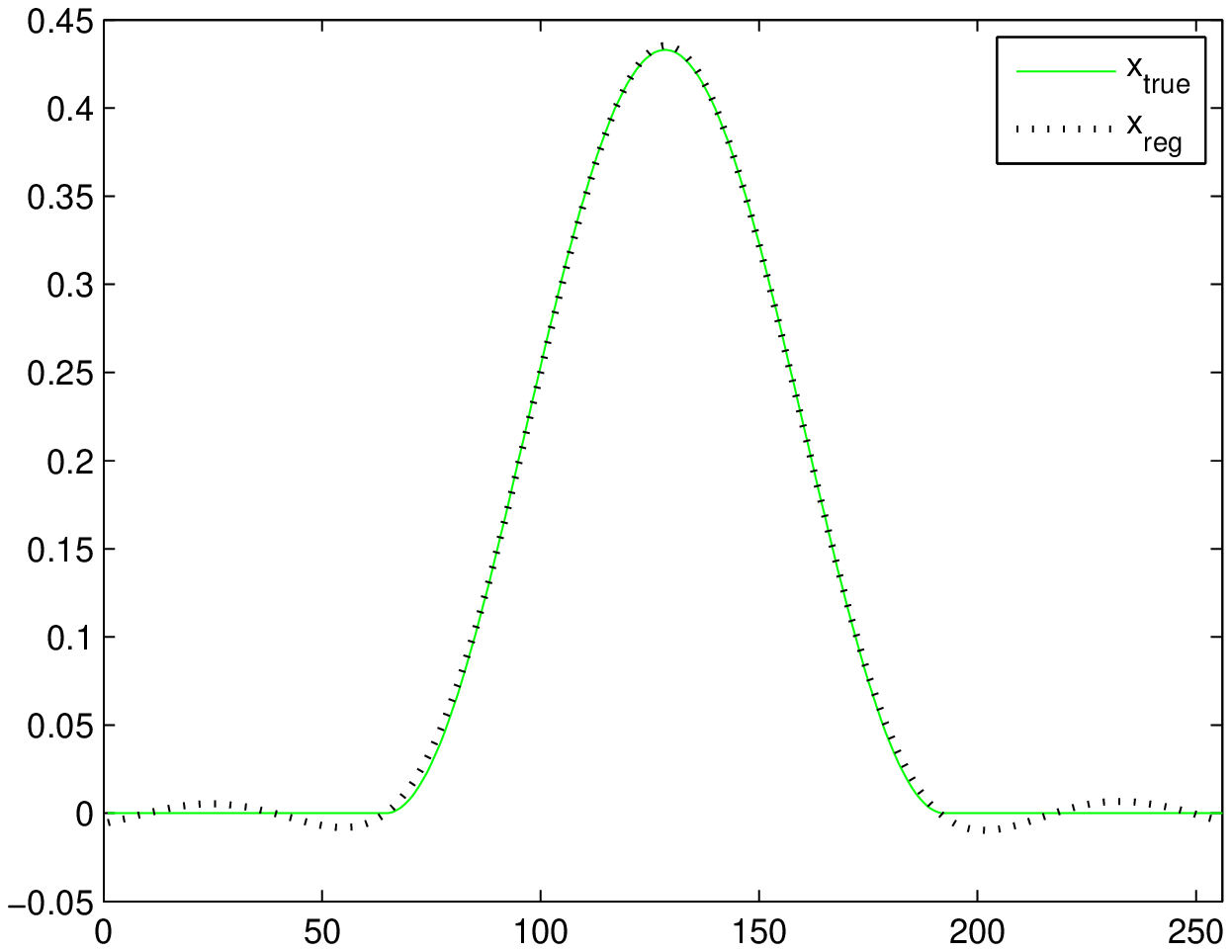}}
  \centerline{(d)}
\end{minipage}
\caption{(a)-(b): The relative errors $\|x^{(k)}-x_{true}\|/\|x_{true}\|$
by LSQR and the hybrid LSQR for $\varepsilon=10^{-3}$; (c)-(d):
The best possible regularized solutions $x_{reg}$ by LSQR
for $\mathsf{heat}$ (left) and $\mathsf{phillips}$ (right).}
\label{fig6}
\end{figure}

We first have a close look at the severely and moderately ill-posed problems.
Figure~\ref{fig5a} (a)-(b) and Figure~\ref{fig6} (a)-(b) plot
the relative errors of regularized solutions obtained by the two
methods for $\mathsf{shaw,\ wing}$ and
$\mathsf{heat,\ phillips}$. Clearly, we see that for each problem the
relative errors reach the same minimum level. After semi-convergence of LSQR,
the TSVD method applied to projected problems simply stabilizes the
regularized solutions with the minimum error and does not improve
them. This means that LSQR has already found best possible regularized solutions
at semi-convergence and has the full regularization,
and regularization applied to projected problems does not help
and is unnecessary at all. In practice, we simply stop LSQR after its
semi-convergence for severely and moderately ill-posed problems.

For these four problems, for test purposes we choose
$x_{reg}=\arg\min _{k}\|x^{(k)}-x_{true}\|$ for LSQR, which are just
the iterates obtained by LSQR at semi-convergence.
Figure~\ref{fig5a} (c)-(d) and Figure~\ref{fig6} (c)-(d) show that the
regularized solutions $x_{reg}$ are generally excellent approximations to the
true solutions $x_{true}$. The exception is the problem $\mathsf{wing}$ whose
underlying integral equation has a discontinuous solution, which
corresponds to the true solution $x_{true}$ whose
entries have big jumps in the discrete case,
as depicted in Figure~\ref{fig5a} (d). For it, the regularized solution
$x_{reg}$ deviates from $x_{true}$ considerably and the
relative error is not small. This is because all CG type methods applied to
either \eqref{eq1} or $A^TAx=A^Tb$ or $\min\|AA^Ty-b\|$ with $x=A^Ty$
compute smooth regularized solutions. More insightfully, LSQR and CGLS are
equivalent to implicitly solving the Tikhonov regularization problem
\eqref{tikhonov}, and their regularized solutions are
of the filtered form \eqref{eqfilter2}. It is well known that
the regularization term $\lambda^2\|x\|$ in Tikhonov regularization does not
suit for discontinuous solutions. The continuous ill-posed problems with
discontinuous or non-smooth solutions are from numerous important applications,
including linear regression, barcode reading, gravity surveying in geophysics,
image restoration and some others \cite{aster,hansen10,mueller}. For them, a
better alternative is use the 1-norm $\lambda^2\|Lx\|_1$ as the regularization
term, which leads to the Total Variation Regularization
\cite{aster,engl00,hansen10,mueller,vogel02} or Errors-in-Variables Modeling
called in \cite{huffel}, where $L\not=I$ is some
$p\times n$ matrix with no restriction to $p$ and is typically taken to be the
discrete approximation to the first or second derivative
operator \cite[Ch.8]{hansen10}.

\begin{figure}
\begin{minipage}{0.48\linewidth}
  \centerline{\includegraphics[width=6.0cm,height=5cm]{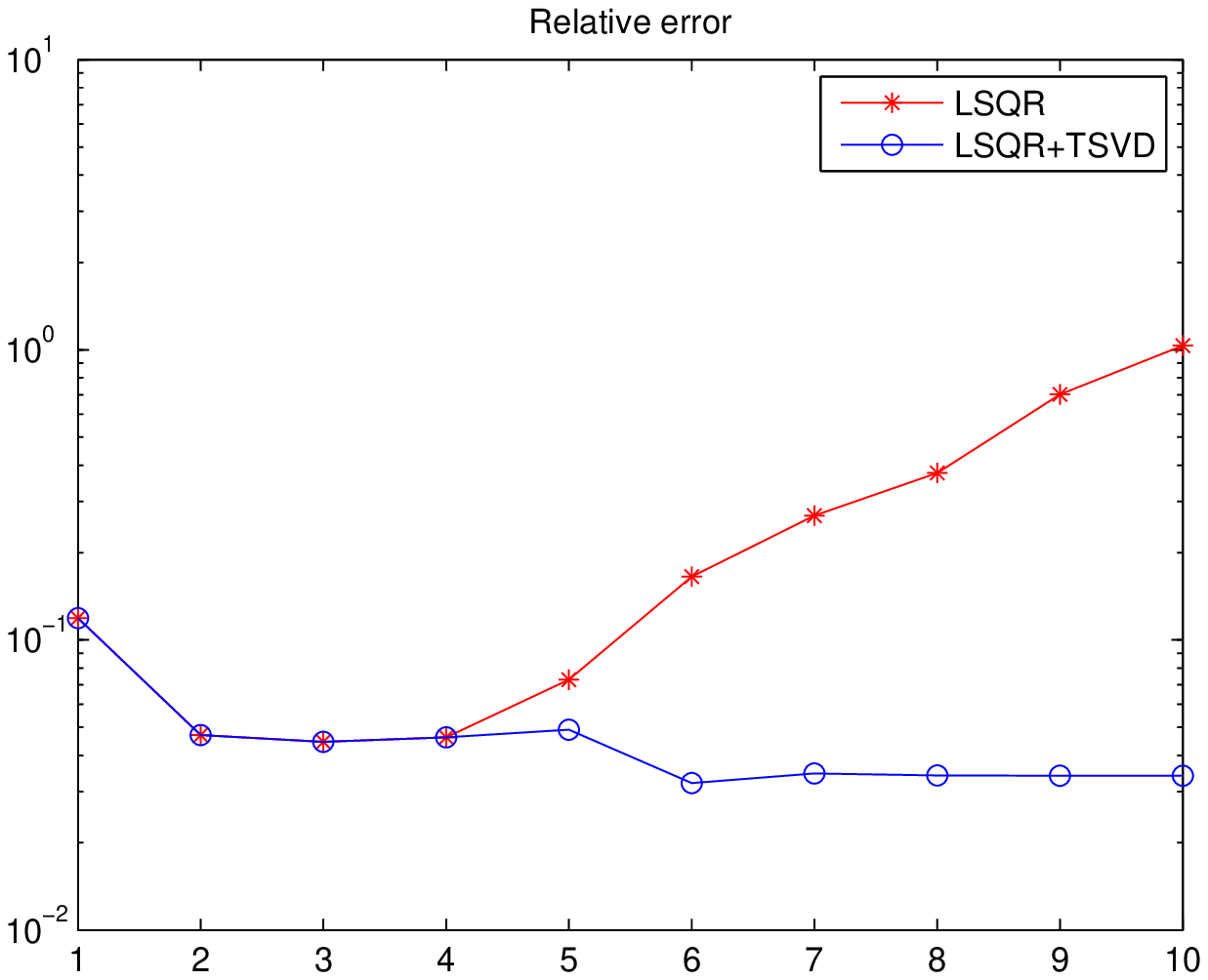}}
  \centerline{(a)}
\end{minipage}
\hfill
\begin{minipage}{0.48\linewidth}
  \centerline{\includegraphics[width=6.0cm,height=5cm]{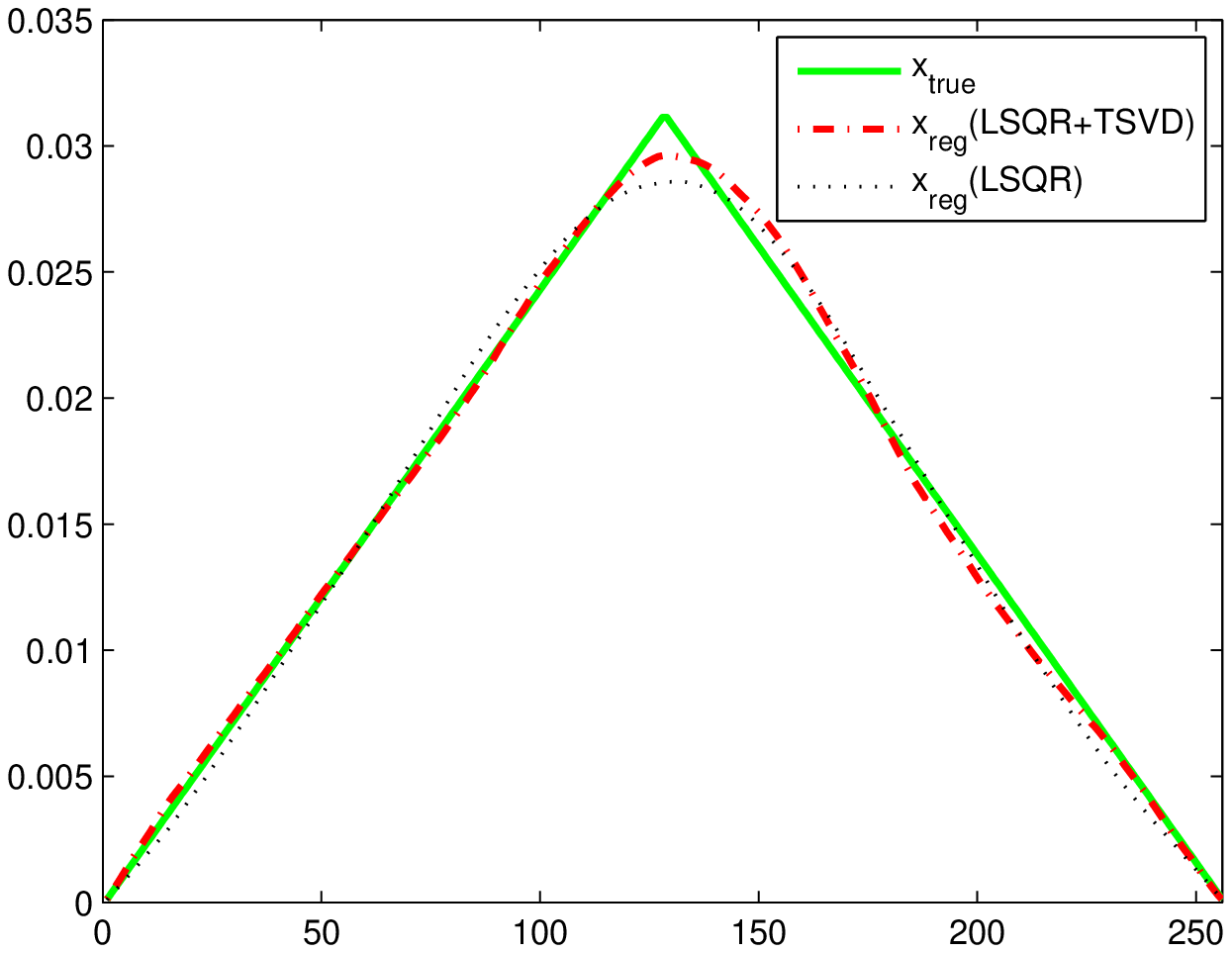}}
  \centerline{(b)}
\end{minipage}
\caption{(a)-(b):
The relative errors $\|x^{(k)}-x_{true}\|/\|x_{true}\|$ and
the regularized solutions $x_{reg}$ by LSQR at semi-convergence and the
best possible regularized solution by the hybrid LSQR for {\sf deriv2}.}
\label{figmild2}
\end{figure}

Now we investigate the behavior of LSQR and the hybrid LSQR for
$\mathsf{deriv2}$. Figure~\ref{figmild2} (a) indicates that
the relative errors of $x^{(k)}$ by the hybrid LSQR reach a considerably
smaller minimum level than those by LSQR, illustrating that LSQR
has only the partial regularization. Precisely, we find that the semi-convergence
of LSQR occurs at iteration $k=4$, but the regularized solution is not acceptable.
The hybrid LSQR uses a larger six dimensional Krylov subspace
$\mathcal{K}_6(A^TA,A^Tb)$ to construct a more accurate regularized solution.
We also choose $x_{reg}=\arg\min _{k}\|x^{(k)}-x_{true}\|$ for LSQR and the
hybrid LSQR, respectively. Figure~\ref{figmild2} (b) indicates that the
best regularized solution by the hybrid LSQR is a considerably better
approximation to $x_{true}$ than that by LSQR, especially in the
non-smooth middle part of $x_{true}$.

\subsection{Decay behavior of $\alpha_k$ and $\beta_{k+1}$}

For the severely ill-posed $\mathsf{shaw, wing}$ and the moderately
ill-posed $\mathsf{heat, phillips}$, we now
illustrate that $\alpha_k$ and $\beta_{k+1}$ decay as fast as the singular
values $\sigma_k$ of $A$. We take the noise level $\varepsilon=10^{-3}$. The results
are similar for $\varepsilon=10^{-2}$ and $10^{-4}$.

\begin{figure}
\begin{minipage}{0.48\linewidth}
  \centerline{\includegraphics[width=6.0cm,height=5cm]{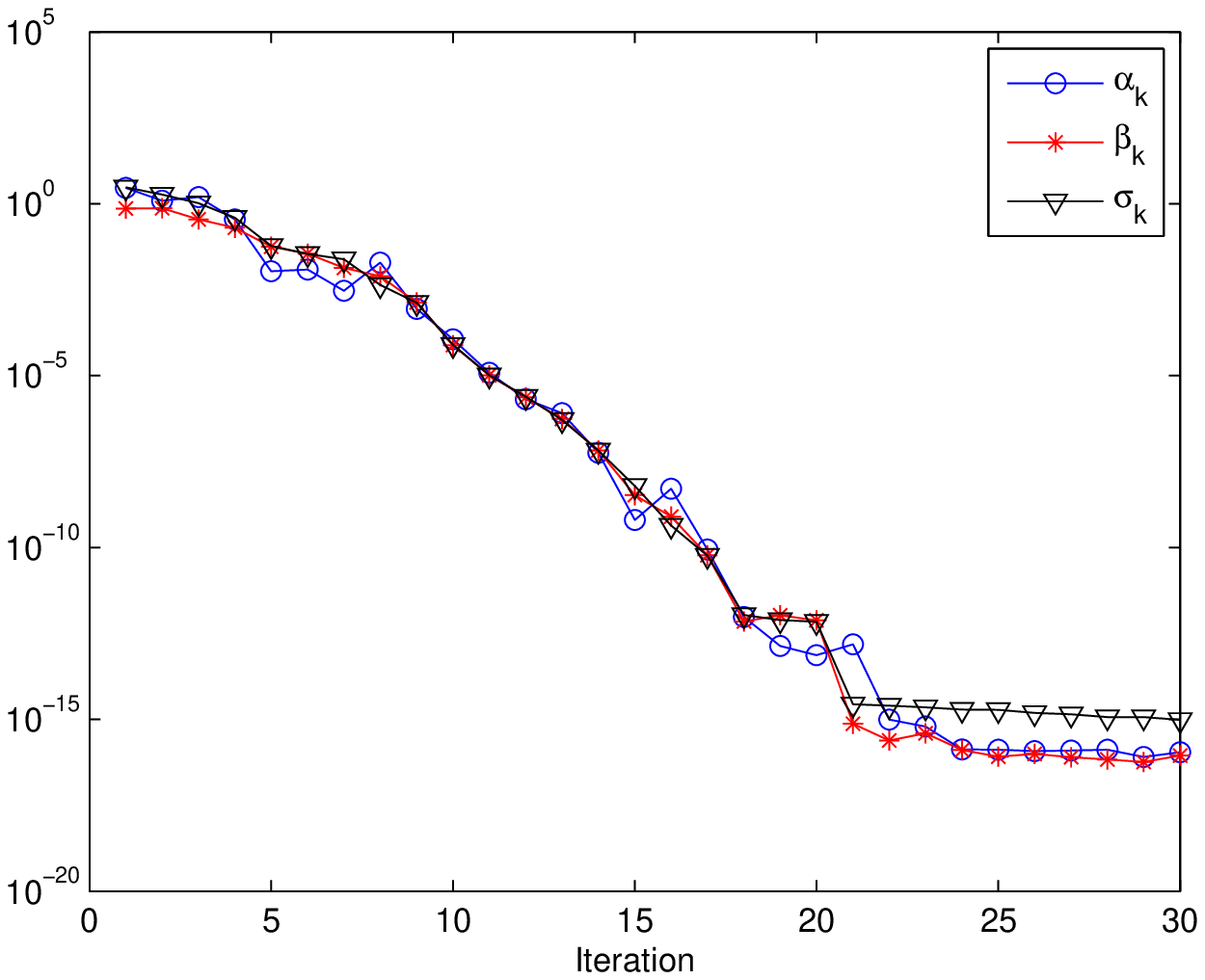}}
  \centerline{(a)}
\end{minipage}
\hfill
\begin{minipage}{0.48\linewidth}
  \centerline{\includegraphics[width=6.0cm,height=5cm]{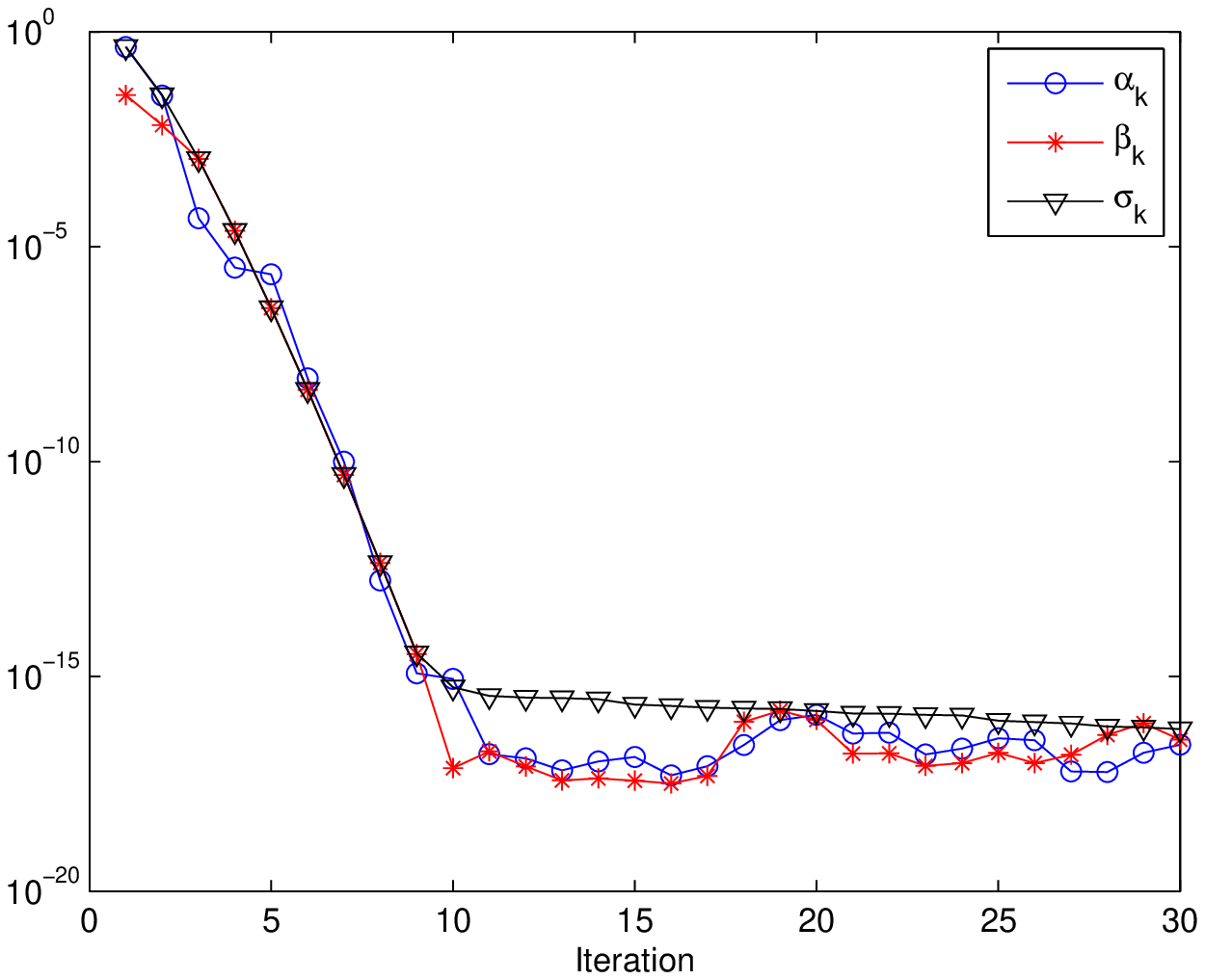}}
  \centerline{(b)}
\end{minipage}
\vfill
\begin{minipage}{0.48\linewidth}
  \centerline{\includegraphics[width=6.0cm,height=5cm]{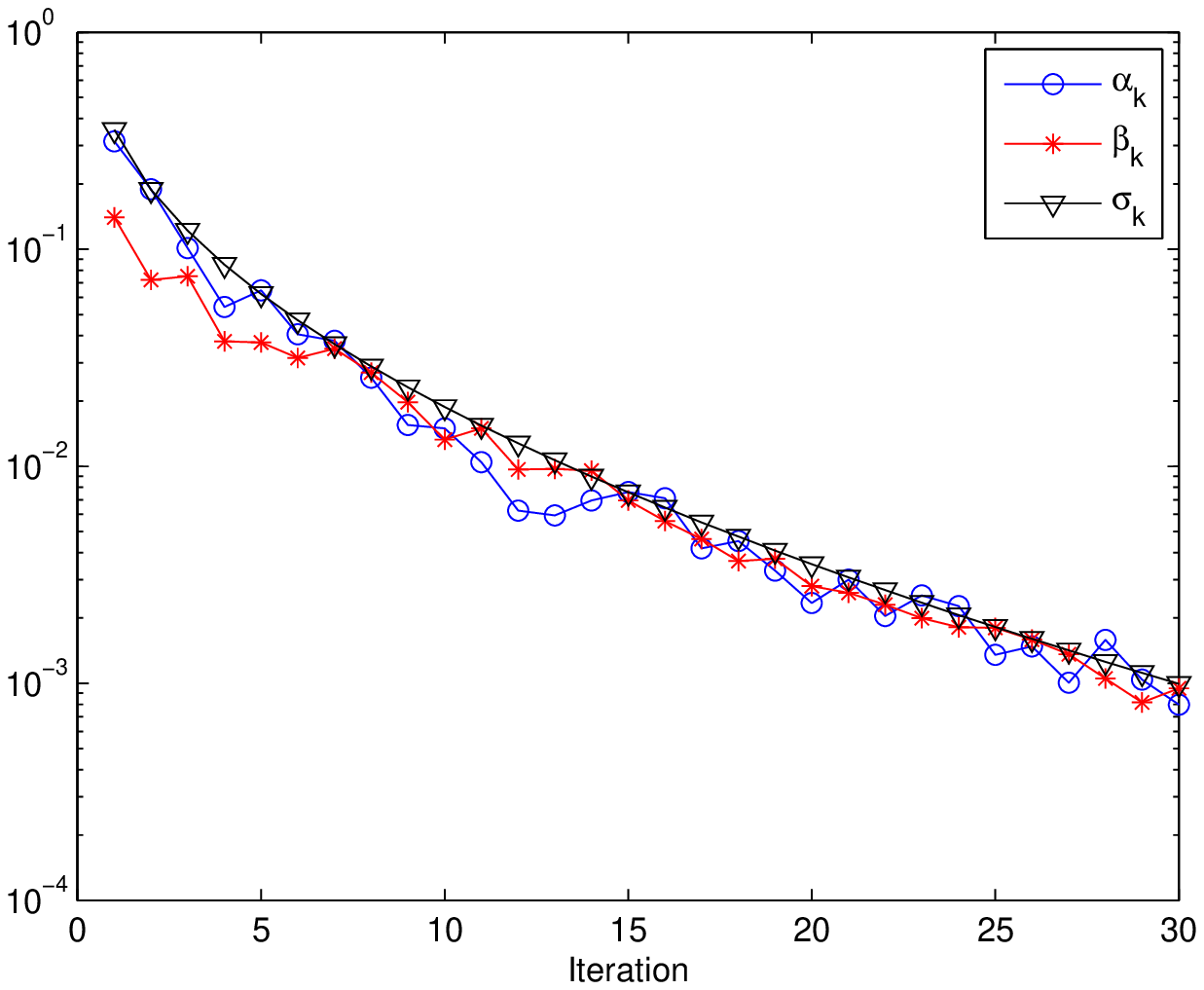}}
  \centerline{(c)}
\end{minipage}
\hfill
\begin{minipage}{0.48\linewidth}
  \centerline{\includegraphics[width=6.0cm,height=5cm]{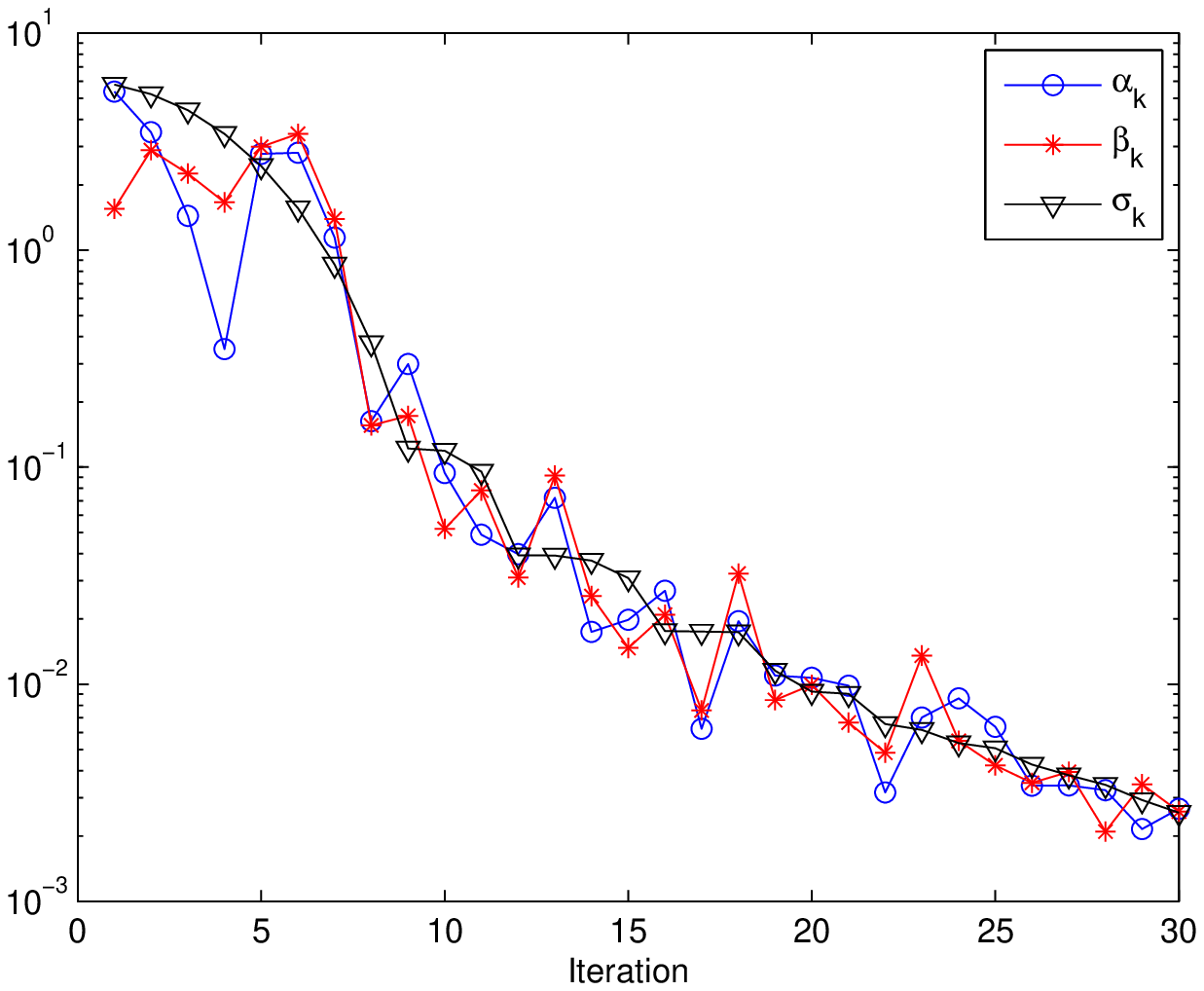}}
  \centerline{(d)}
\end{minipage}
\caption{(a)-(d): Decay curves of the sequences $\alpha_k$,
$\beta_{k+1}$ and $\sigma_k$ for $\mathsf{shaw, wing, i\_laplace}$ and $\mathsf{heat}$
(from top left to bottom right).} \label{fig7}
\end{figure}

Figure~\ref{fig7} illustrates that both $\alpha_k$ and $\beta_{k+1}$
decay as fast as $\sigma_k$, and for $\mathsf{shaw}$ and
$\mathsf{wing}$ all of them decay swiftly and level off at
$\epsilon_{\rm mach}$ due to round-off errors
in finite precision arithmetic. Precisely, they reach
the level of $\epsilon_{\rm mach}$ at $k=22$ and $k=8$ for $\mathsf{shaw}$ and
$\mathsf{wing}$, respectively. Such decay behavior has also
been observed in \cite{bazan14,gazzola14,gazzola-online}, but no theoretical
support was given. These experiments confirm Theorem~\ref{main1} and
Theorem~\ref{main2}, which have proved that $\gamma_k$ decreases as fast as
$\sigma_{k+1}$ and that $\alpha_k$, $\beta_{k+1}$ and $\alpha_k+\beta_{k+1}$
decay as fast as $\sigma_k$.

\subsection{A comparison of LSQR and the TSVD method}

We compare the performance of LSQR and the TSVD method
for the severely ill-posed $\mathsf{shaw,\ wing}$ and moderately ill-posed
$\mathsf{heat,\ phillips}$. We take $\varepsilon=10^{-3}$. For each problem,
we compute the norms of regularized solutions, their relative errors and
the residual norms obtained by the two methods. We plot the L-curves of
the residual norms versus those of regularized solutions in the $\log$-$\log$
scale.

\begin{figure}
\begin{minipage}{0.48\linewidth}
  \centerline{\includegraphics[width=6.0cm,height=5cm]{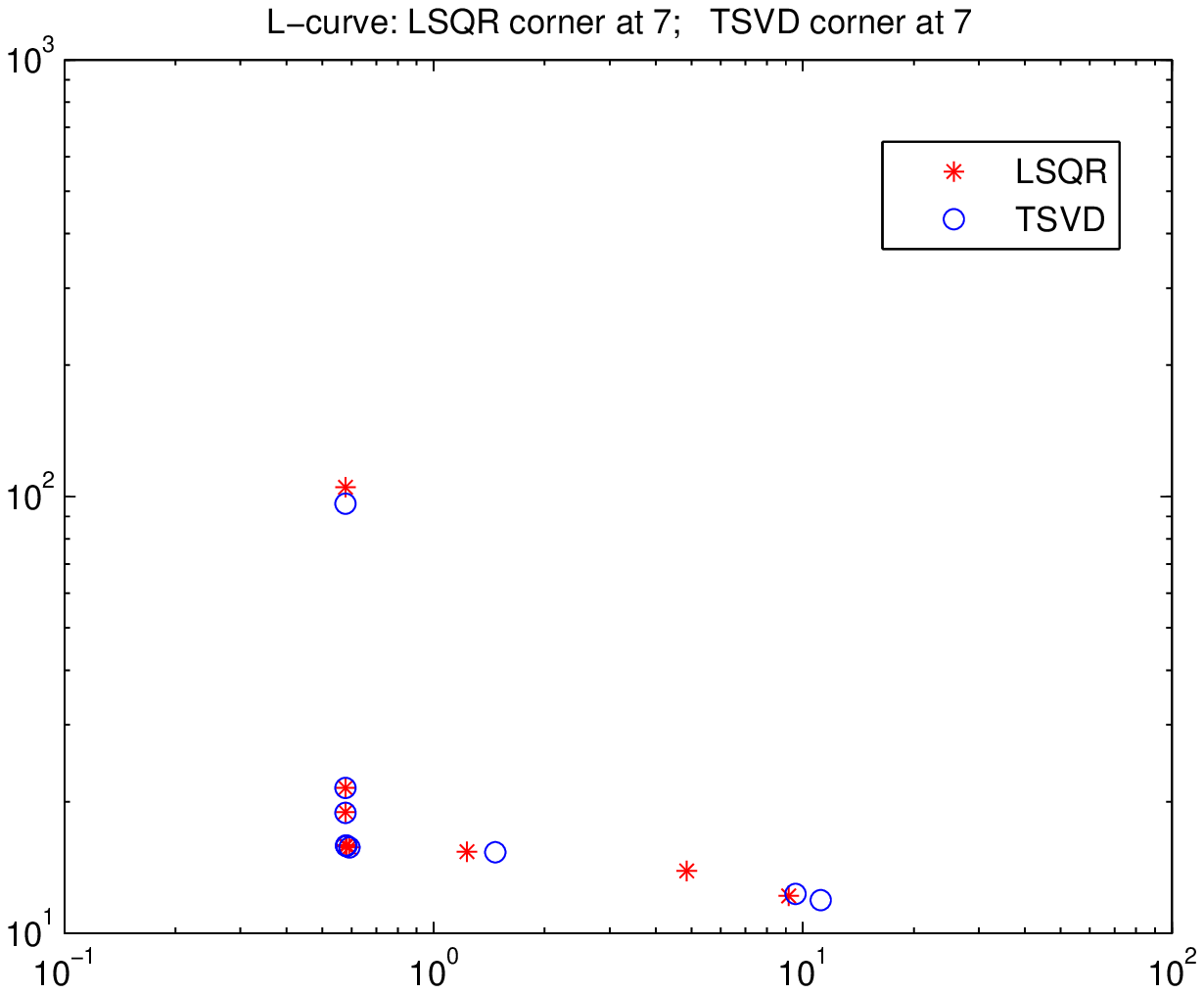}}
  \centerline{(a)}
\end{minipage}
\hfill
\begin{minipage}{0.48\linewidth}
  \centerline{\includegraphics[width=6.0cm,height=5cm]{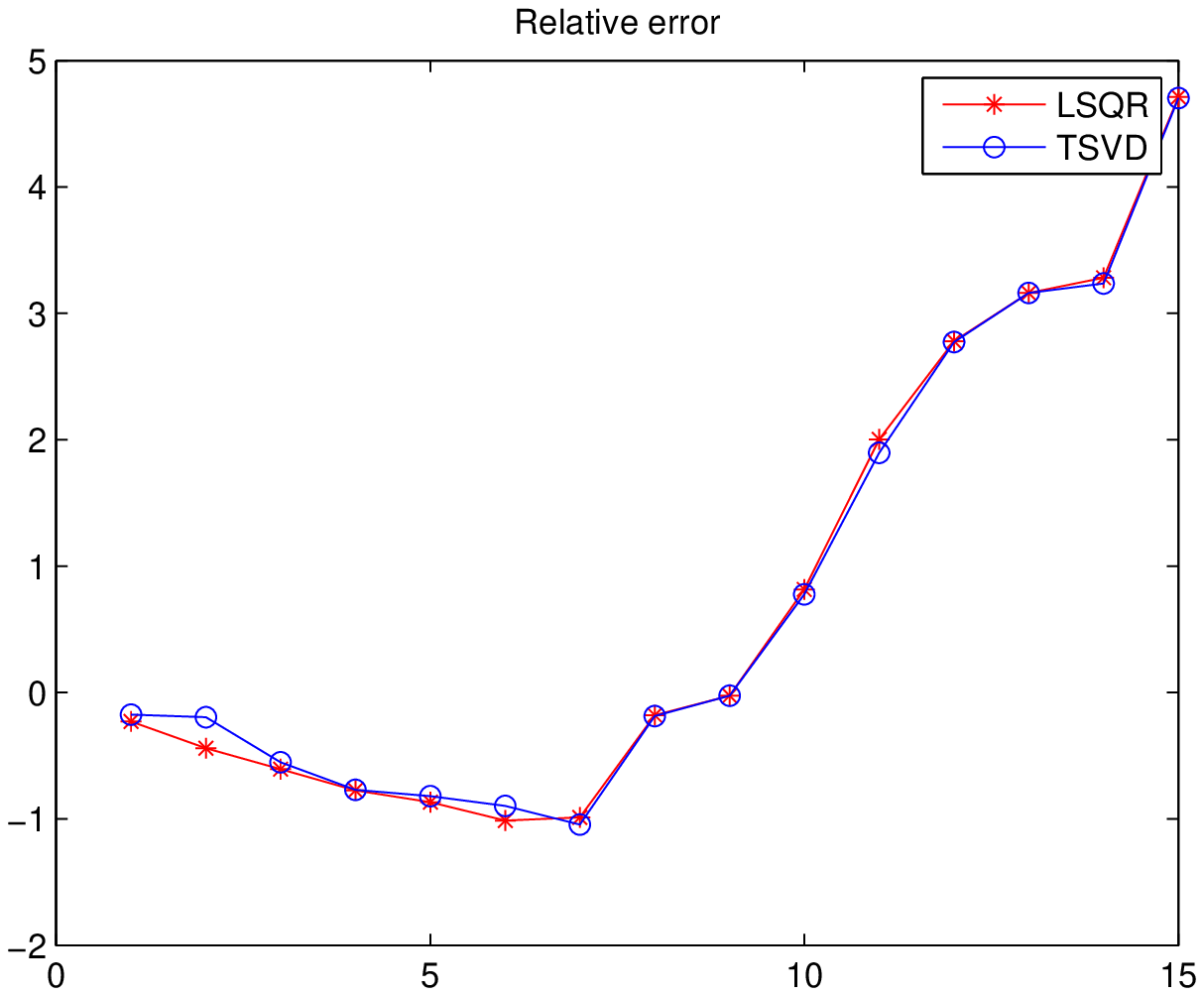}}
  \centerline{(b)}
\end{minipage}
\vfill
\begin{minipage}{0.48\linewidth}
  \centerline{\includegraphics[width=6.0cm,height=5cm]{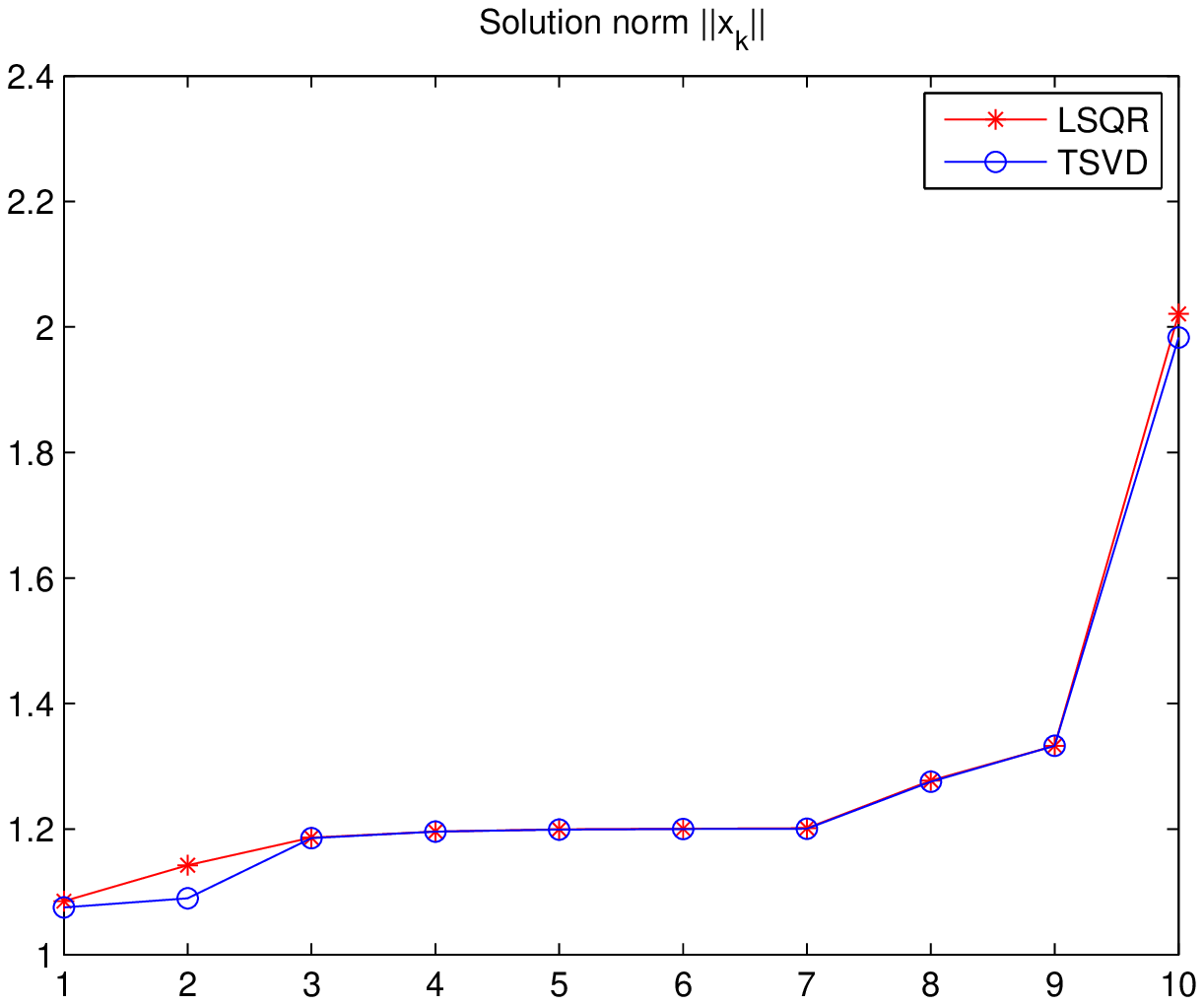}}
  \centerline{(c)}
\end{minipage}
\hfill
\begin{minipage}{0.48\linewidth}
  \centerline{\includegraphics[width=6.0cm,height=5cm]{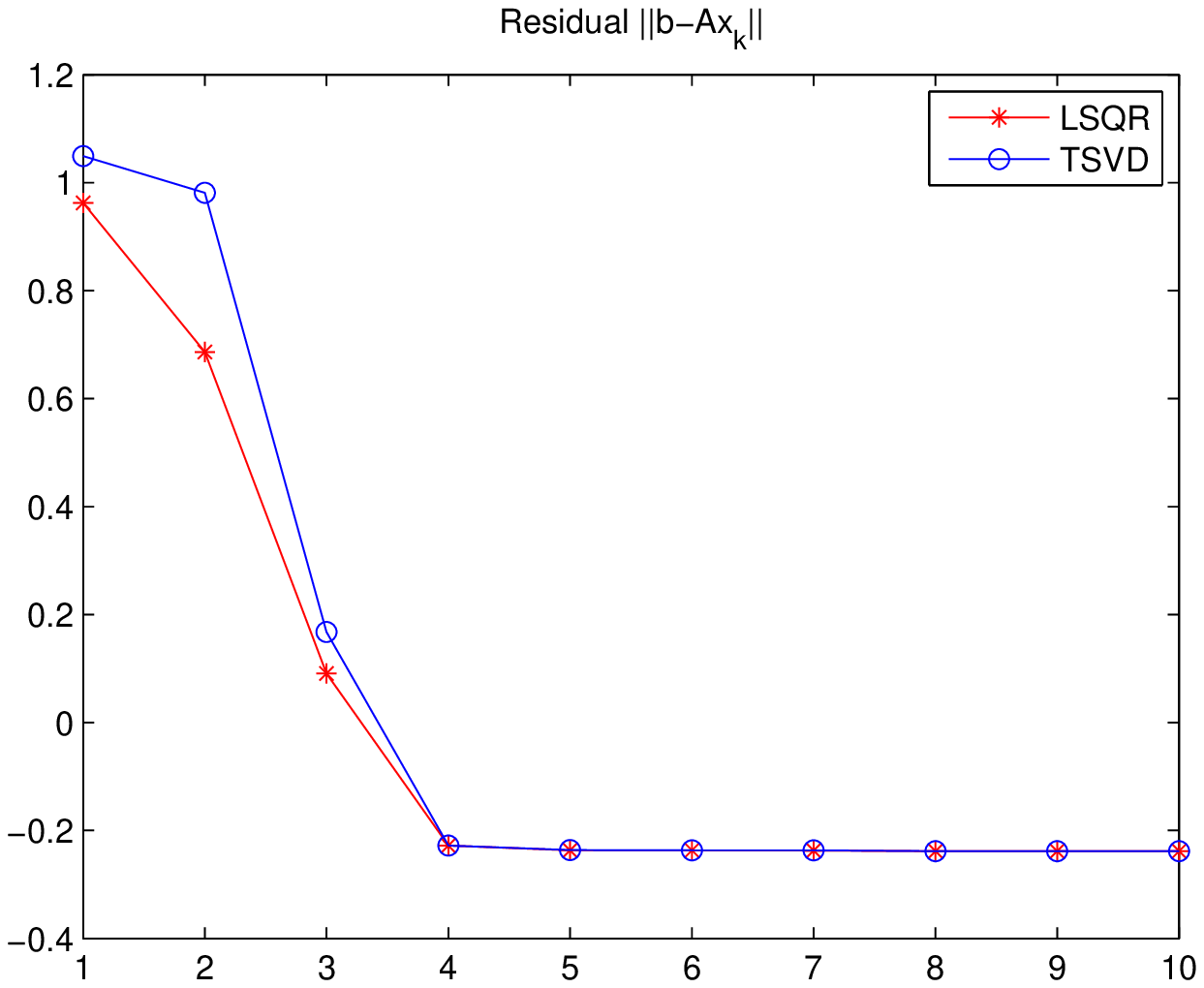}}
  \centerline{(d)}
\end{minipage}
\caption{Results for the severely ill-posed problem $\mathsf{shaw}$.}
\label{lsqrtsvd1}
\end{figure}

\begin{figure}
\begin{minipage}{0.48\linewidth}
  \centerline{\includegraphics[width=6.0cm,height=5cm]{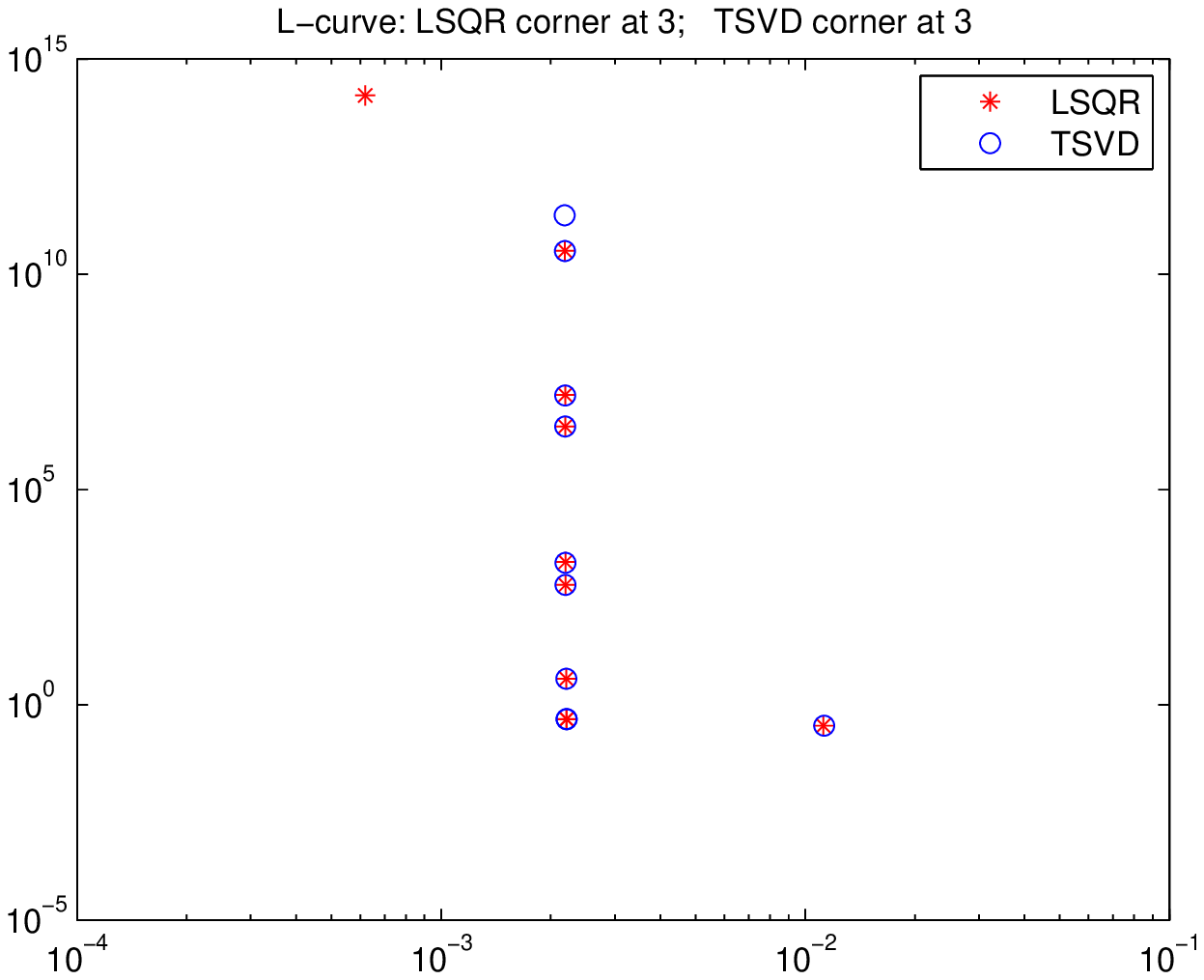}}
  \centerline{(a)}
\end{minipage}
\hfill
\begin{minipage}{0.48\linewidth}
  \centerline{\includegraphics[width=6.0cm,height=5cm]{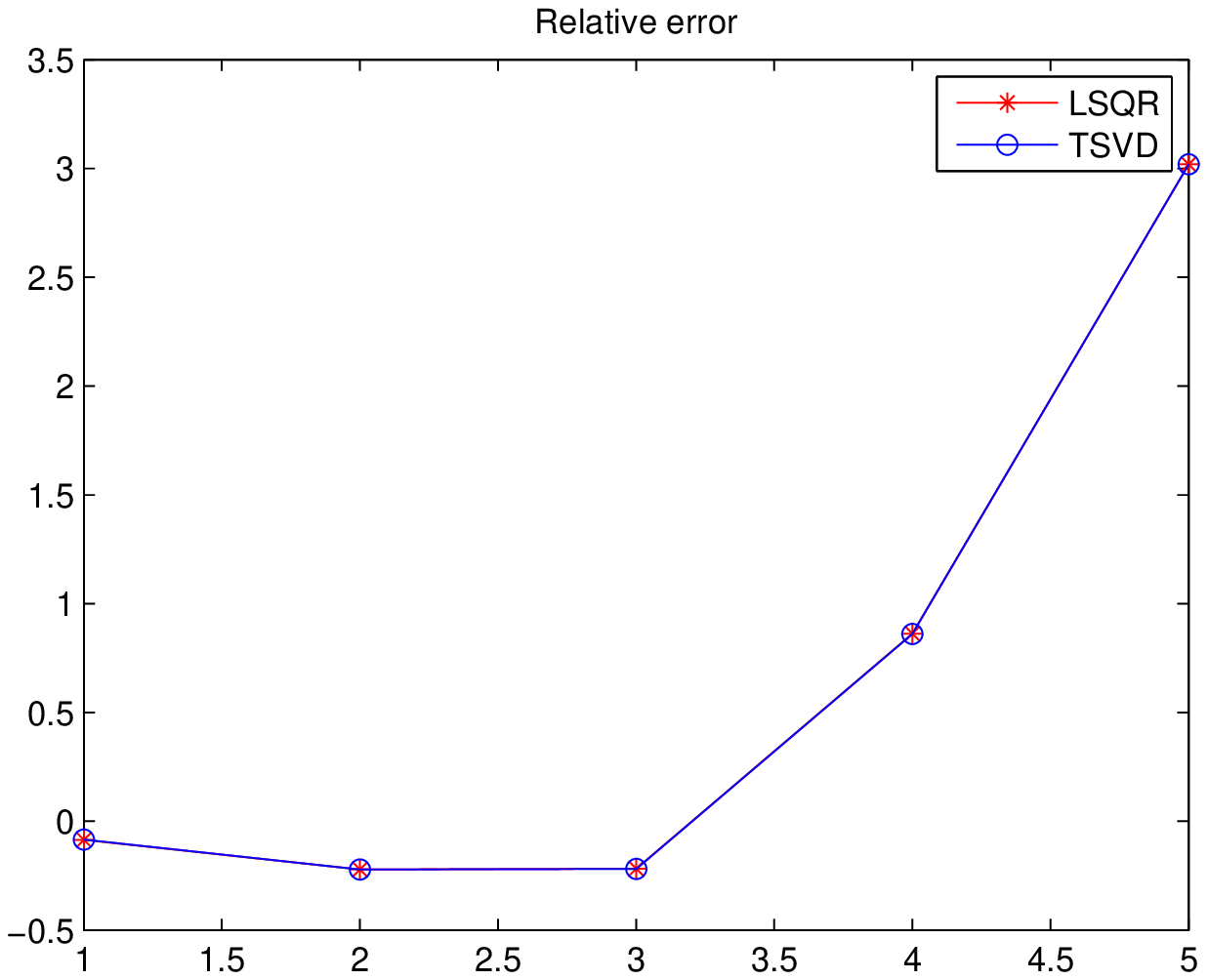}}
  \centerline{(b)}
\end{minipage}
\begin{minipage}{0.48\linewidth}
  \centerline{\includegraphics[width=6.0cm,height=5cm]{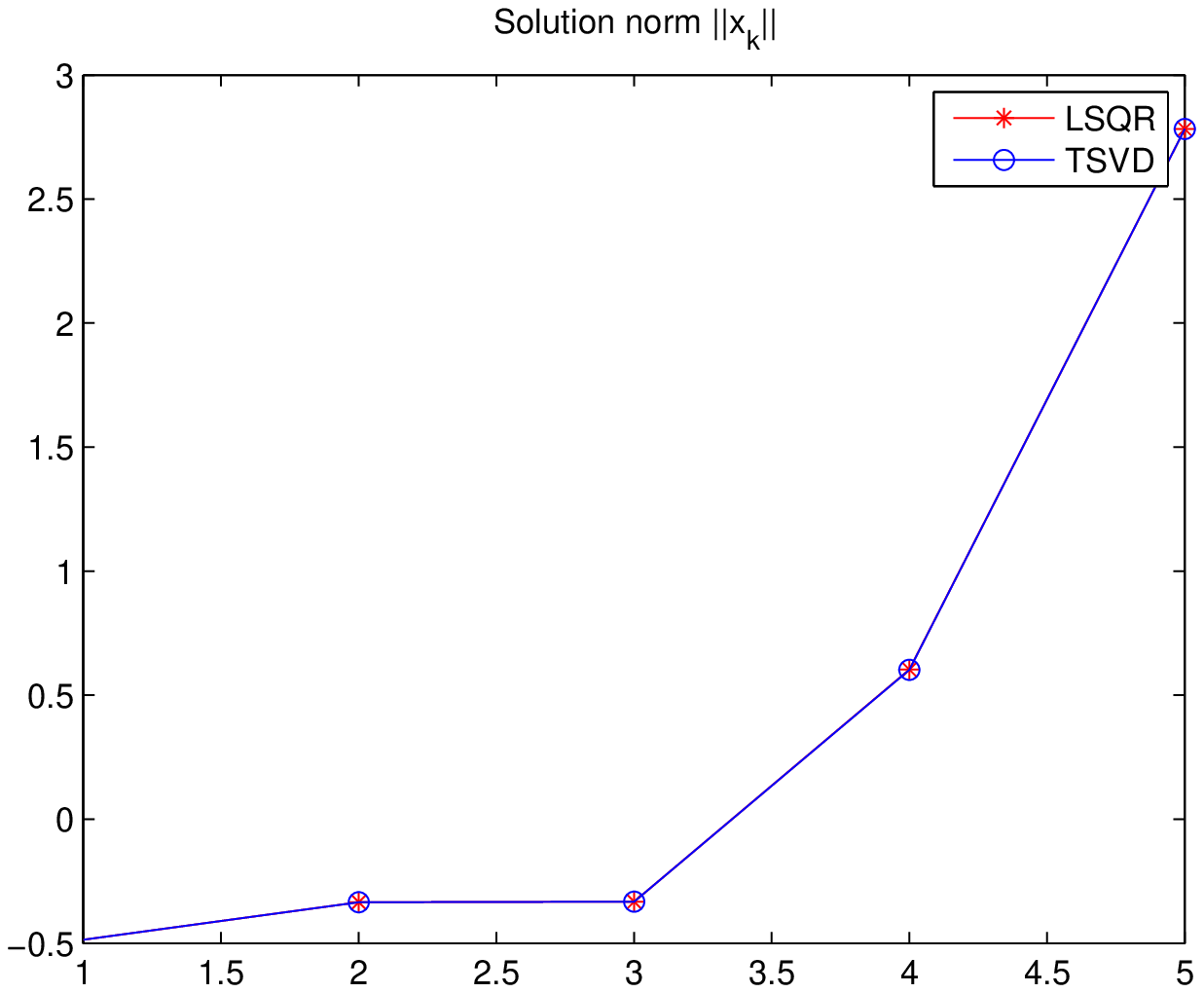}}
  \centerline{(c)}
\end{minipage}
\hfill
\begin{minipage}{0.48\linewidth}
  \centerline{\includegraphics[width=6.0cm,height=5cm]{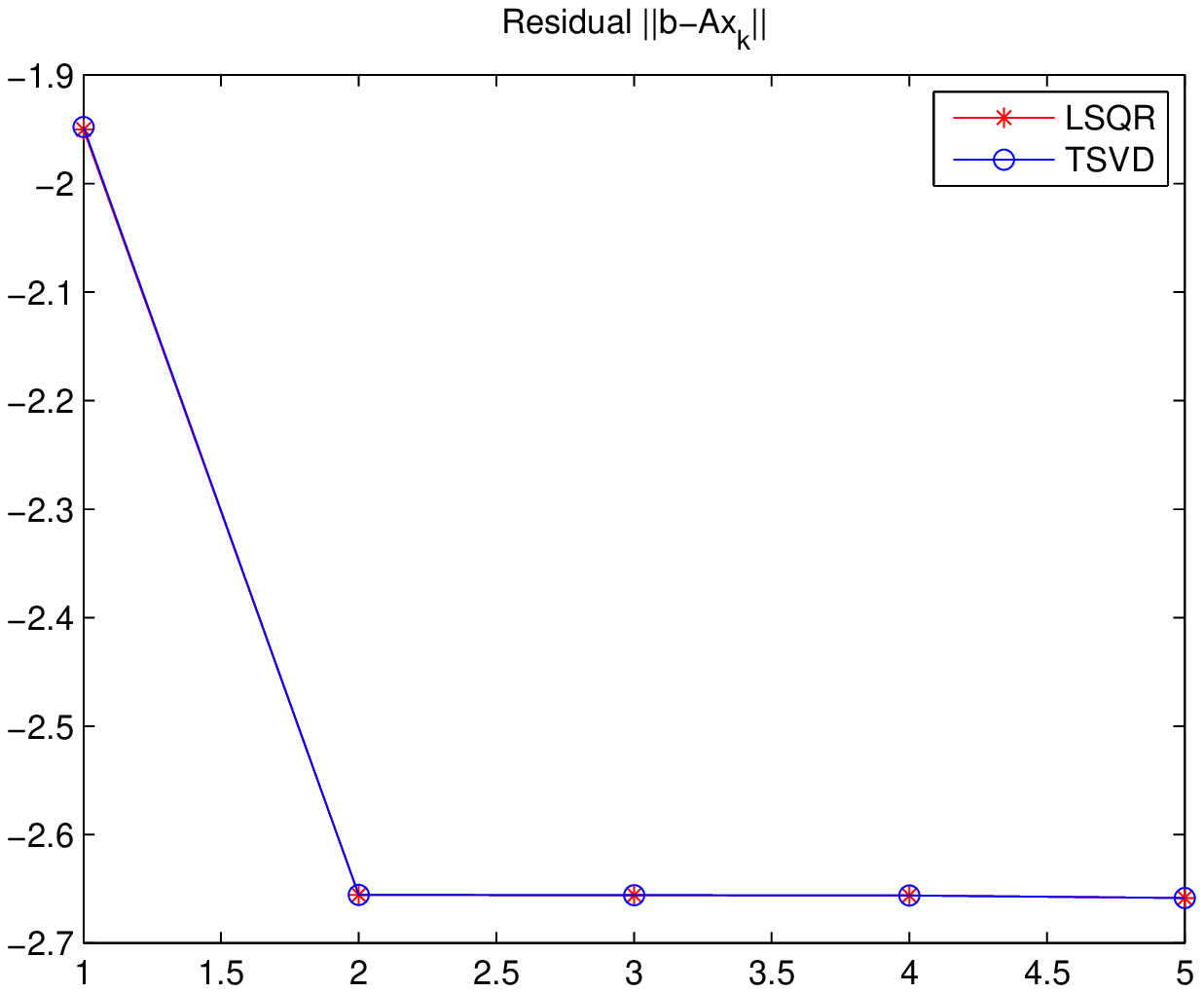}}
  \centerline{(d)}
\end{minipage}
\caption{Results for the severely ill-posed problem $\mathsf{wing}$.}
\label{lsqrtsvd2}
\end{figure}

Figures~\ref{lsqrtsvd1}--\ref{lsqrtsvd2} indicate LSQR and the TSVD method
behave very similarly for $\mathsf{shaw}$ and $\mathsf{wing}$. They
illustrate that, for $\mathsf{wing}$, the norms of approximate
solutions and the relative errors by the two methods
are almost indistinguishable for the same $k$, and, for $\mathsf{shaw}$, the
residual norms by LSQR decreases more quickly than the ones by the TSVD method
for $k=1,2,3$ and then they become almost identical starting from $k=4$.
The L-curves tell us that the two methods obtain the best regularized solutions
when $k_0=7$ and $k_0=3$ for $\mathsf{shaw}$ and $\mathsf{wing}$, respectively.
The values of $k_0$ determined by the L-curves are exactly the ones
at which semi-convergence occurs, as indicated by (b) and (c) in
Figures~\ref{lsqrtsvd1}--\ref{lsqrtsvd2}. These results demonstrate that
LSQR has the full regularization and resembles the TSVD method very much.

\begin{figure}
\begin{minipage}{0.48\linewidth}
  \centerline{\includegraphics[width=6.0cm,height=5cm]{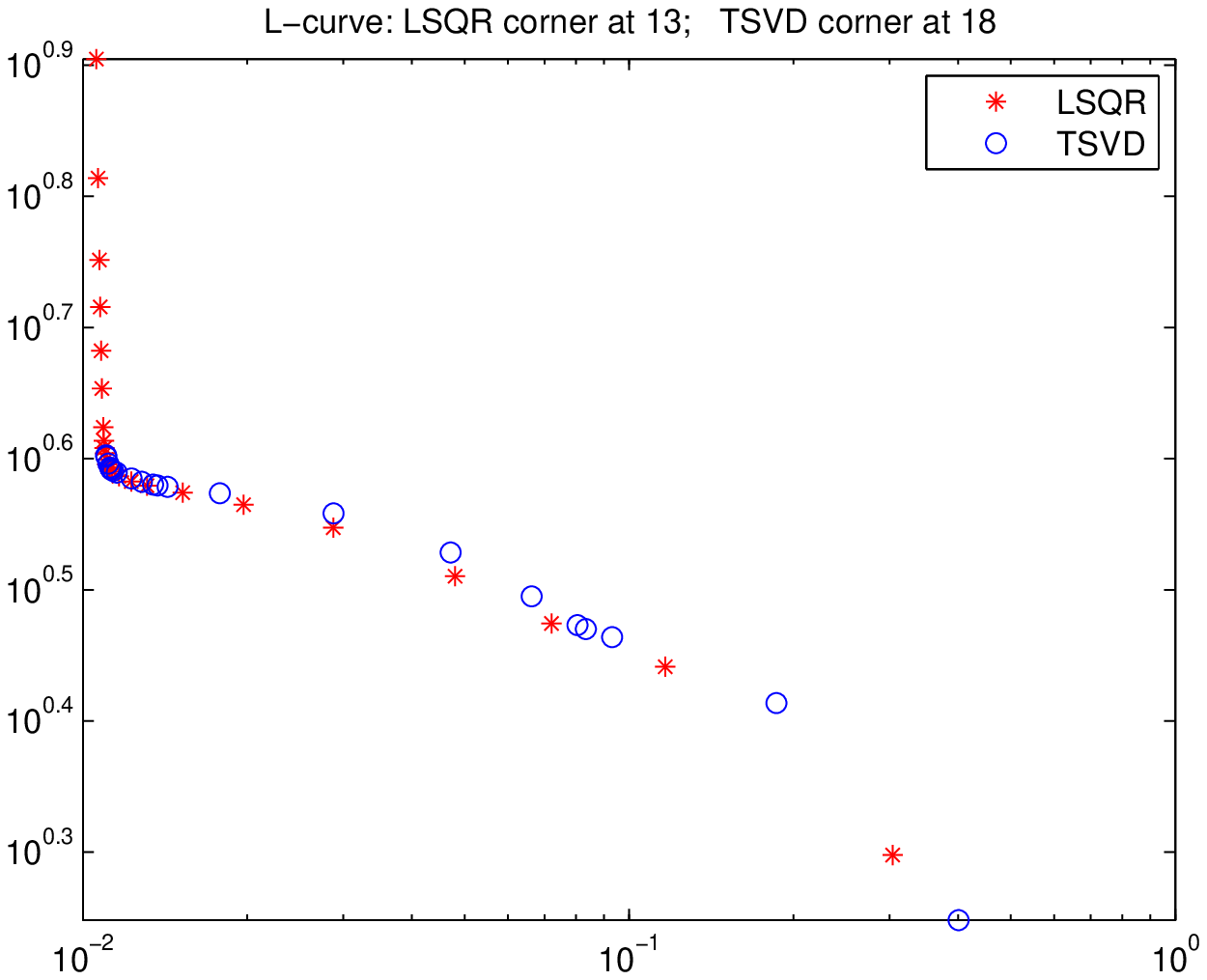}}
  \centerline{(a)}
\end{minipage}
\hfill
\begin{minipage}{0.48\linewidth}
  \centerline{\includegraphics[width=6.0cm,height=5cm]{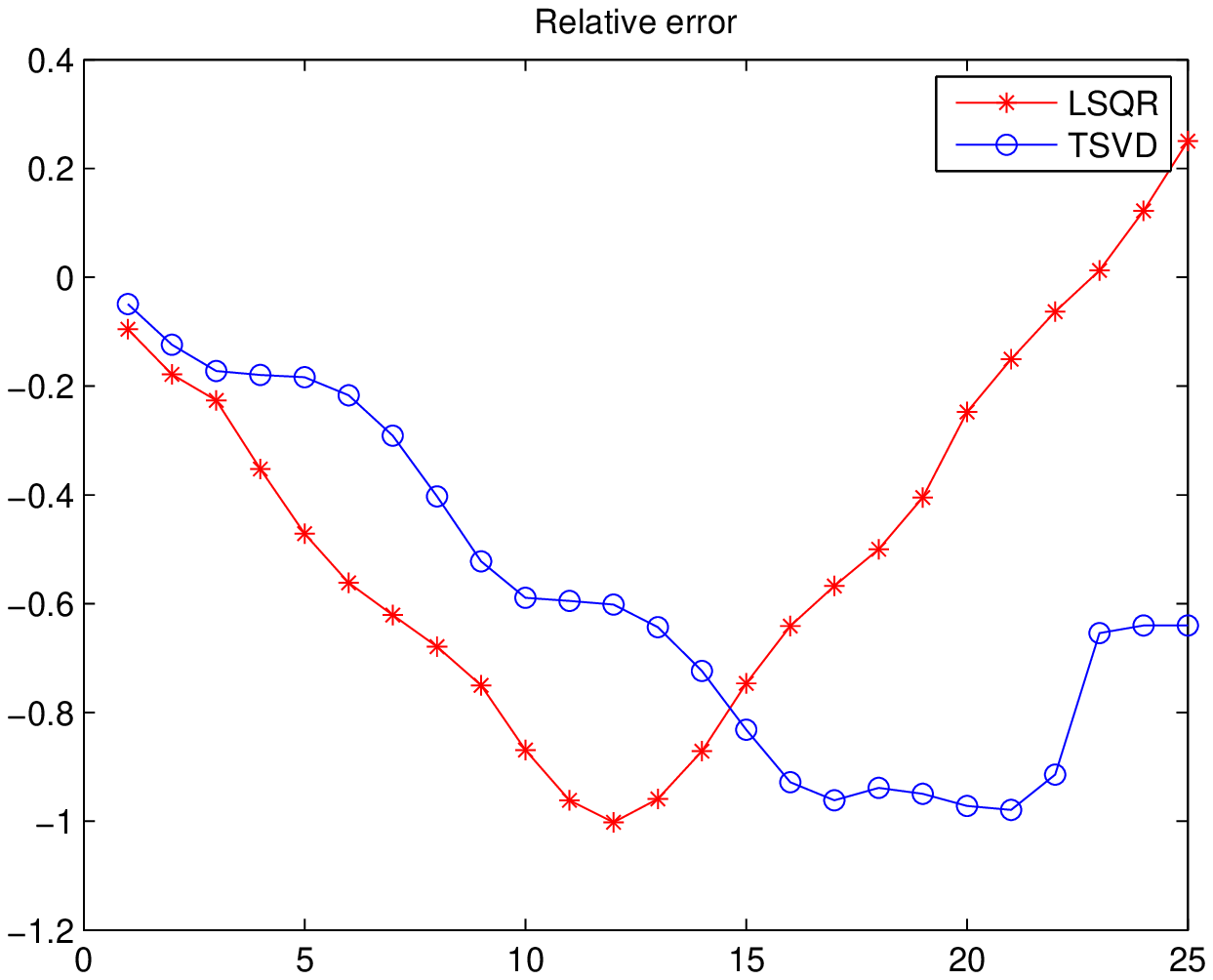}}
  \centerline{(b)}
\end{minipage}

\vfill
\begin{minipage}{0.48\linewidth}
  \centerline{\includegraphics[width=6.0cm,height=5cm]{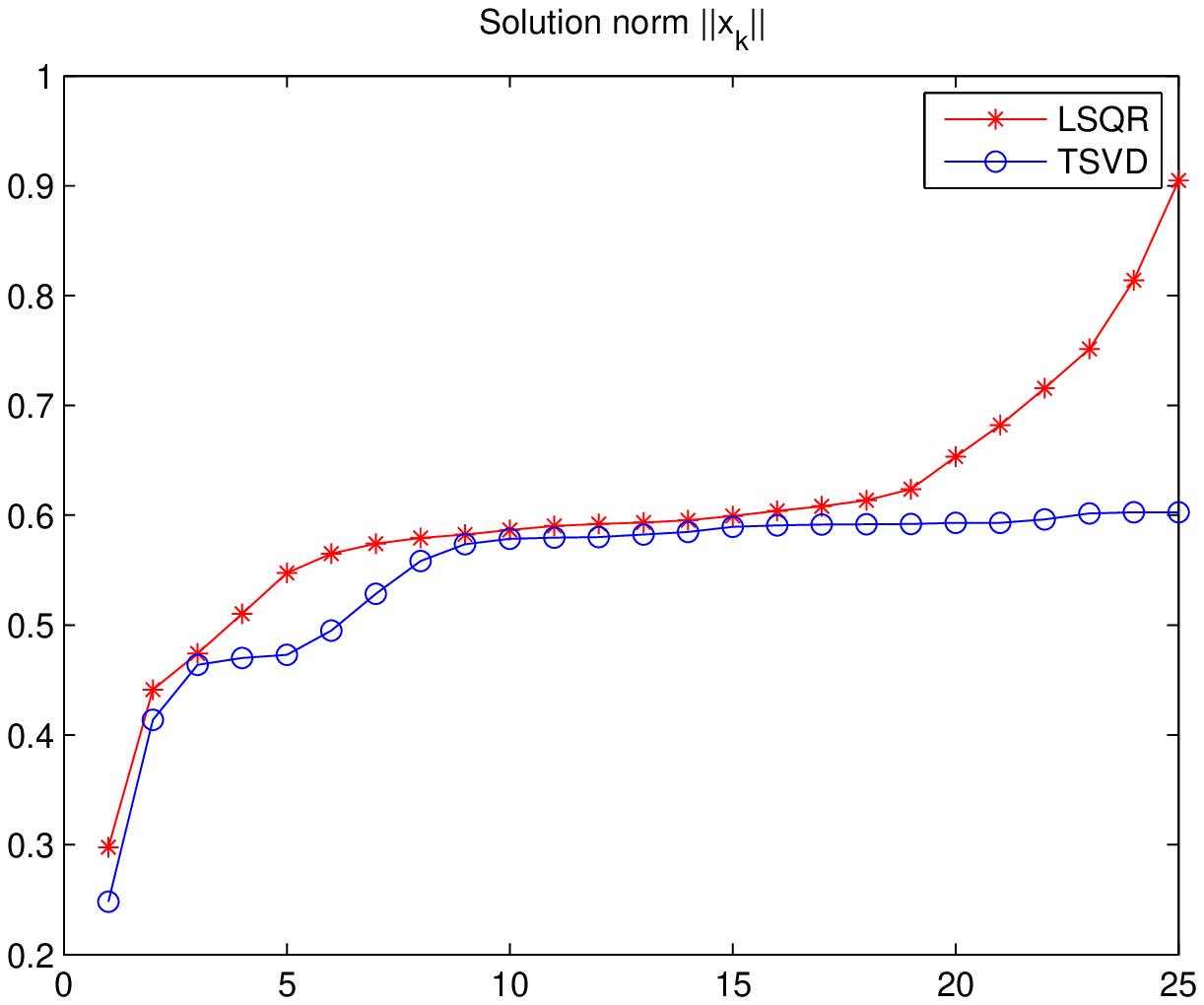}}
  \centerline{(c)}
\end{minipage}
\hfill
\begin{minipage}{0.48\linewidth}
  \centerline{\includegraphics[width=6.0cm,height=5cm]{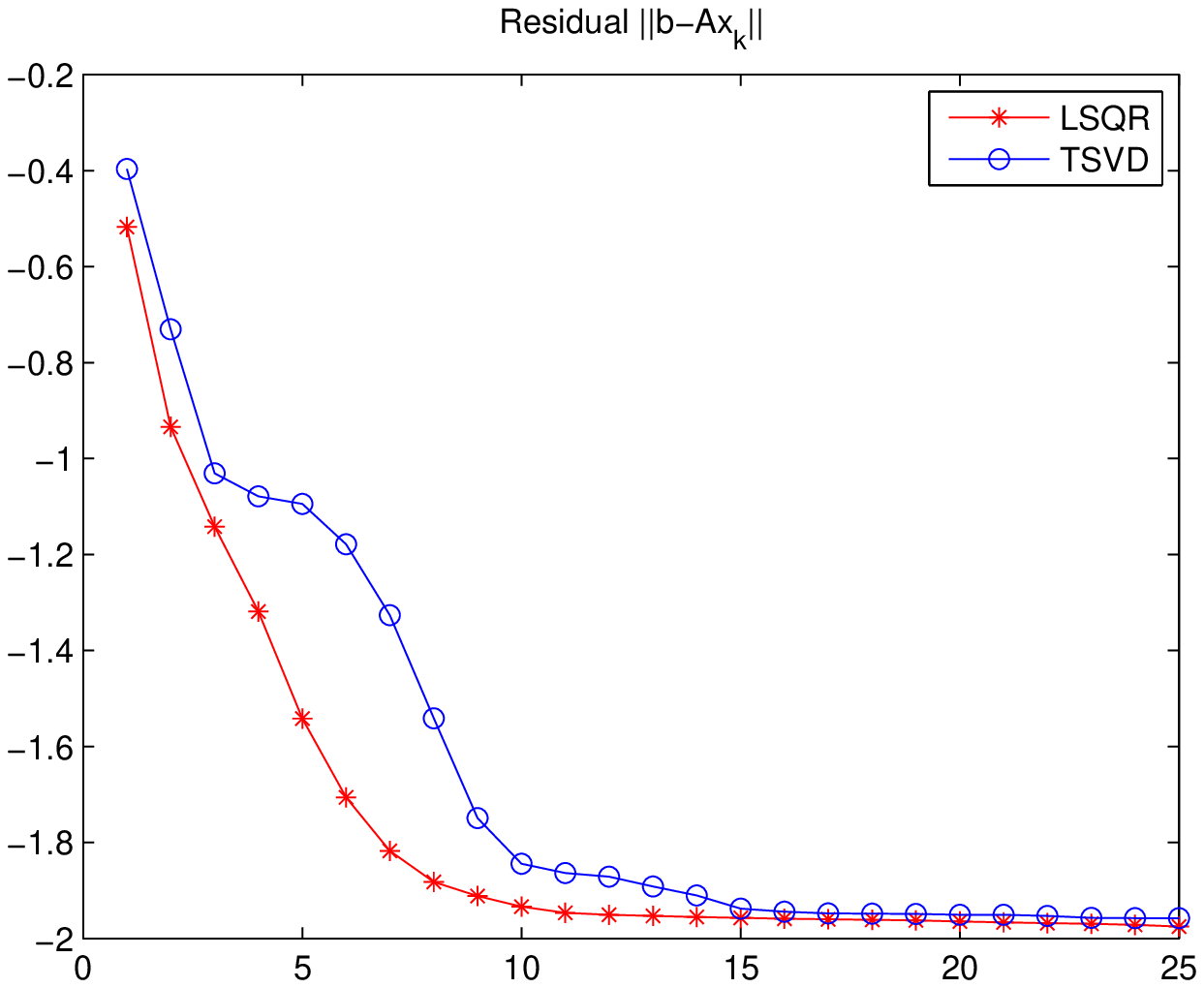}}
  \centerline{(d)}
\end{minipage}
\caption{Results for the moderately ill-posed problem $\mathsf{heat}$.}
\label{lsqrtsvd3}
\end{figure}

\begin{figure}
\begin{minipage}{0.48\linewidth}
  \centerline{\includegraphics[width=6.0cm,height=5cm]{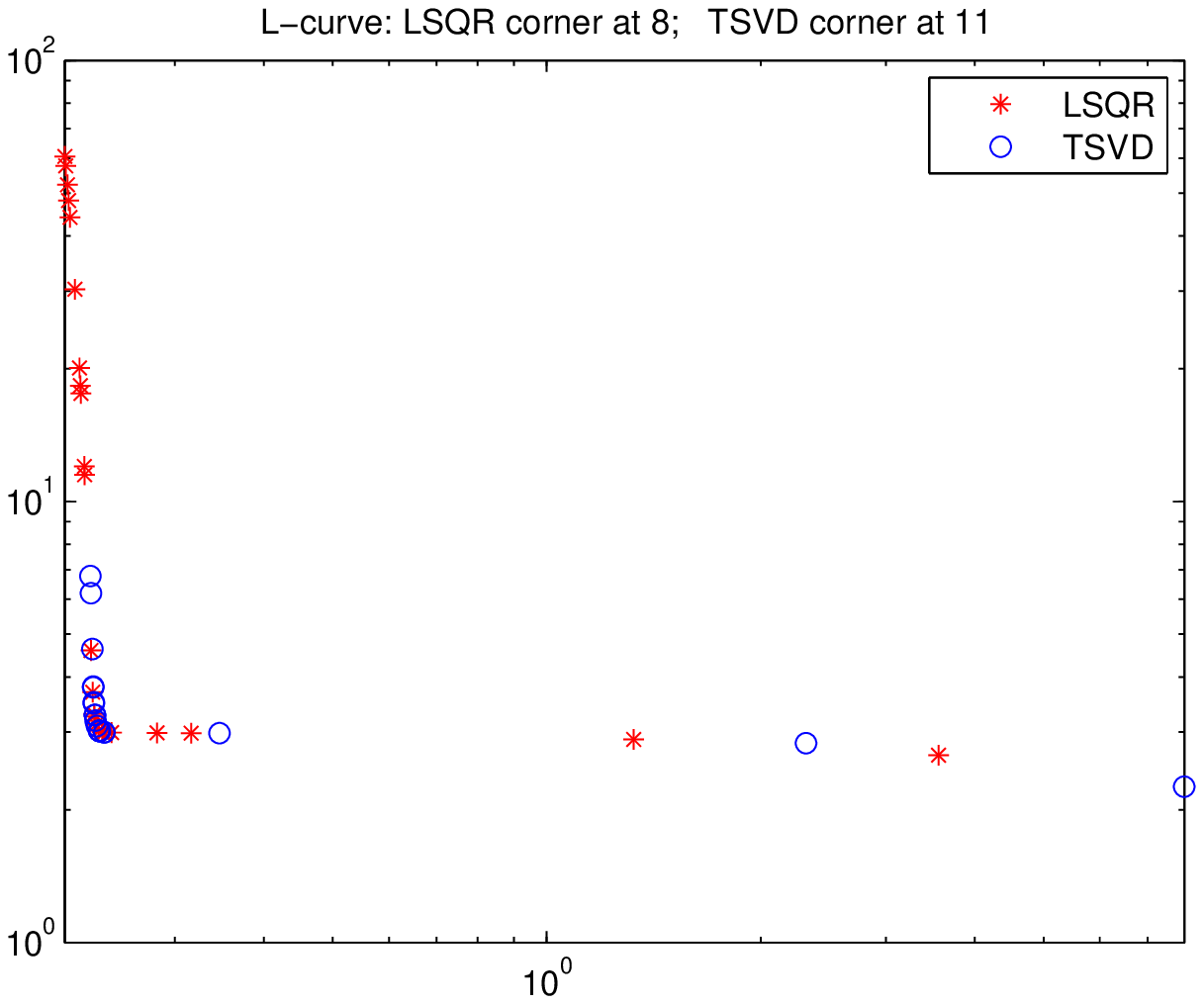}}
  \centerline{(a)}
\end{minipage}
\hfill
\begin{minipage}{0.48\linewidth}
  \centerline{\includegraphics[width=6.0cm,height=5cm]{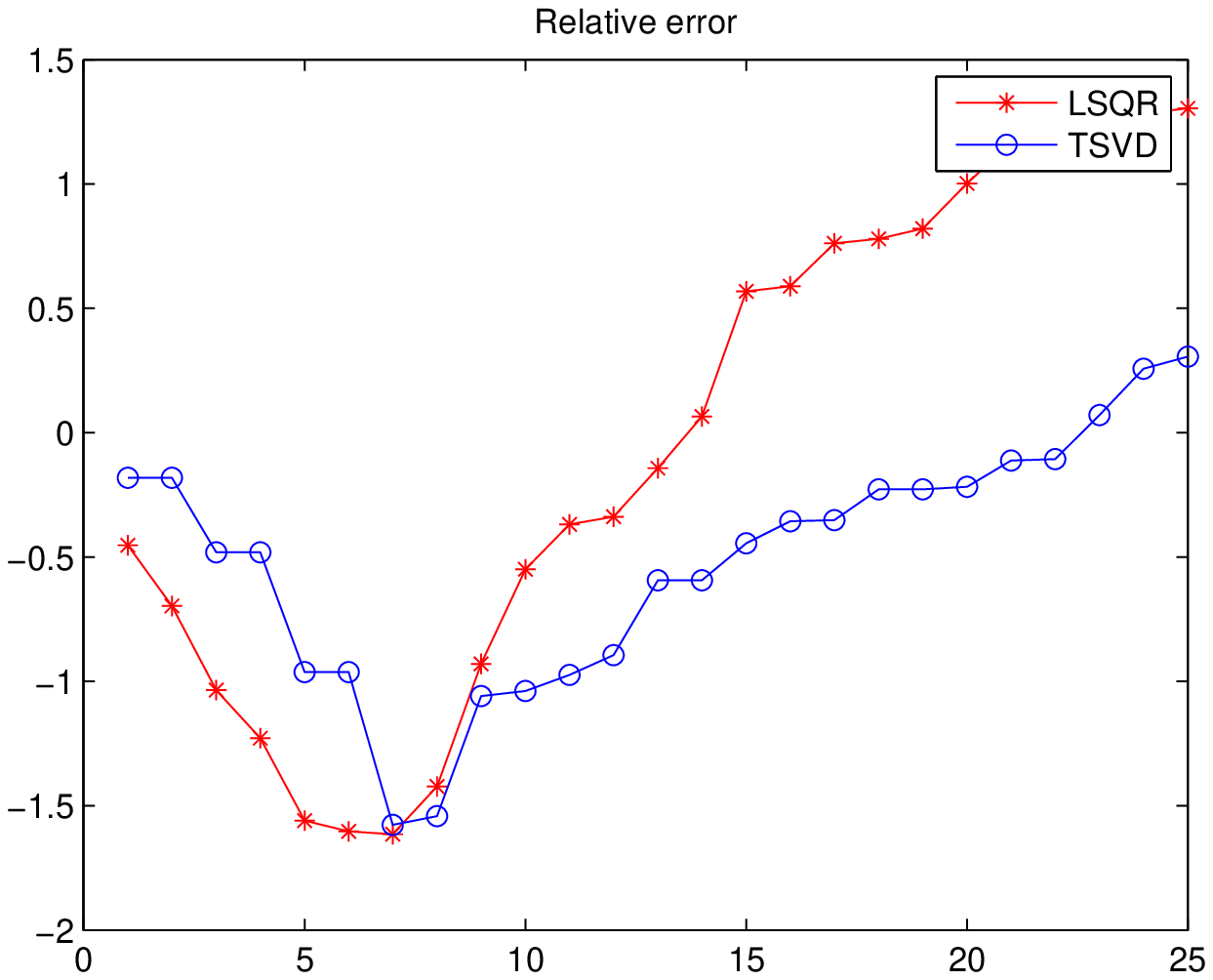}}
  \centerline{(b)}
\end{minipage}

\vfill
\begin{minipage}{0.48\linewidth}
  \centerline{\includegraphics[width=6.0cm,height=5cm]{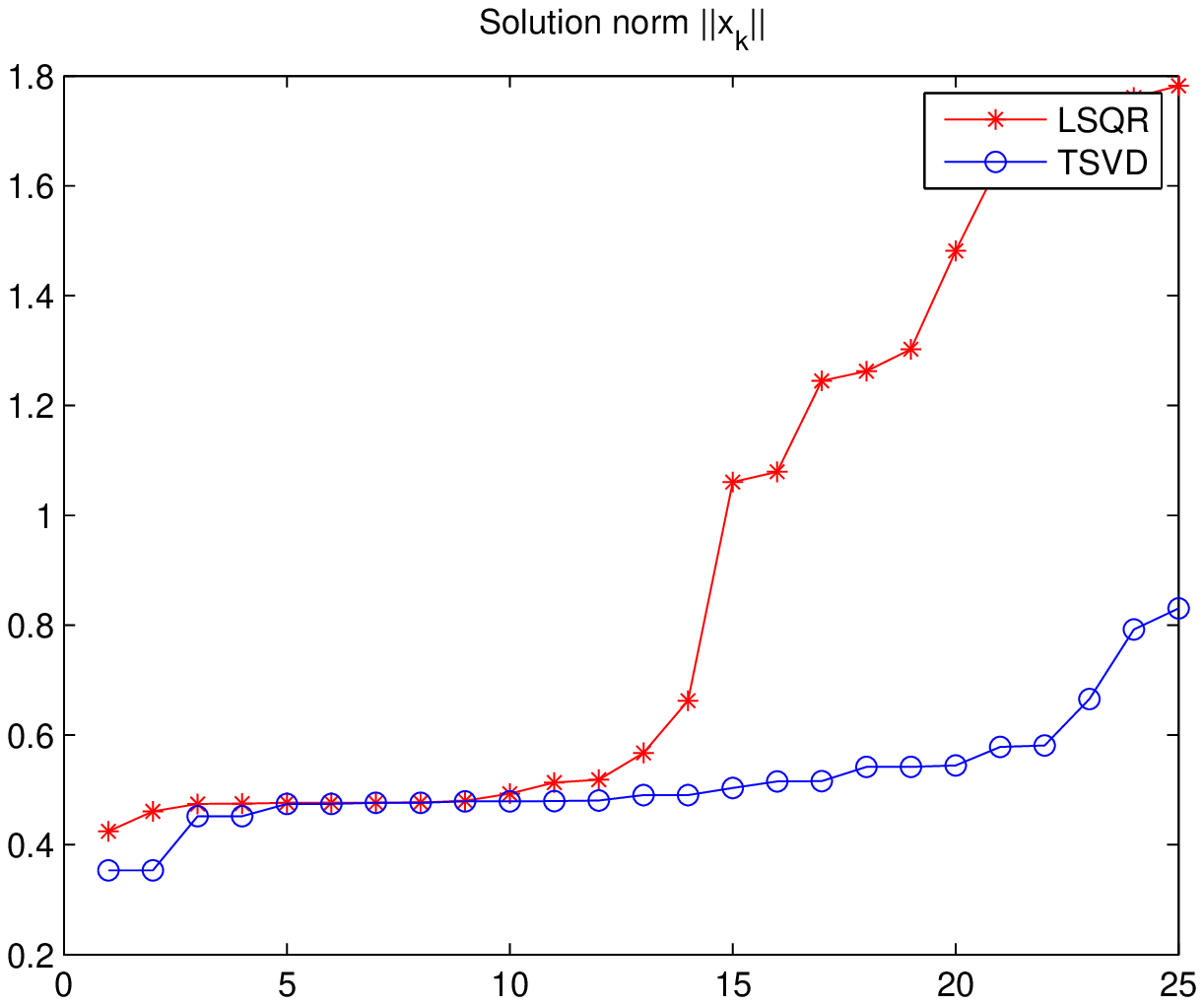}}
  \centerline{(c)}
\end{minipage}
\hfill
\begin{minipage}{0.48\linewidth}
  \centerline{\includegraphics[width=6.0cm,height=5cm]{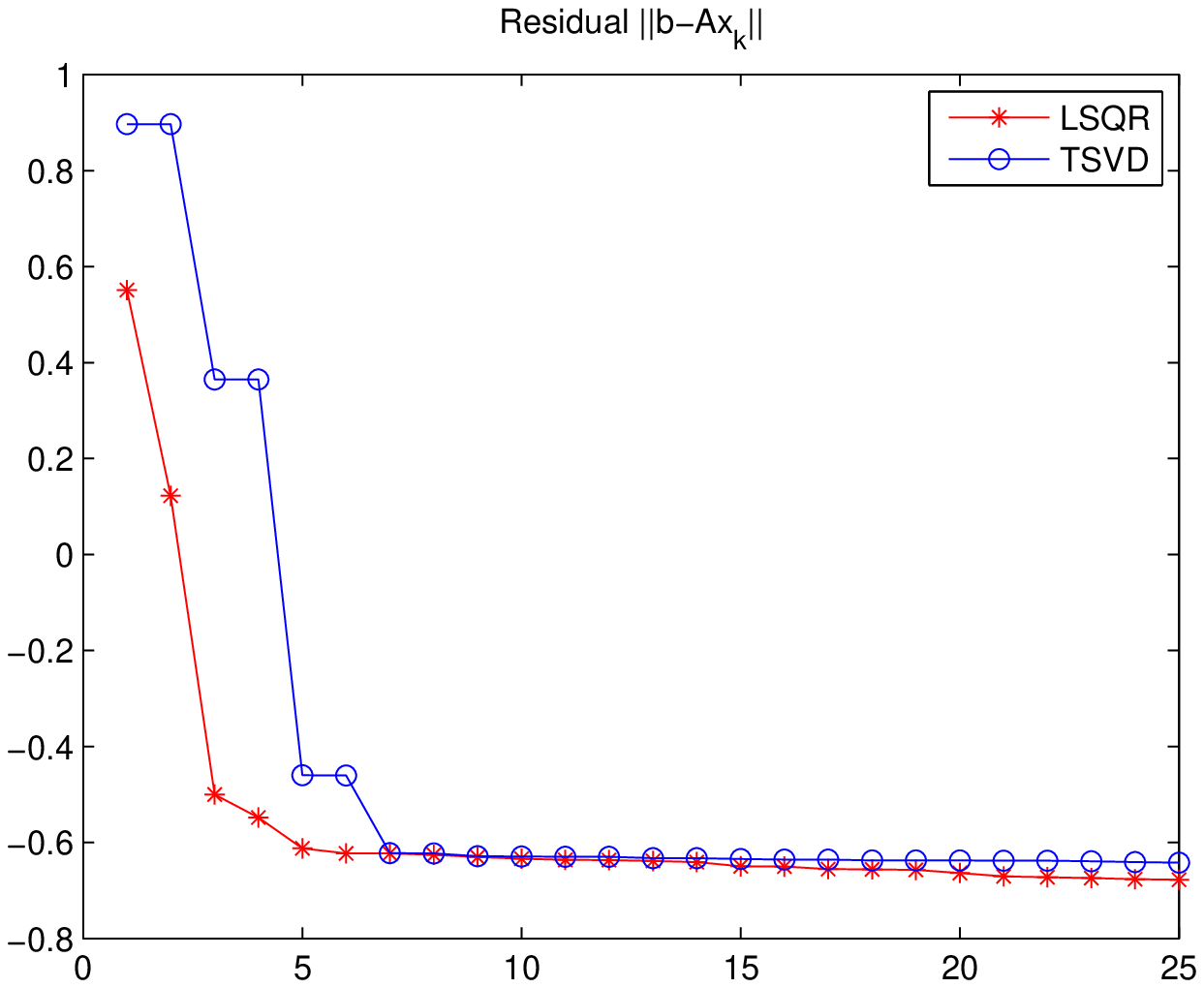}}
  \centerline{(d)}
\end{minipage}
\caption{ Results for the moderately ill-posed problem $\mathsf{phillips}$.}
\label{lsqrtsvd4}
\end{figure}

For each of $\mathsf{heat}$ and $\mathsf{phillips}$,
Figures~\ref{lsqrtsvd3}--\ref{lsqrtsvd4} demonstrate that
the best regularized solution obtained by LSQR is at
least as accurate as, in fact, a little bit more accurate than
that by the TSVD method, and the corresponding
residual norms decreases and drop below at least the same level as those
by the TSVD method. The residual norms by the two methods then stagnate after
the best regularized solutions are found. All these confirm that
LSQR has the full regularization. The fact that the best regularized solutions
by LSQR can be more accurate than the best TSVD solutions is not unusual.
We can explain why. Note that the true solutions $x(t)$ to the
integral equations that generate the problems $\mathsf{heat}$
and $\mathsf{phillips}$ are at least first order
differentiable. It is known that, in the infinite dimensional space setting,
for a linear compact operator equation $Kx=g$,
the TSVD method and standard-form Tiknonov regularization method have been
shown to be order optimal only when the true solution
is continuous or first order differentiable, and they are not order optimal
for stronger smoothness assumptions on the true solution. In contrast,
CGLS is order optimal, and the smallest error of the iterates is of the
same order as the worst-case error for the arbitrarily smooth
true solution, that is, given the same noise level, the smoother the true
solution is, the more accurate the best regularized solution is. In other
words, for the smoother true solution, the best regularized solution by CGLS
is generally more accurate than the counterpart corresponding to
the continuous or first order differentiable true solution;
see, e.g., \cite[p.187-191]{engl00} and \cite[p.13,34-36,40]{kirsch}.
Consequently, for the discrete \eqref{eq1} resulting from such kind of
continuous compact linear equation,  once
the mathematically equivalent LSQR has the full regularization,
its best regularized solution is at least as accurate as and
can be more accurate than the best regularized solution by the TSVD method
or the standard-form Tiknonov regularization method when the true
solution of a continuous compact linear operator
equation is smoother than only continuous or first order differentiable.

From the figures we observe some obvious differences
between moderately and severely ill-posed problems. For $\mathsf{heat}$,
it is seen that the relative errors and
residual norms converge considerably more quickly for the LSQR solutions
than for the TSVD solutions. Figure~\ref{lsqrtsvd3} (b) tells us that
LSQR only uses 12 iterations to find the best regularized
solution and the TSVD method finds the best regularized
solution for $k_0=21$, while the L-curve gives $13$ iterations and
$k_0=18$ iterations, respectively. Similar differences are observed for
$\mathsf{phillips}$, where Figure~\ref{lsqrtsvd4} (b) indicates
that both LSQR and the TSVD method find the
best regularized solutions at $k_0=7$, while the L-curve shows
that $k_0=8$ for LSQR and $k_0=11$ for the TSVD method. Therefore,
unlike for severely ill-posed problems, the L-curve criterion is not
very reliable to determine correct $k_0$ for moderately
ill-posed problems.

We can observe more. Figure~\ref{lsqrtsvd3} shows that
the TSVD solutions improve little and their residual norms decrease
very slowly for the indices $i=4,5,11,12,18,19,20$.  This implies that
the $v_i$ corresponding to these indices $i$ make very little contribution
to the TSVD solutions. This is due to the fact that the Fourier coefficients
$|u_i^T\hat{b}|$ are very small relative to $\sigma_i$ for these indices $i$.
Note that $\mathcal{K}_k(A^TA,A^Tb)$ adapts itself in
an optimal way to the specific right-hand side $b$, while the TSVD method uses
all $v_1,v_2,\ldots,v_k$ to construct a regularized solution,
independent of $b$. Therefore, $\mathcal{K}_k(A^TA,A^Tb)$ picks up only those
SVD components making major contributions to $x_{true}$, such that LSQR uses
possibly fewer $k$ iterations than $k_0$ needed by the TSVD method to capture
those truly needed
dominant SVD components. The fact that LSQR (CGLS) includes fewer SVD components
than the TSVD solution with almost the same accuracy was first noticed
by Hanke \cite{hanke01}.
Generally, for severely and moderately ill-posed problems,
we may deduce that LSQR uses possibly fewer than $k_0$ iterations to compute a
best possible regularized solution if, in practice, some of
$|u_i^T b|$, $i=1,2,\ldots,k_0$ are considerably bigger than
the corresponding $\sigma_i$ and some of them are reverse.
For $\mathsf{phillips}$,
as noted by Hansen \cite[p.32, 123--125]{hansen10}, half of the SVD
components satisfy $u_i^T \hat{b}=v_i^Tx_{true}=0$ for $i$ even, only the
odd indexed $v_1,v_3,\ldots,$ make contributions to $x_{true}$. This is why
the relative errors and residual norms of TSVD solutions
do not decrease at even indices before $x_{k_0}^{tsvd}$ is found.

\subsection{A comparison of LSQR and GMRES, RRGMRES}

GMRES applied to solving \eqref{eq1} with $A$ square computes the iterate
$$
x_k^{g}=\|b\|W_k \bar{H}_k^{\dagger}e_1^{(k+1)}
=W_k \bar{H}_k^{\dagger}W_{k+1}^Tb.
$$
The quantity
\begin{align}
\gamma_k^g&=\|A-W_{k+1}\bar{H}_kW_k^T\| \label{gmresgamma}
\end{align}
measures the accuracy of the rank $k$
approximation $W_{k+1}\bar{H}_kW_k^T$ to $A$, where the columns of
$W_k$ and $W_{k+1}$ are orthonormal bases of $\mathcal{K}_k(A,b)$ and
$\mathcal{K}_{k+1}(A,b)$, respectively,
generated by the Arnoldi process starting with $w_1=b/\|b\|$,
and $\bar{H}_k=W_{k+1}^TAW_k$ is the $(k+1)\times k$ upper Hessenberg
matrix. The size of $\gamma_k^g$ reflects the regularizing effects of
GMRES for solving \eqref{eq1}. We should address that, different from
$\gamma_k$ defined by \eqref{gammak} for LSQR, which has been proved to
decrease monotonically as $k$ increases (cf. \eqref{gammamono}),
mathematically $\gamma_k^g$ has no monotonic property.
Similar to the LSQR iterates $x^{(k)}$ and $\gamma_k$, qualitatively speaking,
if $\gamma_k^g$ decays smoothly in some definitive manner,
then, to some extent, GMRES has regularizing effects;
if they do not decay at all or behave irregularly, then GMRES does not
have regularizing effects and fails to work for \eqref{eq1}.
We test GMRES on the general nonsymmetric $\mathsf{heat}$ and
the following Example 6, and compare it with LSQR.

{\bf Example 6}.
Consider the general nonsymmetric ill-posed problem $\mathsf{i\_laplace}$,
which is severely ill-posed and arises from inverse Laplace transformation.
It is obtained by discretizing the first kind Fredholm integral
equation \eqref{eq2} with $[0,\infty)$ the domains of $s$ and $t$.
The kernel $k(s,t)$, the right-hand side $g(s)$ and the solution $x(t)$ are
given by
\begin{equation*}\label{}
  k(s,t)=\exp(-st),\ \ g(s)=\frac{1}{s+1/2},\ \ x(t)=\exp(-t/2).
\end{equation*}

\begin{figure}
\begin{minipage}{0.48\linewidth}
  \centerline{\includegraphics[width=6.0cm,height=5cm]{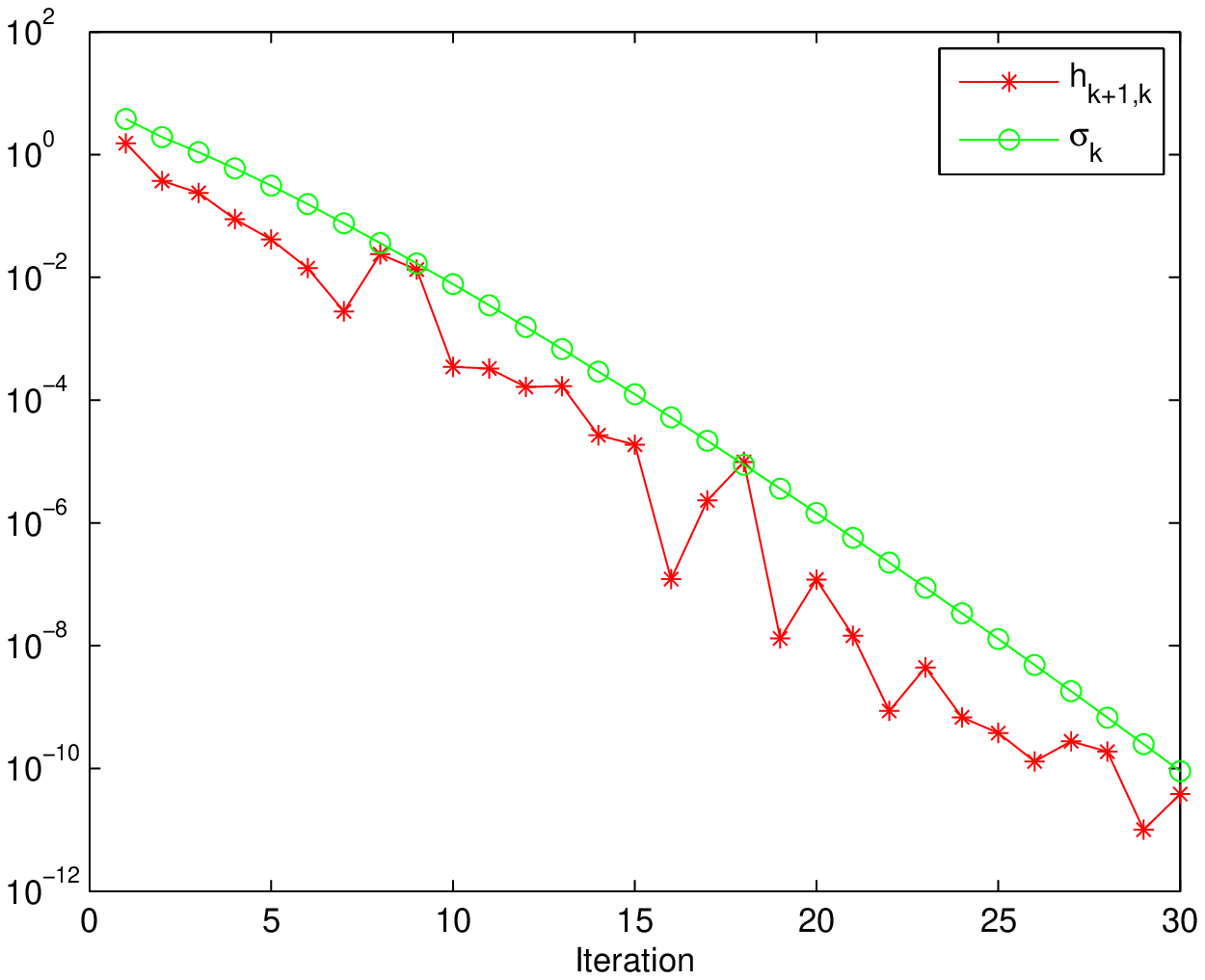}}
  \centerline{(a)}
\end{minipage}
\hfill
\begin{minipage}{0.48\linewidth}
  \centerline{\includegraphics[width=6.0cm,height=5cm]{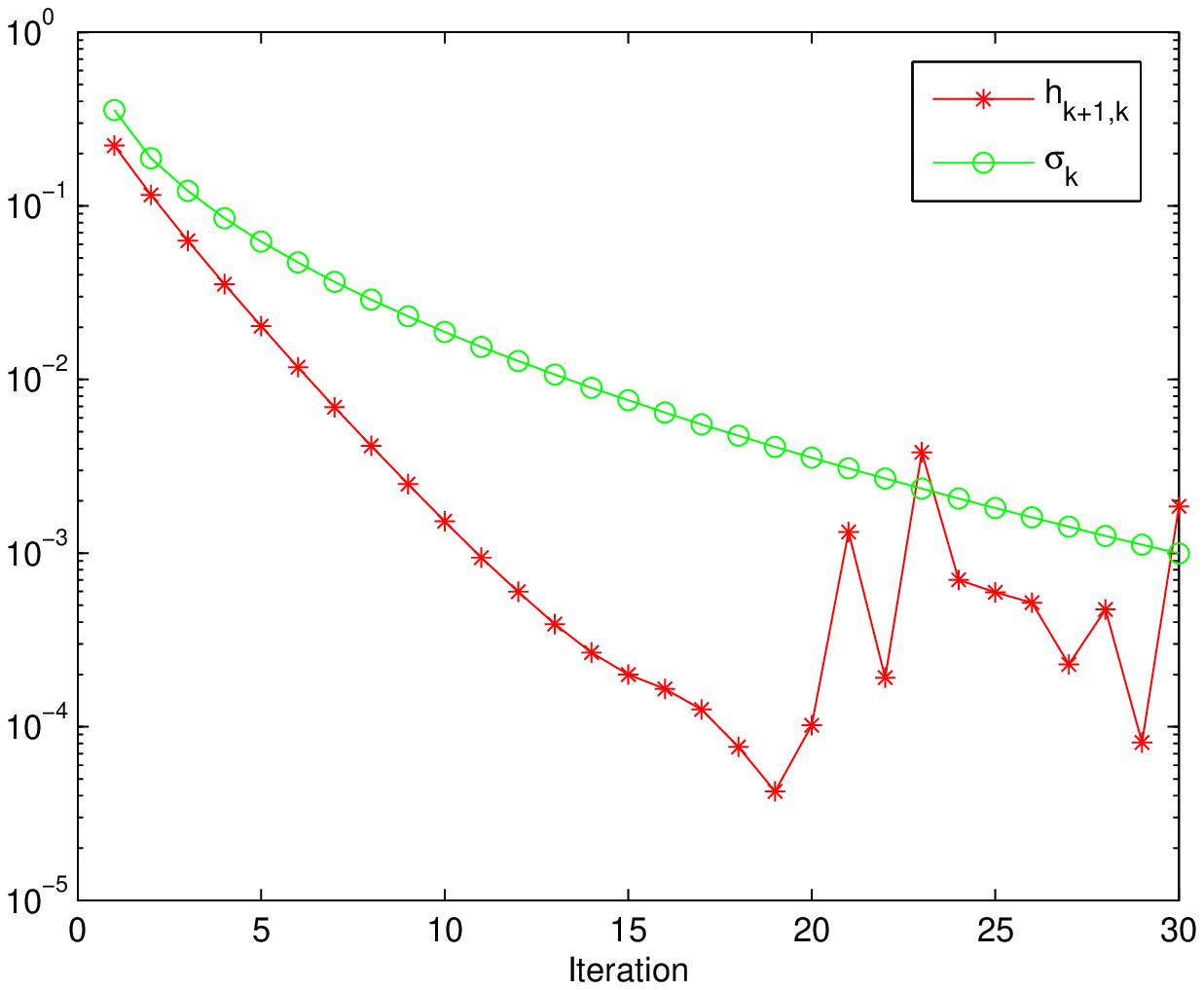}}
  \centerline{(b)}
\end{minipage}
\vfill
\begin{minipage}{0.48\linewidth}
  \centerline{\includegraphics[width=6.0cm,height=5cm]{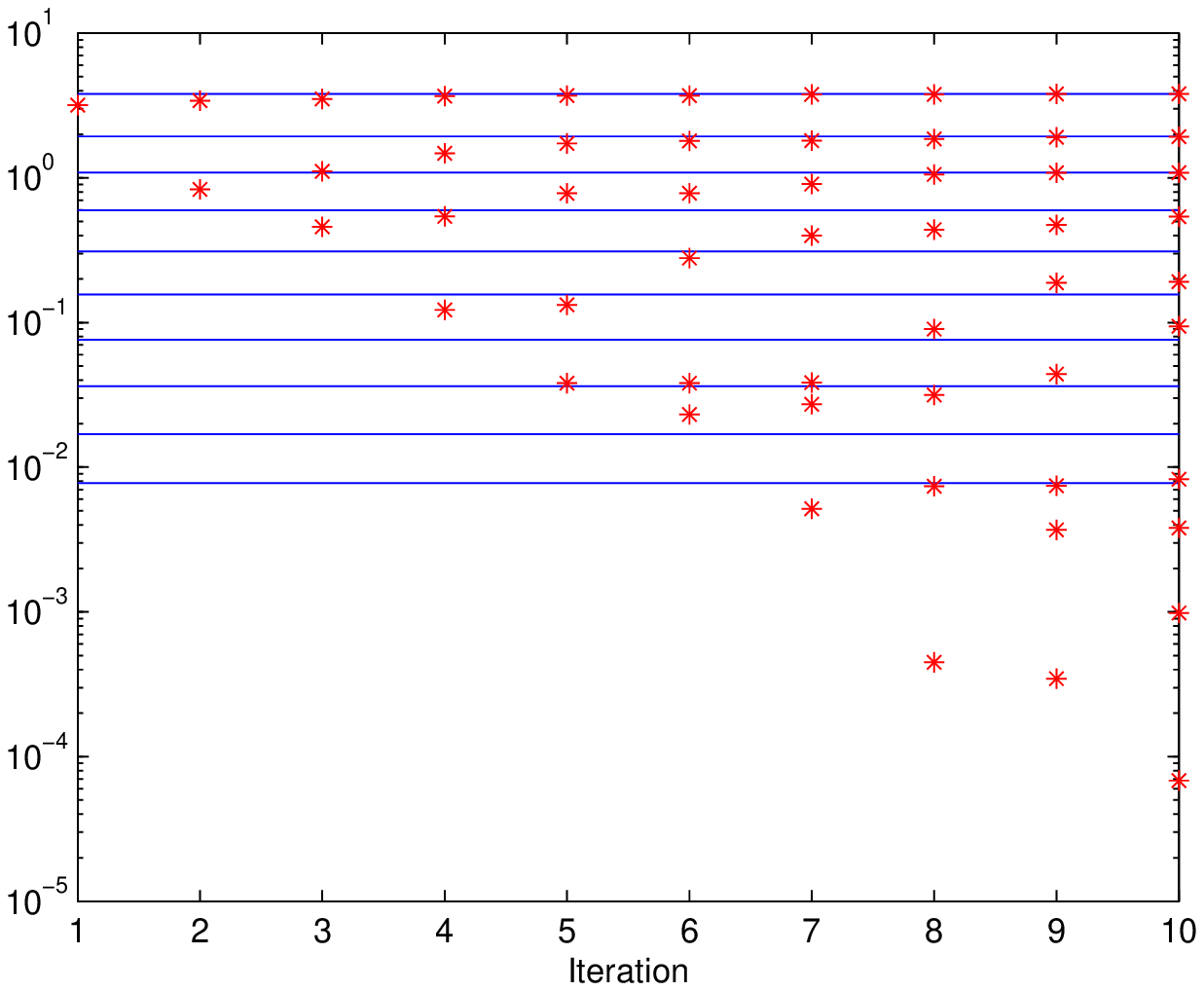}}
  \centerline{(c)}
\end{minipage}
\hfill
\begin{minipage}{0.48\linewidth}
  \centerline{\includegraphics[width=6.0cm,height=5cm]{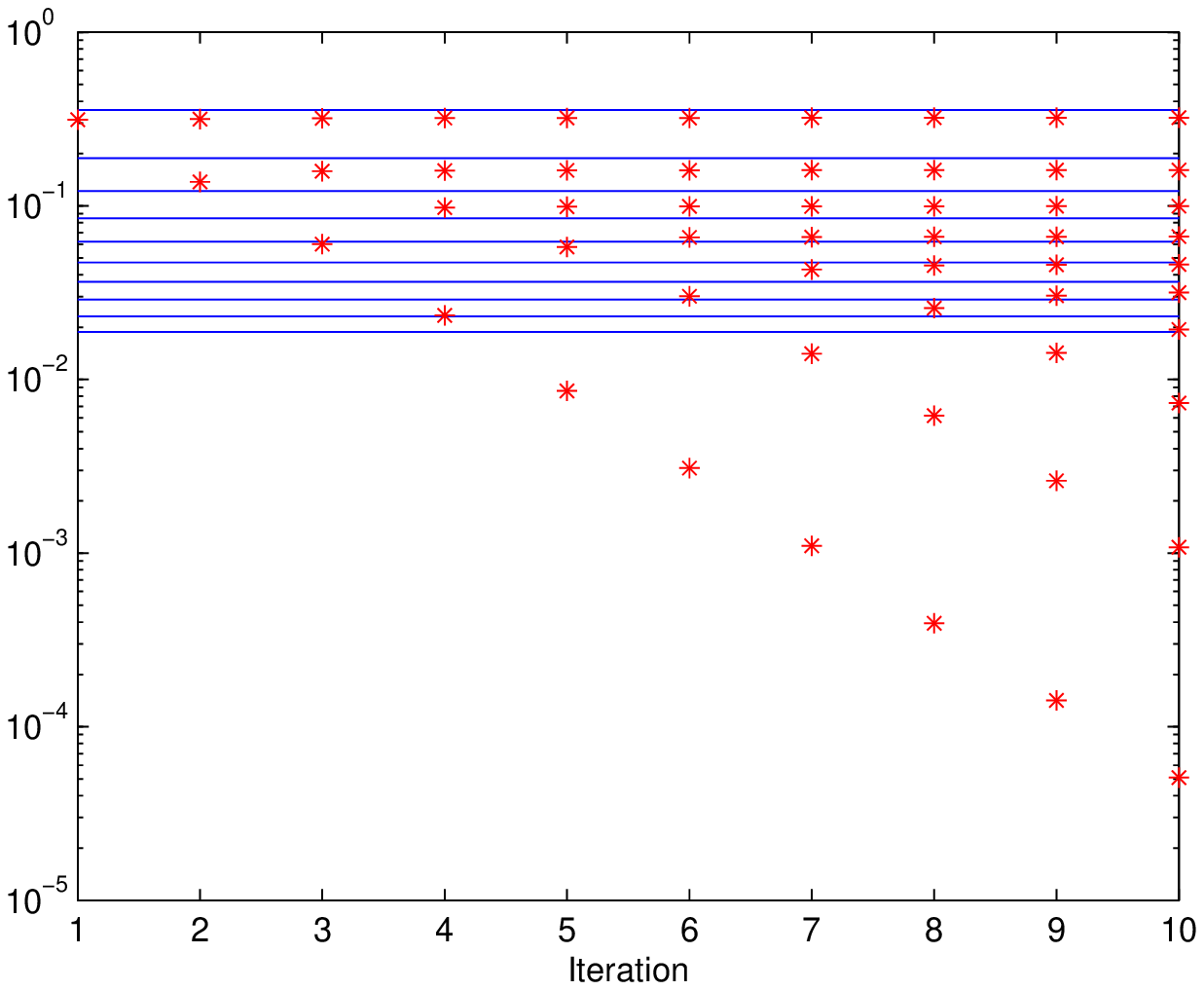}}
  \centerline{(d)}
\end{minipage}
\caption{(a)-(b): Decay curves of the sequences $h_{k+1,k}$
and $\sigma_{k+1}$; (c)-(d): The singular values (star) of
$\bar{H}_k$ and the ones (solid line) of $A$ for
$\mathsf{i\_laplace}$ (left) and $\mathsf{heat}$ (right).}
\label{fig8}
\end{figure}

We investigate the regularizing effects of GMRES with $\varepsilon=10^{-3}$.
Let $H_k=(h_{i,j})\in \mathbb{R}^{k\times k}$ denote the upper
Hessenberg matrix obtained by the $k$-step Arnoldi process. We
observe that the $h_{k+1,k}$ decay quickly with $k$ increasing, generally faster
than $\sigma_k$; see Figure~\ref{fig8} (a)-(b).
This phenomenon may lead to a misbelief that GMRES has
general regularizing effects. However, it is not the case. In fact,
a small $h_{k+1,k}$ exactly indicates that all the eigenvalues of $H_k$
may approximate some $k$ eigenvalues of $A$ well and the Arnoldi method finds an
approximate $k$-dimensional invariant subspace or eigenspace of $A$
\cite{jia01,saad}. We also refer to \cite{jia95,jia98} for a detailed
convergence analysis of the Arnoldi method. Unfortunately,
for a general nonsymmetric matrix $A$, a
small $h_{k+1,k}$ does not mean that the singular values, i.e., the Ritz values,
of $\bar{H}_k$ are also good approximations to some $k$ singular values of
$A$. As a matter of fact, as our analysis in Section \ref{rankapp} and
Section \ref{alphabeta} has indicated, the accuracy of the singular values
of a rank $k$ approximation matrix, here the singular values of $\bar{H}_k$,
as approximations to the $k$ large singular values of $A$ critically relies on
the size of $\gamma_k^g$ defined by \eqref{gmresgamma} other than $h_{k+1,k}$.
Indeed, as indicated by Figure~\ref{fig8} (c)-(d), though
$h_{k+1,k}$ is small, some of the singular values of $\bar{H}_k$ are
very poor approximations to singular values of $A$,
and some of those good approximations are much smaller than $\sigma_{k+1}$
and approximate the singular values of $A$ in disorder
rather than in natural order. It is important to note that, for a
general nonsymmetric or, more rigorously, non-normal $A$,
the $k$-dimensional Krylov subspace $\mathcal{K}_k(A,b)$ that
underlies the Arnoldi process mixes all the left and right singular vectors of
$A$, and the Arnoldi process generally fails to
extract the dominant SVD components and cannot generate a high quality rank
$k$ approximation to $A$, causing that GMRES has no good regularizing effects.

\begin{figure}
\begin{minipage}{0.48\linewidth}
  \centerline{\includegraphics[width=6.0cm,height=5cm]{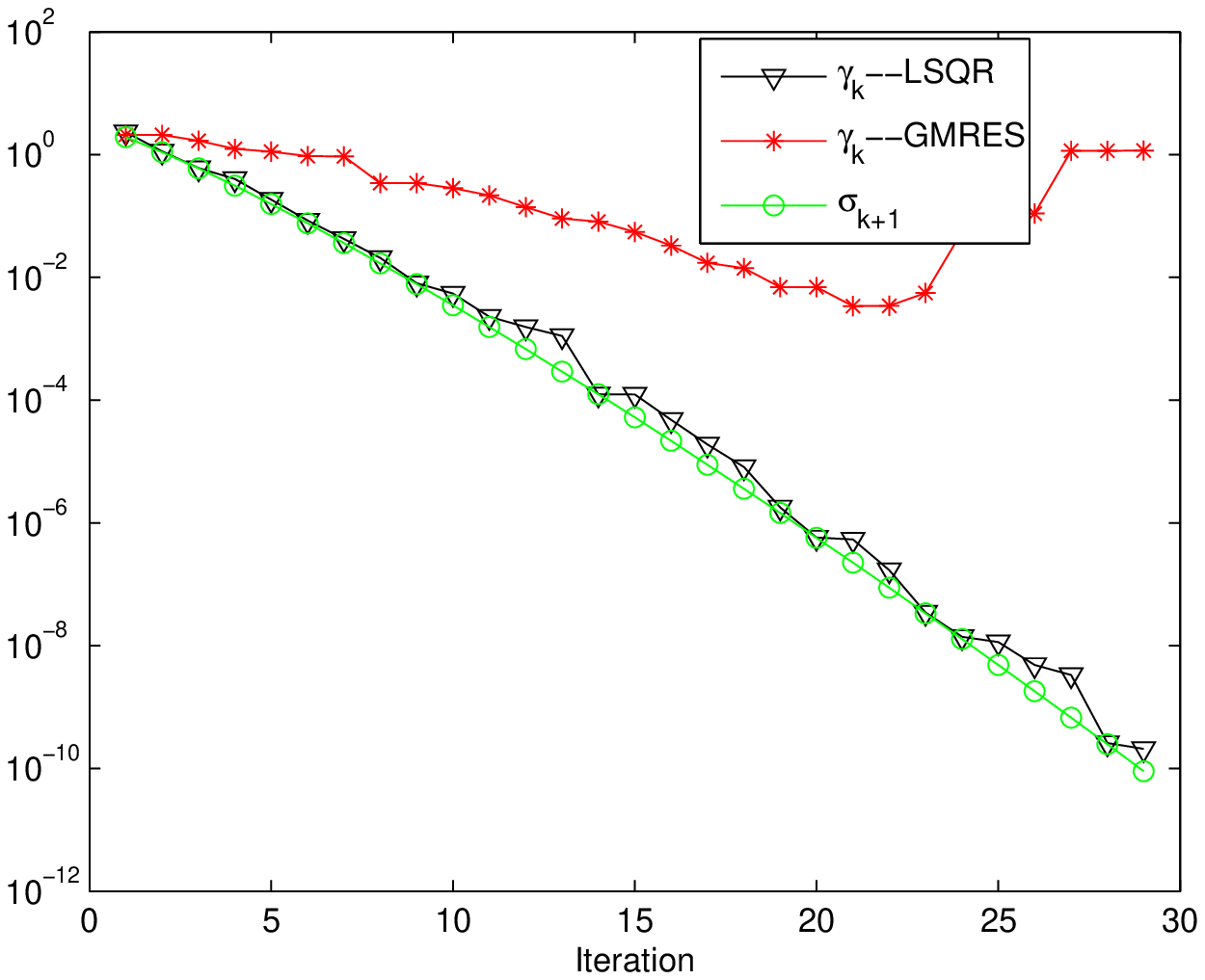}}
  \centerline{(a)}
\end{minipage}
\hfill
\begin{minipage}{0.48\linewidth}
  \centerline{\includegraphics[width=6.0cm,height=5cm]{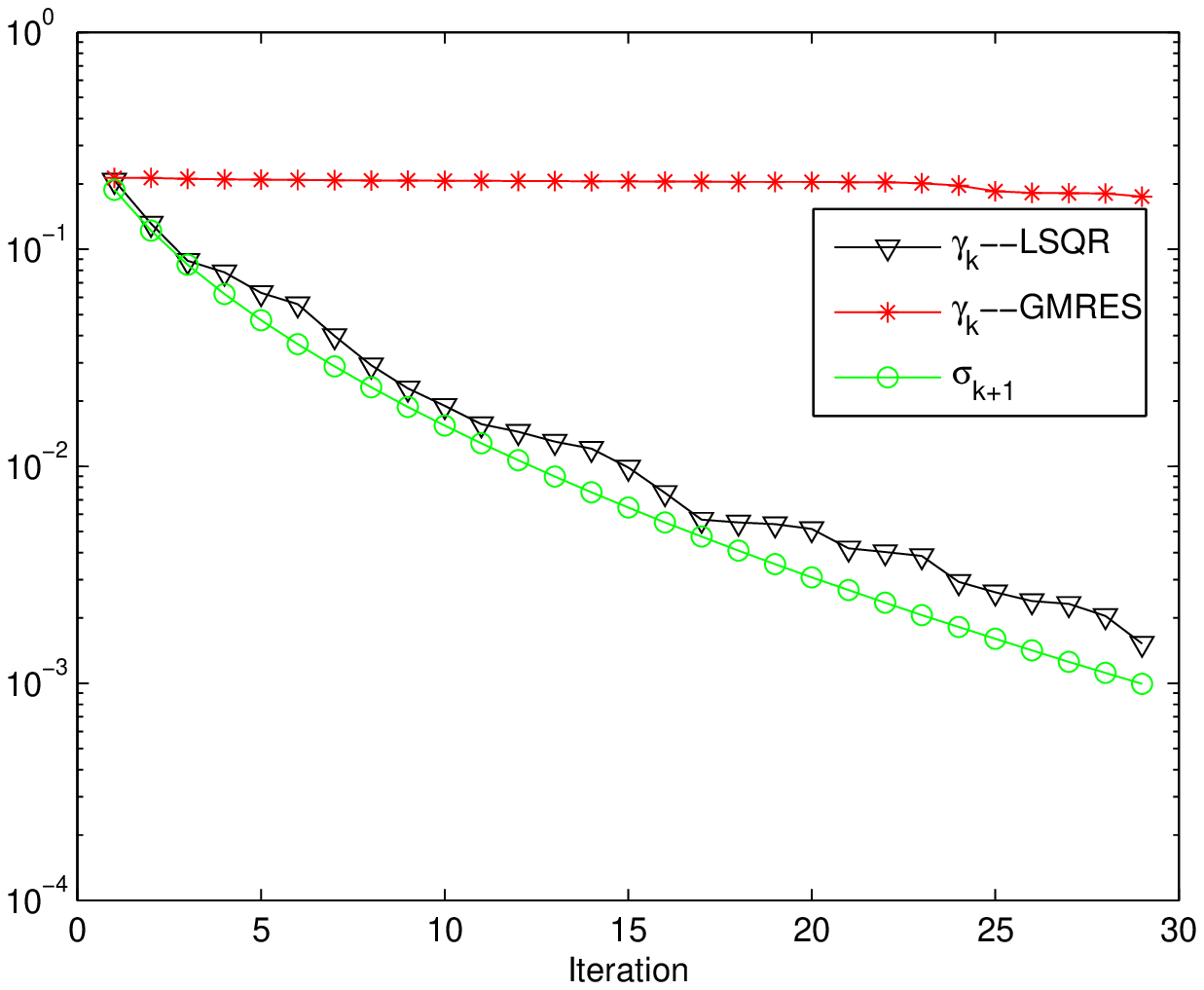}}
  \centerline{(b)}
\end{minipage}
\vfill
\begin{minipage}{0.48\linewidth}
  \centerline{\includegraphics[width=6.0cm,height=5cm]{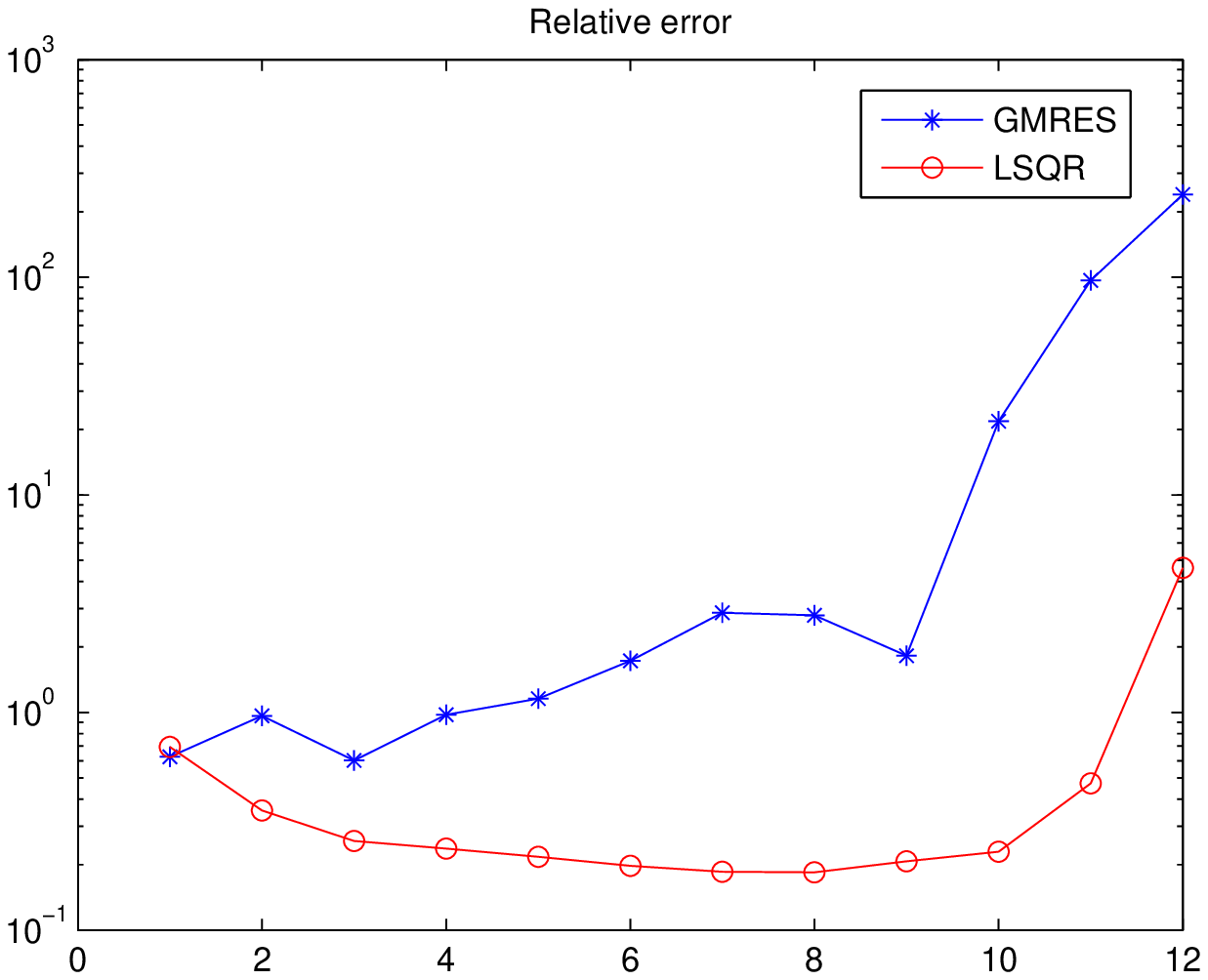}}
  \centerline{(c)}
\end{minipage}
\hfill
\begin{minipage}{0.48\linewidth}
  \centerline{\includegraphics[width=6.0cm,height=5cm]{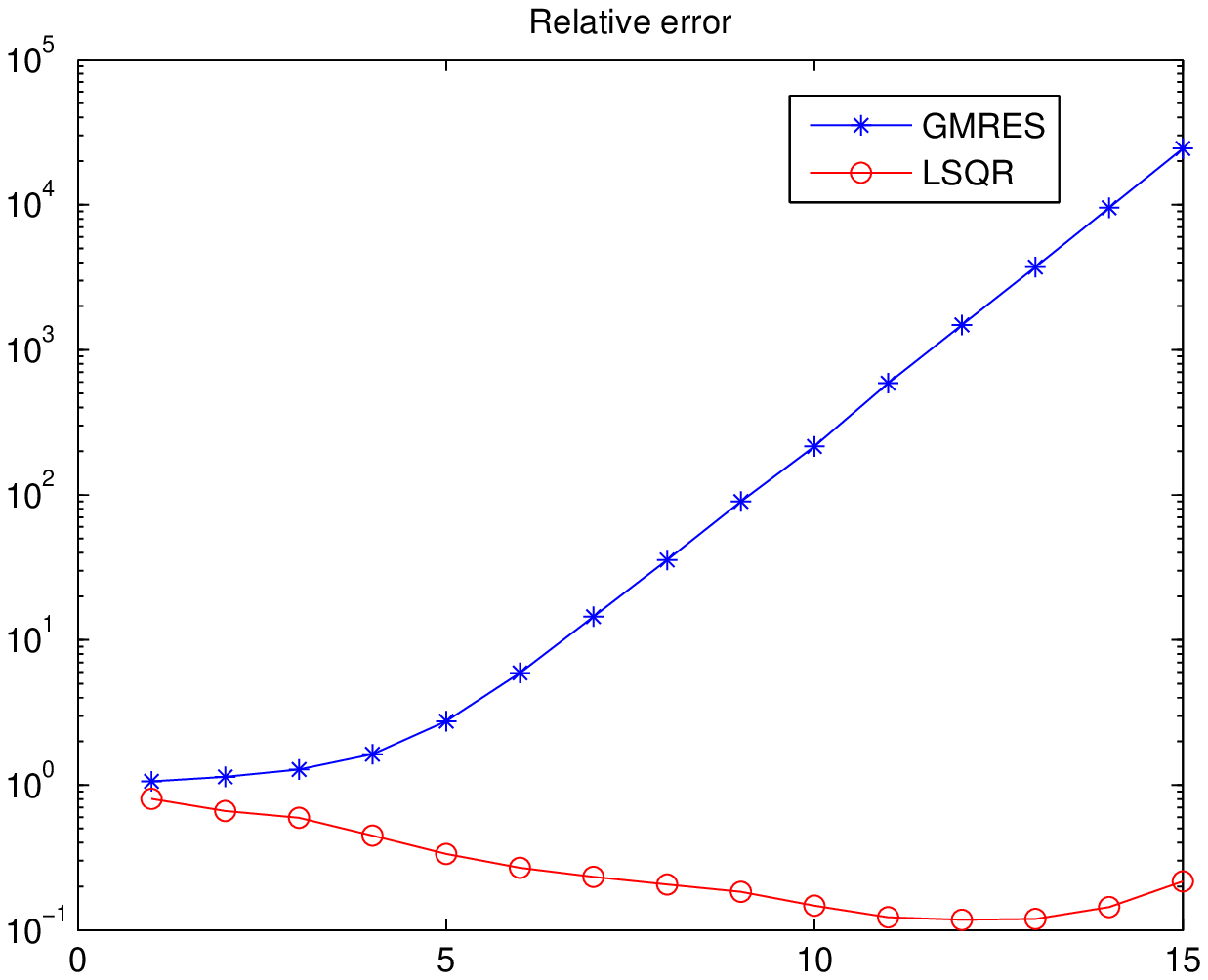}}
  \centerline{(d)}
\end{minipage}
\vfill
\begin{minipage}{0.48\linewidth}
  \centerline{\includegraphics[width=6.0cm,height=5cm]{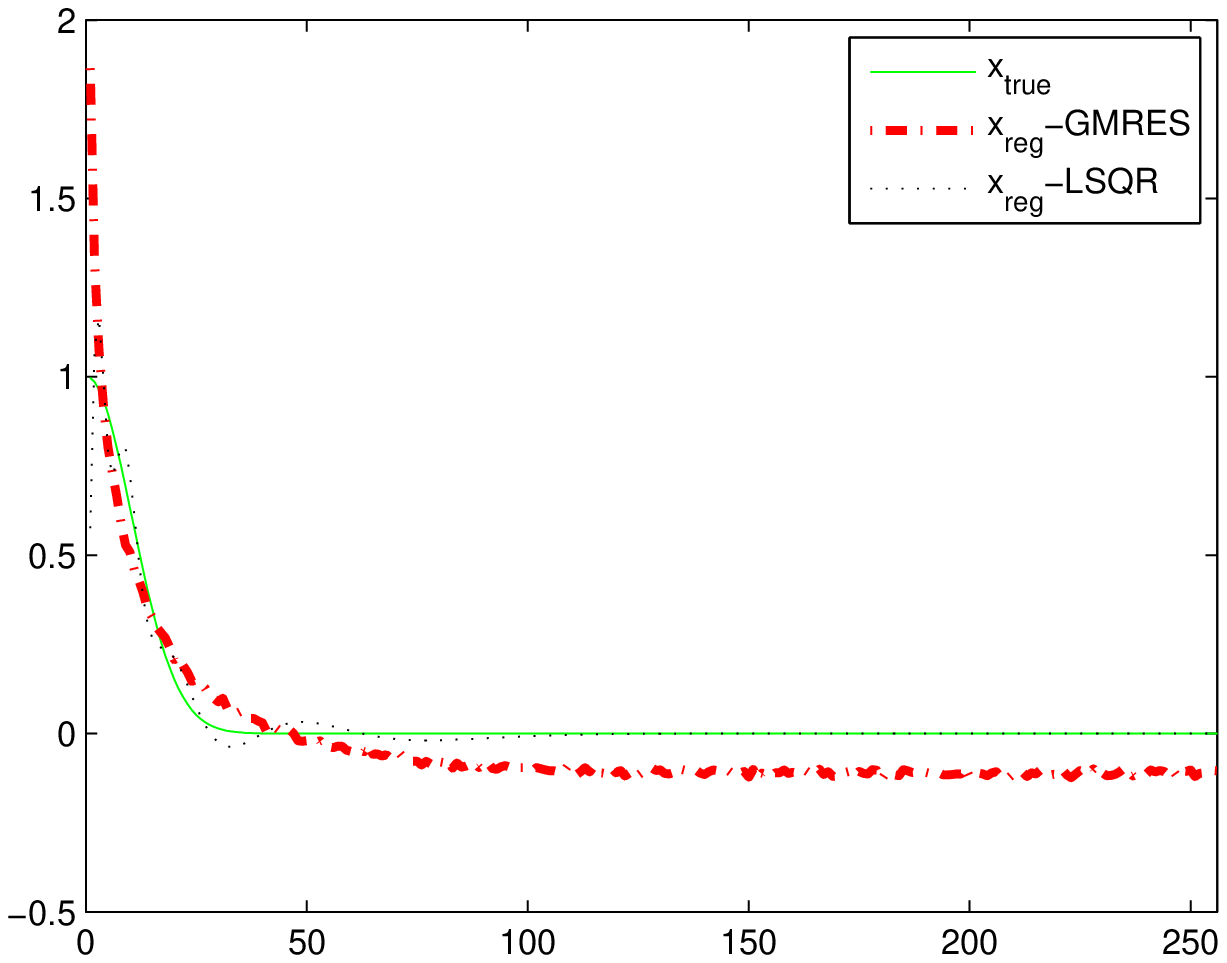}}
  \centerline{(e)}
\end{minipage}
\hfill
\begin{minipage}{0.48\linewidth}
  \centerline{\includegraphics[width=6.0cm,height=5cm]{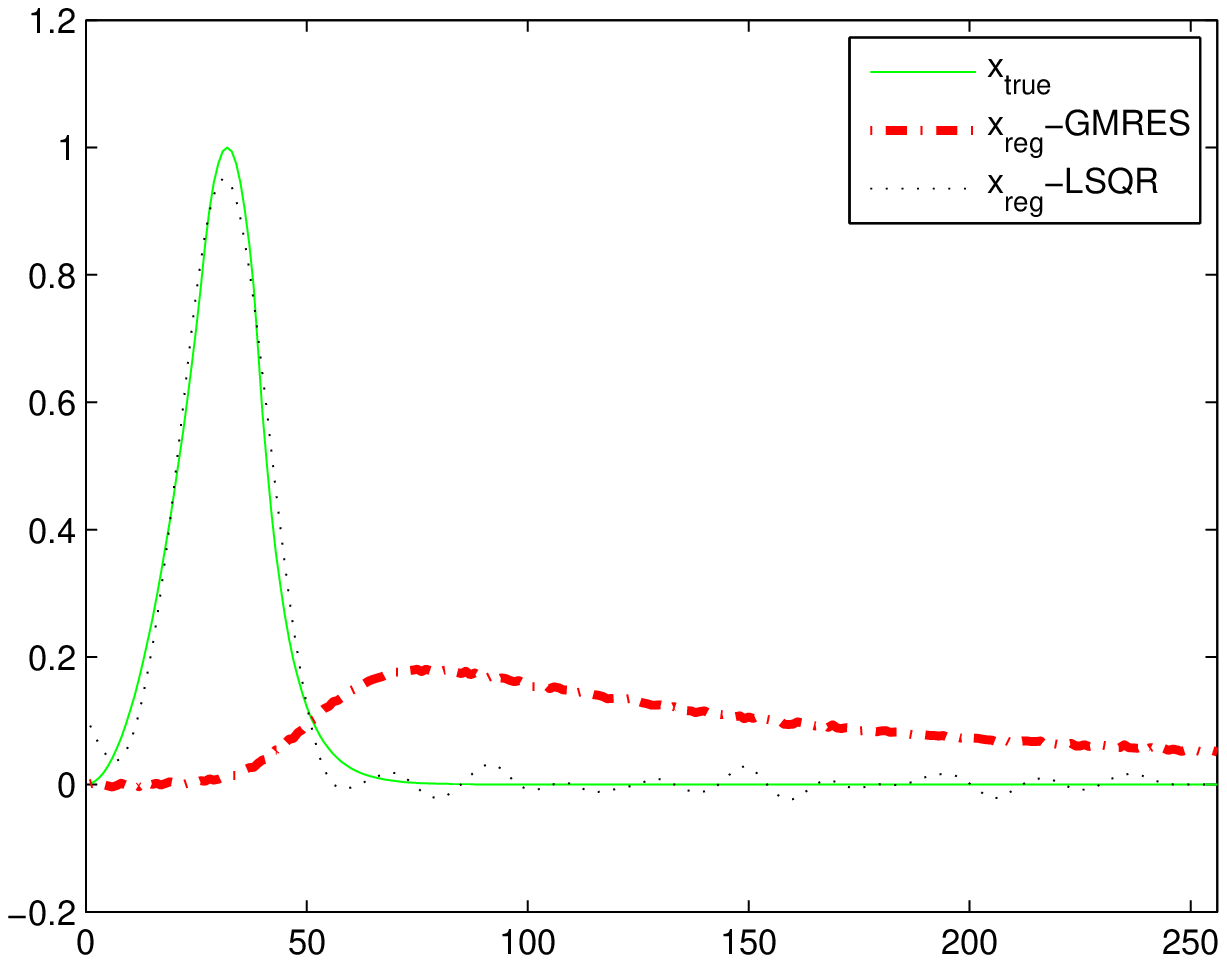}}
  \centerline{(f)}
\end{minipage}
\caption{(a)-(b): Decay curves of the sequences $\gamma_k$, $\gamma_k^g$,
denoted by $\gamma_k$-LSQR and $\gamma_k$-GMRES in the figure,
and $\sigma_{k+1}$; (c)-(d): The relative errors
$\|x^{(k)}-x_{true}\|/\|x_{true}\|$; (e)-(f): The regularized
solutions $x_{reg}$ obtained by LSQR and GMRES for
$\mathsf{i\_laplace}$ (left) and $\mathsf{heat}$ (right).} \label{fig5}
\end{figure}

Figure~\ref{fig5} (a)-(b) gives more justifications.
We have a few important observations: For the two test problems,
the quantities $\gamma_k$ decay as fast as the $\sigma_{k+1}$ for LSQR,
while the $\gamma_k^g$ diverge quickly from the $\sigma_{k+1}$ for GMRES and
do not exhibit any regular decreasing tendency. For $\mathsf{i\_laplace}$,
the $\gamma_k^g$ decrease very slowly until $k=19$, then basically stabilize
for three iterations followed, and finally start to
increase from $k=22$ onwards. As for $\mathsf{heat}$, the
$\gamma_k^g$ are almost constant from beginning to end.
Since all the $\gamma_k^g$ are not small, they
indicate that the Arnoldi process cannot generate any reasonable and meaningful
rank $k$ approximations to $A$ for $k=1,2,\ldots,30$. This is especially true
for $\mathsf{heat}$. Consequently, we are sure that GMRES
fails and does not have regularizing effects
for the two test problems.

We plot $\|x^{(k)}-x_{true}\|/\|x_{true}\|$
by LSQR and $\|x_k^{g}-x_{true}\|/\|x_{true}\|$ by GMRES in
Figure~\ref{fig5} (c)-(d). Obviously, LSQR exhibits semi-convergence,
but GMRES does not and the relative errors obtained by it even increase from
the beginning; see Figure~\ref{fig5} (d). This again demonstrates that
GMRES cannot provide meaningful regularized solutions for these two problems.
Let $x_{reg}=\arg\min _{k}\|x^{(k)}-x_{true}\|$. Figure~\ref{fig5} (e) and
(f) show that LSQR obtains excellent regularized
solutions, while GMRES fails. It is known that
MR-II \cite{fischer96} for $A$ symmetric and
RRGMRES \cite{calvetti00} for $A$ nonsymmetric work on the subspace
$\mathcal{K}_k(A,Ab)$. They were originally designed to solve singular or
inconsistent systems, restricted to
a subspace of range of $A$, and compute the minimum-norm least squares
solutions when the ranges of $A$ and $A^T$ are identical. However,
for the preferred RRGMRES \cite{neuman12}, we have observed phenomena
similar to those for GMRES, illustrating that RRGMRES does not have
regularizing effects for the test problems.
From these typical experiments, we conclude that GMRES and
RRGMRES are susceptible to failure for general nonsymmetric ill-posed
problems and they are not general-purpose regularization methods.
In fact, as addressed in \cite[p.126]{hansen10} and \cite{jensen07},
GMRES and RRGMRES may only work well when either
the mixing of SVD components is weak or the Krylov basis vectors are
just well suited for the ill-posed problem, as addressed in \cite{hansen10}.

\section{Conclusions}\label{concl}

For the large-scale ill-posed problem \eqref{eq1}, iterative solvers
are the only viable approaches. Of them, LSQR and CGLS are most popularly
used for general purposes, and CGME and LSMR are also choices.
They have general regularizing effects and exhibit semi-convergence.
However, if semi-convergence occurs before it
capture all the needed dominant SVD components, then best possible
regularized solutions are not yet found and the solvers have only the
partial regularization. In this case, their hybrid variants have often
been used to compute best possible regularized solutions.
If semi-convergence means that they have already found
best possible regularized solutions, they have the full regularization,
and we simply stop them after semi-convergence.

We have considered the fundamental open question in depth: Do LSQR, CGLS,
LSMR and CGME have the full or partial regularization for severely, moderately
and mildly ill-posed problems? We have first considered the case that
all the singular values of $A$ are simple. As a key and indispensable step, we
have established accurate bounds for the 2-norm
distances between the underlying $k$ dimensional Krylov subspace and the
$k$ dimensional dominant right singular subspace for the three kinds of
ill-posed problems under consideration. Then we have
provided other absolutely necessary background and ingredients. Based on them,
we have proved that, for severely or moderately ill-posed problems
with $\rho>1$ or $\alpha>1$ suitably, LSQR has the full regularization.
Precisely, for $k\leq k_0$
we have proved that a $k$-step Lanczos bidiagonalization produces a near best
rank $k$ approximation of $A$ and the $k$ Ritz values approximate
the first $k$ large singular values of $A$ in natural order, and
no small Ritz value smaller than $\sigma_{k_0+1}$ appears before LSQR
captures all the needed dominant SVD components, so that
the noise $e$ in $b$ cannot deteriorate regularized
solutions until a best possible regularized solution has been found.
We have shown that LSQR resembles the TSVD method for these
two kinds of problems. For mildly ill-posed problems, we have proved
that LSQR generally has only the partial regularization
since a small Ritz value generally appears before all
the needed dominant SVD components are captured. Since CGLS is mathematically
equivalent to LSQR, our assertions on the full or partial regularization
of LSQR apply to CGLS as well.

We have derived bounds for the diagonals and subdiagonals of bidiagonal
matrices generated by Lanczos bidiagonalization. Particularly,
we have proved that they decay as fast as the singular values of $A$
for severely ill-posed problems or moderately ill-posed problems
with $\rho>1$ or $\alpha>1$ suitably and decay more slowly
than the singular values of $A$ for mildly ill-posed problems.
These bounds are of theoretical and practical importance, and they
can be used to identify the degree of ill-posedness without extra cost
and decide the full or partial regularization of LSQR.

Based on some of the results established for LSQR, we have derived accurate
estimates for the accuracy of the rank $k$ approximations to $A$ and
$A^TA$ that are involved in CGME and LSMR, respectively. We have analyzed the
behavior of the smallest singular values of the projected matrices associated
with CGME and LSMR. Using these results, we have shown that
LSMR has the full regularization for severely and moderately ill-posed problems
with $\alpha>1$ and $\alpha>1$ suitably, and it generally has
only the partial regularization for mildly ill-posed probolems. In the meantime,
we have shown that the regularization of CGME has indeterminacy
and is inferior to LSQR and LSMR for each of three kinds of ill-posed problems.
In addition, our results have indicated
that the rank $k$ approximations to $A$ generated by Lanczos
bidiagonalization are substantially more accurate than those obtained
by standard randomized algorithms \cite{halko11}
and the strong RRQR factorizations \cite{gu96}.

With a number of nontrivial modifications and reformulations,
we have shown how to extend all the results obtained
for LSQR, CGME and LSMR to the case that $A$ has multiple singular
values.

We have made detailed and illuminating numerical experiments and confirmed our
theory on LSQR. We have also compared LSQR with GMRES and RRGMRES,
showing that the latter two methods do not general regularizing
effects and fail to deliver regularized solutions for general
nonsymmetric ill-posed problems. Theoretically, this is due to the fact that
GMRES and RRGMRES may work and have regularizing effects only
for (nearly) symmetric or, more generally,
(nearly) normal ill-posed problems, for which the left and right singular
vectors are (nearly) identical to the eigenvectors of $A$.

Our analysis approach can be adapted to MR-II for symmetric ill-posed problems,
and similar results
and assertions are expected for three kinds of symmetric ill-posed problems.
Using a similar approach to that in \cite{huangjia},
the authors \cite{huang15} has made an initial regularization analysis on
MR-II and derived the corresponding $\sin\Theta$ bounds,
which are too large overestimates. Our approach are applicable
to the preconditioned CGLS (PCGLS) and LSQR (PLSQR) \cite{hansen98,hansen10}
by exploiting the transformation technique originally proposed
in \cite{bjorck79,elden82} and advocated in \cite{hanke92,hanke93,hansen07}
or the preconditioned MR-II \cite{hansen10,hansen06},
all of which correspond to a general-form Tikhonov regularization involving the
matrix pair $\{A,L\}$, in which the regularization term $\|x\|^2$ is replaced by
$\|Lx\|^2$ with some $p\times n$ matrix $L\not=I$. It
should also be applicable to the mathematically equivalent LSQR
variant \cite{kilmer07} that is based on a joint bidiagonalization of
the matrix pair $\{A,L\}$ that corresponds to the above general-form Tikhonov
regularization. In this setting, the Generalized SVD (GSVD)
of $\{A,L\}$ or the mathematically equivalent SVD of $AL_A^{\dagger}$
will replace the SVD of $A$ to play a central role in analysis, where
$L_A^{\dagger}=\left(I-\left(A(I-L^{\dagger}L)^{\dagger}A\right)\right)^{\dagger}
L^{\dagger}$ is call the {\em $A$-weighted generalized inverse of $L$}
and $L_A^{\dagger}=L^{-1}$ if $L$ is square and invertible;
see \cite[p.38-40,137-38]{hansen98} and \cite[p.177-183]{hansen10}.

Finally, we highlight on hybrid Krylov iterative solvers and
make some remarks, which deserve particular and enough attention in our opinion.
Because of lack of a complete regularized theory
on LSQR, in order to find a best possible regularized solution for a
given \eqref{eq1}, one has commonly been using some
hybrid LSQR variants without considering the degree of ill-posedness
of \eqref{eq1}; see, e.g., \cite{aster,hansen98,hansen10} and the related papers
mentioned in the introduction. The hybrid CGME \cite{hanke01} and
CGLS \cite{aster,hansen10} have also been used.
However, Bj\"{o}rck \cite{bjorck94} has addressed that the hybrid LSQR variants
are mathematically complicated, and pointed out that it is hard to find
reasonable regularization parameters and tell when to stop them reliably.

For a hybrid LSQR variant, or more generally,
for any hybrid Krylov solver that first projects and then
regularizes \cite{hansen98,hansen10}, the situation is
more serious than what has been realized. It has long commonly accepted that
the approach of {\em "first-regularize-then-project"}
is equivalent to the approach of {\em "first-project-then-regularize"} and
they produce the same solution; see Section 6.4 and Figure 6.10
of \cite{hansen10}. This equivalence seems natural. Unfortunately, they
are {\em not} equivalent when solving \eqref{eq1}. Their equivalence requires
the assumption that the same regularization
parameter $\lambda$ in Tikhonov regularization is {\em used}, so that
both of them solve the same problem and
compute the same regularized solution. However, as far as
regularization methods are concerned, the fundamental point is
that each of the two approaches {\em must determine} its own optimal
regularization parameter $\lambda$ which is {\em unknown} in advance.
Mathematically,
for the approach of "first-regularize-then-project", there is an optimal
$\lambda$ since \eqref{eq1} satisfies the Picard condition, though its
determination is generally costly and may not be computationally viable
for a large \eqref{eq1}. On the contrary, for the approach of
"first-project-then-regularize", one must determine its optimal $\lambda$
for each projected problem, so one will have a sequence of optimal
$\lambda$'s. Whether or not they converge to the optimal regularization
parameter of \eqref{tikhonov} is unclear and lacks theoretical evidence.
For discrete regularization parameters in the TSVD method
for \eqref{eq1} and each of the projected problems, the situation is
similar. Unfortunately, for projected problems, their optimal
regularization parameters and their determination may encounter
insurmountable mathematical and numerical difficulties,
as we will clarify below.

As is well known, the Picard condition is an absolutely necessary condition
for the existence of the squares integrable solution to a linear compact
operator equation; without it, regularization
would be out of the question; see, e.g., \cite{engl00,kirsch,mueller}. This is
also true for the discrete linear ill-posed problem, where
the discrete Picard condition means that $\|x_{true}\|\leq C$ uniformly with some
(not large) constant $C$ such that regularization is
useful to compute a meaningful approximation to it \cite{hansen98,hansen10}.
Nevertheless, to the best of our knowledge, the discrete Picard conditions
for projected problems arising from LSQR or any other Krylov iterative solver
have been paid little attention until very recently \cite{gazzola15}.
Unfortunately, a fatal problem is that {\em the discrete Picard conditions are
not necessarily satisfied for the projected problems}. In \cite{gazzola15},
taking $e=\mathbf{0}$, i.e., $b=\hat{b}$ {\em noise free}, the authors
have proved that the discrete Picard conditions are satisfied
or inherited for the projected problems
under the {\em absolutely necessary assumption} that the
$k$ Ritz values, i.e., the singular values of the projected matrix at
iteration $k$, approximate the $k$ large singular values of $A$ in natural
order, regularization makes sense and can be used to solve the projected problems.
However, as have been stated in \cite{hansen98,hansen10} and
highlighted in this paper, under such assumption, Krylov solvers
themselves will find best possible regularized solutions at
semi-convergence, and there is no need to continue iterating and
regularize the projected problems at all, that is, no hybrid variant
is needed. On the other hand, if the
$k$ Ritz values do not approximate the $k$ large singular values of $A$
in natural order and at least one Ritz value smaller than $\sigma_{k_0+1}$
appears before $k\leq k_0$, the discrete Picard
conditions are {\em essentially not satisfied any longer} for the projected
problems starting from such $k$ onwards. If so, regularization applied to
projected problems is mathematically groundless and numerically may
lead to unavoidable failure.

\begin{figure}[htp]
\begin{minipage}{0.48\linewidth}
  \centerline{\includegraphics[width=6.0cm,height=5cm]{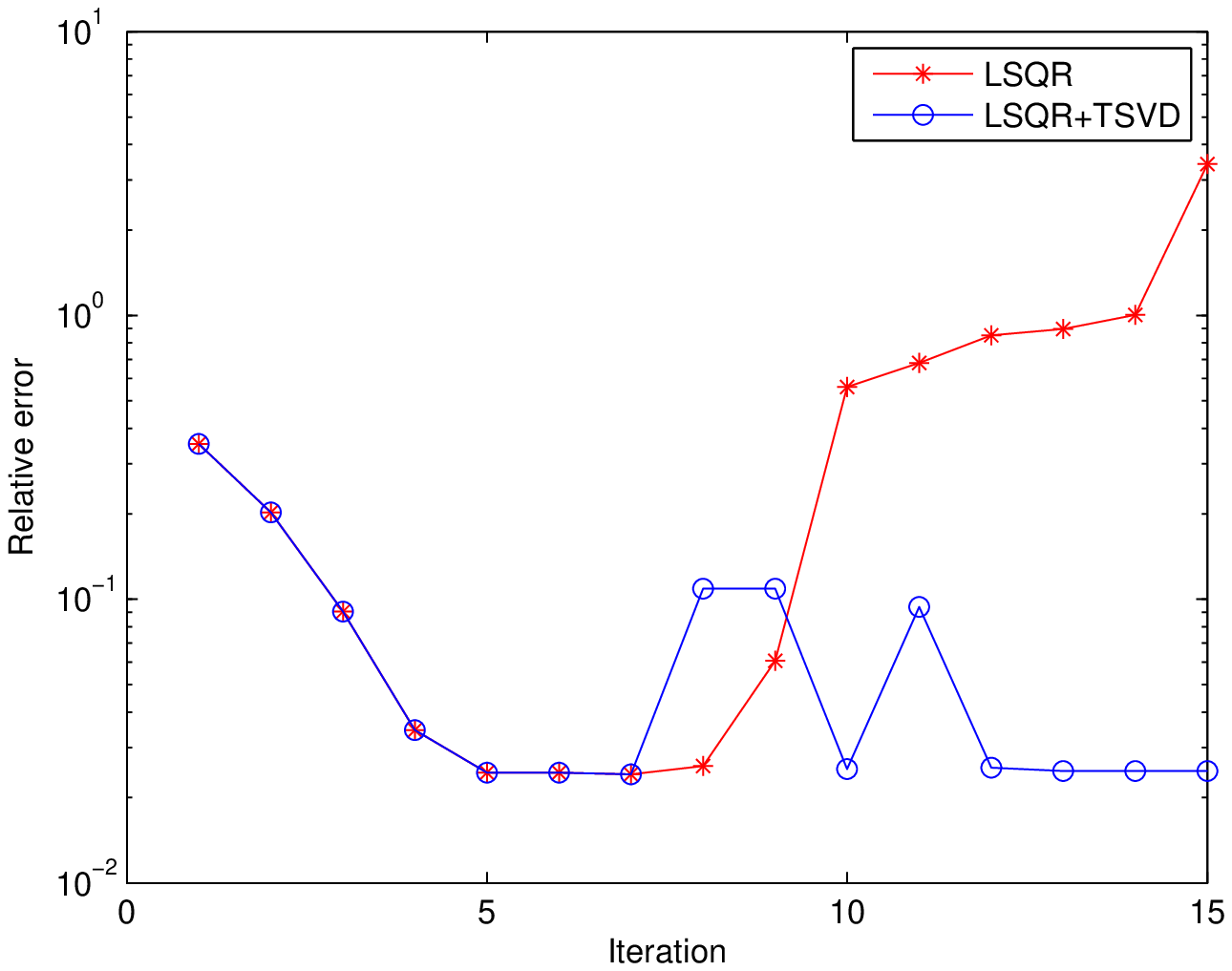}}
  \centerline{(a)}
\end{minipage}
\hfill
\begin{minipage}{0.48\linewidth}
  \centerline{\includegraphics[width=6.0cm,height=5cm]{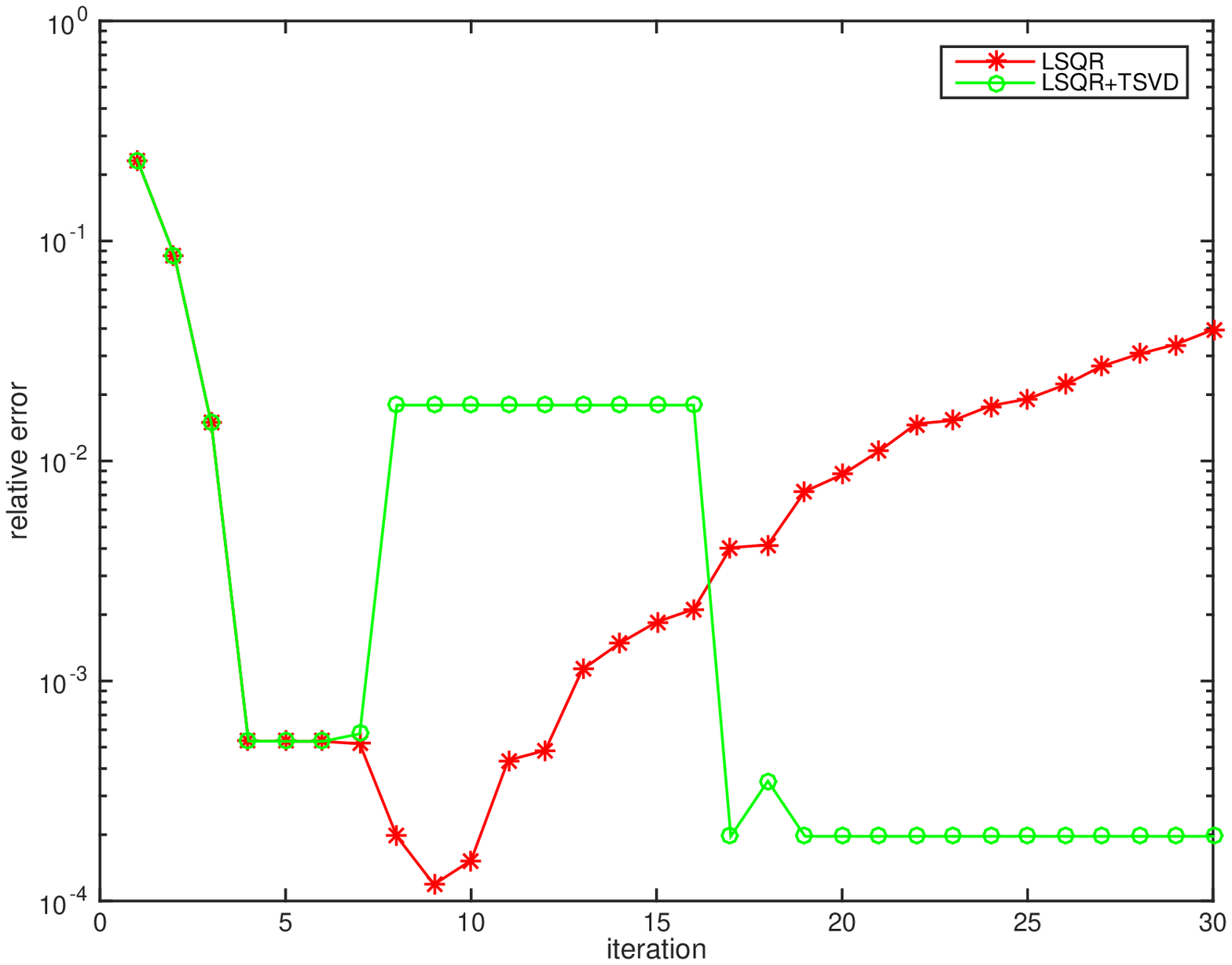}}
  \centerline{(b)}
\end{minipage}
\vfill
\begin{minipage}{0.48\linewidth}
  \centerline{\includegraphics[width=6.0cm,height=5cm]{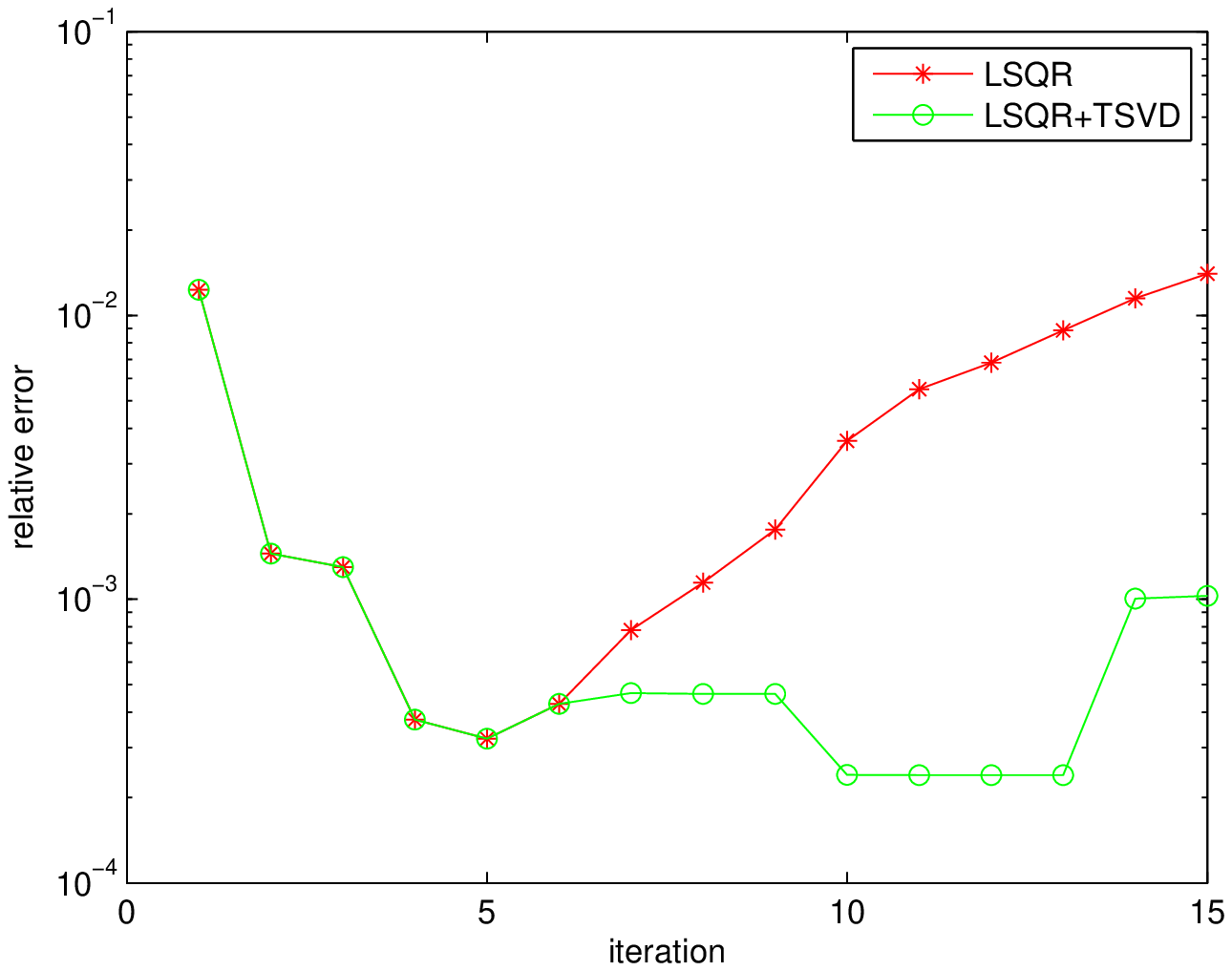}}
  \centerline{(c)}
\end{minipage}
\hfill
\begin{minipage}{0.48\linewidth}
  \centerline{\includegraphics[width=6.0cm,height=5cm]{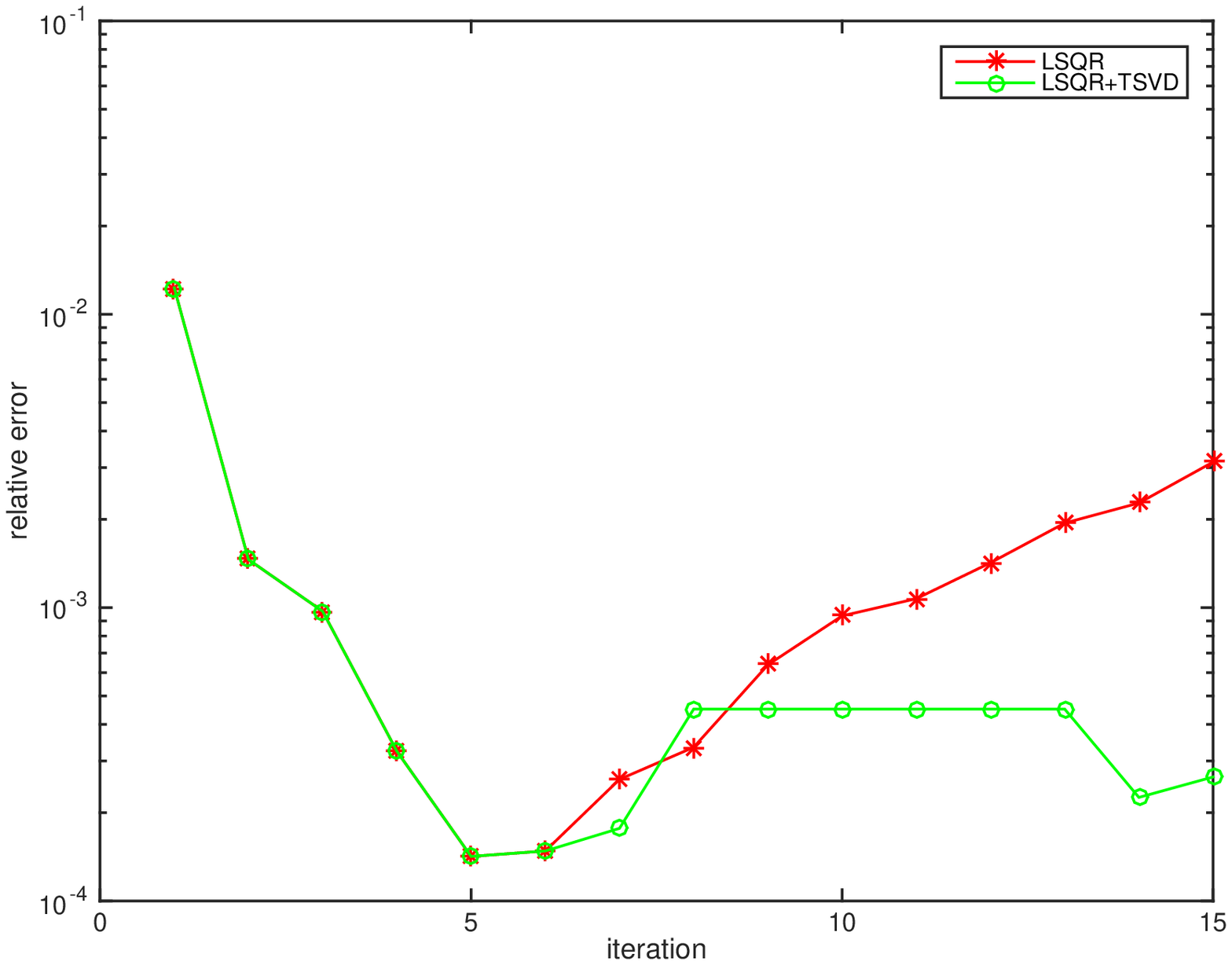}}
  \centerline{(d)}
\end{minipage}
\caption{(a)-(b): The relative errors $\|x^{(k)}-x_{true}\|/\|x_{true}\|$
by LSQR and the hybrid LSQR for {\sf phillips} of $n=1,024$ and $10,240$ with
$\varepsilon=10^{-3}$ and $10^{-4}$; (c)-(d):
The relative errors $\|x^{(k)}-x_{true}\|/\|x_{true}\|$ by LSQR and the hybrid LSQR
for {\sf deriv2} of $n=1,024$ and $10,240$
with $\varepsilon=10^{-3}$.} \label{figerratic}
\end{figure}

We take LSQR as an example for a precise statement on the
discrete Picard conditions for projected problems. Recall that, in the
projected problem \eqref{yk}, the noisy right-hand side is
$\|b\|e_1^{(k+1)}=P_{k+1}^Tb$ with $b=\hat{b}+e$ and the noise-free right-hand
side is $P_{k+1}^T\hat{b}$. Then for
$k=1,2,\ldots,n-1$ and for $n$ arbitrarily large and $\sigma_n\rightarrow 0$
(cf. \cite{gazzola15}), the discrete Picard conditions for
the projected problems are
$$
{\sup}_{k,n}\|B_k^{\dagger}P_{k+1}^T\hat{b}\|\leq C
$$
uniformly with some constant $C$. Numerically, for a given $n$ and
$\sigma_n$ close to zero arbitrarily, once
$\|B_k^{\dagger}P_{k+1}^T\hat{b}\|$ is very large for some $k$, then
the discrete Picard condition actually fails for the corresponding
projected problem. In this case, for the projected
problem, it is hard to apply regularization to the projected problem
and speak of its optimal regularization parameter,
which does not exist at all in the extreme case
that $\|B_k^{\dagger}P_{k+1}^T\hat{b}\|$ is
infinitely unbounded, which amounts to stating that $B_k$ has a
singular value close to zero arbitrarily. As a result,
any regularization applied to it works poorly. Indeed,
for {\sf phillips} and {\sf deriv2} of order $n=1,024$ and $10,240$,
we have observed that the hybrid LSQR exhibits considerable
{\em erratic} other than smooth curves of the errors between the regularized
solutions and $x_{true}$ in the dampening and stabilizing stage,
causing that the hybrid LSQR is unreliable to obtain a best regularized
solution; see Figure~\ref{figerratic}. Actually,
the regularized solutions obtained by the hybrid LSQR after its stabilization
are considerably less
accurate than those by the pure LSQR itself. For {\sf deriv2} of order $n=1,024$,
similar phenomena have also been observed for the hybrid MINRES and
MR-II \cite{huang15}.

The above phenomena are exactly due to the actual failure of the discrete
Picard conditions for the projected problems
because each of the projected matrices starts to have at least one
singular value {\em considerably} smaller than $\sigma_{k_0+1}$ from some
iteration $k\leq k_0$ onwards, which and whose corresponding (left
and right) Ritz vectors does not approximate any singular
triplet of $A$ well. A consequence of such actual
failure is that it is hard to reliably stop the hybrid variants at
right iteration in order to ultimately find a best regularized solution.
Therefore, for the mildly ill-posed problems and moderately ill-posed problems
with $\alpha>1$ not enough, it is appealing to seek other mathematically solid
and computationally viable variants of LSQR, LSMR and MR-II so that
best possible regularized solutions can be found.

\section*{Acknowledgements} I thank Dr. Yi Huang and Mrs. Yanfei Yang for
running the numerical experiments. I am grateful to Professors
\AA. Bj\"{o}rck, P. C. Hansen, L. Reichel and D.~P. O'Leary for their
comments and suggestions that helped improve the presentation of this paper.


\begin{thebibliography}{1}

\bibitem{aster} {\sc R.~C. Aster, B. Borchers and C.~H. Thurber},
{\em Parameter Estimation and Inverse Problems}, Second Edition, Elsevier,
New York, 2013.

\bibitem{reichel05} {\sc J. Baglama and L. Reichel},
{\em Augmented implicitly restarted Lanczos bidiagonalization methods},
SIAM J. Sci. Comput., 27 (2005), pp.~19--42.

\bibitem{baglama07} \sameauthor, {\em Augmented GMRES-type methods},
Numer. Linear Algebra Appl., 14 (2007), pp.~337--350.

\bibitem{bai} {\sc Z. Bai, J. Demmel, J. Dongarra, A. Ruhe and
H.~A. van der Vorst}, {\em Templates for the
Solution of Algebraic Eigenvalue Problems: A Practical Guide},
SIAM, Philadelphia, PA, 2000.

\bibitem{bauer11} {\sc F. Bauer and M.~A. Lukas}, {\em Comparing parameter
choice methods for regularization of ill-posed problems}, Math. Comput. Simul.,
81 (2011), pp.~1795--1841.

\bibitem{bazan10} {\sc F.~S.~V. Baz\'{a}n and L.~S. Borges}, {\em GKB-FP:
an algorithm for large-scale discrete ill-posed problems},
BIT Numer. Math., 50 (2010), pp.~481--507.

\bibitem{bazan14} {\sc F.~S.~V. Baz\'{a}n, M.~C.~C. Cunha and L.~S. Borges},
{\em Extension of GKB-FP algorithm to large-scale
general-form Tikhonov regularization}, Numer. Linear Algebra Appl., 21 (2014),
pp.~316--339.

\bibitem{berisha} {\sc S~Berisha and J.~G. Nagy}, {\em Restore Tools:
Iterative methods for image restoration}, 2012. Available from
{\sf http://www.mathcs.emory.edu/$^\sim$nagy/RestoreTools}.

\bibitem{bjorck88}{\sc \AA. Bj\"{o}rck}, {\em A bidiagonalization algorithm
 for solving large and sparse ill-posed systems of linear equations},
 BIT Numer. Math., 28 (1988), pp.~659--670.

\bibitem{bjorck96} \sameauthor, {\em Numerical Methods for Least Squares
Problems}, SIAM, Philadelphia, PA, 1996.

\bibitem{bjorck15} \sameauthor, {\em Numerical Methods in Matrix Computations},
Texts in Applied Mathematics 59, Springer, 2015.

\bibitem{bjorck79} {\sc \AA. Bj\"{o}rck and L. Eld\'{e}n},
{\em Methods in numerical algebra for ill-posed problems}, Report LiTH-R-33-1979,
Dept. of Mathematics, Link\"{o}ping Univeristy, Sweden, 1979.
Proceedings of the International Symposium on Ill-posed Problems:
Theory and Practice, University of Delaware, Newark, Delaware,
Oct. 2--6, 1979.

\bibitem{bjorck94} {\sc \AA. Bj\"orck, E.~Grimme and P. Van Dooren},
{\em An implicit shift bidiagonalization algorithms for ill-posed problems},
BIT Numer. Math., 34 (1994), pp.~510--534.

\bibitem{calvetti99} {\sc D. Calvetti, G.~H. Golub and L. Reichel},
{\em Estimation of the L-curve via Lanczos bidiagonalization},
BIT Numer. Math., 39 (1999), pp.~603--619.

\bibitem{calvetti02b}{\sc D. Calvetti, P. C. Hansen and L. Reichel},
{\em L-curve curvature bounds via Lanczos bidiagonalization},
Electr. Trans. Numer. Anal., 14 (2002), pp.~20--35.

\bibitem{calvetti00} {\sc D. Calvetti, B. Lewis and L. Reichel},
{\em GMRES-type methods for inconsistent systems}, Linear Algebra Appl.,
316 (2000), pp.~157--169.

\bibitem{calvetti01} \sameauthor,
{\em On the choice of subspace for iterative methods for
linear ill-posed problems}, Int. J. Appl. Math. Comput. Sci.,
11 (2001), pp.~1069--1092.

\bibitem{calvetti02} \sameauthor,
{\em GMRES, L-curve, and discrete ill-posed problems},
BIT Numer. Math., 42 (2002), pp.~44--65.

\bibitem{calvetti02c} \sameauthor, {\em On the regularizing properties of the
GMRES method}, Numer. Math., 91 (2002), pp.~605--625.

\bibitem{calvetti00b} {\sc D. Calvetti, S. Morigi, L. Reichel and F. Sgallari},
{\em Tikhonov regularization and the L-curve for large discrete ill-posed problems},
J. Comput. Appl. Math., 123 (2000), pp.~423--446.

\bibitem{calvetti03} {\sc D. Calvetti and L. Reichel},
{\em Tikhonov regularization of large linear problems},
BIT Numer. Math., 43 (2003), pp.~263--283.

\bibitem{craig} {\sc E.~J. Craig}, {\em The n-step iteration procedures},
J. Math. Phys., 34 (1955), pp.~64--73.

\bibitem{chung08} {\sc J. Chung, J.~G. Nagy and D.~P. O'Leary},
{\em A weighted GCV method for Lanczos hybrid regularization},
Electr. Trans. Numer. Anal., 28 (2008), pp.~149--167.

\bibitem{eicke} {\sc B. Eicke, A.~K. Lious and R. Plato},
{\em The instability of some gradient methods for ill-posed problems},
Numer. Math., 58 (1990), pp.~129--134.

\bibitem{elden82} {\sc L.~Eld\'{e}n}, {\em A weigthed pseudoinverse, generalized
singular values and constrained least squares problems}, BIT, 22 (1982), pp.~487--501.

\bibitem{engl93} {\sc H.~W. Engl}, {\em Regularization methods for the
stable solution of inverse problems}, Surveys Math. Indust., 3 (1993),
pp.~71--143.

\bibitem{engl00} {\sc H.~W. Engl, M. Hanke and A. Neubauer},
{\em Regularization of Inverse problems}, Kluwer Academic Publishers,
2000.

\bibitem{firro97} {\sc R.~D. Fierro, G.~H. Golub, P.~C. Hansen and D.~P. O'Leary},
{\em Regularization by the truncated total least squares},
SIAM J. Sci. Comput., 18 (1997), pp.~1223--1241.

\bibitem{fischer96} {\sc B. Fischer, M. Hanke and M. Hochbruck},
{\em A note on conjugate-gradient type methods for
indefinite and/or inconsistent linear systems}, Numer. Algor., 11
(1996), pp.~181--187.

\bibitem{fong} {\sc D.~C.~-L. Fong and M.~Saunders},
{\em LSMR: an iterative algorithm for sparse least-squares problems},
SIAM J. Sci. Comput., 33 (2011), pp.~2950--2971.

\bibitem{gazzola14} {\sc S. Gazzola},
{\em Regularization techniques based on Krylov methods for
ill-posed linear systems}, Ph. D. thesis, Department of Mathematics,
University of Padua, Italy, 2014.

\bibitem{gazzola13} {\sc S. Gazzola and P. Novati},
{\em Multi-parameter Arnoldi-Tikhonov methods}, Electr.
Trans. Numer. Anal., 40 (2013), pp.~452--475.

\bibitem{gazzola15} \sameauthor, {\em Inheritance of the discrete
Picard condition in Krylov subspace methods}, BIT Numer. Math.,
56 (2016), pp.~893--918.

\bibitem{gazzola14b} {\sc S. Gazzola, P. Novati and M. R. Russo},
{\em Embedded techniques for choosing the parameter in
Tikhonov regularization}, Numer. Linear Algebra Appl., 21 (2014), pp.~796--812.

\bibitem{gazzola-online}{\sc S. Gazzola, P. Novati and M. R. Russo},
{\em On Krylov projection methods and Tikhonov regularization},
Electr. Trans. Numer. Anal., 44 (2015), pp.~83---123.

\bibitem{gazzola16} {\sc S.~Gazzola, E.~Onunwor, L.~Reichel and G.~Rodriguez},
{\em On the Lanczos and Golub-Kahan reduction methods applied to
discrete ill-posed problems}, Numer. Linear Algebra Appl.,
23 (2016). pp.~187--204.

\bibitem{gilyza86} {\sc S.~F. Gilyazov},
{\em Regularizing algorithms based on the conjugate-gradient method}, U.S.S.R.
Comput. Maths. Math. Phys. 26 (1986), pp.~8--13.

\bibitem{gilyazov} {\sc S. F. Gilyazov and N. L. Gol'dman},
{\em Regularization of Ill-Posed Problems by Iteration Methods},
Kluwer Academic Publishers, Boston, 2010.

\bibitem{golub79} {\sc G.~H. Golub, M.~T. Heath and G. Wahba},
{\em Generalized cross-validation as a method for choosing a good ridge
parameter}, Technometrics, 21 (1979), pp.~215--223.

\bibitem{golub89} {\sc G.~H. Golub and D.~P. O'Leary},
{\em Some history of the conjugate gradient method and the Lanczos
algorithms: 1948--1976}, SIAM Rev., 31 (1989), pp.~50--102.

\bibitem{gu2015} {\sc M. Gu}, {\em Subspace iteration randomization and
singular value problems}, SIAM J. Sci. Comput., 37 (2015), pp.~A1139--A1173.

\bibitem{gu96} {\sc M. Gu and S.~C. Eisenstat},
{\em Efficient algorithms for computing a strong rank-revealing QR
factorization}, SIAM J. Sci. Comput., 17 (1996), pp.~848--869.

\bibitem{halko11} {\sc N. Halko, P.~G. Martinsson
J.~A. Tropp}, {\em Finding structure with randomness:
probabilistic algorithms for constructing approximate
matrix decompositions}, SIAM Rev., 53 (2011), pp.~217--288.

\bibitem{hanke92} {\sc M. Hanke}, {\em Regularization with differential
operators: An iterative approach}, Numer. Func. Anal. Opt., 13 (1992),
pp.~523--540.

\bibitem{hanke95} {\sc M. Hanke}, {\em Conjugate Gradient Type Methods for
Ill-Posed Problems}, Longman, Essex, 1995.

\bibitem{hanke96a} \sameauthor, {\em Limitations of the L-curve method in
ill-posed problems}, BIT Numer. Math., 36 (1996), pp.~287--301.

\bibitem{hanke01} \sameauthor, {\em On Lanczos based methods for the
regularization of discrete ill-posed problems}, BIT Numer. Math.,
41 (2001), pp.~1008--1018.

\bibitem{hanke93} {\sc M. Hanke and P.~C. Hansen},
{\em Regularization methods for large-scale problems}, Surveys Math. Indust.,
3 (1993), pp.~253--315.

\bibitem{hanke96} {\sc M. Hanke and J.~G. Nagy},
{\em Restoration of atmospherically blurred images by
symmetric indefinite conjugate gradient techniques}, Inverse Probl.,
12 (1996), pp.~157--173.

\bibitem{hansen90} {\sc P.~C. Hansen}, {\em The discrete Picard condition
for dicrete ill-posed problems}, BIT Numer. Math., 30 (1990), pp.~658--672.

\bibitem{hansen90b} \sameauthor, {\em
Truncated singular value decomposition solutions to discrete ill-posed problems
with ill-determined numerical rank}, SIAM J. Sci. and Stat.,
Comput., 11 (1990), pp.~503--518.

\bibitem{hansen92} \sameauthor, {\em Analysis of discrete ill-posed
problems by means of the L-curve}, SIAM Rev. 34 (1992), pp.~561--580.

\bibitem{hansen98} \sameauthor, {\em Rank-Deficient and Discrete Ill-Posed
Problems: Numerical Aspects of Linear Inversion}, SIAM, Philadelphia, PA, 1998.

\bibitem{hansen07} \sameauthor, {\em Regularization tools version 4.0
for Matlab 7.3}, Numer. Algor., 46 (2007), pp.~189--194.

\bibitem{hansen08} \sameauthor, {\em Regularization tools:
A Matlab package for
analysis and solution of discrete ill-posed problems
version 4.1 for Matlab 7.3}, 2008. Available from
{\sf www.netlib.org/numeralgo}.

\bibitem{hansen10} \sameauthor, {\em Discrete Inverse Problems:
Insight and Algorithms}, SIAM, Philadelphia, PA, 2010.

\bibitem{hansen06} {\sc P.~C. Hansen and T.~K. Jensen},
{\em Smoothing-norm preconditioned for regularizing minimum-residual
methods}, SIAM J. Matrix Anal. Appl., 29 (2006), pp.~1--14.

\bibitem{hansen93} {\sc P.~C. Hansen and D.~P. O'Leary},
{\em The use of the L-curve in the regularization of discrete ill-posed
problems}, SIAM J. Sci. Comput., 14 (1993), pp.~1487--1503.

\bibitem{hansen13} {\sc P.~C. Hansen, V.~Pereyra and G.~Scherer},
{\em Least Squares Data Fitting with Applications}, The Johns Hopkins
University Press, Baltimore, 2013.

\bibitem{hestenes} {\sc M.~R. Hestenes and E. Stiefel},
{\em Methods of conjugate gradients for solving linear systems},
J. Res. Nat. Bur. Stand., 49 (1952), pp.~409--436.

\bibitem{hps09} {\sc M. R. Hn\v{e}tynkov\'{a},
M. Ple\v{s}inger and Z. Strako\v{s}},
{\em The regularizing effect of the Golub-Kahan iterative bidiagonalization
and revealing the noise level in the data}, BIT Numer. Math., 49 (2009),
pp.~669--696.

\bibitem{hofmann86} {\sc B. Hofmann}, {\em Regularization for
Applied Inverse and Ill-posed Problems}, Teubner, Stuttgart, Germany, 1986.

\bibitem{hong92} {\sc Y.~T. Hong and C.~T. Pan},
{\em Rank-revealing QR decompositions and the singular
value decomposition}, Math. Comput., 58 (1992), pp.~213--232.

\bibitem{huangjia} {\sc Y. Huang and Z. Jia}, {\em Some results on the regularization
of LSQR for large-scale ill-posed problems},
Science China Math., doi: 10.1007/s11425-015-0568-4, 2016.

\bibitem{huang15} \sameauthor, {\em On regularizing effects of MINRES and MR-II
for large-scale symmetric discrete ill-posed problems},
arXiv: math.NA/1503.03936, 2015.

\bibitem{ito15} {\sc K.~Ito and B.~Jin}, {\em Inverse Problems: Tikhonov Theory and
Algorithms}, World Scientific Publishing, 2015.

\bibitem{jensen07} {\sc T.~K. Jensen and P.~C. Hansen}, {\em Iterative
regularization with minimum-residual methods}, BIT Numer. Math.,
47 (2007), pp.~103--120.

\bibitem{jia95} {\sc Z. Jia}, {\em The convergence of generalized Lanczos methods
for large unsymmetric eigenproblems}, SIAM J. Matrix Anal. Appl., 16 (1995),
pp.~843--862.

\bibitem{jia98} \sameauthor, {\em Generalized block Lanczos methods
for large unsymmetric eigenproblems}, Numer. Math., 80 (1998), pp.~239--266.

\bibitem{jia05} \sameauthor, {\em The convergence of harmonic Ritz values,
harmonic Ritz vectors and refined harmonic Ritz vectors}, Math. Comput.,
74 (2005), pp.~1441--1456.

\bibitem{jia03} {\sc Z. Jia and D. Niu},
{\em An implicitly restarted
bidiagonalization Lanczos method for computing a partial singular value
decomposition}, SIAM J. Matrix Anal. Appl., 25 (2003), pp.~246--265.

\bibitem{jia10} \sameauthor, {\em A refined
harmonic Lanczos bidiagonalization method and
an implicitly restarted algorithm for computing the smallest
singular triplets of large matrices}, SIAM J. Sci. Comput.,
32 (2010), pp.~714--744.

\bibitem{jia01} {\sc Z. Jia and G.~W. Stewart}, {\em  An analysis of the
Rayleigh--Ritz method for approximating eigenspaces}, Math. Comput.,
70 (2001), pp,~637--647.

\bibitem{johnsson} {\sc C. Johnsson},
{\em On finite element methods for optimal control problems}, Tech. Report
79-04 R, Dept. of Computer Science, University of Gothenburg, 1979.

\bibitem{kaipio} J.~Kaipio and E. Somersalo, {\em Statistical and
Computational Inverse Problems}, Applied Mathematical Sciences 160,
Springer, 2005.

\bibitem{kern} {\sc M.~Kern}, {\em Numerical Methods for Inverse Problems},
John Wiley \& Sons, Inc., 2016.

\bibitem{kilmer07} {\sc M.~E. Kilmer, P.~C. Hansen and
M.~I. Espa\~{n}ol},
{\em A projection-based approach to general-form Tikhonov regularization},
SIAM J. Sci. Comput., 29 (2007), pp.~315--330.

\bibitem{kilmer03} {\sc M.~E. Kilmer and D.~P. O'Leary},
{\em Choosing regularization parameters in iterative methods for
ill-posed problems},
SIAM J. Matrix Anal. Appl., 22 (2001), pp.~1204--1221.

\bibitem{kilmer99} {\sc M.~E. Kilmer and G.~W. Stewart},
{\em Iterative regularization and MINRES}, SIAM J. Matrix Anal.
Appl., 21 (1999), pp.~613--628.

\bibitem{kindermann}{\sc S.~Kindermann}, {\em Convergence analysis of minimization-based
noise level-free parameter choice rules for linear ill-posed problems},
Electr. Trans. Numer. Math., 38 (2011), pp.233--257.

\bibitem{kirsch} {\sc A. Kirsch}, {\em An Introduction to the Mathematical
Theory of Inverse Problems}, Second Edition, Applied Mathematical Sciences 120,
Springer, 2011.

\bibitem{kythe} {\sc P.~K. Kythe and P.~Puri},
{\em Computational Methods for Linear Integral Equations},
Birkh\"{a}user, Boston/Basel/Berlin, 2002.

\bibitem{lanczos} {\sc C.~C. Lanczos}, {\em An iteration method for the
solution of the eigenvalue problem of linear differential and integral
operators}, J. Res. Nat. Bur. Stand., 45 (1950), pp.~255--282.

\bibitem{lawson} {\sc R.~A. Lawson and R.~J. Hanson},
{\em Solving Least Squares Problems}, Prentice-Hall, Englewood Cliffs, NJ,
1974; reprinted by SIAM, Philadelphia, PA, 1995.

\bibitem{lewis09} {\sc B. Lewis and L. Reichel},
{\em Arnoldi-Tikhonov regularization methods}, J. Comput. Appl. Math.,
226 (2009), pp.~92--102.

\bibitem{meurant} {\sc G. Meurant}, {\em The Lanczos and Conjugate Gradient
Algorithms: From Theory to Finite Precision Computations},
SIAM, Philadelphia, PA, 2006.

\bibitem{miller} {\sc K. Miller}, {\em Least squares methods for ill-posed
problems with a prescribed bound}, SIAM J. Math. Anal., 1 (1970), pp.~52--74.

\bibitem{morozov} {\sc V.~A. Morozov}, {\em On the solution of functional
equations by the method of regularization}, Soviet Math. Dokl., 7 (1966),
pp.~414--417.

\bibitem{mueller} {\sc J.~L. Mueller and S. Siltanen}, {\em
Linear and Nonlinear Inverse Problems with Practical Applications},
SIAM, Philadelpha, PA, 2012.

\bibitem{natterer} {\sc F. Natterer}, {\em The Mathematics of Computerized
Tomography}, Reprinted version of the 1986 edition published by Wiley and Teubner,
SIAM, Philadelphia, PA, 2001.

\bibitem{nemi} {\sc A. S. Nemirovskii}, {\em The regularizing properties
of the adjoint gradient method in ill-posed problems}, U.S.S.R. Comput.
Maths. Math. Phys., 26 (1986), pp.~7--16.

\bibitem{neumaier98} {\sc A.~Neumaier}, {\em Solving ill-conditioned and singular
linear systems: a tutorial on regularization}, SIAM Rev., 40 (1998), pp.~636--666.

\bibitem{neuman12} {\sc A. Neuman, L. Reichel and H. Sadok},
{\em Algorithms for range restricted iterative methods for linear
dicrete ill-posed problems}, Numer. Algor., 59 (2012), pp.~325--331.

\bibitem{nolet} {\sc G. Nolet}, {\em Solving or resolving inadequate and noisy
tomographic systems}, J. Comput. Phys., 61 (1985), pp.~463--482.

\bibitem{novati13} {\sc P. Novati and M.~R. Russo}, {\em A GCV based
Arnoldi-Tikhonov regularization method}, BIT Numer. Math.,
54 (2014), pp.~501--521.

\bibitem{oleary81} {\sc D.~P. O'Leary and J.~A. Simmons},
{\em A bidiagonalization-regularization procedure for large
scale discretizations of ill-posed problems},
SIAM J. Sci. Statist. Comput., 2 (1981), pp.~474--489.

\bibitem{paige75}{\sc C.~C. Paige and M.~A. Saunders},
{\em Solution of sparse indefinite systems of linear equations}.
SIAM J. Numer. Anal., 12 (1975), pp.~617--629.

\bibitem{paige82} \sameauthor,
{\em LSQR: an algorithm for sparse linear equations and
sparse least squares}, ACM Trans. Math. Softw., 8 (1982), pp.~43--71.

\bibitem{paige06} {\sc C.~C. Paige and Z.~Z. Strako\v{s}},
{\em Core problems in linear algebraic systems}, SIAM J. Matrix Anal.
Appl., 27 (2006), pp.~861--875.

\bibitem{parlett} {\sc B.~N. Parlett}, {\em The Symmetric Eigenvalue
Problem}, SIAM, Philadelpha, PA, 1998.

\bibitem{phillips} {\sc D.~L. Phillips},
{\em A technique for the numerical solution of certain integral equations
of the first kind}, J. ACM, 9 (1962), pp.~84--97.

\bibitem{reichel13} {\sc L. Reichel and G. Rodriguez}, {\em Old and new
parameter choice rules for discrete ill-posed problems},
Numer. Algor., 63 (2013), pp.~65--87.

\bibitem{renaut} {\sc R.~A. Renaut, S.~Vatankhah and V.~E. Ardestani},
{\em Hybrid and iteratively reweighted regularization by
unbiased predictive risk and weighted GCV}, arXiv: math.NA/1509.00096v1, 2015.

\bibitem{saad} {\sc Y. Saad}, {\em Numerical Methods for Large Eigenvalue
Problems}, Second Edition, SIAM, Philadelphia, PA, 2011.

\bibitem{scales} {\sc J.~A. Scales and A. Gerztenkorn}, {\em Robust methods
in inverse theory}, Inverse Probl., 4 (1988), pp.~1071--1091.

\bibitem{squire} {\sc W. Squire}, {\em The solution of ill-conditioned linear
systems arising from Fredholm equations of the first kind by steepest descents
and conjugate gradients}, Int. J. Numer. Meth. Eng., 10 (1976), pp.~607--617.

\bibitem{stewart98} {\sc G.~W. Stewart}, {\em Matrix Algorithms I:
Basic Decompositions}, SIAM, Philadelphia, PA, 1998.

\bibitem{stewart01} \sameauthor, {\em Matrix Algorithms II:
Eigensystems}, SIAM, Philadelphia, PA, 2001.

\bibitem{stewartsun} {\sc G.~W. Stewart and J.-G Sun}, {\em
Matrix Pertubation Theory}, Academic Press, Boston, 1990.

\bibitem{tal} {\sc A.~A. Tal}, {\em Numerical solution of Fredholm integral
equations of the first kind}, TR-66-34, Computer Science Center,
University of Maryland, College Park, MD, 1966.

\bibitem{tikhonov63} {\sc A.~N. Tikhonov}, {\em Solution of incorrectly
formulated problems and the regulariza-tion method}, Dokl. Akad. Nauk.
SSSR, 151 (1963), pp.~501--504. Soviet Math. Dokl., 4 (1963), pp.~1035--1038.

\bibitem{tikhonov77} {\sc A.~N. Tikhonov and V.~Y. Arsenin},
{\em Solutions of Ill-Posed Problems}, Winston \& Sons, Washington, D.C., 1977.

\bibitem{vorst86} {\sc A. Van der Sluis and H.~A. Van der Vorst},
{\em The rate of convergence of conjugate gradients},
Numer. Math., 48 (1986), pp.~543--560.

\bibitem{vorst90} \sameauthor,
{\em SIRT- and CG-type methods for iterative solution of sparse
linear least squares problems}, Linear Algebra Appl., 130 (1990), pp.~257--302.

\bibitem{vorst02} {\sc H.~A. Van der Vorst}, {\em Computational Methods for
Large Eigenvalue Problems}, In: P.~G. Ciarlet and F.
Cucker (eds.), Handbook of Numerical Analysis, vol. VIII, pp.~3--179. North
Holland Elsevier, Amsterdam (2002).

\bibitem{huffel} {\sc S. Van Huffel and P. Lemmerling} (eds.),
{\em Total Least Squares and Errors-in-Variables Modeling:
Analysis, Algorithms and Applications}, Kluwer Academic Publishers,
Boston, 2002.

\bibitem{varah79} {\sc J.~M. Varah}, {\em A practical examination of some
numerical methods for linear discrete ill-posed problems}, SIAM Rev., 21 (1979),
pp.~100--111.

\bibitem{vogel96} {\sc C.~R. Vogel}, {\em Non-convergence of the L-curve
regularization parameter selection method},
Inverse Probl., 12 (1996), pp.~535--547.

\bibitem{vogel02} \sameauthor, {\em Computational Methods for Inverse
Problems}, SIAM, Philadelphia, PA, 2002.

\bibitem{wahba} {\sc G. Wahba}, {\em Practical approximate solutions to linear
operator equations when the data are noisy}, SIAM J. Numer. Anal.,
14 (1977), pp.~651--667.

\bibitem{wilkinson} {\sc J.~H. Wilkinson}, {\em The Algebraic Eigenvalue Problem},
Clarendon Press, Oxford, 1965; reprinted in 2004.
\end{thebibliography}
\end{document}